\documentclass[oneside,a4paper]{amsart}
\usepackage[T1]{fontenc}
\usepackage[utf8]{inputenc}
\usepackage[english]{babel}
\usepackage[autostyle, english=american]{csquotes}
\MakeOuterQuote{"}

\usepackage{xcolor}
\definecolor{reference}{rgb}{0.20,0.36,0.74}
\definecolor{citation}{rgb}{0,.40,.80}
\usepackage[colorlinks,linktocpage,urlcolor=blue,linkcolor=reference,citecolor=citation]{hyperref}
\calclayout

\usepackage{mathrsfs}
\usepackage{mathtools}
\usepackage{latexsym}
\usepackage{graphicx}
\usepackage{amscd,amssymb,amsmath,amsbsy,amsfonts,amsthm}
\usepackage[mathscr]{eucal}
\usepackage[all, 2cell]{xy}
\UseTwocells
\xyoption{2cell}{\UseTwocells}
\usepackage{verbatim}
\usepackage{enumerate}
\usepackage{tikz}
\usepackage{enumitem}%[wide, topsep=\parskip, itemsep=\parskip, parsep=\parskip]
\setlist[enumerate,1]{label={(\arabic*)}}

\usepackage[backend=biber,style=alphabetic,giveninits=true,maxbibnames=99,maxalphanames=99]{biblatex}

\usepackage[capitalise]{cleveref}
\crefformat{equation}{(#2#1#3)}

\makeatletter
\def\cref@thmoptarg[#1]#2#3#4{%
    \ifhmode\unskip\unskip\par\fi%
    \normalfont%
    \trivlist%
    \let\thmheadnl\relax%
    \let\thm@swap\@gobble%
    \thm@notefont{\fontseries\mddefault\upshape}%
    \thm@headpunct{.}% add period after heading
    \thm@headsep 5\p@ plus\p@ minus\p@\relax%
    \thm@space@setup%
    #2% style overrides
    \@topsep \thm@preskip               % used by thm head
    \@topsepadd \thm@postskip           % used by \@endparenv
    \def\@tempa{#3}\ifx\@empty\@tempa%
      \def\@tempa{\@oparg{\@begintheorem{#4}{}}[]}%
    \else%
      \refstepcounter[#1]{#3}%  <<< cleveref modification
      \@namedef{cref@#3@alias}{#1}% added
      \def\@tempa{\@oparg{\@begintheorem{#4}{\csname the#3\endcsname}}[]}%
    \fi%
    \@tempa}%
\makeatother

\usepackage{relsize}
\usepackage[bbgreekl]{mathbbol}
\DeclareSymbolFontAlphabet{\mathbb}{AMSb} %to ensure that the meaning of \mathbb does not change
\DeclareSymbolFontAlphabet{\mathbbl}{bbold}

\newlength{\bibparskip}
\let\oldthebibliography\thebibliography
\renewcommand\thebibliography[1]{%
  \oldthebibliography{#1}%
  \setlength{\parskip}{\bibitemsep}%
  \setlength{\itemsep}{\bibparskip}%
}

\numberwithin{equation}{subsection}

\DeclareTextFontCommand{\emdef}{\bfseries}

\newcommand{\lasttwoofyear}{
  \expandafter\getlasttwo\number\numexpr\year\relax\relax
}
\def\getlasttwo#1#2#3#4\relax{#3#4}

\newtheorem{thm}{Theorem}

\newtheorem*{thm*}{Theorem}
\newtheorem{lem}[thm]{Lemma}

\newtheorem{cor}[thm]{Corollary}

\newtheorem*{cor*}{Corollary}
\newtheorem{prop}[thm]{Proposition}

\newtheorem*{prop*}{Proposition}

\newtheorem*{conj*}{Conjecture}

\newtheorem{thmx}{Theorem}

\theoremstyle{definition}
\newtheorem{defn}[thm]{Definition}

\newtheorem{construction}[thm]{Construction}
\newtheorem{notation}[thm]{Notation}

\newtheorem{assumption}[thm]{Assumption}

\newtheorem{conv}[thm]{Conventions}

\theoremstyle{remark}
\newtheorem{rem}[thm]{Remark}

\newtheorem{ex}[thm]{Example}

\newtheorem{warn}[thm]{Warning}

\newcommand{\todo}[1]{\textcolor{red}{TODO: #1}}

\definecolor{note_color}{rgb}{0.0,0.7,0.0}

\newcommand{\eset}{\varnothing}
\newcommand{\fcat}{\mathscr}
\newcommand{\inj}{\hookrightarrow}
\newcommand{\surj}{\twoheadrightarrow}
\newcommand{\areq}{\mathbin{{\xrightarrow{\,\sim\,}}}}

\newcommand{\tto}{{\xymatrix{\ar[r]&}}}

\renewcommand{\phi}{\varphi}
\renewcommand{\epsilon}{\varepsilon}

\DeclareMathOperator{\Vect}{Vect}

\DeclareMathOperator{\Sym}{Sym}
\newcommand{\Mod}{\mathrm{Mod}}

\newcommand{\RMod}{\mathrm{RMod}}

\DeclareMathOperator{\rk}{rk}

\DeclareMathOperator{\Prim}{Prim}

\DeclareMathOperator{\Rep}{Rep}

\DeclareMathOperator{\Coind}{coInd}

\renewcommand{\Im}{\mathrm{Im}}

\newcommand{\Top}{{\mathrm{Top}}}
\newcommand{\Orb}{{\mathrm{Orb}}}

\newcommand{\Glo}{{\mathrm{glo}}}
\DeclareMathOperator{\LocSys}{LS}

\DeclareMathOperator{\Spec}{Spec}
\DeclareMathOperator{\Spf}{Spf}

\DeclareMathOperator{\QCoh}{QCoh}

\newcommand{\sm}{{\mathrm{sm}}}

\DeclareMathOperator{\Nm}{Nm}

\newcommand{\Type}{ {\fcat{S}} }

\newcommand{\Cat}{ {\fcat C\mathrm{at}} }

\newcommand{\Prs}{ {\mathrm{Pr}} }

\DeclareMathOperator{\Ab}{Ab}
\DeclareMathOperator{\PreAb}{PreAb}
\DeclareMathOperator{\Grp}{Grp}
\DeclareMathOperator{\Hom}{Hom}

\newcommand{\HHom}{ {\fcat Hom} }

\DeclareMathOperator{\Map}{Map}
\DeclareMathOperator{\Fun}{Fun}

\DeclareMathOperator{\Sp}{Sp}

\DeclareMathOperator{\Ind}{Ind}
\DeclareMathOperator{\coInd}{coInd}

\DeclareMathOperator*{\fib}{\mathsf{fib}}

\DeclareMathOperator*{\cofib}{\mathsf{cofib}}

\DeclareMathOperator{\Alg}{Alg}
\DeclareMathOperator{\CAlg}{CAlg}
\DeclareMathOperator{\cCAlg}{cCAlg}

\DeclareMathOperator{\Lan}{Lan}
\DeclareMathOperator{\Ran}{Ran}

\DeclareMathOperator{\ev}{ev}

\newcommand{\prolim}[1][]{\lim\limits_{\xleftarrow[#1]{}}}
\newcommand{\indlim}[1][]{\lim\limits_{\xrightarrow[#1]{}}}

\DeclareMathOperator{\triv}{triv}

\DeclareMathOperator{\Id}{Id}
\DeclareMathOperator{\Rex}{\mathrm{Rex}}

\newcommand{\Shv}{{\mathcal S\mathrm{hv}}}
\newcommand{\op}{{\mathrm{op}}}

\newcommand{\lax}{{\mathrm{lax}}}
\newcommand{\oplax}{{\mathrm{oplax}}}
\newcommand{\llax}{{\mathrm{l.lax}}}
\newcommand{\rlax}{{\mathrm{r.lax}}}

\newcommand{\res}{{\mathrm{Res}}}
\newcommand{\stable}{{\mathrm{st}}}
\newcommand{\loc}{\mathrm{loc}}

\newcommand{\cSpec}{\mathrm{cSpec}}
\newcommand{\Hyp}{\mathrm{Hyp}}
\newcommand{\Dist}{\mathrm{Dist}}

\newcommand{\Ar}{{\mathrm{Ar}}}

\newcommand{\Unst}{{\mathrm{Un}}}
\newcommand{\coUnst}{{\mathrm{Un}^\mathrm{co}}}

\newcommand{\norm}{\mathrm{Nm}}

\newcommand{\lc}{\mathrm{lc}}

\let\lim\relax
\DeclareMathOperator*{\lim}{\mathsf{lim}}

\let\colim\relax
\DeclareMathOperator*{\colim}{\mathsf{colim}}

\DeclareMathOperator*{\llaxlim}{\mathsf{l.\!laxlim}}
\DeclareMathOperator*{\rlaxlim}{\mathsf{r.\!laxlim}}
\DeclareMathOperator*{\pllaxlim}{\mathsf{p.\!l.\!laxlim}}

\newcommand{\yoneda}{\text{\usefont{U}{min}{m}{n}\symbol{'107}}}
\DeclareFontFamily{U}{min}{}
\DeclareFontShape{U}{min}{m}{n}{<-> dmjhira}{}

\mathchardef\mdef="2D

\DeclareMathOperator{\const}{const}

\usepackage{tikz-cd}
\bibliography{ref}

\newcommand{\gitq}{/\!/}

\newcommand{\ab}{{\mathrm{ab}}}
\newcommand{\Mell}{\mathcal M_{\mathrm{Ell}}}

\newcommand{\ori}{\mathrm{or}}

\newcommand{\LL}{{\mathrm{L}}}
\newcommand{\RR}{{\mathrm{R}}}

\newcommand{\ptemp}{{\mathrm{ptemp}}}
\newcommand{\temp}{{\mathrm{temp}}}

\newcommand{\rep}{{{\mathrm{rep}}}}

\newcommand{\TopStk}{{\mathrm{TStk}}}

\newcommand{\TopPStk}{{\mathrm{TPStk}}}

\newcommand{\gen}{{\mathrm{gen}}}
\newcommand{\ngen}{{\mathrm{ngen}}}

\newcommand{\tors}{\mathrm{tors}}
\newcommand{\nnull}{\mathrm{null}}
\newcommand{\llocal}{\mathrm{loc}}

\newcommand{\complete}{\mathrm{cpl}}
\newcommand{\SpStk}{\mathrm{SpStk}}

\newcommand{\PreOr}{{\mathrm{Pre}}}
\newcommand{\qcs}{{\mathrm{qcs}}}

\numberwithin{thm}{subsection}

\begin{document}
\title{Local coefficients for genuine equivariant cohomology I}
\author{Nikolai Konovalov, Artem Prikhodko}
\date{}

\begin{abstract}
%We introduce the category of genuine equivariant local systems which categorify genuine equivariant cohomology, including Lurie’s tempered cohomology and elliptic cohomology of Grojnowski, Lurie, and Gepner--Meier, similarly to how the category of ordinary local systems categorifies singular cohomology. For a compact Lie group $G$ we identify genuine local systems on the classifying global space $BG$ with the genuine equivariant category of $G$-spectra. We extend Atiyah--Segal's completion theorem to tempered cohomology and use it to construct the category of tempered local systems. For a finite group $G$ this recovers Lurie's construction. We show that tempered category is a smashing localization of the genuine one. 

We develop a theory of equivariant local systems which categorifies genuine equivariant cohomology theories -- such as Atiyah--Segal $K$-theory, Lurie's tempered cohomology, and the equivariant elliptic cohomology of Grojnowski, Greenlees, and Gepner--Meier -- analogously to how the category of ordinary local systems categorifies singular cohomology.

More precisely, we introduce, for a global space $X$ and a coefficient system $\fcat A \colon \Orb^\op \to \Prs^\LL$, the categories $\LocSys^\Glo(X,\fcat A)$ and $\LocSys^\gen(X,\fcat A)$ of globally equivariant and genuine local systems, the latter generalizing the genuine stable category $\Sp^G$ of a compact Lie group $G$ to non-constant coefficients. When the coefficients come from an oriented abelian group stack $A$ over a locally complex periodic base, we also construct the category $\LocSys^\temp(X,A)$ of tempered local systems, extending Lurie's theory beyond finite groups; the defining condition is justified by a tempered form of the Atiyah--Segal completion theorem. We show that global sections induce an equivalence $\LocSys^\temp(BU(1),A) \simeq \QCoh(A)$ and we prove that $\LocSys^\temp(X,A)$ is a smashing localization of $\LocSys^\gen(X,A)$ if $X$ is an orbispace.
\end{abstract}

\maketitle

\setcounter{tocdepth}{2}
\tableofcontents

%\begingroup
%\let\clearpage\relax

\newcommand{\nc}{\mathrm{nc}}
\newcommand{\SpDM}{\mathrm{SpDM}}
\newcommand{\SpPStk}{\mathrm{SpPStk}}
\newcommand{\Spet}{\operatorname{\mathrm{Sp\acute et}}}
\newcommand{\fpqc}{\mathrm{fpqc}}

\section{Introduction}
Let $X\in \Type$ be a homotopy type and let $\fcat{C}$ be a category. The category of \emph{local systems} on $X$ with coefficients in $\fcat{C}$ is the functor category $\LocSys(X,\fcat{C}) := \Fun(X,\fcat{C})$. If $\fcat{C} = \Sp$ and $R \in \Sp$ is a spectrum, then the $R$-cohomology of $X$ is computed as the global sections
$$C^*(X,R) \simeq \lim_X \underline{R}$$
of the constant local system with value $R$. The category $\LocSys(X,\Sp)$ is, however, a far richer invariant than the cohomology it computes: it is presentably symmetric monoidal, it is functorial in $X$ with $*$-pullbacks and $*$-pushforwards satisfying base change and projection formulas, it contains the Thom local systems $\mathbbl{1}^E$ of vector bundles $E$ on $X$, and its global sections compute cohomology twisted by an arbitrary local system. In this sense $\LocSys(X,\Sp)$ is a \emph{categorification} of the cohomology of $X$, which governs Poincar\'e duality, Thom isomorphisms, and the rest of the standard properties of the cohomology groups. 

% In this sense $\LocSys(X,\Sp)$ \emph{categorifies} the cohomology of $X$, and it is this categorification, rather than the cohomology groups themselves, which supports Poincar\'e duality, Thom isomorphisms, and the rest of the standard toolkit.

Passing to the equivariant setting, one immediately observes that there are several candidates for a stable equivariant homotopy theory. Indeed, fix a (non-trivial) compact Lie group~$G$. By Elmendorf's theorem~\cite{Elmendor_SFPS}, the category $\Type^G$ of $G$-homotopy types is equivalent to the presheaf category $\Fun(\Orb_G^\op, \Type)$, and one can consider the \emph{naive} stabilization
$$\Sp^{nG} := \Fun(\Orb_G^{\op},\Sp),$$
which represents equivariant cohomology theories. However, the category $\Sp^{nG}$ lacks several crucial properties of the category $\Sp$. For example, the suspension spectrum $\Sigma^\infty_+ M$ of a closed $G$-manifold $M$ is \emph{not} a dualizable object in $\Sp^{nG}$ and the category is \emph{not} generated by dualizable objects. To resolve this issue, Lewis, May, and Steinberger~\cite{LMS86} introduced the \emph{genuine} equivariant stable category
$$\Sp^G := \Sp^{nG}[\{(\mathbb{S}^V)^{-1}\}_{V}]$$
obtained by $\otimes$-inverting the representation spheres. Only after this inversion is $\Sigma^\infty_+ M$ dualizable and $\Sp^G$ has equivariant Atiyah duality, transfers, and Wirthm\"uller's isomorphism. Moreover, the most interesting equivariant invariants such as Atiyah--Segal equivariant $K$-theory~\cite{Segal68} or equivariant elliptic cohomology~\cite{Greenlees_rational_elliptic},~\cite{GM},~\cite{Lur_Ell3} are  represented by objects of~$\Sp^G$. 

This suggests the question which the present paper addresses among other things: what is the category of local systems which categorifies a \emph{genuine} equivariant cohomology theory? We also point out here that there is another approach to cohomology groups using \emph{sheaves} instead of local systems and we will develop the corresponding equivariant theory in the forthcoming paper~\cite{localcoeff_2}.

%The category $\Sp^{nG}$ could be used to represent Borel equivariant cohomology theories (which are represented by the full subcategory $\Sp^{BG}:=\Fun(BG,\Sp)$), but the most interesting equivariant invariants like the Atiyah--Segal equivariant $K$-theory~\cite{Segal68} or the equivariant elliptic cohomology~\cite{Greenlees_rational_elliptic},~\cite{GM},~\cite{Lur_Ell3} are  \emph{not} of this kind. Namely, they are \emph{genuine}, i.e. they are represented by objects of the genuine equivariant stable category
%$$\Sp^G := \Sp^{nG}[\{(\mathbb{S}^V)^{-1}\}_{V}]$$
%of~, obtained by $\otimes$-inverting the representation spheres. Only after this inversion does one have transfers, the Wirthm\"uller isomorphism, and equivariant Atiyah duality. This suggests the question which the present paper addresses among other things: what is the category of local systems which categorifies a \emph{genuine} equivariant cohomology theory? We also point out here that there is another approach to cohomology groups using \emph{sheaves} instead of local systems and we will develop the corresponding equivariant theory in the forthcoming paper~\cite{localcoeff_2}.

For a finite group, and for complex periodic coefficients given by an oriented $p$-divisible group, an answer was given by Lurie in the theory of \emdef{tempered local systems}~\cite[Section~5]{Lur_Ell3}. The aim is then to construct such a theory in a setting which (i) treats all compact Lie groups on the same footing -- in particular the positive-dimensional ones; (ii) is functorial in the equivariant space in a way which supports the usual formalism of pullbacks, pushforwards, and duality; and (iii) works for not necessarily complex oriented coefficients. In this paper, we fulfill the first two goals and partially the third one.

\iffalse
The aim is then to construct such a theory in a setting which (i) treats all compact Lie groups on the same footing -- in particular the positive-dimensional ones; %, which are indispensable for elliptic phenomena, cf.~\cite{Greenlees_rational_elliptic, GM}; 
(ii) works for not necessarily complex oriented coefficients; and (iii) is functorial in the equivariant space in a way which supports the usual formalism of pullbacks, pushforwards, and duality. In this paper, we answer the first and the third questions.
\fi

Before we demonstrate our main results, let us briefly explain inputs for our constructions, which are \emph{global spaces and coefficient systems}.
In order to treat all compact Lie groups simultaneously, we work over the category of \emdef{global spaces}
$$\Type^{\Glo} := \Fun(\Orb^{\op},\Type)$$
of Gepner--Henriques and Rezk~\cite{GepnerHenriques07},~\cite{Rezk_GlobalHomotopy}. Here the global orbit category $\Orb$ has as objects the classifying stacks $BG$ of compact Lie groups $G$, see \cref{definition: orbit category}. The category $\Type^\Glo$ contains $BG$ for every compact Lie group $G$, as well as the (non-full) subcategory of \emdef{orbispaces} $\Type^{\Orb}$, which is related to $G$-equivariant homotopy theory by the equivalence $\Type^{\Orb}_{/BG}\simeq \Type^G$ of \cref{proposition: G-spaces and orbispace}; under this equivalence a $G$-space $X$ corresponds to its global quotient $X\gitq G$. Two classes of morphisms in $\Orb$ play distinct and complementary roles throughout: the \emdef{full} morphisms, which induce surjections on isotropy groups, and the \emdef{faithful} ones, which induce injections, see \cref{section: faithful morphisms}.

Recall that a generalized cohomology theory is completely determined by its values on a point by the Eilenberg--Steenrod theorem. In particular, given an abelian group $A$, one constructs the cohomology theory $H^*(-,A)$. A feature of equivariant homotopy theory is that $\Type^G$ has several basic building blocks indexed by $G$-orbits. This suggests that instead of a single abelian group~$A$, one has to consider a functor
$$A\colon \Orb_G^\op \tto \Ab.$$
Such functors are called \emph{coefficient systems} and the resulting equivariant cohomology theories were studied by Bredon~\cite{Bredon67}.

That being said, we consider a functor
$$\fcat A \colon \Orb^\op \tto \Prs^\LL$$
as an input for our categorification, that is, a coefficient category $\fcat{A}(BG)$ for each compact Lie group $G$. We continue to call such functors \emdef{coefficient systems} (see~\cref{definition: coefficient system}) and one may view $\fcat A$ as a global $2$-space in a certain sense, cf~\cite{GLP26}. In order to get a reasonable theory we usually impose certain assumptions on $\fcat A$ as in \cref{assumpt_limit_prerserving_cs}. Two main examples are the following
\begin{itemize}
\item A constant coefficient system with values in a category $\fcat C$, to be denoted $\underline{\fcat{C}}$, see \cref{example: constant system}.

\item A coefficient system $\fcat A$ coming from a preoriented abelian group stack $A$ over a spectral base $S$, with
$$\fcat A(BG) := \QCoh(A[\widehat G]),$$
where $\widehat G$ denotes the Pontryagin dual of $G$, see \cref{example: coefficient system from preoriented abelian group}.
\end{itemize}

%The equivariant setting also affords more freedom in the choice of coefficients. Rather than fixing a single coefficient category, we allow a \emdef{coefficient system}
%$$\fcat{A}\colon \Orb^{\op}_{\fcat{T}} \tto \CAlg(\Prs^\LL),$$
%that is, a coefficient category $\fcat{A}(BG)$ for each orbit of a global family $\fcat T$ of compact Lie groups, functorially in $\Orb_{\fcat T}$; see \cref{definition: coefficient system}.
%The examples of interest to us are of geometric origin, and follow~\cite{GM}: a preoriented abelian group stack $A$ over a non-connective spectral stack $S$ determines the coefficient system
%$$\fcat{A}(BG) := \QCoh(A[\widehat{G}]), \qquad \widehat{G} = \Hom(G,\mathbb{T}),$$
%see \cref{example: coefficient system from preoriented abelian group}. 
Taking $A = \mathbb{G}_{a}$ over $\mathbb{Q}[\beta^{\pm 1}]$ (or more generally, any $1$-dimensional rational group scheme) produces the rational theories of~\cite{Greenlees_rational_elliptic}; taking $A = \mathbb{G}_{m, KU}$ produces equivariant $K$-theory; and taking $A$ to be an oriented spectral elliptic curve over the moduli stack $\Mell^{\ori}$ produces equivariant topological modular forms, see respectively \cref{example: additive group}, \cref{example: greenlees}, \cref{example: multiplicative group}, and \cref{example: oriented elliptic curve}. %In this geometric situation the coefficients are genuinely non-affine, and consequently the naive global sections functor loses information. We therefore work throughout with an \emdef{enhanced global sections} functor $\Gamma_{\fcat A}(X,-)$ taking values in $\fcat A(X) \simeq \QCoh(A(X))$, the quasi-coherent sheaves on the stack $A(X)$ of $A$-tempered cohomology of $X$ in the sense of~\cite[Section~6]{GM}; see \cref{section: global section}.

\subsection{Main constructions and results}\label{section: main results}
For a global space $X$ and a coefficient system $\fcat{A}$  we construct two categories, which are related via the adjunction
$$\Sigma^\infty_X\colon \xymatrix{\LocSys^\Glo(X,\fcat A) \ar@<0.5ex>[r] & \ar@<0.5ex>[l] \LocSys^\gen(X,\fcat A)} \colon \Omega^\infty_X,$$
and, under a certain orientation hypothesis, a third one $\LocSys^\temp(X,A)$.

The first one is the category of \emdef{globally equivariant local systems} $\LocSys^\Glo(X,\fcat A)$ from \cref{definition: global local systems}. This category is defined as a partial left lax limit of $\fcat A$ over the category $\Orb_{/X}$ of orbits over $X$, with the full morphisms marked. Concretely, such a local system assigns to each point $x\colon BG \to X$ an object $\mathcal{L}^{BG}\in \fcat A(BG)$, together with comparison maps which are required to be equivalences only along full morphisms, see \cref{remark: informal global locsys}. If $X=BG$ and the coefficient system $\fcat{A}=\underline{\Sp}$ is constant, then $\LocSys^\Glo(X,\fcat A)\simeq \Sp^{nG}$. So, the category $\LocSys^\Glo(X,\fcat A)$ may be thought of as a generalization of one of the main constructions from~\cite{LMS86} to the case of non-constant coefficients. Although, the basic functoriality (like $\#$-pushforwards along faithful morphisms or base change formulas) for the constant coefficient system follows directly from the pointwise formulas for Kan extensions, a non-constant coefficient systems require more careful abstract treatment as done in \cref{appendix: cart_lax_limits} and \cref{appendix: right adjoint lax limits}. %\todo{\cref{definition: global local systems} is followed by a remark unwinding the definition; a \texttt{\textbackslash label} on it would let us point the reader there. IGNORE this comment.} This category plays the role of $\Sp^{nG}$ and is, for us, an auxiliary construction; but it is at this level that we establish the basic functoriality --- $\#$-pushforwards along faithful morphisms, base change, projection and K\"unneth formulas --- which is then inherited by the other two.

The second category is the category $\LocSys^\gen(X,\fcat{A})$ of \emdef{genuine local systems}. This is obtained by $\otimes$-inverting Thom local systems. Namely, on an orbit $BG$ we set
$$\LocSys^{\gen}(BG,\fcat{A}) := \LocSys^\Glo(BG,\fcat A)[\{(\mathbbl 1^V)^{-1}\}_{V\in \Vect_{BG}}]$$
in evident analogy with $\Sp^G$, and then we \emph{right Kan extend} to all global spaces, see \cref{def_genuine_LS} and \cref{construction: genuine local system}. The right Kan extension may seem ad hoc as one is tempted to simply $\otimes$-invert all Thom local systems over $X$. However, a general global space (or even an orbispace) carries too few vector bundles for the naive $\otimes$-inversion to be of any use, see \cref{ex_LS_ngen_for_N}. Fortunately, for global quotients of $G$-spaces these two constructions agree, see \cref{prop_LS_gen_vs_ngen_global_quotients}.

We also note that if $\fcat{A}$ is induced by a preoriented abelian group stack $A$ over $S$ and both $A$ and $S$ are affine, then the theory of globally equivariant (resp. genuine) local systems $\LocSys^{\Glo}(-,\fcat{A})$ (resp. $\LocSys^{\gen}(-,\fcat{A})$) can be reduced to the case of a constant coefficient systems, as in \cref{example: locsysglo geometric type}. However, non-constant coefficients are crucial for equivariant cohomology theories arising from a \emph{non-affine} abelian group stack, e.g. from the spectral elliptic curve.

Our first result collects the basic properties of genuine local systems, and in particular says that our construction is a correct generalization of $\Sp^G$.

\begin{thmx}\label{theorem: intro genuine}
Let $\fcat{A}$ be a strongly continuous symmetric monoidal coefficient system of Beck--Chevalley type (see \cref{definition: coefficient system} and \cref{assumpt_limit_prerserving_cs}). Then the assignment $$X \mapsto \LocSys^{\gen}(X,\fcat A)$$ extends to a functor
$$\LocSys^{\gen}(-,\fcat A)\colon (\Type^{\Glo})^{\op} \tto \CAlg(\Prs^{\LL})$$
which preserves small limits (i.e. satisfies descent) and has the following properties.
\begin{enumerate}
\item If $\fcat{A} = \underline{\Sp}$ is constant, then there is a symmetric monoidal equivalence $\LocSys^{\gen}(BG,\underline{\Sp}) \simeq \Sp^G$ for every compact Lie group $G$.
\item If $f\colon X\to Y$ is a faithful morphism of global spaces, then $f^{\gen,*}$ admits a left adjoint $f^{\gen}_{\#}$, and the base change and projection formulas hold, see \cref{proposition: genuine sharp} and \cref{corollary: projection formula for genuine sharp}.
\item If $\fcat{A}$ is of geometric type, then $\LocSys^{\gen}(BG,\fcat A)$ is a rigid, $\fcat{A}(BG)$-atomically generated $\fcat A(BG)$-algebra, with $\fcat{A}(BG)$-atomic generators $p^{\gen}_{\#}\mathbbl{1}_{BH}$ indexed by faithful morphisms $p\colon BH \to BG$, see \cref{corollary: genuine BG rigid over ABG}.
\item There is a family of symmetric monoidal, colimit preserving and jointly conservative geometric fixed point functors 
$$\Phi^{x}_{\fcat A}\colon \LocSys^{\gen}(X,\fcat A) \tto \fcat{A}(BH),$$
indexed by the points $x \colon BH \to X$, and for $X = BG$ these assemble into a recollement generalizing the classical isotropy separation of $\Sp^G$, see \cref{construction: geometric fixed points}, \cref{corollary: gfp are jointly conservative global}, and \cref{corollary: isotropy separation for gen}.
\end{enumerate}
\end{thmx}

Genuine local systems, however, do not yet solve the problem we started from. If the coefficient system is of geometric type (\cref{definition: coefficient system is of geometric type}), then
$$\LocSys^{\Glo}(BG,\fcat{A}) \simeq \Mod_{\fcat{A}_G}\left(\Sp^{nG}\otimes \fcat{A}(BG)\right) \;\; \text{and} \;\; \LocSys^{\gen}(BG,\fcat{A}) \simeq \Mod_{\Sigma^\infty_{BG}\fcat{A}_G}\left(\Sp^{G}\otimes \fcat{A}(BG)\right),$$
see \cref{example: geom type gen BG}. In other words, the genuine equivariant cohomology theories arising this way are exactly genuine suspensions of naive ones. In particular, for $A = \mathbb{G}_{m,KU}$ the category $\LocSys^\gen(BG,\fcat A)$ does \emph{not} recover genuine equivariant $K$-theory. 

To resolve this issue we assume that the coefficient system $\fcat{A}$ is induced by an \emph{oriented} abelian group stack $A$ over $S$, see \cref{definition: complex oriented abelian group object}. Then, following~\cite{Lur_Ell3}, we introduce an additional category $\LocSys^\temp(X,A)$ of  \emdef{tempered local systems}. Namely, $\LocSys^\temp(X,A)$ is the full subcategory of $\LocSys^{\Glo}(X,A)$ consisting of those $\mathcal{L}$ for which, for every faithful morphism $i\colon BH \to BG$ over $X$ between compact abelian Lie groups, the comparison map
$$(\mathcal{L}^{BG})^{\wedge}_{(H,G)} \tto (i_*\mathcal{L}^{BH})^{hG/H}$$
is an equivalence, see \cref{def_LS_temp}. %The condition is an axiomatization of the Atiyah--Segal completion theorem, and it is satisfied by the unit precisely because $A$ is assumed oriented, see \cref{lem_TempAS_structure_sheaf}. 
The condition cutting out tempered local systems is justified by the following completion theorem, which is the technical heart of \cref{section: complex periodic coefficients} and a tempered form of the Atiyah--Segal completion theorem~\cite[Theorem~2.1]{AtiyahSegal_Completion}.

\begin{thmx}[\cref{lem_TempAS_structure_sheaf}]\label{theorem: intro AS}
Let $A$ be an oriented non-connective spectral abelian group stack over a locally complex periodic base $S$, and let $i\colon BH \inj BG$ be a faithful morphism of compact abelian Lie groups. Then the canonical map
$$\theta_{H,G}\colon \mathcal{O}_{A[\widehat{G}]} \tto \left(i_*\mathcal{O}_{A[\widehat{H}]}\right)^{hG/H}$$
exhibits the target as the $(H,G)$-completion of the source, see \cref{nota_tors_loc_comp_BH_BG}.
\end{thmx}

Since we work over a base which is only locally complex periodic and over abelian group objects which need not be affine, the proof cannot follow~\cite[Section~4.6]{Lur_Ell3} verbatim. We prove \cref{theorem: intro AS} by reducing to the case of the inclusion of the trivial group into $U(1)$, where the statement holds essentially by the definition of an orientation, see \cref{LStemp_TateInvs_are_local} and \cref{lem_TempAS_structure_sheaf}.

Our next result collects the basic properties of tempered local systems similar to \cref{theorem: intro genuine}.

\begin{thmx}\label{theorem: intro tempered}
Let $A$ be an oriented non-connective spectral abelian group stack over a locally complex periodic base $S$. Then the full subcategory $\LocSys^{\temp}(X,A) \subseteq \LocSys^{\Glo}(X,A)$, $X \in \Type^\Glo$ has the following properties.
\begin{enumerate}
\item It is preserved by $*$-pullbacks, and the functor $\LocSys^{\temp}(-,A)\colon \Type^{\Glo,\op}\to \Prs^\LL$ preserves small limits, see \cref{pb_preserves_LStemp} and \cref{proposition: temp descent}.
\item The inclusion admits a left adjoint $L_X$ whose kernel is the $\otimes$-ideal of null local systems.  Consequently $\LocSys^{\temp}(X,A)$ carries a unique symmetric monoidal structure making $L_X$ symmetric monoidal, and $\LocSys^{\temp}(-,A)$ lifts to $\CAlg(\Prs^\LL)$, see \cref{LStemp_admits_left_adj_global} and \cref{thm_CMon_str_on_LStemp}.
\item If $T$ is a connected compact abelian Lie group, then the global sections functor
$$\Gamma_A(BT,-)\colon \LocSys^{\temp}(BT,A) \tto \QCoh(A(BT))$$
is an equivalence, see \cref{LStemp_on_tori}. Thus for $T \cong U(1)^{\times r}$ one has $\LocSys^\temp(BT,A)\simeq \QCoh(A^{\times r})$.
\end{enumerate}
\end{thmx}

Part (3) is the most important here as it shows that the categories $\LocSys^{\temp}(-,A)$ are categorifications we are looking for. This assertion implies that the $U(1)$-tempered cohomology of a point is the structure sheaf $\mathcal{O}_A$ as one expects. We point out that connectedness is essential -- for $G$ non-connected the functor $\Gamma_A(BG,-)$ need not be conservative, see \cref{remark: non-conservative in general}.

Our next theorem asserts that, if $X$ is an orbispace, then the passage from genuine to tempered local systems is rather mild. 

\begin{thmx}[\cref{theorem: temp are modules in gen} and \cref{corollary: temp in genuine}]\label{theorem: intro comparison} Let $A$ be an oriented non-connective spectral abelian group stack over a locally complex periodic base $S$. Let $X$ be an orbispace. Then the adjunction
$$L^{\gen}_X \colon \xymatrix{\LocSys^{\gen}(X,A) \ar@<0.5ex>[r] & \ar@<0.5ex>[l] \LocSys^\temp(X,A)}  \colon \beta_X$$
exhibits $\LocSys^\temp(X,A)$ as a smashing localization of $\LocSys^\gen(X,A)$. More precisely, $\beta_X$ is fully faithful, the commutative algebra $\beta_X(\mathbbl{1})$ is idempotent, and there is a natural equivalence
$$\LocSys^{\temp}(X,A) \simeq \Mod_{\beta_X(\mathbbl 1)}\LocSys^{\gen}(X,A).$$
Moreover, a genuine local system $\mathcal{L}$ lies in the essential image of $\beta_X$ if and only if $\Phi^{x}\mathcal{L}$ is $(\fcat{P},G)$-local for every faithful point $x\colon BG \to X$ with $G$ abelian, see \cref{definition: proper localization}.
\end{thmx}

In other words, $\LocSys^\temp(X,A)$ is uniquely determined by a single idempotent algebra $\beta_X(\mathbbl{1})\in\LocSys^\gen(X,A)$. %, see \cref{theorem: temp are modules in gen}. 
Therefore, in the case of an orbispace, the three constructed categories are related as follows

\[\xymatrix@C=3.0cm@R=3.2pc{
\LocSys^{\Glo}(X, A) \ar@<0.7ex>[r]^-{\Sigma^\infty_X} \ar@/^2.2pc/[rr]^-{L_X} & \LocSys^{\gen}(X, A) \ar@<0.7ex>[r]^-{L^\gen_X\simeq (-)\otimes \beta_X(\mathbbl{1})} \ar@<0.7ex>[l]^-{\Omega^\infty_X} & \LocSys^{\temp}(X, A).\ar@{^{(}->}@<0.7ex>[l]^-{\beta_X} \ar@{^{(}->}@/^2.2pc/[ll]^-{\iota_X}
}\]

%Throughout, $\fcat{A}$ denotes a strongly continuous symmetric monoidal coefficient system of Beck--Chevalley type (see \cref{definition: coefficient system} and \cref{assumpt_limit_prerserving_cs}), and $A$ denotes an abelian group object in non-connective spectral stacks over a base $S$.

%Part~(3) is the precise sense in which our tempered local systems are a geometric theory: on the classifying space of a torus they are simply quasi-coherent sheaves on a power of $A$. 
%The comparison between the two theories is our main result.

The hypothesis that $X$ be an orbispace cannot be dropped: for $X=\Sigma BC_p$ -- the suspension of $BC_p$ inside $\Type^\Glo$ -- the functor $\beta_X$ is neither fully faithful nor continuous, see \cref{example: temp is not local}. 
\iffalse
\todo{Do we want this statement in intro?} Specializing \cref{theorem: intro comparison} to a point of the moduli of coefficients, we obtain the statement which motivated the whole discussion.

\begin{thmx}\label{theorem: intro corollary}
Let $G$ be a compact abelian Lie group and let $A$ be nc-affine. Then there is a genuine equivariant $E_\infty$-ring $R_G \in \CAlg(\Sp^G)$ and identifications
$$\LocSys^{\temp}(BG,A) \simeq \Mod_{R_G}(\Sp^G), \qquad \LocSys^{\gen}_{\ab}(BG,A)\simeq \Mod_{\Sigma^{\infty}_{BG}\Omega^{\infty}_{BG}R_G}(\Sp^G),$$
under which $R_G$ is an idempotent algebra over $\Sigma^{\infty}_{BG}\Omega^{\infty}_{BG}R_G$. Moreover, $\LocSys^{\temp}(BG,A)$ is a rigid, $\QCoh(A[\widehat G])$-atomically generated $\QCoh(A[\widehat G])$-algebra, with atomic generators $p^{\temp}_{\#}\mathbbl 1_{BH}$ indexed by faithful morphisms $p\colon BH \to BG$; see \cref{cor_rigidity_of_LStemp}.
\end{thmx}

In other words: although the category of genuine local systems does not by itself recover genuine equivariant cohomology theories such as equivariant $K$-theory or equivariant elliptic cohomology, a smashing localization of it does, at least in the complex periodic case. \todo{The introduction of \cref{section: comparison of tempered and genuine local systems} states \cref{theorem: intro corollary} for an arbitrary compact Lie group $G$, but \cref{example: geom type gen BG} seems to require $G \in \fcat{T} = \Orb_\ab$. We have stated it for $G$ abelian --- please confirm, or tell us how the non-abelian case follows.}
\fi

\subsection{Relation to other work}\label{section: relation to other work} Our definitions are inspired by many previous constructions and may be considered as their generalizations. Here we give a brief comparison with selected papers; a detailed account is given in the main body of the paper.

The area of equivariant homotopy theory where all groups are studied simultaneously is usually called \emph{global homotopy theory}, a term introduced by Schwede in~\cite{Schwede}. He constructed the category of global spectra $\Sp^{\mathrm{gl}}$ and showed that many natural equivariant cohomology theories, such as Atiyah--Segal $K$-theory, are represented by global spectra. As was shown in~\cite{LNP25}, the category $\Sp^{\mathrm{gl}}$ is a partial lax limit over $\Orb$ of genuine $G$-equivariant homotopy theories $\Sp^G$. A similar strict limit over $\Orb^\rep$ is the category $\LocSys^\gen(\mathcal{N},\underline{\Sp})$ for the \emph{normal subgroup classifier} $\mathcal{N}\in\Type^\Glo$, see \cref{example: normal subgroup classifier} and \cref{ex_LS_ngen_for_N}.

The theory of genuine local systems has recently attracted attention from the point of view of \emph{parametrized homotopy theory} (see e.g.~\cite{para_exposeI} and~\cite{Shah23}) and \emph{twisted ambidexterity} (see e.g.~\cite{Cnossen23}). In fact,~\cite{Cnossen23} already contains a construction for the category $\LocSys^\gen(X,\underline{\fcat{C}})$ with \emph{constant coefficients} and an analog of \cref{theorem: intro genuine}. We extend Cnossen's results to the case of \emph{non-constant} coefficients; this is not a formal exercise, as it requires a finer understanding of lax limits and their functoriality. 

Our construction of the category $\LocSys^\Glo(X,\fcat{A})$ of globally equivariant local systems is inspired by the definition of \emph{pretempered local systems} $\LocSys^\ptemp(X, A)$ from~\cite[Construction~5.1.3]{Lur_Ell3}. In fact, $\LocSys^\ptemp(X, A)$ is a special case of $\LocSys^\Glo(X,\fcat{A})$ for a suitable coefficient system~$\fcat{A}$. For this reason, some arguments in \cref{section: globally equivariant local systems} and in \cref{section: complex periodic coefficients} mimic the arguments from~\cite{Lur_Ell3}. We also note that we obtain the category $\LocSys^{\gen}(X,\fcat{A})$ from $\LocSys^\Glo(X,\fcat{A})$ by using the formal $\otimes$-inversion of~\cite{Robalo15} or, more precisely, its parametrized refinement from~\cite{Cnossen23}.

Of course, the theory of tempered local systems originates in Lurie's work on tempered cohomology~\cite{Lur_Ell3}, where it is developed for finite groups, taking as input an oriented $p$-divisible group. Our \cref{def_LS_temp} is modelled on~\cite[Section~5]{Lur_Ell3}, the difference being that we require as input an honest abelian group object $A$ in spectral stacks, which is what allows us to treat compact Lie groups of positive dimension. For finite abelian groups the two notions agree by definition.%, see \cref{temp_LS_on_BCp} and \cref{remark: tempered local system comp gen finite}. %\todo{Do we want to include a general comparison statement with~\cite{Lur_Ell3} for finite groups, or is the case of $BC_p$ together with the remarks enough? The abstract currently claims the recovery of Lurie's category for an arbitrary finite group.}

The geometric input for our coefficient systems, together with the tempered cohomology stack $A(X)$, is taken from Gepner--Meier~\cite{GM}, where equivariant topological modular forms are constructed. Our description of $\LocSys^\temp(BT,A)$ for a torus $T$ should be compared with the rational $S^1$-equivariant elliptic cohomology of Greenlees~\cite{Greenlees_rational_elliptic} (see also~\cite{Grojnowski_delocalized}), which is recovered by \cref{example: greenlees is complex oriented}. %The global homotopy theory we use is that of Gepner--Henriques~\cite{GepnerHenriques07} and Rezk~\cite{Rezk_GlobalHomotopy}; the partially lax limit which appears in \cref{definition: global local systems} is closely related to the description of global spectra given by Linskens, Nardin and Pol~\cite{LNP25}. Finally, the formal $\otimes$-inversion we use follows Robalo~\cite{Robalo15} and its parametrized refinement follows Cnossen~\cite{Cnossen23}, and the rigidity statements in \cref{theorem: intro genuine} and \cref{cor_rigidity_of_LStemp} rest on the theory of (locally) rigid categories of Ramzi~\cite{Ramzi24_dualizable, Ramzi26}.

The connection between genuine local systems and tempered ones has already been discussed in~\cite{Lur_Ell3} without a formal proof. For finite abelian groups, \cref{theorem: intro comparison} was obtained by Gepner, Linskens and Pol in~\cite{GLP26}, in the language of global $2$-rings and their genuine refinements. Their result is in one respect more general than ours, in that the coefficients are allowed to come from an oriented divisible group which need not be induced by an abelian group object. To the best of our knowledge, neither the construction of tempered local systems nor their comparison with genuine ones for positive-dimensional compact Lie groups has been studied before. The coefficient systems are closely related to their global $2$-rings and we compare their pregenuine and rigid global $2$-rings with our assumptions on coefficient systems in \cref{remark: geometric and pregenuine}.

Importantly, we do \emph{not} cover any ambidexterity results in this paper -- i.e.\@ the comparison of $\#$- and $*$-pushforwards along non-faithful morphisms. For the case of finite groups, the interested reader may look at e.g.~\cite{Lur_Ell3},~\cite{BDL26}, or~\cite{BenMoshe24}. It seems that new phenomena occur for positive-dimensional compact Lie groups, see e.g. \cref{remark: ambidexterity for elliptic cohomology}. A detailed account will appear elsewhere. 

\subsection{Future directions}\label{section: future work}
%This work arose as an offshoot of another project, where for a topological stack $X$ we construct categories of genuine equivariant and tempered sheaves and establish $6$-functor formalism of these constructions. Inside these categories there are natural subcategories of locally constant objects, which, for a locally contractible in an appropriate sense $X$, are equivalent to the purely homotopy theoretic constructions of this work. The theories of local systems and sheaves are parallel to some extent, but at some point it became too cumbersome to keep developing them under the same roof, so we decided to separate it into different projects. We also used this opportunity to distribute technical stuff concerning lax limits between two manuscripts. The theory of genuine and tempered sheaves will be the content of the sequel of this paper.

This work arose as an offshoot of another project, where for a topological stack $X$ we construct the categories $\Shv^\gen(X,\fcat{A})$  and $\Shv^\temp(X,A)$ of genuine equivariant and tempered sheaves respectively. The main goal is to establish a $6$-functor formalism for $\Shv^\gen(X,\fcat{A})$ and $\Shv^\temp(X,A)$. By restricting to the subcategories $\Shv^{\gen,\lc}(X,\fcat{A})$ and $\Shv^{\temp,\lc}(X,A)$ of locally constant objects, we show the equivalences
$$\Shv^{\gen,\lc}(X,\fcat{A}) \simeq \LocSys^\gen(\Pi^\Glo(X),\fcat{A}) \;\; \text{and} \;\; \Shv^{\temp,\lc}(X,A) \simeq \LocSys^\temp(\Pi^\Glo(X), A), $$
where~$X$ is a sufficiently nice topological stack and $\Pi^\Glo(X) \in \Type^\Glo$ is the underlying global space of~$X$, see \cref{example: topological stacks}. In other words, we compare the ``topological'' category of equivariant sheaves with the ``homotopical'' category of equivariant local systems. The theory of local systems runs parallel to the theory of sheaves to some extent, but it is also a precursor to it, which is our main reason for treating it in a standalone manuscript. The theory of genuine and tempered sheaves will be the content of the sequel~\cite{localcoeff_2} to the present paper.

Let $A$ be an affine oriented abelian group object, e.g.\@ $A = \mathbb G_{m,KU}$. Then there is an associated genuine $U(1)$-equivariant $E_\infty$-ring spectrum $\fcat{A}_{U(1)} \in \CAlg(\Sp^{U(1)})$ such that
$$\LocSys^\gen(BU(1), A) \simeq \Mod_{\Sigma^\infty_{BU(1)} \Omega^\infty_{BU(1)} \fcat{A}_{U(1)}} \Sp^{U(1)}.$$
It follows that $\LocSys^\gen(BU(1), A)$ a priori loses some information about $\fcat{A}_{U(1)}$ and its modules in $\Sp^{U(1)}$. In this case the situation is rescued by the category of tempered local systems, which is equivalent to $\Mod_{\fcat{A}_{U(1)}} \Sp^{U(1)}$. One may then wonder whether there is a variant $\LocSys^{\gen,\prime}(X,\fcat{A})$ of genuine local systems such that $\LocSys^{\gen,\prime}(BU(1),\fcat{A})$ identifies with $\Mod_{\fcat{A}_{U(1)}} \Sp^{U(1)}$ for an arbitrary genuine equivariant ring $\fcat{A}_{U(1)}$ (and a suitable coefficient system $\fcat{A}$)? In particular, such a construction $\LocSys^{\gen,\prime}(X,A)$ should be equivalent to tempered local systems $\LocSys^{\temp}(X,A)$ without any localizations, if one starts with an oriented abelian group object $A$. The main difficulty is that one has to lift a coefficient system
$$\fcat A\colon \Orb^\op \tto \CAlg(\Prs^\LL)$$
to some sort of ``genuine equivariant commutative $2$-ring''. At the time of writing, we do not know what the precise definition of a categorification of a genuine ring should be. %We plan to return to this question in a future work.

It may be the case that the question posed in the previous paragraph is too vague. However, the following special case looks rather promising. Let $A$ be an oriented abelian group object over $\Spec(R)$ and let $A_{\geq 0}$ be its connective cover over $\Spec(R_{\geq 0})$. Then $A_{\geq 0}$ is no longer oriented, but it is still preoriented. So, we consider the coefficient system 
$$\fcat{A}_{\geq 0} \colon \Orb^\op \tto \CAlg(\Prs^\LL),\;\; BG \mapsto \QCoh(A_{\geq 0}[\widehat{G}]) $$
as before. We expect that there is a variant of genuine local systems $\LocSys^{\gen,\prime}(X,\fcat{A}_{\geq 0})$ such that 
$$\LocSys^{\gen,\prime}(BG,\fcat{A}_{\geq 0}) \simeq \Mod_{\fcat{A}_G} \LocSys^{\gen}(BG,\fcat{A}_{\geq 0})$$
for some $\fcat{A}_{G} \in \CAlg \LocSys^{\gen}(BG,\fcat{A}_{\geq 0})$ and
$$\LocSys^{\gen,\prime}(BU(1),\fcat{A}_{\geq 0}) \simeq \QCoh(\mathrm{Bl}_Z A_{\geq 0}),$$
where $\mathrm{Bl}_Z A_{\geq 0}$ is a \emph{blow-up} of $A_{\geq 0}$ at a certain closed spectral subscheme $Z \subset A_{\geq 0}$ replacing the condition that the geometric fixed points be local. If $A=\mathbb{G}_{m,KU}$, then the ring $\fcat{A}_{G}$ seems to be closely related to \emph{equivariant connective $K$-theory} of~\cite{Greenlees04}; if $A$ is a spectral elliptic curve, this opens a route to the definition of \emph{equivariant connective elliptic cohomology}.

\subsection{Outline of the paper}\label{section: outline}
Here we give a brief outline of the manuscript; more detail is given at the beginning of each section. \Cref{section: global spaces} is a review of global homotopy theory: the global orbit category $\Orb$, the classes of full and faithful morphisms, the category of orbispaces and its relation to $G$-spaces, and global families of compact Lie groups. We do not claim any original results here; we include this material to fix notation, which may differ slightly from one work to another.

\Cref{section: globally equivariant local systems} introduces coefficient systems and the category $\LocSys^{\Glo}_{\fcat T}(X,\fcat A)$ of globally equivariant local systems. After discussing the coefficient systems attached to preoriented abelian group stacks in \cref{ssect_coefs_from_PreAbStk}, we show that $\LocSys^{\Glo}_{\fcat T}(-,\fcat A)$ is right Kan extended from orbits, that it is much better behaved on orbispaces than on general global spaces (\cref{section: local systems on orbispaces}), and that $*$-pullbacks along faithful morphisms admit left adjoints satisfying base change and the projection formula (\cref{subsection: functoriality along faithful morphisms}). \Cref{section: global section} constructs the enhanced global sections functor and an isotropy separation decomposition of $\LocSys^{\Glo}_{\fcat T}(BG,\fcat A)$.

\Cref{section: genuine equivariant local system} develops genuine local systems. We introduce Thom local systems in \cref{section: thom local systems}, define $\LocSys^{\gen}_{\fcat T}(-,\fcat A)$ in \cref{section: geuine local systems on global spaces} and establish its functoriality, compare it with the naive $\otimes$-inversion for global quotients and prove the rigidity statement in \cref{section: LSgen on G-spaces}, and construct the geometric fixed point functors in \cref{section: geometric fixed points}.

\Cref{section: complex periodic coefficients} turns to complex periodic coefficients and tempered local systems. After the necessary preliminaries on formal completions in spectral algebraic geometry (\cref{section: formal completions}), complex orientability (\cref{section: complex oriented coefficient systems}) and oriented abelian group objects (\cref{section: oriented abelian group objects}), we prove \cref{theorem: intro AS} in \cref{section: as comparison}, define tempered local systems in \cref{section: tempered locsys}, construct their monoidal structure in \cref{section: monoidal structure}, and compute the tempered local systems on classifying stacks of tori in \cref{section: tempered local of tori}.

\Cref{section: comparison of tempered and genuine local systems} contains the proof of \cref{theorem: intro comparison}. We first show in \cref{section: thom is tensor-invertible} that tempered Thom local systems are $\otimes$-invertible, so that $L_X$ factors through $L^{\gen}_X$, and then in \cref{section: smashing localization} we compute the geometric fixed points of the objects in the essential image of the right adjoint $\beta_X$ and deduce that $L^{\gen}_X$ is smashing. In \cref{section: dexterity}, we study the category $\LocSys^{\temp}(X,A)$ as a module in $\Prs^\LL$ over the symmetric monoidal category $\LocSys^{\Glo}_\ab(X,A)\in \CAlg(\Prs^\LL)$.

Four appendices collect material of an independent nature: \cref{appendix: sag} is a reminder on (non-connective) spectral algebraic geometry; \cref{appendix: cart_lax_limits} on Cartesian fibrations, lax limits and lax Kan extensions; \cref{appendix: right adjoint lax limits} on adjoint functors between right lax limits; and \cref{section: object inversion} on the $\otimes$-inversion of a collection of objects in a symmetric monoidal category. %\todo{ObjInv.tex has no \texttt{\textbackslash label} on its \texttt{\textbackslash section}; please add one (e.g.\@ \texttt{appendix: object inversion}) so that it can be cross-referenced here.}

\subsection{Acknowledgements} The entire project started in a collaboration with Ivan Perunov; both authors are deeply indebted for his insights and his help. The authors are also grateful to Chris Brav, Jack Morgan Davies, Sil Linskens, and Tomer Schlank for many fruitful conversations. We also thank Dmitry Kubrak, Kaif Hilman, and Sanath Devalapurkar for their interest to the paper and helpful discussions.

The first author is also grateful to the University of Chicago for excellent working condition; in particular, for providing access to Claude Opus 5. LLMs were used to prepare the introduction and improve the writing. All mathematical content as well as any errors are due to the authors.

The second author gratefully acknowledges support from the Basic Research Program of HSE University (HSE-BR-2025-060).

\section{Global spaces}\label{section: global spaces}
Let $G$ be a compact Lie group and let us denote  by $\Top^G$ the category of topological spaces equipped with continuous $G$-action and $G$-equivariant maps between them. One of the first observations of equivariant homotopy theory is that the naive localization of $\Top^G$ at a class of $G$-equivariant continuous maps which are non-equivariant weak equivalences loses too much information. E.g.\@ one can recover Borel equivariant (generalized) cohomology this way, but already Atiyah--Segal's equivariant $K$-theory does not factor through this localization. Instead one considers the localization of $\Top^G$ with respect to the finer class of \emph{weak $G$-equivariant homotopy equivalences}. Then, by Elmendorf's theorem \cite{Elmendor_SFPS}, the restricted Yoneda embedding
$$\Type^G:= \Top^G[\{\text{weak $G$-homotopy equivalences}\}^{-1}] \tto \Fun(\Orb_G^\op, \Type)$$
is an equivalence, where $\Orb_G$ denotes the full $\infty$-subcategory of the left hand side spanned by spaces of the form $G/H$ for all closed subgroups $H\subseteq G$.

Now, if one wants to study equivariant phenomena which are not specific to a fixed $G$, it is convenient to combine all $\Type^G$ into a single category. Such a category of \emdef{global homotopy types} or \emdef{global spaces $\Type^\Glo$} was introduced in \cite{GepnerHenriques07} and thoroughly studied in \cite{Rezk_GlobalHomotopy}. In this section we review the relevant constructions and results about $\Type^\Glo$.

\subsection{Category \texorpdfstring{$\Orb$}{Orb}}\label{section: category orb}
Similarly to the case of $\Type^G$, the category $\Type^\Glo$ is defined as presheaves on the \emdef{global orbit category $\Orb$} (\Cref{definition: global spaces}). It is convenient for us to give slightly different construction of $\Orb$ than in~\cite{GepnerHenriques07} based on classifying topological stacks.
\begin{defn}\label{defn: stacks}
For a regular cardinal $\kappa$ let us denote by
$$\TopPStk_\kappa := \Fun(\Top_\kappa^\op, \Type)$$
the category of presheaves of homotopy types on $\Top_\kappa$, the full subcategory of the category of topological spaces $\Top$ spanned by topological spaces with cardinality of the underlying set less than $\kappa$. The category $\TopStk_\kappa$ is defined as a full subcategory of $\TopPStk_\kappa$ spanned by sheaves with respect to the Grothendieck topology on $\Top_\kappa$ which covering families are disjoint unions of jointly surjective open embeddings.

For $\lambda \ge \kappa$ the inclusion $\Top_\kappa \inj \Top_\lambda$ induces fully faithful embeddings
$$\TopPStk_\kappa \inj \TopPStk_\lambda, \qquad \TopStk_\kappa \inj \TopStk_\lambda.$$
The category of \emdef{topological prestacks $\TopPStk$ (resp.\@ stacks $\TopStk$)} is defined as a union\footnote{In fact, for this work it is enough to consider $\TopStk_\kappa$ for any regular $\kappa > \mathfrak{c}$, where $\mathfrak{c}$ is the cardinal of the set of real numbers.}
$$\TopPStk := \indlim \TopPStk_\kappa, \qquad \TopStk := \indlim \TopStk_\kappa.$$
\end{defn}

\begin{construction}\label{construction_Pi}
Let $\Pi\colon \Top \to \Type$ denote the functor mapping a topological space $X$ to its underlying homotopy type $\Pi(X)$ (also called the fundamental $\infty$-groupoid of $X$). We will denote by
$$\Pi\colon \TopPStk \tto \Type$$
the left Kan extension of the functor above along the Yoneda embedding $\Top \inj \TopPStk$. Equivalently, this is the unique colimit preserving functor from $\TopPStk$ to $\Type$ which coincides with the fundamental groupoid functor when restricted to $\Top$. If $X$ is a topological prestack, then we will call $\Pi(X)$ the \emdef{underlying homotopy type} or the \emdef{fundamental groupoid of $X$}.
\end{construction}

\begin{rem}\label{Pi_preeserves_fprods}
Note that $\Pi\colon \TopPStk \tto \Type$ preserves finite products. Indeed, $\Pi\colon \Top \to \Type$ preserves finite products and finite products in $\Type$ preserve colimits in each argument.
\end{rem}

\begin{ex}\label{example: classifying spaces}
Let $G$ be a topological group and let $BG$ denote the topological classifying stack of $G$ (e.g.\@ defined as the delooping of $G$ in the topos $\TopPStk$). Then, since $\Pi\colon \TopPStk \tto \Type$ preserves colimits and finite products, $\Pi(BG)$ is naturally equivalent to the usual classifying space of the grouplike $E_1$-space $\Pi(G)$.
\end{ex}

Note that the category $\TopPStk$ is closed symmetric monoidal. In particular, the category $\TopPStk$ is enriched over itself, see \cite[Corollary~7.4.10]{GepnerHauseng_Enriched}. If $X,Y \in \TopPStk$, we will write $\Map(X,Y) \in \TopPStk$ for the internal $\Hom$ in the category $\TopPStk$.

\begin{defn}\label{definition: TopStk-enriched orbit category}
The \emdef{enriched orbit category} $\Orb_{\TopStk}$ is the enriched full subcategory of $\TopPStk$ spanned by the classifying stacks $BG$ for compact Lie groups.
\end{defn}

Since the functor $\Pi\colon \TopPStk \to \Type$ is symmetric monoidal, we can transfer the enrichment along the functor $\Pi$ in order to obtain the ($\infty$-)category, see e.g. \cite[Definition~C.3.1]{HM24}.

\begin{defn}\label{definition: orbit category}
The \emdef{orbit category} $\Orb$ is the ($\infty$-)category $\tau_{\Pi}(\Orb_{\TopStk})$ obtained from the enriched orbit category $\Orb_{\TopStk}$ by the transfer of enrichment along the symmetric monoidal functor $\Pi\colon \TopPStk \to \Type$.
\end{defn}

%\begin{defn}
%The \emdef{orbit category} $\Orb$ is the $\infty$-category associated to the topological category whose objects are compact Lie groups and whose mapping spaces are given by
%$$\Map_{\Orb}(H,G) = |\Map_{\Grp}(H,G)\gitq G|,$$
%the geometric realization of the action groupoid of $G$ acting on the space of continuous group homomorphisms $\Map_{\Grp}(H, G)$ by conjugation. Composition is induced by the composition of group homomorphisms.
%\end{defn}

\begin{rem}\label{remark: orbit category references}
Note that the definition of $\Orb$ agrees with the definition given in our main reference \cite[Section~2.2]{Rezk_GlobalHomotopy}, see also \cite[Definition~6.1]{LNP25}. Note as well that loc.\@ cit.\@ use the notation $\mathrm{Glo}$ for $\Orb$. Our definition of $\Orb$ also agrees with the definition given in \cite[Section 4]{GepnerHenriques07} restricted to compact Lie groups.
\end{rem}

%\begin{rem}\label{remark: orbit category gepner}
%We repeat \cite[Remark~6.2]{LNP25}. The definition of $\Orb$ also agrees with the definition given in Section 4 of \cite{GepnerHenriques07} restricted to compact Lie groups, up to one difference. \cite{LNP25} apply thin geometric realization to the action groupoids to obtain a topologically enriched category, while the original definition uses fat geometric realization. Up to a technical condition, the two conventions define Dwyer-Kan equivalent topological categories. See \cite[Remark~3.10]{korschgen} for a more detailed discussion.
%\end{rem}

\begin{rem}\label{remark: objects of orb}
We think about objects of $\Orb$ as \emph{connected topological groupoids}. Therefore, we will usually denote an object of $\Orb$ as $BG$ if we want to emphasize the isotropy group, or by a single letter, like $T$, if the isotropy group is not important.
\end{rem}

\begin{ex}\label{example: mapping spaces in Hom}
Let $G, H$ be a pair of topological groups. Then the mapping stack $\Map(BH, BG)$ between classifying stacks $BH,BG$ is equivalent to the quotient stack $[\Map_{\Grp}(H, G)/G]$, where $\Map_{\Grp}(H, G)$ denotes the space of continuous group homomorphisms and $G$ acts on %$\Map_{\Grp}(H, G)$
it by conjugation. In particular, for compact Lie groups $H, G$, by \cite[Proposition~7.1]{Rezk18}, 
$$\Hom_{\Orb}(BH,BG) \simeq  \coprod_{[\alpha]} BC_G(\Im(\alpha)),$$
where $[\alpha]$ runs over the set of $G$-conjugacy classes in $\Map_{\Grp}(H,G)$ and $C_G(\Im(\alpha))$ denotes the centralizer of the image of $\alpha$ in $G$.
\end{ex}

\begin{rem}\label{remark: geometric realization}
The functor $\Pi\colon \TopPStk \to \Type$ induces the canonical functor $\Pi\colon \Orb \to \Type$. By \cite[Theorem~1.1]{Rezk18}, the map
\begin{equation}\label{equation: geom is ff}
\Hom_{\Orb}(BH,BG) \tto \Hom_{\Type}(\Pi(BH),\Pi(BG))
\end{equation}
is an equivalence if a compact Lie group $G$ is $1$-truncated, i.e. $G$ is an extension of a finite group by a compact torus. However, for arbitrary compact Lie groups $H$ and $G$, the map~\eqref{equation: geom is ff} is highly non-trivial. For instance, let $H=G=SU(2)$ be the special unitary group. Then, by \cref{example: mapping spaces in Hom}, we have $$\Hom_{\Orb}(BH,BG)\simeq BC_2 \sqcup BSU(2).$$ However, $\Pi(BH)=\Pi(BG)\simeq \mathbb{HP}^\infty$ and $\Hom_\Type(\mathbb{HP}^\infty, \mathbb{HP}^\infty)$ has infinitely many connected components, see e.g. \cite{Mislin87} and \cite{DZ87}.
\end{rem}

\begin{defn}[Section~3.1 in \cite{Rezk_GlobalHomotopy}]\label{definition: global spaces}
The category $\Type^\Glo$ of \emdef{global spaces} is the category $\Fun(\Orb^{\op},\Type)$ of presheaves on $\Orb$.
\end{defn}

\begin{ex}[Example~3.1.8 in \cite{Lur_Ell3}]\label{example: constant global spaces}
Let $X\in \Type$ be a homotopy type. Then the \emdef{constant global space} $\underline{X} \in \Type^{\Glo}$ is defined by the formula $\underline{X}(T)=X$. Since the category $\Orb$ has a final object $\ast\simeq Be$ (i.e. the classifying stack of a trivial group), the functor
$$\Type \tto \Type^{\Glo}, \;\; X \mapsto \underline{X} $$
is fully faithful. In particular, a global space $X$ is constant if and only if $X$ is equivalent to the left Kan extension along the embedding $\ast \xhookrightarrow{Be} \Orb$.
\end{ex}

The rich class of global spaces is provided by topological (pre)stacks.

\begin{ex}\label{example: topological stacks}
Note that the restricted Yoneda embedding induces the functor
$$\yoneda\colon \TopPStk \tto \Fun_{\TopPStk}(\Orb^{\op}_{\TopStk}, \TopPStk), \;\; X \mapsto \Map(-,X)$$
where $\Fun_{\TopPStk}(\Orb^{\op}_{\TopStk}, \TopPStk)$ is the category of $\TopPStk$-enriched functors. The transfer of enrichment along the symmetric monoidal functor $\Pi\colon \TopPStk \to \Type$ induces the functor
$$\Fun_{\TopPStk}(\Orb^{\op}_{\TopStk}, \TopPStk) \tto \Fun(\Orb^{\op}, \Type) \simeq \Type^{\Glo},\;\; F \mapsto \Pi \circ F,$$
see e.g. \cite[Section~C.3]{HM24}. By taking the composite, we obtain the functor
$$\Pi^{\Glo}\colon \TopPStk \tto \Fun(\Orb^{\op},\Type) $$
which send a topological prestack $X\in \TopPStk$ to a global space $\Pi^{\Glo}(X)$ whose $T$-points $(\Pi^{\Glo}X)(T)$ are defined by the formula 
$$(\Pi^{\Glo} X)(T) \simeq \Pi\Map(T, X).$$
The global space $\Pi^{\Glo}(X)$ is called the \emdef{underlying global space} of the topological prestack~$X$. The functor $\Pi^{\Glo}$ preserves Cartesian products and coproducts, however may not preserve pullbacks or arbitrary colimits.
\end{ex}

\begin{ex}[Example~3.1.6 in \cite{Lur_Ell3}]\label{example: yoneda global spaces}
Let $X\in \Type$ be a homotopy type. Let $X^{-} \in \Type^{\Glo}$ denote the global space defined by the formula $X^{T}= \Hom_{\Type}(\Pi(T),X)$, $T\in\Orb$. More precisely, $X^{-}$ is the composite
$$X^{-} \colon \Orb^{\op} \xrightarrow{\Pi} \Type^{\op} \xrightarrow{\Hom(-,X)} \Type. $$
Note that $X^{-}$ is the underlying global space of the constant topological prestack $\underline{X}\in \TopPStk$.
\end{ex}

\subsection{Faithful morphisms}\label{section: faithful morphisms}
\begin{defn}\label{definition: faithful morphism representable}
A morphism $f\colon BH \to BG \in \Orb$ is \emdef{faithful} if the induced homomorphism $H \to G$ is \emph{injective} for some (equivalently, for any) choice of a base point of $BH$. We write $\Orb^\rep$ for the wide subcategory of $\Orb$ whose $\Hom$-spaces are spanned by faithful morphisms.
\end{defn}

\begin{defn}\label{definition: full morphism}
A morphism $f\colon BH \to BG \in \Orb$ is \emdef{full} if the induced homomorphism $H \to G$ is \emph{surjective} for some (equivalently, for any) choice of a base point of $BH$. 
\end{defn}

\begin{notation}\label{notation: full and faithful}
We will write $f\colon T' \twoheadrightarrow T$ if $f\in \Orb$ is full and $f\colon T' \hookrightarrow T$ if $f$ is faithful.
\end{notation}

We record for later applications that the class of full morphisms is closed under pullbacks in $\Orb$.

\begin{prop}\label{proposition: full morphisms are closed under pullbacks}
Let $f\colon G' \to G$ and $g\colon H \to G$ be (continuous) homomorphisms of compact Lie groups, and let $H'=H\times_G G'$. Suppose that $f$ is surjective. Then the canonical map 
$$BH' \tto BH \times_{BG} BG' \in \Type^\Glo$$
is an equivalence.
\end{prop}

\begin{proof}
See \cite[Proposition~6.1.1]{Rezk_GlobalHomotopy}.
\end{proof}

\begin{notation}\label{notation: surjective marking}
For later applications it will be convenient to mark all the full morphisms in the category $\Orb$; we denote this marked category by $(\Orb,\mathfrak{S})$.
\end{notation}

It is clear that full morphisms are left orthogonal (see \cite[Definition~5.2.8.1]{Lur_HTT}) to faithful morphisms in $\Orb$. Moreover, full and faithful morphisms form a \emph{factorization system} on $\Orb$, see \cite[Definition~5.2.8.8]{Lur_HTT}.

\begin{defn}[Section~3.4 in \cite{Rezk_GlobalHomotopy}]\label{definition: faithful morphism}
A morphism $f\colon X\to Y$ of global spaces is \emdef{faithful} if $f$ is right-orthogonal to any full morphism $T' \to T \in \Orb$ between representable presheaves.
\end{defn}

\begin{rem}\label{remark: faithful morphism}
Unwinding the definition, we see that a morphism $f\colon X \to Y$ is faithful if and only if the commutative square
\[\xymatrix{
X(BG/N) \ar[r] \ar[d]^-f & X(BG) \ar[d]^f \\
Y(BG/N) \ar[r] & Y(BG)
}\]
is fibered for all compact Lie groups $G$ and all closed normal subgroups $N \subseteq G$.
\end{rem}

We record the following properties of faithful morphisms, whose proofs are immediate from the definition.

\begin{prop}[Proposition~3.4.5 in \cite{Rezk_GlobalHomotopy}]\label{proposition: basic properties of faithful}
$\quad$
\begin{enumerate} 
\item The pullback of a faithful morphism is faithful.
\item If $g\colon Y \to Z$ is faithful, then $f\colon X \to Y$ is faithful if and only if $g\circ f$ is faithful.
\item Let $f_\bullet \colon X_\bullet \to Y_\bullet$ be a natural transformation of functors $X_\bullet, Y_\bullet \colon \fcat{I} \to \Type^\Glo$ such that
each $f_i\colon X_i \to Y_i$ is faithful. Then $\lim_{\fcat{I}} f_\bullet$ is faithful.
\end{enumerate}
\end{prop}

\subsection{Orbispaces}\label{section: orbispaces}
In this section we review an important \emph{non}-full subcategory of $\Type^\Glo$, the so-called category of orbispaces $\Type^\Orb$. For any compact Lie group $G$ the corresponding global space $BG$ is contained in  $\Type^\Orb$, and the category of $G$-spaces $\Type^G$ is naturally equivalent to $\Type^\Orb_{/BG}$ (\Cref{proposition: G-spaces and orbispace}). Some of the results about genuine and tempered local systems are either proved by reduction to the case of orbispaces, or just do not hold for an arbitrary global space. This justifies our interest in $\Type^\Orb$ in the present work.

\begin{defn}[Section~4.5 in \cite{Rezk_GlobalHomotopy}]\label{definition: orbispaces}
The category $\Type^\Orb$ of \emdef{orbispaces} is the category $\Fun(\Orb^{\rep,\op},\Type)$ of presheaves on $\Orb^{\rep}$.
\end{defn}

We write $\iota\colon \Orb^{\rep} \to \Orb$ for the evident functor. Let $\iota_!\colon \Type^\Orb \to \Type^\Glo$ denote the left Kan extension along $\iota^\op$.

\begin{prop}\label{proposition: lan orbispaces}
Let $X \in \Type^{\Orb}$ be an orbispace. Then $$(\iota_!X)(T)\simeq \coprod_{T\twoheadrightarrow T_0}X(T_0).$$
\end{prop}

\begin{proof}
This assertion is a particular case of a more general statement about the category of presheaves on an arbitrary category $\fcat{C}$ equipped with a factorization system $(\fcat{C}_L,\fcat{C}_R)$. The proof below works in the general setting and one obtains the particular case by assuming $\fcat{C}=\Orb$, $\fcat{C}_R=\Orb^{\rep}$, and $\fcat{C}_L$ is the wide subcategory of $\fcat{C}$ whose $\Hom$-spaces are spanned by full morphisms. %Recall that $(\fcat{C}_L,\fcat{C}_R)$ is a factorization system on $\fcat{C}$.

We write $\Ar_L(\fcat{C})$ for the full subcategory of the arrow category $\Ar(\fcat{C}):=\fcat{C}^{[1]}$ spanned by the edges in $\fcat{C}_L$. Let $\widetilde{\Ar}(\fcat{C})$ denote the pullback of the following diagram
$$
\xymatrix{
\widetilde{\Ar}(\fcat{C}) \ar[r] \ar[d] & \Ar_L(\fcat{C}) \ar[d]^{t} \\
\fcat{C}_R \ar[r] & \fcat{C}.
}
$$
Informally, an object of $\widetilde{\Ar}(\fcat{C})$ is an edge $T' \twoheadrightarrow T$ and a morphism $(T'\twoheadrightarrow T) \to (T'_0\twoheadrightarrow T_0)\in \widetilde{\Ar}(\fcat{C})$ is a commutative square
$$\xymatrix{
T' \ar[r]\ar@{->>}[d] & T'_0 \ar@{->>}[d]  \\
T \ar@{^{(}->}[r]& T_0
}$$
in the category $\fcat{C}$. Let $s\colon \widetilde{\Ar}(\fcat{C}) \to \fcat{C}$ and $t\colon \widetilde{\Ar}(\fcat{C}) \to \fcat{C}_R$ denote the source and target projections, respectively. 

Note that the functor $\iota \colon \fcat{C}_R \to \fcat{C}$ factors as the composite
$$\iota \colon \fcat{C}_R \xrightarrow{c} \widetilde{\Ar}(\fcat{C}) \xrightarrow{s} \fcat{C}, $$
where the functor $c\colon \fcat{C}_R \to \widetilde{\Ar}(\fcat{C})$ maps an object $X \in \fcat{C}_R$ to the identity edge $\Id_X \in \widetilde{\Ar}(\fcat{C})$. Therefore, $$\iota_! \simeq s_! \circ c_!.$$
We observe that the target projection $t\colon \widetilde{\Ar}(\fcat{C}) \to \fcat{C}_R$ is the left adjoint to $c\colon \fcat{C}_R \to \widetilde{\Ar}(\fcat{C})$. Hence, $c_! \simeq t^*$.

As in \cite[Proposition~6.7]{LNP25}, the source projection $s\colon \widetilde{\Ar}(\fcat{C}) \to \fcat{C}$ is a Cartesian fibration. We note that the strict fiber $\widetilde{\Ar}(\fcat{C})_T$ of the source projection $s$ is equivalent to the discrete category whose objects are full edges $f\colon T \twoheadrightarrow T_0$ for an object $T\in \fcat{C}$. Therefore, we have
$$s_!F(T)\simeq \coprod_{T\twoheadrightarrow T_0} F(T\twoheadrightarrow T_0)$$
for a functor $F\colon \widetilde{\Ar}(\fcat{C})^\op \to \Type$ and $T \in \fcat{C}$.
Finally, we deduce
\[(\iota_! X)(T) \simeq (s_! t^* X)(T) \simeq \coprod_{T\twoheadrightarrow T_0} t^*(X)(T\twoheadrightarrow T_0) \simeq \coprod_{T\twoheadrightarrow T_0} X(T_0).\qedhere\]
\end{proof}

\begin{ex}\label{example: normal subgroup classifier}
Let $\mathcal{N}\in \Type^\Glo$ denote the left Kan extension $\iota_!(\ast)$ of the final object $\ast \in \Type^\Orb$. Then, by \cref{proposition: lan orbispaces}, we have that the space $\mathcal{N}(BG)$ is discrete and $$\mathcal{N}(BG) \simeq \{\text{normal subgroups $N$ of $G$}\}$$ for every compact Lie group $G$. In other words, $\mathcal{N}$ is the \emdef{normal subgroup classifier}, see \cite[Section~4.1]{Rezk_GlobalHomotopy}.
\end{ex}

The computation of \cref{proposition: lan orbispaces} together with \cref{remark: faithful morphism} implies the following corollary.

\begin{cor}\label{corollary: lan and faithful}
Let $f\colon X\to Y$ be a morphism in $\Type^\Orb$. Then $\iota_!(f)$ is faithful. \qed
\end{cor}

So far, we have seen that if a global space $X$ is obtained from an orbispace, then there exists a morphism $X\to \mathcal{N}$. We recall that such a morphism is essentially unique.

\begin{prop}\label{proposition: uniqueness of orbi}
Let $X \in \Type^\Glo$ be a global space. Then the space $\Hom^\rep(X,\mathcal{N})$ of faithful morphisms to $\mathcal{N}$ is a $(-1)$-type. In other words, $\Hom^\rep(X,\mathcal{N})$ is either contractible, or empty.
\end{prop}

\begin{proof}
See \cite[Proposition~4.3.2]{Rezk_GlobalHomotopy}.
\end{proof}

Let $(\Type^{\Glo}_{/ \mathcal{N}})^{\rep} \subset \Type^{\Glo}_{/ \mathcal{N}}$ denote the full subcategory of the slice category $\Type^{\Glo}_{/ \mathcal{N}}$ whose objects are exactly the faithful morphisms $X\to \mathcal{N}$. By \cref{example: normal subgroup classifier} and \cref{corollary: lan and faithful}, the left Kan extension $\iota_!\colon \Type^{\Orb}\to \Type^{\Glo}$ induces the functor
$$\overline{\iota}_! \colon \Type^{\Orb} \tto (\Type^{\Glo}_{/ \mathcal{N}})^{\rep}.$$

\begin{thm}\label{theorem: orbi and global}
The functor $\overline{\iota}_! \colon \Type^{\Orb} \to (\Type^{\Glo}_{/ \mathcal{N}})^{\rep}$ is an equivalence.
\end{thm}

\begin{proof}
See \cite[Proposition~4.6.1]{Rezk_GlobalHomotopy}.
\end{proof}

\begin{rem}\label{remark: orbi not equivalence}
We warn the reader that the functor $\iota_!\colon \Type^{\Orb}\to \Type^{\Glo}$ is \emph{neither} fully faithful, \emph{nor} essentially surjective.
\end{rem}

\begin{rem}\label{remark: orbi inverse}
One can show that the inverse functor $(\overline{\iota}_!)^{-1}\colon (\Type^{\Glo}_{/ \mathcal{N}})^{\rep} \to \Type^\Orb$ maps a faithful morphism $X\to \mathcal{N}$ to the orbispace $Y$ defined by the formula $Y(BG)\simeq X(BG)\times_{\mathcal{N}(BG)}\{e\}$.
\end{rem}

\begin{rem}\label{remark: orbi abuse}
In the sequel, we will abuse notation by saying that a global space $X\in \Type^\Glo$ is an orbispace if there exists a faithful morphism $X\to \mathcal{N}$. Indeed, by \cref{proposition: uniqueness of orbi}, existence of such a morphism is a property and not additional data.
\end{rem}

In the rest of the section, we recall the relation between $G$-spaces and orbispaces, where $G$ is a compact Lie group.

\begin{defn}\label{definition: enriched G-orbit}
Let $G$ be a compact Lie group. The \emdef{enriched $G$-orbit category} $\Orb_{G,\Top}$ is the full topological subcategory of the category of compact (smooth) $G$-manifolds and equivariant maps spanned by the $G$-manifolds whose $G$-action is transitive.
\end{defn}

\begin{defn}\label{definition: G-orbit}
Let $G$ be a compact Lie group. The \emdef{$G$-orbit category} $\Orb_{G}$ is the underlying $\infty$-category of the enriched $G$-orbit category $\Orb_{G,\Top}$.
\end{defn}

\begin{ex}\label{example: G-orbits}
Recall that an object of the category $\Orb_G$ has the form $G/H$, where $H$ is a \emph{closed} subgroup of~$G$. Moreover, if $H,K$ are two closed subgroups of $G$, then the $\Hom$-space $\Hom_{\Orb_G}(G/H,G/K)$ is equivalent to the underlying homotopy type $\Pi((G/K)^H)$ of the $H$-fixed points $(G/K)^H$.
\end{ex}

Note that we have an enriched functor 
$$q\colon \Orb_{G,\Top} \to \Orb_{\TopStk}, \;\;G/H \mapsto [(G/H)/G]\simeq BH.$$
The transfer of enrichment along the symmetric monoidal functor $\Pi\colon \TopPStk \to \Type$ induces the functor
$$q\colon \Orb_{G} \tto \Orb. $$
We observe that the functor $q$ factors through the functor $\overline{q}\colon \Orb_{G} \to \Orb^{\rep}_{/BG}$ which maps a $G$-orbit $G/H$ to a faithful morphism $BH \to BG$. The next result is well-known, see e.g.\@ \cite{GepnerHenriques07}, \cite[Proposition~3.5.1]{Rezk_GlobalHomotopy}, or \cite[Lemma~6.12]{LNP25} for a detailed account.

\begin{prop}\label{proposition: two G-orbit}
The functor $\overline{q}\colon \Orb_{G} \to \Orb^{\rep}_{/BG}$ is an equivalence. \qed 
\end{prop}

\begin{defn}\label{definition: G-spaces}
Let $G$ be a compact Lie group. The category $\Type^G$ of \emdef{$G$-spaces} is the category $\Fun(\Orb^\op_G,\Type)$ of presheaves on $\Orb_G$.
\end{defn}

\begin{rem}\label{remark: G-spaces}
Equivalently, by \cref{proposition: two G-orbit}, the category $\Type^G$ of $G$-spaces  is equivalent to the category $\Fun(\Orb^{\rep,\op}_{/BG},\Type)$ of presheaves on $\Orb^{\rep}_{/BG}$.
\end{rem}

Abusing notation, we write $q\colon \Orb^{\rep}_{/BG} \to \Orb^{\rep}$ for the projection. Let $q_!\colon \Type^G \to \Type^{\Orb}$ denote the left Kan extension along $q$. Since left Kan extensions preserve representable presheaves, we have a functor $\overline{q}_!\colon \Type^G \to \Type^{\Orb}_{/BG}$. The next assertion is formal now, see also \cite[Proposition~3.5.1]{Rezk_GlobalHomotopy}.

\begin{prop}\label{proposition: G-spaces and orbispace}
The functor $\overline{q}_!\colon \Type^G \to \Type^{\Orb}_{/BG}$ is an equivalence. Moreover, the composite $$\overline{\iota}_!\circ \overline{q}_!\colon \Type^G \tto \Type^{\Orb}_{/BG} \tto (\Type^{\Glo}_{/BG})^{\rep}$$
is an equivalence as well. Here $(\Type^{\Glo}_{/BG})^{\rep}$ is the full subcategory of $\Type^{\Glo}_{/BG}$ spanned by the faithful morphisms $X\to BG$.
\end{prop}

\begin{proof}
The first assertion follows by e.g. \cite[Proposition~5.1.6.10]{Lur_HTT}. The second one follows by \cref{theorem: orbi and global}.
\end{proof}

\begin{notation}\label{notation: gitq}
We will write $-\gitq G\colon \Type^G \to \Type^{\Glo}$ for the left Kan extension along the forgetful functor $q\colon \Orb_G \to \Orb$.
\end{notation}
\begin{rem}
Let $X \in \Type^G$ be a $G$-space. We claim that the space of maps
$$\Hom_{\Type^\Orb}(BH, X\gitq G)$$
is naturally equivalent to
$$\coprod_{H \inj G / \sim} (X(G/H))_{hC_G(H)},$$
where $H \inj G / \sim$ denotes the set of group embeddings up to a $G$-conjugation. To see this, by construction, we have to compute the left Kan extension along a coCartesian fibration $\Orb^{\rep,\op}_{/BG} \to \Orb^{\rep, \op}$. By the colimit formula for Kan extensions, this is equivalent to
$$\colim_{\Hom_{\Type^\Orb}(BH, BG)} X(G/H).$$
However, by \Cref{example: mapping spaces in Hom}, we have
$$\Hom_{\Type^\Orb}(BH, BG) \simeq \coprod_{H \inj G/\sim} BC_G(H).$$
Hence,
$$\colim_{\Hom_{\Type^\Orb}(BH, BG)} X(G/H) \simeq \coprod_{H\inj G/\sim} (X(G/H))_{hC_G(H)}$$
as claimed.

By combining the last observation with \Cref{proposition: lan orbispaces}, we deduce that
$$(X\gitq G)(BH) \simeq \coprod_{\alpha \colon H \to G/\sim} (X(G/{\Im(\alpha)}))_{hC_G(\Im(\alpha))}.$$
\end{rem}

\begin{ex}\label{example: normal subgroup classifier as a colimit}
Let $\mathcal N = \iota_!(*) \in \Type^\Glo$ be the normal subgroup classifier, see \cref{example: normal subgroup classifier}. Since $\mathcal N$ is a presheaf on $\Orb$, $\mathcal{N}$ is a colimit of representable global spaces $BG\in \Orb$. We claim that $\mathcal{N}$ is a colimit of representable global spaces and \emph{faithful} morphisms between them. More precisely, let
$$\yoneda\colon \Orb^{\rep} \tto \Fun(\Orb^{\rep,\op}, \Type) = \Type^{\Orb} $$
be the Yoneda embedding. Then, the colimit
$$\colim_{T \in \Orb^{\rep}} \yoneda(T) \simeq \ast $$
is the final object (this is a general fact about any Yoneda embedding). Therefore,
$$\mathcal{N} \simeq \iota_!(*) \simeq \colim_{T\in \Orb^{\rep}}\iota_!(\yoneda(T)) \simeq \colim_{T\in \Orb^{\rep}} T $$
and all morphisms in the colimit are faithful by \cref{corollary: lan and faithful}.
\end{ex}

\begin{ex}\label{example: any orbispace is a colimit of G-spaces}
Let $X\in \Type^{\Glo}$ be an orbispace considered as a global space (see \cref{remark: orbi abuse}). Then $X \simeq X \times_{\mathcal{N}}\mathcal{N}$. Therefore, by \cref{example: normal subgroup classifier as a colimit}, we obtain
$$X\simeq X\times_{\mathcal{N}}\left(\colim_{BG\in \Orb^{\rep}} BG\right) \simeq \colim_{BG\in \Orb^{\rep}}X\times_{\mathcal{N}}BG.$$
By \cref{proposition: basic properties of faithful}, all morphisms in the colimit are faithful. Moreover, for each compact Lie group $G$, the orbispace $X\times_{\mathcal{N}}BG$ can be viewed as a $G$-space by \cref{proposition: G-spaces and orbispace}. So, informally, we presented any orbispace as a colimit of quotient spaces and faithful maps between them.
\end{ex}

\subsection{Families of subgroups}
Some genuine equivariant cohomology theories are first constructed only for spaces equivariant with respect to a specific family of compact Lie groups and then are formally extended for all groups. E.g.\@ Lurie's theory of tempered cohomology is defined for spaces equipped with a finite abelian group action, and Gepner--Meier's elliptic cohomology is constructed for spaces equivariant for compact abelian Lie groups. In this section we introduce some notions related to families of compact Lie groups and corresponding subcategories of global spaces and orbispaces. These subcategories are the natural domains of definition for the corresponding genuine equivariant cohomology and associated categories of local systems, which we develop in this work.
\begin{defn}\label{notat_glob_family_of_grps}
An \emdef{orbifamily} $\fcat{T}$ is a collection of compact Lie groups closed under isomorphisms and passage to subgroups. An orbifamily $\fcat{T}$ is a \emdef{global family} if $\fcat{T}$ is also closed under quotients. Finally, a global family is \emdef{multiplicative} if it is closed under Cartesian products.
\end{defn}

\begin{notation}\label{notation: families and subcategory}
Let $\fcat{T}$ be an orbifamily. Then we write $\Orb^{\rep}_{\fcat T}$ for the full subcategory of $\Orb^{\rep}$ spanned by the objects of the form $BG$, $G\in \fcat{T}$. Similarly, if $\fcat{T}$ is a global family, we write $\Orb_{\fcat{T}}$ for the full subcategory of $\Orb$ spanned by the objects $BG$, $G\in \fcat{T}$.
\end{notation}

\begin{ex}\label{example: examples of index categories}
Here are typical examples of (global) families:
\begin{itemize}
\item The family consisting of the trivial group $\{e\}$.

\item The family of finite (resp.\@ finite abelian) groups.

\item The family $\fcat{T}$ of all compact abelian Lie groups. We will denote $\Orb_{\fcat T}$ (resp. $\Orb^{\rep}_{\fcat T}$) by $\Orb_\ab$ (resp. $\Orb^{\rep}_{\ab}$).

\item The family of all compact Lie groups.
\end{itemize}
\end{ex}

\begin{ex}\label{example: orbifamily associated to G}
Let $G$ be a compact Lie group. We will denote by $\leq \!\! G$ (resp.\@ by $\fcat{P}_G$) the orbifamily of compact Lie groups which are isomorphic to a subgroup (resp.\@ a proper subgroup) of $G$. We will omit the subscript and simply write $\fcat{P}$ if the group $G$ is clear from the context.
\end{ex}

\begin{notation}\label{notation: T-orbispaces}
If $\fcat{T}$ is an orbifamily, then we will write $\Type^{\Orb}_{\fcat{T}}$ for the category of presheaves $\Fun(\Orb^{\rep,\op}_{\fcat{T}},\Type)$. If $\fcat T$ is a global family, we denote the presheaf category $\Fun(\Orb^{\op}_{\fcat{T}},\Type)$ by $\Type^\Glo_{\fcat T}$.
\end{notation}

Let $\delta_{\fcat{T}}\colon \Orb^{\rep}_{\fcat{T}}\to \Orb^{\rep}$ be the evident embedding functor. Since $\delta_{\fcat{T}}$ is fully faithful, the left Kan extension 
$$\delta_{\fcat{T},!}\colon \Type^{\Orb}_{\fcat{T}} \to \Type^{\Orb}$$
is also fully faithful.

\begin{defn}\label{definiton: classfier subgroups_family}
Let $\fcat{T}$ be an orbifamily of compact Lie groups. The \emdef{$\fcat{T}$-normal subgroup classifier} $\mathcal{N}_{\fcat{T}} \in \Type^\Glo$ is the left Kan extension $\iota_!\delta_{\fcat{T},!}(*)$ of the final object $*\in \Type^{\Orb}_{\fcat T}$. 
\end{defn}

\begin{rem}\label{remark: family_classfier is left kan} 
As in \cref{example: normal subgroup classifier}, we have
$$\mathcal{N}_\fcat{T}(BG) = \{\text{normal subgroups $N$ of $G$\;|\; $G/N \in \fcat{T}$}\}$$
for all compact Lie groups $G$. In particular, $\mathcal{N}_{\fcat{T}}\times_{\mathcal{N}}\mathcal{N}_{\fcat{T}} \simeq \mathcal{N}_{\fcat{T}}.$ 
\end{rem}

Note that $\mathcal{N}_{\fcat{T}}\in \Type^{\Glo}$ is an orbispace, so there is a (unique) faithful morphism $\mathcal{N}_{\fcat{T}} \to \mathcal{N}$. If $X\in \Type^{\Glo}$ is another orbispace, then we will write
$$\pi_{\fcat{T}}\colon X_{\fcat{T}}= X\times_{\mathcal N} \mathcal{N}_{\fcat{T}} \tto X$$
for the canonical projection. In particular, $\pi_{\fcat{T}}$ is a faithful morphism. 
\begin{rem}\label{remark: faithful points family}
Let $G$ be a compact Lie group. Then a faithful point $f\colon BG \to X$ lifts to (a necessarily faithful) point of $X_{\fcat{T}}$ if and only if $G\in \fcat{T}$. Moreover, any such lift is unique. In other words, $X_{\fcat{T}}$ has less faithful points in a way controlled by the family~$\fcat{T}$.
\end{rem}

The next assertion is straightforward.

\begin{prop}\label{proposition: family restriction}
Let $\fcat{T} $ be an orbifamily and let $X\in \Type^{\Glo}$ be an orbispace. Then there is a natural equivalence $X_{\fcat{T}} \simeq \iota_!\delta_{\fcat{T},!}\delta^*_{\fcat{T}}\iota^*(X)$ such that the projection $\pi_{\fcat{T}}$ is equivalent to the counit map of the adjunction. \qed
\end{prop}

\begin{ex}\label{example: orbifamily colimit}
Let $\fcat{T}$ be an orbifamily and $X$ be an orbispace. Then we have 
$$X\simeq \colim_{T\in \Orb^{\rep}_{/X}} T \;\; \text{and}\;\; X_{\fcat{T}} \simeq \colim_{T\in \Orb^{\rep}_{\fcat{T}/X}} T.$$
\end{ex}

\begin{ex}\label{example: constant global space}
Let $X$ be an orbispace and let $\fcat{T}$ be the trivial family. Then $\delta^*_{\fcat{T}}\iota^*(X) \simeq X(e)$ and $X_{\fcat{T}}\simeq \underline{X(e)}$ is the constant global space, see \cref{example: constant global spaces}. In particular, the canonical projection 
$$\pi_{\fcat{T}}\colon \underline{X(e)} \tto X $$
is faithful.
\end{ex}

\begin{ex}\label{example: G-spaces}
Let $G$ be a compact Lie group and $X\to BG$ be a faithful morphism in $\Type^{\Glo}$, i.e. $X$ is a $G$-space. Then $$X_{\leq G}\simeq X\times_{BG} BG \times_{\mathcal{N}} \mathcal{N}_{\leq G} \simeq X\times_{BG} BG \simeq X,$$
see \cref{example: orbifamily associated to G}.
\end{ex}

\section{Globally equivariant local systems}\label{section: globally equivariant local systems}
In this section we construct and study the category of globally equivariant local systems. It serves a somewhat technical role, but the results we prove in this section are the basis for defining and proving similar assertions for genuine equivariant and tempered local systems. Recall that non-equivariantly for a space $X$ and a category $\fcat A$ one can consider the category $\LocSys(X, \fcat A)$ of local systems on $X$ with coefficients in $\fcat A$. Equivariantly, we have an additional flexibility to choose as an input for our construction a coefficient category for each global orbit, i.e.\@ a functor
$$\fcat A \colon \Orb^\op_{\fcat T} \tto \Cat,$$
where $\fcat T$ is some family of compact Lie groups. We call such a functor a \emdef{coefficient system}. We discuss this notion and consider several examples in \Cref{section: coefficient systems}. In \Cref{ssect_coefs_from_PreAbStk} we study an important class of coefficient systems attached to a preoriented abelian group spectral stack. These coefficient systems are used in subsequent sections to categorify tempered and elliptic cohomology.

The construction of $\LocSys^\Glo_{\fcat T}$ is introduced in \Cref{section: local systems on global spaces}. By definition it is contravariantly functorial in global spaces, and as a first important technical result we show that in fact it is right Kan extended from the subcategory of $\fcat T$-orbits. We also establish several other basic properties of $\LocSys^\Glo$ in this section. Moreover, for formal reasons the pullback functor admits a right adjoint, and we show that it is well behaved if the coefficient system is limit preserving.

Unfortunately, the coefficient systems we are most interested in  are not of this type in general. But it turns out that even under milder assumptions on $\fcat A$, $\LocSys^\Glo_{\fcat T}(-, \fcat A)$ is still well-behaved when restricted to the non-full subcategory of orbispaces and faithful morphisms. In \Cref{section: local systems on orbispaces} we give a more compact description $\LocSys^\Glo_{\fcat T}(X, \fcat A)$ for $X \in \Type^\Orb$ and deduce several nice properties. In particular, we show that the pullback functor along a faithful morphism admits a left adjoint $\#$-pushforward functor, such that the base change and projection formula hold. In \Cref{subsection: functoriality along faithful morphisms} we extend these results to faithful morphisms between general global spaces.

Finally, in \Cref{section: global section} for a global space $X$ we construct an adjunction
$$\mathbbl 1_X \otimes - \colon \fcat A(X) \xymatrix{\ar@{^(->}@<0.5ex>[r] & \ar@<0.5ex>[l]} \LocSys^\Glo_{\fcat T}(X, \fcat A) \colon \Gamma_{\fcat A}(X, -),$$
where the right hand side refines the usual global sections functor. Moreover, for a compact Lie group $G \in \fcat T$ we show that $\LocSys^\Glo_{\fcat T}(BG, \fcat A)$ admits a semi-orthogonal decomposition (a.k.a. recollement) generalizing the classical isotropy separation.

\subsection{Coefficient systems}\label{section: coefficient systems}
For a homotopy type $X$ and a category $\fcat A$ one can consider the category of local systems on $X$ with coefficients in $\fcat A$. It turns out that passing to the equivariant setting, one obtains an additional flexibility of choosing a coefficient category $\fcat A(BG)$ for each compact Lie group $G$, which are functorial with respect to the morphisms in $\Orb$. We call such a choice a \emdef{coefficient system}, see \Cref{definition: coefficient system}. Moreover, for the corresponding theory of local systems to be reasonable, one needs to put some additional assumptions on $\fcat A$. In this section we introduce the corresponding notions of coefficient systems and construct several examples, including the one coming from an oriented abelian group object. For the rest of this section $\fcat T$ denotes a global family of compact Lie groups, see \Cref{notat_glob_family_of_grps}.

\begin{defn}\label{definition: coefficient system}
A \emdef{coefficient system $\fcat A$} is a functor
$$\fcat A \colon \Orb^\op_{\fcat{T}} \tto \Prs^\LL.$$
For a morphism $\phi \colon T_0 \to T$, we will denote by $\phi_{\fcat A}^*$ the induced functor
$$\phi_{\fcat A}^*\colon \fcat A(T) \tto \fcat A(T_0).$$
We say that the coefficient system is stable (resp. $E_n$-monoidal, etc.) if we are given a factorization of $\fcat A$ through a forgetful functor $\Prs^{\LL}_{\mathrm{st}} \to \Prs^\LL$ (resp. $\Alg_{E_n}(\Prs^\LL) \to \Prs^\LL$, etc.)
\end{defn}

\begin{rem}
The notation $\phi^*_\fcat{A}$ is not always convenient for us, so we will also use other notation like $\fcat{A}^*(\phi)$ for the same transition functor in the sequel. We hope this ambiguity will not confuse the reader too much.
\end{rem}

\begin{notation}
We will write $\phi_{\fcat{A},*}$ or $\fcat{A}_*(\phi)$ for the right adjoint to the transition functor $\phi^*_{\fcat{A}}$ in the coefficient system $\fcat{A}$.
\end{notation}

\begin{defn}\label{assumpt_limit_prerserving_cs}
A coefficient system $\fcat A\colon \Orb^{\op}_{\fcat{T}} \to \Prs^\LL$ is 
\begin{enumerate}
\item
\emdef{strongly continuous} if the transition functor $\phi_{\fcat{A}}^*$ in $\fcat A$ is strongly continuous for every faithful morphism $\phi \in \Orb_{\fcat{T}}$, i.e. the right adjoint $\phi_{\fcat{A},*}$ preserves small colimits;
\item
\emdef{limit-preserving} if the transition functors $\phi_{\fcat{A}}^*$ in $\fcat A$ preserve limits for all morphisms $\phi \in \Orb_{\fcat{T}}$;
\item
\emdef{of Beck--Chevalley type} if, for every full morphism $\phi\colon T^\prime \surj T$ and any morphism $\psi\colon T_0 \to T$, the commutative square
\[\xymatrix{
\fcat A(T_0^\prime) & \ar[l]_{(\phi^\prime_{\fcat A})^*} \fcat A(T_0) \\
\fcat A(T^\prime) \ar[u]^-{(\psi^\prime_{\fcat A})^*} & \ar[l]_-{\phi^*_{\fcat A}}\ar[u]_-{\psi^*_{\fcat A}} \fcat A(T),
}\]
where $T_0^\prime := T^\prime \times_T T_0$ (see \cref{proposition: full morphisms are closed under pullbacks}), is vertically right adjointable. That is the natural transformation 
$$\psi_{\fcat{A},*}\phi_{\fcat{A}}^* \tto (\psi^\prime_{\fcat{A}})_*(\phi^\prime_{\fcat{A}})^* $$
is an equivalence.
\end{enumerate}
\end{defn}

\begin{defn}\label{definition: coefficient system is of geometric type}
An $E_\infty$-monoidal coefficient system $\fcat A\colon \Orb^{\op}_{\fcat{T}} \to \Prs^\LL$ is \emdef{of geometric type} if $\fcat{A}$ is a strongly continuous coefficient system of Beck--Chevalley type such that the right adjoint $$\phi_{\fcat{A},*} \colon \fcat{A}(T') \tto \fcat{A}(T)$$ is conservative and strictly $\fcat{A}(T)$-linear for every faithful morphism $\phi\colon T'\to T$. That is $\phi_{\fcat{A},*}$ satisfies the projection formula with respect to $\phi_{\fcat{A}}^*$.
\end{defn}

\begin{rem}\label{remark: geometric and pregenuine}
Let $\fcat{A}\colon \Orb^{\op}_{\fcat{T}} \to \CAlg(\Prs^\LL)$ be an $E_\infty$-monoidal coefficient system. Then the right Kan extension
$$\overline{\fcat{A}}\colon (\Type^{\Glo}_{\fcat{T}})^{\op}\tto \CAlg(\Prs^\LL)$$
along the Yoneda embedding $\Orb_{\fcat{T}} \hookrightarrow \Type^{\Glo}_{\fcat{T}}$ is a \emph{global 2-ring} in the sense of~\cite[Definition~6.4]{GLP26}. In~\cite{GLP26}, the authors develop the notions of \emph{pregenuine} (see~\cite[Definition~10.1]{GLP26}) and \emph{rigid} (see~\cite[Definition~13.4]{GLP26}) global 2-rings. We warn the reader that if a global 2-ring $\overline{\fcat{A}}$ is rigid, then the coefficient system $\fcat{A}$ is strongly continuous, but not vice versa. Also, if $\overline{\fcat{A}}$ is rigid, then $\fcat{A}$ is not of geometric type, since the right adjoints $\phi_{\fcat{A},*}$ may not be conservative, although they satisfy the projection formula. 

Also, the Beck--Chevalley property in~\cite[Definition~10.1(1)]{GLP26} is different from the one in \cref{assumpt_limit_prerserving_cs} since we only consider pullbacks along the \emph{full} morphism and pushforwards along \emph{all} morphisms. Whereas~\cite[Definition~10.1(1)]{GLP26} deals with pullbacks along \emph{all} morphisms and pushforwards only along \emph{faithful} maps.
\end{rem}

\begin{rem}\label{remark: geometric type and modules}
If $\fcat{A}\colon \Orb^\op_{\fcat{T}} \to \Prs^\LL$ is a coefficient system of geometric type and $H\subset G$ is a closed subgroup in a compact Lie group $G\in \fcat{T}$, then $\fcat{A}(BH) \simeq \Mod_R\fcat{A}(BG)$ is the category of $R$-modules in $\fcat{A}(BG)$ for some $R\in \CAlg(\fcat{A}(BG))$. However, the categories $\fcat{A}(BH)$ and $\fcat{A}(BG)$ still can be rather complicated.
\end{rem}

\begin{ex}\label{example: constant system}
Let $\fcat{C} \in \CAlg(\Prs^{\LL})$ be a presentably symmetric monoidal category. We denote by
$$\underline{\fcat{C}}\colon \Orb^\op_{\fcat{T}} \tto \Prs^\LL $$ 
the \emdef{constant coefficient system} valued at $\fcat{C}$. Clearly, the constant coefficient system $\underline{\fcat{C}}$ is of geometric type.
\end{ex}

\begin{ex}\label{example: itterated coeff system}
Let $\fcat{C} \in \CAlg(\Prs^{\LL}_{\mathrm{st}})$ be a stable presentably symmetric monoidal category. The coefficient system
$$\fcat{C}^{\flat}\colon \Orb^\op_{\ab} \xrightarrow{\Pi} \Type^{\op} \xrightarrow{\LocSys(-,\fcat{C})} \Prs^\LL $$
is a coefficient system of geometric type. Indeed, let $BH \to BG$ be a faithful morphism, then we have the adjoint pair
$$\res^G_H\colon \xymatrix{\fcat{C}^\flat(BG)\simeq \LocSys(BG,\fcat{C}) \ar@<0.5ex>[r] & \ar@<0.5ex>[l] \LocSys(BH,\fcat{C})\simeq \fcat{C}^\flat(BH) } \colon \Coind^G_H.$$
Note that the right adjoint $\Coind^G_H$ satisfies the projection formula and $$\res^G_H\Coind^G_H(-)\simeq [\Sigma^{\infty}_+G/H,-].$$ Since $G$ is abelian, $\Sigma^{\infty}_+G/H$ is equivalent to a wedge of sphere spectra with the trivial $H$-action. So, $\Coind^G_H$ is also conservative. This argument proves the main property from~\cref{definition: coefficient system is of geometric type} and we will leave it to the reader to check the strong continuity and the Beck--Chevalley property from~\cref{assumpt_limit_prerserving_cs}.
\end{ex}

\begin{ex}\label{example: cochains} Let $R\in \CAlg(\Sp)$ be an $E_\infty$-ring and let 
$$R^{(-)}\colon \Type^{\op} \tto \CAlg(\Sp)$$ be the functor which sends a space $X$ to its $R$-cochains $R^X$. Then the composite
$$\fcat{A}_R\colon \Orb^\op \xrightarrow{\Pi} \Type^{\op} \xrightarrow{R^{(-)}} \CAlg(\Sp) \xrightarrow{\Mod_{(-)}} \Prs^\LL $$
is a coefficient system. Suppose additionally that $R$ is complex orientable, then the restricted coefficient system
$$\fcat{A}_R\colon \Orb^\op_{\ab} \xrightarrow{\Pi} \Type^{\op} \xrightarrow{R^{(-)}} \CAlg(\Sp) \xrightarrow{\Mod_{(-)}} \Prs^\LL $$
is of geometric type. Indeed, it is clear that the coefficient system $\fcat{A}_R$ is strongly continuous and it satisfies the main property of \cref{definition: coefficient system is of geometric type}. So, we only have to check the Beck--Chevalley property. However, by~\cite[Proposition~4.6.11]{Lur_Ell2}, if
\[\xymatrix{
BH^\prime \ar[d] \ar@{->>}[r] & BH \ar[d] \\
BG^\prime \ar@{->>}[r] & BG,
}\]
is a fibered square in $\Orb_{\ab}$ (see \cref{proposition: full morphisms are closed under pullbacks}), then $$R^{BH'} \simeq R^{BH}\otimes_{R^{BG}} R^{BG'} \in \CAlg(\Sp).$$
This implies the assertion.
\end{ex}

Suppose that $X=BT\simeq \mathbb{CP}^\infty$, then the $E_\infty$-ring 
$$R^{BT}\simeq R^{\mathbb{CP}^\infty}\simeq \lim R^{\mathbb{CP}^n}$$ is the inverse limit of ``nilpotent'' $E_\infty$-rings. Therefore, it is more natural to consider the coefficient system $\widehat{\fcat{A}}_R$ which assigns to a topological groupoid $BG$ the category $\Mod_{R^{BG}}^{\complete}$ of ``complete'' $R^{BG}$-modules, rather than the category $\Mod_{R^{BG}}$ of all modules as in \cref{example: cochains}. We will make this construction precise in \cref{example: complete cochains}.

\begin{ex}\label{example: complete cochains} Let $R\in \CAlg(\Sp)$ be an $E_\infty$-ring, let $X \in \Type$ be a space and let $e\colon \{*\}\to X$ be a point. Then $e^*\colon R^X \to R$ is a map of $E_\infty$-rings, i.e. $R$ is an $R^X$-$E_\infty$-ring. Recall from \cite[Definition~2.15]{MNN17} that an $R^X$-module $M$ is \emph{$R$-acyclic} if $M\otimes_{R^X} R\simeq 0$ and \emph{$R$-complete} if it is right-orthogonal to all $R$-acyclic modules. We write $\Mod_{R^X}^{\complete}$ for the full subcategory of $\Mod_{R^X}$ spanned by $R$-complete modules. Note that, if $X$ is connected, the full subcategory $\Mod_{R^X}^{\complete}$ does not depend on the choice of the basepoint $e$. Moreover, by the paragraph above~\cite[Definition~2.19]{MNN17}, the inclusion
$$\Mod_{R^X}^{\complete} \hookrightarrow \Mod_{R^X} $$
admits a left adjoint $L^R_X\colon \Mod_{R^X} \to \Mod_{R^X}^{\complete}$, which is the \emph{$R$-completion functor}. So the assignment $X\mapsto \Mod_{R^X}^{\complete}$ glues into the functor\footnote{Here $\Type_{\geq 1} \subset \Type$ stands for the full subcategory of \emph{connected} homotopy types.}
$$\Mod_{R^{(-)}}^{\complete}\colon \Type^{\op}_{\ge 1} \tto \Prs^\LL$$
and we obtain the refinement
$$\widehat{\fcat{A}}_{R}\colon \Orb^\op \xrightarrow{\Pi} \Type^{\op}_{\ge 1} \xrightarrow{\Mod_{R^{(-)}}^{\complete}} \Prs^\LL $$
of the coefficient system $\fcat{A}_R$ from \cref{example: cochains}.
\end{ex}

As in \cref{example: cochains}, one expects that the coefficient system $\widehat{\fcat{A}}_R$ is well-behaved if $R$ is complex orientable and $\widehat{\fcat{A}}_R$ is restricted to \emph{abelian} topological groupoids.

\begin{lem}\label{lemma: faithful morphisms completion}
Let $R$ be a complex orientable $E_\infty$-ring and $f\colon BH \to BG$ be a morphism in $\Orb_{\ab}$. Then the right adjoint 
$$f_*\colon \Mod_{R^{BH}} \to \Mod_{R^{BG}}$$
preserves the full subcategories of $R$-complete objects. Moreover, if $f$ is faithful, then $f_*$ commutes with $R$-completion, i.e. the natural transformation
$$\theta_f\colon L^R_{BG} f_* M \tto f_*L^R_{BH} M $$
is an equivalence for any $M\in \Mod_{R^{BH}}$.
\end{lem}

\begin{proof}
The first part is clear, since the left adjoint $f^*(-)\simeq (-)\otimes_{R^{BG}}R^{BH}$ preserves $R$-acyclic objects. 

We will prove the second part. Suppose that $f$ is faithful and let $g\colon BG \to BT$ be a faithful morphism such that $T$ is a (high-dimensional) torus. Since $g_*$ is conservative, $\theta_f$ is an equivalence if and only if $\theta_g$ and $\theta_{g\circ f}$ are equivalences. So, without loss of generality, we can assume that $G=T$ is a torus.

Suppose that $G=T$ is a torus. Then it is enough to show that $f_*$ preserves $R$-acyclic modules. Let $M\in \Mod_{R^{BH}}$ be an $R$-acyclic module, then 
$$f_*(M)\otimes_{R^{BT}}R\simeq M\otimes_{R^{BH}} \left(R^{BH} \otimes_{R^{BT}} R\right). $$
Since $T$ is a torus, we have $R^{BH} \otimes_{R^{BT}} R \simeq R^{T/H}$ by~\cite[Proposition~4.6.11]{Lur_Ell2}. Note that $T/H$ is again a torus, so the $R$-cochains $R^{T/H}$ is a direct sum of several copies of $R$ as an $R^{BH}$-module. This implies that  $f_*(M)\otimes_{R^{BT}}R\simeq 0$, which implies the assertion.
\end{proof}

\begin{cor}[Quillen formal group]\label{corollary: Quillen formal group} Let $R\in \CAlg(\Sp)$ be a complex orientable $E_\infty$-ring. The composite
$$\widehat{\fcat{A}}_{R}\colon \Orb^\op_{\ab} \xrightarrow{\Pi} \Type^{\op}_{\ge 1} \xrightarrow{\Mod_{R^{(-)}}^{\complete}} \Prs^\LL $$
is a coefficient system of geometric type.
\end{cor}

\begin{proof}
The assertion follows immediately from \cref{example: cochains} and \cref{lemma: faithful morphisms completion}. 
\end{proof}

\begin{rem}\label{remark: flat and complete}
Let $R$ be a complex orientable $E_\infty$-ring. Then there is a natural transformation
$$(\Mod_R)^\flat \tto \widehat{\fcat{A}}_R $$
of coefficient systems. This natural transformation is \emph{not} an equivalence, but it is an equivalence when restricted to the full subcategory of $\Orb_{\ab}$ spanned by classifying groupoids of tori, see~\cite[Theorem~7.43]{MNN17}.
\end{rem}

\begin{rem}\label{remark: two completeness}
Let $R$ be an $E_\infty$-ring and let $G$ be a compact Lie group. Then we will write $I_G \subset \pi_0(R^{BG})$ for the \emph{augmentation ideal}, i.e.
$$I_G=\ker(\pi_0(R^{BG}) \tto \pi_0(R)).$$
Suppose that $R$ is \emph{complex periodic}, see~\cref{definition: complex periodic}, and $G$ is abelian. We note that an $R^{BG}$-module $M$ is $R$-complete if and only if $M$ is $I_G$-complete, see~\cite[Definition~7.3.1.1]{Lur_SAG} and \cref{section: formal completions}.
\end{rem}

\subsection{Coefficient system associated with abelian group object}\label{ssect_coefs_from_PreAbStk}
In this section we construct a family of examples of geometric origin which is most important for this work. First we need to review a few geometric preliminaries.
\begin{notation}\label{notat_SpDmNcFlatQcs}
Let $S$ be a non-connective spectral stack. We denote by $\SpDM_S^{\nc,\flat,\qcs}$ the full subcategory of $\SpStk_{/S}^\nc$ spanned by objects $X$ such that for every non-connective $E_\infty$-ring $R$ and a map $\Spec R \to S$ the pullback $X_R := \Spec R \times_S X$ is representable by a non-connective spectral Deligne--Mumford stack (\cite[Definition 1.4.4.2]{Lur_SAG}) which is flat over $\Spec R$ and the connective cover $X_{R,\ge 0}$ is quasi-compact and separated (in the strong sense of \cite[Definition 3.2.0.1]{Lur_SAG}, i.e.\@ such that the relative diagonal is a closed embedding).
\end{notation}
\begin{rem}
Note that a separated spectral Deligne--Mumford stack is an algebraic space. Hence, in the previous definition, $X_{R, \ge 0}$ is represented by a separated spectral algebraic space.
\end{rem}
\begin{defn}\label{definition: ab nc stack}
By abuse of notation we define \emdef{non-connective spectral abelian group stack $A$ over $S$} to be an abelian group object in $\SpDM_S^{\nc,\flat,\qcs}$. We will write $\Ab(\SpDM_S^{\nc,\flat,\qcs})$ for the category of abelian group non-connective spectral stacks.
\end{defn}

\begin{ex}
Let $R$ be a non-connective $E_\infty$-ring. Then a strict abelian variety over $R$ in the sense of \cite[Definition 1.5.1]{Lur_Ell1} satisfies the conditions of \Cref{definition: ab nc stack}. In fact, in loc.~cit.\@ Lurie additionally requires $A_{\ge 0}$ to be proper, locally almost of finite presentation, geometrically reduced, and geometrically connected. Similarly, if $S$ is a non-connective spectral Deligne--Mumford stack and $E$ is an elliptic curve over $S$ in the sense of \cite[Definition 5.6]{GM}, then $E$ is also an example of a non-connective spectral abelian group stack over $S$.
\end{ex}

\begin{defn}[Definition~3.6 in \cite{GM}]\label{definition: preab nc stack}
Let $\fcat{X}$ be a category with finite products. A \emdef{preoriented abelian group object} in $\fcat{X}$ is an abelian group object $A\in \Ab(\fcat{X})$ equipped with a morphism $BU(1) \to \Hom_{\fcat{X}}(*,A)$ of abelian group objects in $\Type$. We define the category $\PreAb(\fcat{X})$ of preoriented abelian group objects in $\fcat{X}$ as $\Ab(\fcat{X})\times_{\Ab(\fcat{\Type})} \Ab(\fcat{\Type})_{BU(1)/}$, where $\Ab(\fcat{X}) \to \Ab(\Type)$ is corepresented by the terminal object $* \in \fcat{X}$.
\end{defn}

In particular, an abelian group non-connective spectral stack $A \in \Ab(\SpDM_S^{\nc,\flat,\qcs})$ is \emdef{preoriented} if $A$ is equipped with the morphism $BU(1) \to A(S)$ of abelian group objects in~$\Type$. We will write $\PreAb(\SpDM_S^{\nc,\flat,\qcs})$ for the category of preoriented abelian group non-connective spectral stacks.

\begin{rem}\label{rem_preor_induces_U_1_eq_lift}
Let $A$ be an abelian group object in $\fcat X$ with the unit $e\colon * \to A$. Then a preorientations of $A$ in particular induces a lift of $e$ to an $U(1)$-equivariant map, where both the source and the target are equipped with the trivial $U(1)$-action.
\end{rem}
\begin{rem}\label{pre_or_depends_only_formal_cpl}
Compare with~\cite[Section~4.3]{Lur_Ell2} and~\cite[Section~5]{GM}. Let $A$ be an abelian group object over $S$ and let
$$\widetilde e \colon BU(1) \times S \tto A$$
be a preorientation, where we consider $BU(1)$ as a constant spectral stacks. Since $BU(1)$ is a connected homotopy type, the map $\widetilde e$ factors through the formal completion $\widehat A$ of $A$ at the unit section (see \Cref{constr_compl_and_comp}). In particular, the restriction map
$$\PreOr(\widehat A) \tto \PreOr(A)$$
is an equivalence.
\end{rem}

We recall that if $G$ is a compact abelian Lie group, then its Pontryagin dual group $\widehat{G} =\Hom(G,U(1))$ is a (discrete) finitely generated abelian group. By~\cite[Construction~3.8]{GM}, the Pontryagin duality extends to the functor
$$\widehat{(-)} \colon \Orb_{\ab}^{\op} \tto \PreAb(\Type) $$
such that $\widehat{BG} \simeq \Hom_{\Orb}(BG,BU(1))\simeq \widehat{G}\times BU(1)$ and the preorientation of $\widehat{BG}$ is induced by the canonical map $BG \to \{*\}$. 

Let $M$ be a discrete finitely generated abelian group. Then we will consider $M$ as a \emph{constant} abelian group non-connective spectral stack over $S$. If $A\in \Ab(\SpDM_S^{\nc,\flat,\qcs})$, then we will write $A[M]\in \SpDM_S^{\nc,\flat,\qcs}$ for the non-connective spectral stack of \emdef{$M$-torsion points} of $A$, that is
$$A[M]=\Map_{\Ab(\SpDM_S^{\nc,\flat,\qcs})}(M,A).$$
Similarly, if $G$ is a compact abelian Lie group and $A\in \PreAb(\SpDM_S^{\nc,\flat,\qcs})$ is preoriented, then we will write $A[\widehat{BG}]$ for the stack $\Map_{\PreAb(\SpDM_S^{\nc,\flat,\qcs})}(\widehat{BG},A)$. We note that $A[\widehat{BG}]\simeq A[\widehat{G}]$, but the first formula is functorial over \emph{all} maps $BH \to BG$ in $\Orb_{\ab}$.

\begin{ex}\label{example: coefficient system from preoriented abelian group}
Let $A \in \PreAb(\SpDM_S^{\nc,\flat,\qcs})$ be a preoriented abelian group non-connective spectral stack over $S$. Then \cite[Construction~3.13]{GM} gives a canonical stable symmetric monoidal coefficient system $\fcat A \colon \Orb^{\op}_{\ab} \to \CAlg(\Prs^{\LL}_{\mathrm{st}})$ as the following composite
$$\fcat{A} \colon \Orb^{\op}_{\ab} \xrightarrow{\widehat{(-)}} \PreAb(\Type) \xrightarrow{\Map_{\PreAb(\SpDM_S^{\nc,\flat,\qcs})}(-,A)} \SpStk^{\op}_{/S} \xrightarrow{\QCoh} \CAlg(\Prs^{\LL}_{\mathrm{st}}).$$
In other words, $\fcat A(BG) \simeq \QCoh(A[\widehat{G}])$ for every compact abelian Lie group $G$.
\end{ex}
Our goal for the rest of this section is to show that the coefficient system from the previous example is of geometric type. %provided the structure $p\colon A \to S$ is flat.

\begin{lem}\label{lem_grp_embeding_induces_aff_morph}
Let $A$ be a preoriented non-connective abelian group spectral stack over $S$ and let $H \inj G$ be an inclusion of compact abelian Lie groups. Then the induced morphism
$$A[\widehat{H}] \tto A[\widehat{G}]$$
is non-connective affine, and the induced morphism of connective covers is a closed embedding.

\begin{proof}
In fact, the proof of a similar assertion for spectral elliptic curves \cite[Proposition 6.3]{GM} works in our setting as well. We review the argument here for reader's convenience. Let $Q:= G/H$. By \cref{proposition: full morphisms are closed under pullbacks}, we have a fibered square
\[\xymatrix{
A[\widehat{H}] \ar[r] \ar[d] & A[\widehat{G}] \ar[d] \\
S \ar[r]^-{e} & A[\widehat{Q}],
}\]
where the lower horizontal map $e$ classifies the unit in $A[\widehat{Q}]$. Since nc-affine maps are closed under pullbacks, it is enough to show that $e$ is nc-affine.

Since the property of a morphism to be nc-affine is local on target, we can assume without loss of generality that $S = \Spec R$ is nc-affine. In this case, since $A[\widehat{Q}]$ is flat over $S$, the morphism $e$ is a base change along $S \to S_{\ge 0} := \Spec R_{\ge 0}$ of the unit map
$$e_{\ge 0} \colon S_{\ge 0} \tto A[\widehat{Q}]_{\ge 0}$$
of the abelian group spectral algebraic spaces $A[\widehat{Q}]_{\ge 0}$ over $S_{\ge 0}$, see \cite[Proposition 2.8.2.10]{Lur_SAG}. We claim that $e_{\ge 0}$ is affine (in fact, a closed embedding). More generally, let $p\colon X \to S_{\ge 0}$ be a separated morphism of spectral Deligne--Mumford stacks, and let $s\colon S_{\ge 0} \to X$ be a section of $p$. We claim that $s$ is a closed embedding. Indeed, $s$ fits into a fibered square
\[\xymatrix{
S_{\ge 0} \ar[r] \ar[d]^s & X \ar[d]^{\Delta} \\
X \ar[r]^-{(s, \Id)} & X \times_{S_{\ge 0}} X,
}\]
where the relative diagonal $\Delta$ is a closed embedding by the separatedness assumption on $p$. Since closed embeddings are stable under base change, the result follows.
\end{proof}
\end{lem}

In order to verify the other conditions of geometric type coefficient systems, we will use the following non-connective spectral analog of Nisnevich-local decomposition of quasi-compact quasi-separated algebraic spaces studied in \cite[Sections 2.5 and 3.4]{Lur_SAG}.
\begin{defn}
A morphism $f\colon X \to Y$ of non-connective spectral stacks is called \emdef{relatively scalloped} if, for any non-connective affine $U$ mapping to $Y$, the fiber product $U \times_Y X$ is a non-connective spectral Deligne--Mumford stack admitting a scallop decomposition (see \cite[Definition 2.5.3.1]{Lur_SAG}).
\end{defn}
\begin{rem}\label{rem_scalloped_in_connective_case}
By \cite[Theorem 3.4.2.1]{Lur_SAG}, a morphism of \emph{connective} spectral stacks is relatively scalloped if and only if it is representable in quasi-compact quasi-separated spectral algebraic spaces.
\end{rem}

The following result can be proved exactly as for derived stacks (see e.g.\@ \cite[Proposition 2.2.2]{GaitsRozI}).
\begin{lem}\label{lem_push_is_nice_for_rep_qcqs_alg_sp_morphs}
Let $f\colon X \to Y$ be relatively scalloped morphism of non-connective spectral stacks. Then:
\begin{enumerate}
\item The pushforward functor
$$f_*\colon \QCoh(X) \tto \QCoh(Y)$$
preserves colimits and is $\QCoh(Y)$-linear, that is, the natural projection formula comparison morphism
$$\mathcal F \otimes f_*(\mathcal G) \tto f_*(\mathcal F \otimes f^*(\mathcal G))$$
is an equivalence.

\item Let $g\colon Y^\prime \to Y$ be an arbitrary map of non-connective spectral stacks, and denote $X^\prime := Y^\prime \times_Y X$. Then the commutative square of categories
\[\xymatrix{
\QCoh(X^\prime) & \ar[l]_-{g^{\prime*}} \QCoh(X) \\
\QCoh(Y^\prime) \ar[u]_{f^{\prime*}} & \ar[l]_-{g^*}\ar[u]_{f^*} \QCoh(Y)
}\]
is vertically right adjointable, that is, the natural map
$$g^* \circ f_* \tto f^\prime_* \circ g^{\prime *}$$
is an equivalence.
\end{enumerate}

\begin{proof}
We first prove the second assertion. By presenting $Y$ as a colimit of non-connective affines and by using that adjointable squares are stable under small limits (\cite[Corollary 4.7.4.18]{Lur_HA}), we can assume that $Y$ is nc-affine. Moreover, by restricting the base-change comparison map to any nc-affine mapping to $Y^\prime$, we can assume that $Y^\prime$ is also nc-affine. So, by assumption, we have a fiber square of non-connective spectral Deligne--Mumford stacks, in which case the result follows from \cite[Proposition 2.5.4.5]{Lur_SAG}.

By the previously proven base change property, the continuity of $f_*$ and the projection formula follow formally from the case where $X$ is a non-connective spectral Deligne--Mumford stack admitting a scallop decomposition and $Y$ is nc-affine. The continuity in this case is treated in \cite[Proposition 2.5.4.3]{Lur_SAG}, and the projection formula holds, since $\QCoh(Y)$ is generated under shifts and colimits by the monoidal unit.
\end{proof}
\end{lem}

\begin{cor}\label{corollary: preab is geom type}
Let $p\colon A \to S$ be as in \Cref{definition: ab nc stack}. Then the coefficient system $\fcat A(T) := \QCoh(A[\widehat{T}])$ from \Cref{example: coefficient system from preoriented abelian group} is of geometric type.

\begin{proof}
By \Cref{lem_grp_embeding_induces_aff_morph}, the last condition of the geometric coefficient system is satisfied. By \Cref{lem_push_is_nice_for_rep_qcqs_alg_sp_morphs} and \cref{proposition: full morphisms are closed under pullbacks}, to verify that $\fcat A$ is strongly continuous and is of Beck--Chevalley type, it suffices to show that for any morphism $BH \to BG$ the induced morphism of non-connective spectral stacks
$$f\colon A[\widehat{H}] \tto A[\widehat{G}]$$
is relatively scalloped. Since the desired property is local on the target, we can assume that $S$ is non-connective affine.

By assumption, the stack $A[\widehat{T}]$ is flat over $S$ and its connective cover $A[\widehat{T}]_{\ge 0}$ is a quasi-compact separated spectral algebraic space for each $T \in \Orb_\ab$. Since quasi-compact quasi-separated spectral algebraic spaces are closed under fibered products, by using~\cite[Theorem~3.4.2.1]{Lur_SAG}, we deduce that the morphism 
$$f_{\ge 0}\colon A[\widehat{H}]_{\ge 0} \tto A[\widehat{G}]_{\ge 0}$$
is relatively scalloped. By flatness, we know that $f$ is the base change of $f_{\ge 0}$ along the nc-affine morphism $A[\widehat{G}] \to A[\widehat{G}]_{\ge 0}$, hence $f$ is also relatively scalloped.
\end{proof}
\end{cor}

\begin{ex}[Additive group]\label{example: additive group}
Compare with \cite[Section~1.6.4]{Lur_Ell2}. Let $\mathbb{G}_{a,\mathbb{Q}}$ denote the affine line $\mathbb{A}^1_{\mathbb{Q}} = \Spec(\mathbb{Q}[t])$ over $\Spec(\mathbb{Q})$. We regard $\mathbb{G}_{a,\mathbb{Q}}$ as a commutative group scheme over~$\mathbb{Q}$, with addition law given by the comultiplication
$$\Delta \colon \mathbb{Q}[t] \tto \mathbb{Q}[t] \otimes_{\mathbb{Q}} \mathbb{Q}[t], \;\; t\mapsto t\otimes 1 + 1\otimes t.$$
As such $\mathbb{G}_{a,\mathbb{Q}}$ is a (strict) abelian group object in spectral prestacks over $\mathbb{Q}$. For a rational $E_\infty$-algebra~$R$, we define $\mathbb{G}_{a,R}$ as $\mathbb{G}_{a,\mathbb{Q}}\times_{\mathbb{Q}}\Spec(R)$. So, $\mathbb{G}_{a,R}$ is again a strict abelian group object in spectral prestacks over $R$ corepresented by the \emph{free} $R$-$E_\infty$-algebra $R\otimes \Sigma^{\infty}_+\mathbb{N}$. Note that this construction does \emph{not} work over a non-rational base, see~\cite[Proposition~1.6.20]{Lur_Ell2}.

By~\cite[Definition~3.6]{GM}, the preorientations of $\mathbb{G}_{a,R}$ correspond to the maps 
$$\Spec R \times BT \tto \mathbb{G}_{a,R}$$
of strict abelian group objects in spectral prestacks over $R$. The latter is equivalent to the maps 
$$\Spec R\times S^2 \tto \mathbb{A}^1_{R}$$
of spectral prestacks over $R$. Since the stack $\mathbb{A}^1_{R}$ is nc-affine, we obtain that the preorientations of $\mathbb{G}_{a,R}$ correspond one-to-one to the elements of $\pi_2(\mathbb{A}^1_{R}(R))\cong\pi_2(\Omega^{\infty} R)$. Moreover, suppose that $\mathbb{G}_{a,R}$ is preoriented by the element $e\in \pi_2(\Omega^{\infty} R)$. Then, by construction, the image of the coordinate $t$ under the natural composite
\begin{equation}\label{equation: preorientation of affine line and coordinate}
\pi_0\Gamma(\mathbb{G}_{a,R},\mathcal{O}) \tto \pi_0\Gamma(\Spec R\times BT,\mathcal{O}) \tto \pi_0\Gamma(\Spec R\times S^2,\mathcal{O}) \cong \pi_2(\Omega^\infty R)
\end{equation}
is again $e$.

In the sequel, if $R=\mathbb{Q}[\beta,\beta^{-1}]$, $|\beta|=2$, then we fix the preorientation of $\mathbb{G}_{a,R}$ given by the class $\beta \in \pi_2(\mathbb{A}^1_{R}(R))\cong \pi_2(\Omega^{\infty} R)$, see the paragraph after~\cite[Definition~3.6]{GM} or \cite[Definition~4.3.1]{Lur_Ell2}.
\end{ex}

\begin{ex}\label{example: greenlees}
The next example is inspired by the main construction of \cite{Greenlees_rational_elliptic}. Let $\mathbb{G}$ be a connected smooth $\mathbb{Q}$-group scheme (\emph{not} necessarily affine) of dimension $1$. We will write $\widehat{\mathbb{G}}$ for the corresponding formal group, i.e. the formal completion of $\mathbb{G}$ around the identity. Let us fix a \emph{coordinate} $t$ of $\widehat{\mathbb{G}}$, that is, an isomorphism 
$$t\colon \widehat{\mathbb{G}} \cong \widehat{\mathbb{G}}_a.$$ 
As in \cref{example: additive group}, let $R=\mathbb{Q}[\beta,\beta^{-1}]$, $|\beta|=2$. We define $S=\Spec(R)$ and $\mathbb{G}_R = \mathbb{G}\times_{\mathbb{Q}}R$. Then, $\mathbb{G}_R$ is a (strict) abelian group object in non-connective spectral stacks, i.e. $\mathbb{G}_R \in \Ab(\SpDM_S^{\nc,\flat,\qcs})$. We recall that the choice of the coordinate $t$ endows $\mathbb{G}_R$ with the canonical preorientation. Indeed, by \Cref{pre_or_depends_only_formal_cpl} we have equivalences
$$\mathrm{Pre}(\mathbb{G}) \simeq \mathrm{Pre}(\widehat{\mathbb{G}})\xrightarrow{t} \mathrm{Pre}(\widehat{\mathbb{G}}_a)$$
between the spaces of preorientations. Finally, we have the canonical preorientation of $\widehat{\mathbb{G}}_a$ from~\cref{example: additive group}.
%then $\widehat{\mathbb{G}}=\Spf(\mathcal{O}(\widehat{\mathbb{G}}))$ and the ring $\mathcal{O}(\widehat{\mathbb{G}})$ is local. Let $m_{\mathbb{G}}$ We fix a \emph{coordinate} $t\in $  
\end{ex}

\begin{ex}[Strict multiplicative group]\label{example: multiplicative group} Compare with~\cite[Construction~1.6.10]{Lur_Ell2} or \cite[Section~4]{GM}. We define a (strict) abelian group object $\mathbb{G}_m$ in spectral prestacks over $\mathbb{S}$ by
$$\mathbb{G}_m(R) = \Hom_{\CAlg}(\Sigma^{\infty}_+\mathbb{Z},R) \simeq \Hom_{\mathrm{CMon}(\Type)}(\mathbb{Z},\Omega^{\infty} R). $$
More precisely, the group structure is given by considering $\mathbb{G}_m(R)$ as a functor from $\mathrm{Lat}\to \Type$, sending $L$ to $\Hom_{\CAlg}(\Sigma^{\infty}_+ L,R).$ By the definition, $\mathbb{G}_m$ is corepresented by $\Sigma^{\infty}_+ \mathbb{Z}$. For an $E_\infty$-algebra~$R$, we define $\mathbb{G}_{m,R}$ as $\mathbb{G}_{m}\times \Spec(R)$.

If $R=KU$ is a (periodic) complex $K$-theory, then we fix a preorientation of $\mathbb{G}_{m,KU}$ as in \cite[Theorem~4.1]{GM} or, equivalently, as in \cite[Remark~4.3.8]{Lur_Ell2}.
\end{ex}

\begin{ex}[Deformation of $\mathbb{G}_m$ to $\mathbb{G}_a$]\label{example: deformation of Gm to Ga}
The next example seems especially interesting as it provides a deformation of $\mathbb{G}_m$ from \cref{example: multiplicative group} to $\mathbb{G}_a$ from \cref{example: additive group} over a rational base. Let $R_0=\mathbb{Q}[\varepsilon]$, $|\varepsilon|=0$ and consider the ring $H=R_0[t,(1-\varepsilon t)^{-1}]$. We equip the ring $H$ with the following structure of a (cocommutative) Hopf algebra. The comultiplication 
$$\nabla\colon H\tto H\otimes_{R_0}H =R_0[t_1,t_2,(1-\varepsilon t_1)^{-1},(1-\varepsilon t_2)^{-1}]$$
is induced by the formula $\nabla(t)= t_1+t_2-\varepsilon t_1 t_2$ and the inverse $\iota\colon H\to H  $ is given by 
$$\iota (t)=-t(1-\varepsilon t)^{-1}.$$ Let $\mathbb{G}=\Spec(H)$ denote the group scheme over $\Spec(R_0)$ given by the Hopf algebra $H$. Note that the fiber $\mathbb{G}_{0}$ over $\varepsilon=0$ is the additive group $\mathbb{G}_a$ and  $\mathbb{G}_{1}$ over $\varepsilon=1$ is the multiplicative group $\mathbb{G}_m$. Finally, we note that there is a canonical isomorphism $\widehat{\mathbb{G}}\cong \widehat{\mathbb{G}}_a$ as $\mathbb{G}$ is endowed with the distinguished coordinate $t$.

Similarly to \cref{example: additive group}, let $R=R_0[\beta,\beta^{-1}]$, $|\beta|=2$. We define $S=\Spec(R)$ and $\mathbb{G}_R = \mathbb{G}\times_{R_0}R$. Then, $\mathbb{G}_R$ is a (strict) abelian group object in non-connective spectral stacks, i.e. $\mathbb{G}_R \in \Ab(\SpDM_S^{\nc,\flat,\qcs})$. Again, the coordinate $t$ endows $\mathbb{G}_R$ with the canonical preorientation. Indeed, we have equivalences
$$\mathrm{Pre}(\mathbb{G}) \simeq \mathrm{Pre}(\widehat{\mathbb{G}})\xrightarrow{t} \mathrm{Pre}(\widehat{\mathbb{G}}_a)$$
between the spaces of preorientations. Finally, we have the canonical preorientation of $\widehat{\mathbb{G}}_a$ from~\cref{example: additive group} given by $\beta \in \pi_2(\mathbb{G}_a(R))\cong \pi_2(R)$.
\end{ex}

\subsection{Local systems on global spaces}\label{section: local systems on global spaces}
In this section for a global space $X$ and a coefficient system $\fcat A$ we construct a category $\LocSys^\Glo(X, \fcat A)$ of globally equivariant local systems on $X$ and study its properties. It serves a somewhat technical role similar to the category of naive $G$-spectra
$$\Sp^{nG} := \Fun(\Orb_G^\op, \Sp)$$
in equivariant homotopy theory, and is used as an intermediate step to construct and to study the categories of genuine equivariant and tempered local systems. Again, for the rest of this section $\fcat T$ denotes a global family of compact Lie groups.

We recall the notion of partial lax limits first, see e.g.~\cite{Berman24} or~\cite{LNP25} for a detailed account.

\begin{defn}\label{defintion: marked category} A \emdef{marked category} is a category $\fcat{I}$ along with a collection of edges $\fcat{I}^\dagger \subset \Map(\Delta^1, \fcat{I})$ which contains all equivalences and which is stable under composition. The edges from $\fcat{I}^\dagger$ are called \emdef{marked}. 
\end{defn}

\begin{defn}\label{definition: partial lax limit}
Let $F \colon (\fcat{I}, \fcat{I}^\dagger) \to \Cat$ be a diagram of categories indexed by a marked category $(\fcat{I},\fcat{I}^\dagger)$. We define its \emdef{partial left lax limit} 
$$\pllaxlim \left(\xymatrix{(\fcat{I},\fcat{I}^\dagger) \ar[r]^-{F} & \Cat}\right) $$
as the full subcategory of the left lax limit $\llaxlim_{\fcat{I}}F$ (see \cref{lax_limit}) spanned by the objects $(c_i \in F(i))_{i\in \fcat{I}}$ such that the structure map
$$\alpha_\phi \colon F(\phi)(c_i) \tto c_j$$
is an equivalence for each marked edge $\phi\colon i \to j$, i.e. $\phi \in \fcat{I}^\dagger$. The partial right lax limit is defined similarly.
\end{defn}

Now we are ready to give the main definition of the present paper.

\begin{defn}\label{definition: global local systems}
Let $X$ be a global homotopy type and let $\fcat A$ be a coefficient system on $\Orb_{\fcat T}$. We define the category of \emdef{globally equivariant local systems on $X$ with coefficients in $\fcat A$} as a partial left lax limit
$$\LocSys^\Glo_{\fcat T}(X, \fcat A) = \pllaxlim \left(\xymatrix{((\Orb_{\fcat T /X})^\op,\mathfrak{S}) \ar[r] & \Orb_{\fcat T}^\op \ar[r]^-{\fcat{A}} & \Prs^\LL}\right),$$
where $((\Orb_{\fcat T/X})^\op,\mathfrak{S})$ is a marked category where an arrow $\phi\colon T' \to T$ over $X$ is marked if and only if $\phi$ is full, i.e. $\phi$ induces a surjective homomorphism of isotropy groups for any choice of a point in $T'$.
\end{defn}
\begin{notation}\label{notation_for_LSGlo_allfamily}
We will omit the subscript and simply write  $\LocSys^\Glo(X, \fcat A)$ instead of $\LocSys^\Glo_{\fcat{T}}(X, \fcat A)$ if the family $\fcat{T}$ is clear from the context. In particular, in this section and \cref{section: genuine equivariant local system}, the omitted subscript means that we consider globally equivariant local systems with respect to the family of \emph{all} compact Lie groups. %we will denote $\LocSys^\Glo_{\fcat T}(X, \fcat A)$ simply by $\LocSys^\Glo(X, \fcat A)$.
\end{notation}

\begin{rem}\label{remark: informal global locsys}
Unwinding the definitions, one finds that an object of $\LocSys^\Glo_{\fcat T}(X, \fcat A)$ consists of the following data:
\begin{itemize}
\item an object $\mathcal L(x) \in \fcat A(BG)$ for every compact Lie group $G \in \fcat T$ and a map $x\colon BG \to X$;

\item a map
$$\alpha_{\phi}\colon \phi_{\fcat A}^*(\mathcal L(x)) \tto \mathcal L(x \circ \phi)$$
for a morphism $\phi\colon BH \to BG$ such that $\alpha_{\phi}$ is an equivalence if the morphism $\phi$ is full;

\item an equivalence $$\alpha_{\phi\circ \psi} \simeq \alpha_{\psi} \circ \psi_{\fcat A}^*(\alpha_\phi)$$
for a pair of composable morphisms $\psi, \phi \in \Orb_{\fcat T}$;

\item and higher compatibilities.
\end{itemize}
\end{rem}

\begin{rem}
The notation $\mathcal{L}(x)$ is not always convenient, so we will also use other notation like $\mathcal{L}^x$, $\mathcal{L}^{BG}$, or $\ev_x(\mathcal{L})$ for the same object in the sequel. We hope this ambiguity will not confuse the reader too much.
\end{rem}
\begin{rem}\label{rem_LSGlo_T_depends_on_SGlo_T}
Let $X$ be a global space. Since the construction of $\LocSys^\Glo_{\fcat T}(X, \fcat A)$ involves probing $X$ only by orbits from $\Orb_{\fcat T}$, the resulting category depends only on the restriction of $X$ to $\Type^\Glo_{\fcat T}$, see \Cref{notation: T-orbispaces}.
\end{rem}

\begin{prop}\label{LS_glo_colimits_pointwise}
Let $X \in \Type^\Glo$ be a global space. Then $\LocSys^\Glo_{\fcat T}(X, \fcat{A})$ is presentable and for $G \in \fcat T$ and each $x\colon BG \to X$ the evaluation functor
$$\ev_x \colon \LocSys^\Glo_{\fcat T}(X, \fcat{A}) \tto \fcat A(BG), \qquad \mathcal L \xymatrix{\ar@{|->}[r] &} \mathcal L(x)$$
preserves colimits.

\begin{proof}
This holds for a general partial left lax limit of a small diagram of presentable categories and colimit preserving functors. Indeed, by \cref{remark: lax and weighted}, the left lax limit
$$\fcat{W} = \llaxlim \left(\xymatrix{(\Orb_{\fcat T/X})^\op \ar[r] & \Orb_{\fcat T}^\op \ar[r]^-{\fcat{A}} & \Prs^\LL \ar[r] & \Cat}\right) $$
is a presentable category. We will show that the full subcategory $\LocSys_{\fcat T}^{\Glo}(X,\fcat{A})\subset \fcat{W}$ is also presentable and closed under colimits. Let $\phi\colon T \twoheadrightarrow T_0$ be a full morphism over $X$ and let $F_{\phi}\colon\fcat{W} \to \fcat{A}(T)^{\Delta^1}$ be the functor given by the formula
$$F_{\phi}(\mathcal{L}) = \phi^*_{\fcat{A}}(\mathcal{L}^{T_0})\tto \mathcal{L}^{T}. $$
Note that $F_{\phi}$ preserves colimits. Moreover, let $\fcat{C}_{\phi}\subset \fcat{W}$  be the preimage of the full subcategory of constant arrows $\fcat{A}(T)\subset \fcat{A}(T)^{\Delta^1}$ under $F_\phi$. Then $\fcat{C}_\phi$ is a presentable full subcategory of $\fcat{W}$ closed under colimits, see~\cite[Proposition~5.5.3.13]{Lur_HTT}. Finally, since 
$$\LocSys_{\fcat T}^{\Glo}(X,\fcat{A}) \simeq \bigcap_{\text{$\phi$ is full}} \fcat{C}_{\phi} \subset \fcat{W},$$
we observe that $\LocSys_{\fcat T}^{\Glo}(X,\fcat{A})$ is a presentable full subcategory of $\fcat{W}$ which is closed under colimits again by~\cite[Proposition~5.5.3.13]{Lur_HTT}.
\end{proof}
\end{prop}

\begin{construction}
Let $f\colon X \to Y$ be a morphism of global spaces. Then the postcomposition with $f$ induces a map of diagrams
$$\Orb_{\fcat T/X} \tto \Orb_{\fcat T/Y}$$
which after passage to the partial left lax limits induces a colimit-preserving functor
$$f^{\Glo,*}\colon \LocSys_{\fcat T}^\Glo(Y, \fcat A) \tto \LocSys_{\fcat T}^\Glo(X, \fcat A).$$
Note that the functor $f^{\Glo,*}$ is exact (resp. strict $E_n$-monoidal, etc.) if the coefficient system $\fcat{A}$ is stable (resp. $E_n$-monoidal, etc). The assignment 
$$X \mapsto \LocSys_{\fcat T}^\Glo(X, \fcat A), \;\; (f\colon X\to Y) \mapsto f^{\Glo,*}$$
glues into a functor
\begin{equation}\label{equation: global local system_pullback}
\LocSys_{\fcat T}^{\Glo,*}(-, \fcat A) \colon \Type^{\Glo, \op}_{\fcat T} \tto \Prs^\LL.
\end{equation}

By the adjoint functor theorem, the functor $f^{\Glo,*}$ admits a right adjoint
$$f^{\Glo}_*\colon \LocSys_{\fcat T}^\Glo(X, \fcat A) \tto \LocSys_{\fcat T}^\Glo(Y, \fcat A).$$
Note that the functor $f^{\Glo}_*$ is exact (resp. right lax $E_n$-monoidal, etc.) if the coefficient system $\fcat{A}$ is stable (resp. $E_n$-monoidal, etc). Similar to the pullback functor, we obtain the functor
$$\LocSys_{\fcat T*}^\Glo(-, \fcat A) \colon \Type^{\Glo}_{\fcat T} \tto \Prs^\RR, \;\; (f\colon X\to Y) \mapsto f^{\Glo}_*.$$
\end{construction}

\begin{prop}\label{prop_LSglo_preserves_limits}
The functor~\eqref{equation: global local system_pullback}
$$\LocSys_{\fcat T}^{\Glo,*}(-, \fcat A) \colon \Type^{\Glo, \op} \tto \Prs^\LL$$
preserves small limits.

\begin{proof}
Let $X_\bullet \colon \fcat I \to \Type^\Glo$ be a small diagram with colimit $X$. Since the full subcategory $\Orb_{\fcat T} \subset \Type^{\Glo}$ consists of representable objects, we obtain an equivalence
$$\Orb_{\fcat T/X} \simeq \colim_{i\in \fcat{I}} \Orb_{\fcat T/X_i}$$
of small categories. Therefore,
$$\llaxlim_{BG \in (\Orb_{\fcat T/X})^\op} \fcat A(BG) \simeq \lim_{i\in \fcat{I}^{\op}} \left(\llaxlim_{BG \in (\Orb_{\fcat T/X_i})^\op} \fcat A(BG)\right).$$
Hence it suffices to show that the full subcategory $\LocSys_{\fcat T}^\Glo(X, \fcat A)$ of the left hand side identifies with the full subcategory $\lim_i \LocSys_{\fcat T}^\Glo(X_i, \fcat A)$ of the right hand side. 

Indeed, let 
$$\mathcal L \in \llaxlim ((\Orb_{\fcat T/X})^\op \xymatrix{\ar[r]^-{\fcat A} &} \Prs^\LL)$$
be an object such that its restriction to $\llaxlim ((\Orb_{\fcat T/X_i})^\op \xrightarrow{\fcat A} \Prs^\LL)$ lies in the full subcategory $\LocSys_{\fcat T}^\Glo(X_i, \fcat A)$ for each $i \in \fcat I$. We will show that $\mathcal{L}$ lies in $\LocSys_{\fcat T}^\Glo(X, \fcat A)$. Let $\phi\colon BH \twoheadrightarrow BG \to X$ be a full morphism over $X$. Then, by \cref{definition: global local systems}, we need to show that the comparison map
$$\alpha_\phi\colon \phi_{\fcat A}^*\mathcal L(BH) \tto \mathcal L(BG)$$
is an equivalence. However, since $BG\in \Type^{\Glo}_{\fcat T}$ is a representable, the map $BG \to X$ factors through $X_i$ for some $i \in\fcat{I}$, and so, $\alpha_\phi$ can be considered as a comparison map in the lax limit $\llaxlim((\Orb_{\fcat T/X_i})^\op \xrightarrow{\fcat A} \Prs^\LL)$, where $\alpha_{\phi}$ is an equivalence by assumption.
\end{proof}
\end{prop}
In particular, we have the following limit presentation of $\LocSys_{\fcat T}^\Glo(-, \fcat A)$.
\begin{cor}\label{cor_LSGlo_LKan_extended}
Let $X \in \Type^\Glo$ be a global space. Then the natural map
$$\LocSys_{\fcat T}^\Glo(X, \fcat A) \tto \lim_{T \in (\Orb_{\fcat T/X})^\op} \LocSys_{\fcat T}^\Glo(T, \fcat A)$$
is an equivalence.

\begin{proof}
By \Cref{rem_LSGlo_T_depends_on_SGlo_T}, the functor 
$$\LocSys^\Glo_{\fcat T}(-, \fcat A) \colon \Type^\Glo \tto \Prs^\LL$$
is equivalent to the right Kan extension of its restriction to $\Type^\Glo_{\fcat T}$. Since the embedding
$$\Type^\Glo_{\fcat T} \inj \Type^\Glo$$ preserves colimits by construction, the result follows from \Cref{prop_LSglo_preserves_limits} and the co-Yoneda lemma.
\end{proof}
\end{cor}

Further results require some assumptions on coefficient system $\fcat A$, see \Cref{assumpt_limit_prerserving_cs}.
\begin{prop}\label{propostion: *-push for global spaces}
Let $\fcat A\colon \Orb^{\op}_{\fcat{T}} \to \Prs^\LL$ be a limit-preserving coefficient system of Beck--Chevalley type. Then:
\begin{enumerate}
\item The pullback functor
$$f^{\Glo,*}\colon \LocSys_{\fcat T}^\Glo(Y, \fcat A) \tto \LocSys_{\fcat T}^\Glo(X, \fcat A)$$
preserves limits and admits a left adjoint $f^{\Glo}_{\#}\dashv f^{\Glo,*}$ for \emph{any} morphism $f\colon X \to Y$ of global spaces.

\item For a morphism $f\colon X \to Y$ of global spaces and $\mathcal L \in \LocSys_{\fcat T}^\Glo(X, \fcat A)$ the value of $f^{\Glo}_* \mathcal L$ on $y\colon BG \to Y$ is given by the limit
$$\lim_{BH\in \Orb_{\fcat T/BG\times_Y X}^\op} \fcat A_*(\phi)(\mathcal L(BH)),$$
where $\phi\colon BH \to BG\times_Y X \to BG$ denotes the structure map followed by the canonical projection.

\item For a fibered diagram
\[\xymatrix{
X^\prime \ar[r]^-q \ar[d]^g & X \ar[d]^f \\
Y^\prime \ar[r]^-p & Y
}\]
of $\fcat T$-global spaces the base change maps
$$q^{\Glo}_{\#} \circ g^{\Glo,*} \tto f^{\Glo,*} \circ p_{\#} , \qquad p^{\Glo,*} \circ f^{\Glo}_* \tto g^{\Glo}_* \circ q^{\Glo,*}$$
are equivalences.
\end{enumerate}

\begin{proof}
For the first part, we note that the limits in both categories $\LocSys_{\fcat T}^\Glo(X, \fcat A)$ and $\LocSys_{\fcat T}^\Glo(Y, \fcat A)$ are computed pointwise in their partial left lax limits definitions provided the assumption on the coefficient system $\fcat A$. Therefore, the pullback $f^{\Glo,*}$ preserves limits, and so, the left adjoint $f^{\Glo}_{\#}$ exists by the adjoint functor theorem.

Let
$$f^\flat\colon \Orb_{\fcat T/X}^\op \tto \Orb_{\fcat T/Y}^\op$$
denotes the functor induced by the composition with $f$. By \Cref{right_adjoint_llax_limits}, for an object
$$\mathcal L \in \llaxlim (\Orb^\op_{\fcat T/X} \xymatrix{\ar[r]^-{\fcat A} &} \Prs^\LL)$$
the value of the right adjoint to the functor
\begin{equation}\label{eq_global_pw_with_assumptions}
\llaxlim_{BG\in \Orb^\op_{\fcat T/Y}} \fcat A(BG) \tto \llaxlim_{BG \in \Orb^\op_{\fcat T/X}} \fcat A(BG)
\end{equation}
induced by $f^\flat$, is given by the limit
$$\lim_{BH\in \Orb_{\fcat T/BG\times_Y X}^\op} \fcat A_*(\phi)(\mathcal L(BH)).$$
By abuse of notation, let us denote this right adjoint also by $f^{\Glo}_*$. We will show that $f^{\Glo}_*$ preserves the subcategory given by the partial lax limit condition of \Cref{definition: global local systems}. To see this we need to show that the natural map
$$\fcat A(\psi)^*\left((f^{\Glo}_*\mathcal L)({BG})\right) \tto (f^{\Glo}_*\mathcal L)({BG^\prime})$$
is an equivalence for any full morphism $\psi\colon BG^\prime \surj BG$ over $Y$ and any $\mathcal{L} \in \LocSys_{\fcat T}^{\Glo}(X,\fcat{A})$. Under the previous identification of $f^{\Glo}_*$ and using our assumptions on $\fcat A$ we find
\begin{align*}
\fcat A^*(\psi)\left((f^{\Glo}_*\mathcal L)(BG)\right) & \simeq \fcat A^*(\psi)\left(\lim_{BH\in \Orb_{\fcat T/BG\times_Y X}^\op}
\fcat A_*(\phi)(\mathcal L(BH))\right) \\
&\simeq \lim_{BH\in \Orb_{\fcat T/BG\times_Y X}^\op} \fcat A_*(\phi^\prime) \fcat A^*(\psi^\prime)(\mathcal L(BH)) \\
&\simeq \lim_{BH\in \Orb_{\fcat T/BG\times_Y X}^\op} \fcat A_*(\phi^\prime)(\mathcal L(BH^\prime)) \\
&\simeq \lim_{BH^\prime\in \Orb_{\fcat T/BG^\prime\times_Y X}^\op} \fcat A_*(\phi^\prime)(\mathcal L(BH^\prime))\\ &\simeq (f^{\Glo}_*\mathcal L)({BG^\prime}),
\end{align*}
where $BH^\prime = BH\times_{BG} BG^\prime$ and $\phi^\prime\colon BH^\prime \to BG^\prime$ and $\psi^\prime\colon BH^\prime \to BH$ denote the natural projections. Above, in the second equivalence we use that $\fcat A^*(\phi)$ preserves small limits together with the Beck--Chevalley property, in the third equivalence we use that the map $\fcat A^*(\psi^\prime)(\mathcal{L}(BH)) \simeq \mathcal{L}(BH')$ is an equivalence for $\mathcal L \in \LocSys_{\fcat T}^\Glo(X, \fcat A)$, and in the second to the last equivalence we use that the functor
$$-\times_{BG} BG^\prime \colon \Orb_{\fcat T/BG\times_Y X} \tto \Orb_{\fcat T/BG^\prime\times_Y X}$$
is cofinal as the right adjoint to the forgetful functor.

%Then the functor induced by $f^\flat$ on the left lax limits
%\begin{equation}\label{eq_LSGlo_computing_pf_1}
%\lim^\llax_{\Orb_{/Y}^\op} \fcat A(-) \tto \lim^\llax_{\Orb_{/Y}^\op} \fcat A(-)
%\end{equation}
%coincides with the functor induced on the left lax limits of the map of diagrams
%$$\fcat A(-) \tto \Ran_{f^\flat} \fcat A(-).$$
%By \todo{ref} for $BG \in \Orb_{/Y}$ we have a fully faithful embedding
%$$(\Ran^{\llax}_{f^\flat} \fcat A)(BG) \xymatrix{\ar@{^(->}[r]&} \lim^{\llax}_{BH \in \Orb_{BG \times_Y X}^\op} \fcat A(BH).$$
%Note that it admits a right adjoint by the adjoint functor theorem. Moreover, the composition
%$$\xymatrix{\fcat A \ar[r] & \Ran^{\llax}_{f^\flat} \fcat A \ar[r] & \lim^{\llax}_{BH \in \Orb_{- \times_Y X}^\op} \fcat A(BH)}$$
%is a strict transformation of the diagrams \todo{why?} (but in general neither the first nor the second arrows are strict separately). It follows by \todo{ref} that the value on $BG$ of the right adjoint to the functor \eqref{eq_LSGlo_computing_pf_1} can be computed as a right adjoint to the functor
%$$\fcat A(BG) \tto \lim^{\llax}_{BH \in \Orb_{/BG\times_Y X}^\op} \fcat A(BH)$$
%and to prove the assertion we only left to show that this construction sends the full subcategory $\LocSys_{\fcat T}^\Glo(X, \fcat A)$ into $\LocSys_{\fcat T}^\Glo(Y, \fcat A)$ \todo{check}.

For the last assertion, it suffices to prove only the second equivalence, the first one will follow by passing to the left adjoints. By the previous part, the value of $p^{\Glo,*} f^{\Glo}_* \mathcal L$ at a point $y^\prime\colon BG \to Y^\prime$ is equivalent to the limit
$$\lim_{BH \in \Orb_{\fcat T/BG\times_{Y} X}^\op} \fcat{A}_*(\phi)(\mathcal L(BH))$$
and the value of $g^{\Glo}_* q^{\Glo,*} \mathcal{L}$ at the same point is equivalent to the limit
$$\lim_{BH \in \Orb_{\fcat T/BG\times_{Y^\prime} X^\prime}^\op} \fcat{A}_*(\phi)(\mathcal L(BH)).$$
But $BG \times_{Y^\prime} X^\prime \simeq BG \times_{Y^\prime} (Y^\prime \times_Y X) \simeq BG \times_{Y} X$.
\end{proof}
\end{prop}

We have the following variant of the previous proposition with weaker assumptions on the coefficient system $\fcat A$, if we restrict our attention to global orbits.
\begin{prop}\label{corollary: pushforward for orbits}
Let $\fcat A\colon \Orb^{\op}_{\fcat{T}} \to \Prs^\LL$ be a coefficient system of Beck--Chevalley type.  Suppose that $f\colon BG' \twoheadrightarrow BG$ is a surjective morphism of representable global spaces and $\mathcal L \in \LocSys_{\fcat T}^\Glo(BG', \fcat A)$ is a globally equivariant local system. Then the value of $f^{\Glo}_* \mathcal L$ at $y\colon BH \to BG$ is given by 
$$(f^{\Glo}_*\mathcal{L})(y) \simeq \fcat{A}_*(g)(\mathcal{L}(BH'))$$
where $H'=H\times_G G'$ and $g\colon BH' \twoheadrightarrow BH$ is the canonical projection.

\begin{proof}
By \Cref{proposition: full morphisms are closed under pullbacks}, the fiber product $BH\times_{BG} BG'$ of global spaces is equivalent to $BH'$. Therefore, arguing as in the proof of \cref{propostion: *-push for global spaces}, the right adjoint to
$$
\llaxlim_{T\in \Orb^\op_{\fcat T/BG}} \fcat A(T) \tto \llaxlim_{T\in \Orb^\op_{\fcat T/BG'}} \fcat A(T)
$$
is given by $\mathcal{L} \mapsto f^{\Glo}_*\mathcal{L}$, $(f^{\Glo}_*\mathcal{L})(y) \simeq \fcat{A}_*(g)(\mathcal{L}(BH'))$, see \cref{right_adjoint_llax_limits}. Then we only need to prove that $f^{\Glo}_*$ preserves the partial lax limits subcategories. The latter follows by the assumption that the coefficient system $\fcat{A}$ is of Beck--Chevalley type. 
\end{proof}
\end{prop}

\begin{rem}\label{remark: *-pushforward limit}
We note that \cref{corollary: pushforward for orbits} does not require that the coefficient system $\fcat{A}$ is limit-preserving in opposite to \cref{propostion: *-push for global spaces}.
\end{rem}

The following examples show that the pullback functor of globally equivariant local systems does \emph{not} necessarily preserve limits without additional assumptions on the morphism or a coefficient system, which is in contrast to the case of the usual local systems.

\begin{ex}\label{example: Ga does not have left adjoint}
Let $\mathbb{G}_a$ be a strict additive group over $\mathbb{Q}$, see \cref{example: additive group} and let $\fcat{A}\colon \Orb^{\op}_{\ab}\to \CAlg(\Prs^\LL)$ be the associated coefficient system of geometric type from \cref{example: coefficient system from preoriented abelian group}. Then, by \cref{example: locsysglo geometric type} below, we have
$$\LocSys^{\Glo}(BU(1),\fcat{A}) \simeq \Mod_{R}\left(\LocSys^{\Glo}(BU(1),\underline{\Mod_{\mathbb{Q}}})\right)\simeq \Mod_{R}\left(\Fun(\Orb_{U(1)}^{\op},\Mod_{\mathbb{Q}})\right) $$
for an $E_\infty$-ring $R\in \CAlg\left(\Fun(\Orb_{U(1)}^{\op},\Mod_{\mathbb{Q}})\right)$. Moreover, $R^{U(1)}\simeq \mathbb{Q}[t]$. Under this identification, the pullback functor 
$$p^{\Glo,*}\colon \Mod_{\mathbb{Q}} \simeq \LocSys_{\fcat T}^\Glo(*, \fcat A) \tto \LocSys_{\fcat T}^\Glo(BU(1), \fcat A)$$
along the projection $p\colon BU(1) \to *$ is given by
$$M \mapsto \underline M \otimes_{\underline{\mathbb{Q}}} R.$$
Here $\underline M\in\Fun(\Orb_{U(1)}^{\op},\Mod_{\mathbb{Q}})$ is the constant presheaf with value $M$. Since
$$R^{U(1)} \simeq \mathbb{Q}[t] \simeq \bigoplus_{\mathbb N} \mathbb{Q}$$
is not finite over $\mathbb{Q}$, the functor $p^{\Glo,*}$ does not preserve limits. 
\end{ex}

\begin{ex}\label{example: Gm does not have left adjoint}
The similar argument also covers the case of the strict multiplicative group $\mathbb{G}_m$ over $KU$, see \cref{example: multiplicative group}. Let $\fcat{A}\colon \Orb^{\op}_{\ab}\to \CAlg(\Prs^\LL)$ be the associated coefficient system of geometric type from \cref{example: coefficient system from preoriented abelian group}. Again, by \cref{example: locsysglo geometric type} below, we have 
$$\LocSys^{\Glo}(BU(1),\fcat{A}) \simeq \Mod_{KU_{U(1)}} \left(\Fun(\Orb_{U(1)}^\op, \Sp)\right),$$
where $KU_{U(1)}$ is the naive spectrum underlying the genuine $U(1)$-equivariant complex $K$-theory. Under this identification, the pullback functor 
$$p^{\Glo,*}\colon \Mod_{KU} \simeq \LocSys_{\fcat T}^\Glo(*, \fcat A) \tto \LocSys_{\fcat T}^\Glo(BU(1), \fcat A) \simeq \Mod_{KU_{U(1)}} \Fun(\Orb_{U(1)}^\op, \Sp)$$
along the projection $p\colon BU(1) \to *$ is given by
$$M \mapsto \underline M \otimes_{\underline{KU}} KU_{U(1)},$$
where $\underline M$ is the constant presheaf with value $M$. Since
$$(KU_{U(1)})^{U(1)} \simeq \Gamma(\mathbb G_{m, KU}, \mathcal O) \simeq KU[\mathbb Z] \simeq \bigoplus_{\mathbb Z} KU$$
is not finite over $KU$, the functor $p^{\Glo,*}$ does not preserve limits.
\end{ex}

\begin{rem}\label{LSGlo_no_bs_in_general}
Let
\[\xymatrix{
X^\prime \ar[r]^-q \ar[d]^g & X \ar[d]^f \\
Y^\prime \ar[r]^-p & Y
}\]
be a fibered diagram of global spaces. Then the base change map
$$p^{\Glo,*} \circ f^{\Glo}_* \tto g^{\Glo}_* \circ q^{\Glo,*}$$
is \emph{not} an equivalence in general. Indeed, set $X = \coprod_I Y$ and let $f$ be the fold map. Then, by \Cref{prop_LSglo_preserves_limits}, we have
$$\LocSys_{\fcat T}^\Glo(X, \fcat A) \simeq \prod_I \LocSys_{\fcat T}^\Glo(Y, \fcat A)$$
and the pushforward $f_*$ is equivalent to the $I$-indexed product. Therefore, the base change property would imply that the pullback $p^{\Glo,*}$ preserves all $I$-indexed products (and hence all small limits). However, \cref{example: Gm does not have left adjoint} shows that this is not always the case.
\end{rem}

\subsection{Local systems on orbispaces}\label{section: local systems on orbispaces}
The theory of globally equivariant (and even more to it genuine and tempered as we will se below) local systems behaves better when restricted to the full subcategory of orbispaces discussed in \Cref{section: orbispaces}. Mostly, this is due to the technical fact, that for an orbispace $X$, in contrast with the case of general global spaces, one can represent $\LocSys^\Glo(X, \fcat A)$ as a lax limit without any partial conditions, see \Cref{LSGlo_on_orbispaces} below.

First, if $X \in \Type^{\Glo}$ is an orbispace, then the limit presentation of $\LocSys_{\fcat T}^\Glo(X, \fcat A)$ from \Cref{cor_LSGlo_LKan_extended} can be simplified.
\begin{prop}\label{proposition: LSglo for orbispaces}
Let $X\in \Type^{\Glo}$ be an orbispace. Then the natural functor
$$\LocSys_{\fcat T}^\Glo(X, \fcat A) \tto \lim_{BG \in (\Orb_{\fcat T/^\rep X})^\op} \LocSys_{\fcat T}^\Glo(BG, \fcat A)$$
is an equivalence, where $\Orb_{\fcat T/^\rep X}$ denotes the full subcategory of $\Orb_{\fcat T/X}$ spanned by faithful morphisms\footnote{Note that morphism between object of $\Orb_{\fcat T/^\rep X}$ are automatically faithful by the $2$-out-of-$3$ property from \Cref{proposition: basic properties of faithful}(2).}.

\begin{proof}
Recall that the functor $\iota_!\colon \Type^\Orb_{\fcat T} \to \Type^\Glo_{\fcat T}$ is the left Kan extension along a (non-full) embedding $\iota\colon \Orb^\rep_{\fcat T} \to \Orb_{\fcat T}$. In particular, $\iota_!$ preserves colimits. Since $X$ is in the essential image of $\iota_!$, we have
$$X \simeq \colim_{BG \in \Orb_{/^\rep X}} BG,$$
and the assertion following by~\cref{prop_LSglo_preserves_limits}.
\end{proof}
\end{prop}

\begin{ex}\label{ex_LSglo_for_N}
\iffalse
Let $\mathcal N = \delta_!(*) \in \Type^\Glo$ be the normal subgroup classifier given by the rule
$$\mathcal N(BG) = \{\text{groupoid of normal subgroups of $G$ modulo conjugation}\},$$
see \cite[Section 4.1]{Rezk_GlobalHomotopy}.
We claim that
$$\LocSys^\Glo(\mathcal N, \fcat A) \simeq \lim_{BG \in \Orb^{\rep, \op}} \LocSys^\Glo(BG, \fcat A).$$
Indeed, by \cref{proposition: LSglo for orbispaces}, it is enough to show that
$$* \simeq \colim_{BG \in \Orb^\rep} BG \in \Type^{\Orb}.$$ However, $*$ is the final object of $\Type^\Orb$, so the assertion follows from a more general fact: for a small category $\fcat C$ and the Yoneda embedding $y_{\fcat C}\colon \fcat C \inj \Fun(\fcat C^\op, \Type)$, the colimit
$$\colim_{c\in \fcat C} y_{\fcat C}(c) \in \Fun(\fcat C^\op, \Type)$$
is equivalent to the final object.
\fi
Let $\mathcal{N} \in \Type^{\Glo}$ be the normal subgroup classifier from \cref{example: normal subgroup classifier}. Then, by \cref{example: normal subgroup classifier as a colimit}, we have $$\mathcal{N}\simeq \colim_{BG\in \Orb^{\rep}}BG.$$ Therefore, for the constant coefficient system $\underline{\Sp}\colon \Orb^\op \to \Prs^{\LL}$, we have
$$\LocSys^\Glo(\mathcal N, \underline{\Sp}) \simeq \lim_{BG \in \Orb^{\rep, \op}} \Sp^{nG},$$
where $\Sp^{nG} = \Fun(\Orb_G^\op, \Sp)$ are the stabilizations of the categories of $G$-spaces for various compact Lie groups $G$. We note that a similar limit category has been studied before in~\cite{LNP25}.
\end{ex}

Similarly, if $X\in \Type^{\Orb}_{\fcat T}$ is an orbispace, then the category $\LocSys_{\fcat T}^\Glo(X, \fcat A)$ admits a more compact description as a left lax limit.
\begin{prop}\label{LSGlo_on_orbispaces}
Let $X\in \Type^{\Orb}$ be an orbispace and denote by $\Orb_{\fcat T/^\rep X}$ the subcategory of $\Orb_{\fcat T/X}$ spanned by orbits with a faithful map to $X$. Then the restriction along the embedding
$$i\colon \Orb_{\fcat T/^\rep X} \xymatrix{\ar@{^(->}[r] &} \Orb_{\fcat T/X}$$ induces an equivalence
$$\LocSys_{\fcat T}^\Glo(X, \fcat A) \xymatrix{\ar[r]^-\sim&} \llaxlim_{BG \in (\Orb_{\fcat T/^\rep X})^\op} \fcat A(BG).$$
In particular, if the coefficient system $\fcat A \colon \Orb^{\op}_{\fcat T} \to \Prs^{\LL}$ is constant with the value $\fcat A_0$, then
$$\LocSys_{\fcat T}^\Glo(X, \fcat A) \simeq \Fun((\Orb_{\fcat T/^\rep X})^\op, \fcat A_0).$$

\begin{proof}
Note that the functor $i$ admits a left adjoint $p$ which sends $x\colon BG \to X$ to $B(G/I_X(x)) \inj X$. Here, $I_X(x)$ is the inertia group of the point $x$, see \cite[Proposition 4.3.1]{Rezk_GlobalHomotopy}. Hence the restrictions along $i$ and $p$ induce the adjoint pair
$$p^*\colon \xymatrix{\displaystyle\llaxlim_{BG \in (\Orb_{\fcat T/^\rep X})^\op} \fcat A(BG) \ar@<0.5ex>[r] & \ar@<0.5ex>[l] \displaystyle\llaxlim_{BG \in (\Orb_{\fcat T/X})^\op} \fcat A(BG)}  \colon i^*.$$
Moreover, since $i$ is fully faithful, so is $p^*$. We will show that the essential image of $p^*$ is $\LocSys_{\fcat T}^\Glo(X, \fcat A)$.

By  construction, if $\mathcal L \in \llaxlim({(\Orb_{\fcat T/^\rep X})^\op} \xrightarrow{\fcat A} \Prs^\LL)$, then the value of $p^*(\mathcal L)$ on a point $x\colon BG \to X$ is equivalent to the pullback $\phi^*_{\fcat A} \mathcal L(p(x))$ along the natural map 
$$\phi \colon BG \surj BG/I_X(x).$$
However, for a surjection $\psi\colon BH \to BG$ over $X$, the natural quotient map $p(x) \to p(x\circ \psi)$ is an equivalence. It follows that $p^*(\mathcal L)$ lies in $\LocSys_{\fcat T}^\Glo(X, \fcat A)$. The reverse inclusion $\LocSys_{\fcat T}^\Glo(X, \fcat A) \subseteq \Im(p^*)$ can be proven in a similar way.
\end{proof}
\end{prop}

\begin{ex}\label{example: LSGlo of N}
Let $X=\mathcal{N} \in \Type^{\Glo}$ be the normal subgroup classifier from \cref{example: normal subgroup classifier} and let $\fcat{A} = \underline{\Sp}$ be the constant coefficient system. Then the diagram
$$\xymatrix{(\Orb_{/^\rep \mathcal{N}})^\op \ar[r] & \Orb^{\op} \ar[r]^-{\underline{\Sp}} & \Prs^\LL}$$
is constant and $\Orb_{/^\rep \mathcal{N}} \simeq \Orb^{\rep}$. Therefore, by \cref{LSGlo_on_orbispaces}, we have
$$\LocSys^\Glo(\mathcal N, \underline{\Sp}) \simeq \Fun(\Orb^{\rep,\op},\Sp).$$
In particular, together with \cref{ex_LSglo_for_N}, we recover the equivalence
$$\Fun(\Orb^{\rep,\op},\Sp) \simeq \lim_{BG \in \Orb^{\rep,\op}} \Fun(\Orb^{\op}_G,\Sp). $$
\end{ex}

\begin{cor}\label{LSGlo_for_T_orbSp_without_family_is_fine}
Let $X$ be a $\fcat T$-global space. Then there is a natural equivalence
\begin{equation}\label{eq_with_T_or_without}
\LocSys^\Glo(X, \fcat A) \areq \LocSys^\Glo_{\fcat T}(X, \fcat A).
\end{equation}

\begin{proof}
The comparison functor \eqref{eq_with_T_or_without} is the one induced on partial left lax limits by the change of marked indexing diagrams
$$\Orb_{\fcat T/X} \tto \Orb_{/X}.$$
By \Cref{prop_LSglo_preserves_limits}, as functors on $\Type^\Glo_{\fcat T}$, both sided of \eqref{eq_with_T_or_without} are right Kan extended from $\Orb_{\fcat T}$. Hence it is enough to show that
$$\LocSys^\Glo(T, \fcat A) \areq \LocSys^\Glo_{\fcat T}(T, \fcat A)$$
for $T \in \Orb_{\fcat T}$. But in this case by definition of global family
$$\Orb_{\fcat T/^\rep T} \areq \Orb_{/^\rep T}.$$
The result then follows from \Cref{LSGlo_on_orbispaces}.
\end{proof}
\end{cor}
\begin{cor}\label{cor_limits_in_LS_Glo_pointwise_on_orbisp}
Let $X$ be an orbispace. Then, for a faithful point $x\colon BG \to X$, the evaluation functor
$$\ev_x\colon \LocSys_{\fcat T}^\Glo(X, \fcat A) \tto \fcat A(BG), \qquad \mathcal L \xymatrix{\ar@{|->}[r] &} \mathcal L(BG)$$
preserves limits.

\begin{proof}
This holds for a general left lax limit of a small diagram of categories admitting limits (even without the assumption that the transition functors in the diagram preserve limits).
\end{proof}
\end{cor}

By \Cref{cor_LSGlo_LKan_extended}, the categories $\LocSys_{\fcat T}^\Glo(BG,\fcat{A})$, $G$ is compact, are basic building blocks for the theory, therefore the following special case of \cref{LSGlo_on_orbispaces} is particularly important.
\begin{ex}\label{LS_glo_on_BG}
Let $G$ be a compact Lie group. Then the category $\Orb_{\fcat T/^\rep BG}$ is equivalent to $\Orb_{G, \fcat T}$, where the latter is the full subcategory of $\Orb_G$ spanned by $G$-orbits with stabilizers in $\fcat T$. Hence
$$\LocSys_{\fcat T}^\Glo(BG, \fcat A) \simeq \llaxlim_{G/H \in \Orb_{G, \fcat T}^\op} \fcat A(BH).$$
In particular, if the coefficient system $\fcat A=\underline{\fcat{A}}_0$ is constant with the value $\fcat A_0 \in \Prs^\LL$, then
$$\LocSys^\Glo_{\fcat T}(BG, \fcat A) \simeq \Fun(\Orb_{G, \fcat T}^\op, \fcat A_0).$$
In particular, $\LocSys^\Glo(BG, \underline{\Type})$ is the category of $G$-spaces and $\LocSys^\Glo(BG, \underline{\Sp})$ is the naive (or semi-naive) $G$-equivariant stable category.
\end{ex}

\begin{ex}[Change-of-coefficients]\label{example: locsysglo geometric type}
Let $G$ be a compact Lie group from our family $\fcat T$, and let $q\colon \Orb_{G} \to \Orb_{\fcat T}$ be the projection. Suppose that the coefficient system $\fcat{A}\colon \Orb^\op_{\fcat T} \to \CAlg(\Prs^\LL)$ is of geometric type (see \cref{definition: coefficient system is of geometric type}). Then there exists a natural transformation
$$\theta\colon \underline{\fcat{A}(BG)} \tto q^*\fcat{A}$$
of functors from $\Orb_{G, \fcat T}$ to $\CAlg(\Prs^\LL)$, where $\underline{\fcat{A}(BG)}$ is the constant functor with the value $\fcat{A}(BG)$. By \cref{LSGlo_on_orbispaces}, the natural transformation $\theta$ induces the functor
$$\theta^*\colon \LocSys^{\Glo}(BG,\underline{\fcat{A}(BG)}) \tto \LocSys^{\Glo}(BG,\fcat{A})$$
such that $(\theta^*(\mathcal{L}))^{\phi} \simeq \phi^*_{\fcat{A}}(\mathcal{L}^{\phi})$ for $\phi\colon BH \hookrightarrow BG$ and $\mathcal{L}\in \LocSys^{\Glo}_{\fcat T}(BG,\underline{\fcat{A}(BG)})$. 

Moreover, let $\theta^R\colon \fcat{A} \to \underline{\fcat{A}(BG)}$ be the right mate of the natural transformation $\theta$, see \cite{haugseng2021lax}. We note that $\theta^R$ is only a \emph{left lax} natural transformation. Therefore, $\theta^*$ admits a right adjoint
$$\theta_*\colon \LocSys^{\Glo}(BG,\fcat{A}) \tto \LocSys^{\Glo}(BG,\underline{\fcat{A}(BG)})$$
induced by the left lax natural transformation $\theta^R$, i.e. $(\theta_*(\mathcal{L}))^{\phi} \simeq \phi_{\fcat{A},*}(\mathcal{L}^{\phi})$ for $\phi\colon BH \hookrightarrow BG$ and $\mathcal{L}\in \LocSys^{\Glo}_{\fcat T}(BG, \fcat{A})$.

By construction, the functor $\theta^*$ is symmetric monoidal. Since $\fcat{A}$ is strongly continuous, we observe that the right adjoint $\theta_*$ commutes with colimits. Finally, since $\fcat{A}$ is of geometric type and both functors $\theta^*, \theta_*$ are pointwise, the functor $\theta_*$ is conservative and $\LocSys^{\Glo}(BG,\underline{\fcat{A}(BG)})$-linear. Therefore, by e.g.~\cite[Proposition~5.29]{MNN17}, there is a natural equivalence
$$\LocSys^{\Glo}(BG, \fcat{A}) \simeq \Mod_{\fcat{A}_G}\left(\LocSys^{\Glo}(BG,\underline{\fcat{A}(BG)})\right),$$
where $\fcat{A}_G=\theta_*\mathbbl{1}^{\Glo} \in \CAlg \LocSys^{\Glo}(BG,\underline{\fcat{A}(BG)})$. By construction, $\fcat{A}(BH)$ is the category of modules over $(\fcat{A}_G)^{BH}\in \CAlg(\fcat{A}(BG))$.
\end{ex}

\begin{cor}\label{corollary: global local systems G-spaces}
Let $G$ be a compact Lie group and $X\in \Type^G$ be a $G$-space. Then there is a natural equivalence
$$\LocSys^{\Glo}_{\fcat T}(X\gitq G,\fcat{A}) \simeq \llaxlim \left(\xymatrix{\Orb_{G, \fcat T}^{\op} \ar[r]^-{(X,\fcat{A})} & \Type \times \Prs^\LL \ar[rr]^-{\LocSys_{\#} \times \Id} && \Prs^\LL \times \Prs^\LL \ar[r]^-{\otimes} & \Prs^\LL}\right). $$
\end{cor}

\begin{proof}
Set $Y=X\gitq G$ and let $p\colon Y\to BG$ be the canonical faithful morphism, see \cref{remark: G-spaces}. Abusing notation, let 
$$p\colon (\Orb_{\fcat T/^\rep Y})^{\op} \tto (\Orb_{\fcat T/^\rep BG})^{\op}$$
denote the projection and let $F\colon (\Orb_{\fcat T/^\rep Y})^\op \to \Prs^\LL$ denote the composite
$$\xymatrix{(\Orb_{\fcat T/^\rep Y})^\op \ar[r] & \Orb^{\rep, \op}_{\fcat T} \ar[r]^-{\fcat{A}} & \Prs^\LL.}$$
Then, by \cref{LSGlo_on_orbispaces}, we have
$$\LocSys^{\Glo}_{\fcat T}(Y,\fcat{A}) \simeq \llaxlim_{(\Orb_{\fcat T/^\rep Y})^\op} F \simeq \llaxlim_{(\Orb_{\fcat T/^\rep BG})^\op} \Ran^{\llax}_p F.$$
Note that $p$ is a coCartesian fibration. So, by \cref{corollary: right oplax along cocartesian}, we have
\begin{align*}
\Ran^{\llax}_pF(G/H) &\simeq \rlaxlim\left(\xymatrix{\Hom_{/BG}(BH,X) \ar[r] & \Orb^{\rep} \ar[r]^-{\fcat{A}_*} & \Cat}\right)\\
&\simeq \LocSys(X(G/H),\fcat{A}(BH)).
\end{align*}
Finally, by loc.cit., the functoriality is given by the left adjoints to the coCartesian lifts. 
\end{proof}

\begin{ex}\label{example: locsysglo geom type general}
As in \cref{example: locsysglo geometric type}, let $G$ be a compact Lie group from $\fcat T$. Suppose $X$ is a $G$-space. Note that in this case by \Cref{LSGlo_for_T_orbSp_without_family_is_fine}
$$\LocSys^\Glo(X//G, \fcat B) \simeq \LocSys^\Glo_{\fcat T}(X//G, \fcat B)$$
for any coefficient system $\fcat B$. We will write $\LocSys_{\#}\circ X(-)\colon \Orb^\op_G \to \Prs^\LL$ for the composite
$$\LocSys_{\#}\circ X(-)\colon \xymatrix{\Orb^\op_G \ar[r]^-{X} & \Type \ar[r]^-{\LocSys_{\#}} & \Prs^\LL.}$$
Let $\fcat{A}$ be a coefficient system of geometric type and let $q\colon \Orb_G \to \Orb_{\fcat T}$ denote the projection. Since $\fcat{A}$ is strongly continuous, the adjoint pair
$$\theta \colon \underline{\fcat{A}(BG)} \xymatrix{\ar@<+0.5ex>[r] & \ar@<+0.5ex>[l]} q^*\fcat{A} \colon \theta^R$$
in $\Fun(\Orb_G^\op,\Prs^\LL)^{\llax}$ induces the adjoint pair
$$\Id \otimes \theta \colon \LocSys_{\#}\circ X(-) \otimes \underline{\fcat{A}(BG)} \xymatrix{\ar@<+0.5ex>[r] & \ar@<+0.5ex>[l]} \LocSys_{\#}\circ X(-) \otimes q^*\fcat{A} \colon \Id \otimes \theta^R. $$
So, by \cref{corollary: global local systems G-spaces}, we obtain the adjoint pair
$$\theta_X^*\colon \LocSys^{\Glo}(X \gitq G,\underline{\fcat{A}(BG)}) \xymatrix{\ar@<+0.5ex>[r] & \ar@<+0.5ex>[l]} \LocSys^{\Glo}(X\gitq G,\fcat{A}) \colon \theta_{X,*}.$$
Finally, as in \cref{example: locsysglo geometric type}, we can apply \cite[Proposition~5.29]{MNN17} and we obtain
$$\LocSys^{\Glo}(X\gitq G, \fcat{A}) \simeq \Mod_{\fcat{A}_X}\left(\LocSys^{\Glo}(X \gitq G,\underline{\fcat{A}(BG)})\right),$$
where $\fcat{A}_X=\theta_{X,*}\mathbbl{1}^{\Glo}_X \in \CAlg\left(\LocSys^{\Glo}(X \gitq G,\underline{\fcat{A}(BG)})\right)$. By construction, there is a canonical equivalence $\fcat{A}_X \simeq p^{\Glo,*}\fcat{A}_G$, where $p\colon X \gitq G\to BG$ is the canonical faithful morphism.
\end{ex}

Our next goal is to show that the assignment $X \mapsto \LocSys_{\fcat T}^\Glo(X, \fcat A)$ has some additional functoriality when restricted to the category of orbispaces and faithful morphisms.
\begin{prop}\label{prop_LSGlo_sharp_pushforward}
Let $f\colon X \to Y$ be a faithful morphism of orbispaces. Then the pullback functor
$$f^{\Glo,*}\colon \LocSys_{\fcat T}^\Glo(Y, \fcat A) \tto \LocSys_{\fcat T}^\Glo(X, \fcat A)$$
admits a left adjoint $f^{\Glo}_{\#}$ and for $G \in \fcat T$ and a faithful $y\colon BG \to Y$ and $\mathcal L \in \LocSys_{\fcat T}^\Glo(X, \fcat A)$ there is a natural equivalence
$$f^{\Glo}_{\#}(\mathcal L)^y \simeq \colim_{x\in \Hom_{/Y}(BG, X)} \mathcal L^x \in \fcat A(BG),$$
where $\Hom_{/Y}(BG, X)$ is considered as a groupoid. Moreover, if the coefficient system $\fcat A$ is $E_n$-monoidal, then the functor $f^{\Glo}_{\#}$ is strictly $\LocSys_{\fcat T}^\Glo(Y, \fcat A)$-linear.

\begin{proof}
By \Cref{LSGlo_on_orbispaces} and the equivalence $\Prs^{\LL,\op} \simeq \Prs^{\mathrm{R}}$, we have
$$\LocSys_{\fcat T}^\Glo(X, \fcat A) \simeq \rlaxlim_{BG \in \Orb_{\fcat T/^\rep X}} \fcat A_*(BG), \qquad \LocSys_{\fcat T}^\Glo(Y, \fcat A) \simeq \rlaxlim_{BG \in \Orb_{\fcat T/^\rep Y}} \fcat A_*(BG),$$
where the functor
$$\fcat A_* \colon \Orb_{\fcat T} \tto \Prs^{\RR}$$
is obtained from $\fcat A$ by passing to the right adjoints of transition functors. Under these identifications, the pullback functor $f^{\Glo,*}$ is induced by the morphism of diagrams
$$f^{\flat}\colon \Orb_{\fcat T/^\rep X} \tto \Orb_{\fcat T/^\rep Y}.$$
Since all the transition functors in the right adjoint diagrams $\fcat A_*$ preserve limits, so is the induced functor $f^{\Glo,*}$. By the adjoint functor theorem, the functor $f^{\Glo,*}$ admits a left adjoint $f^{\Glo}_{\#}$.

Moreover, we have the following description of $f^{\Glo}_{\#}$. We observe that
$$\rlaxlim_{\Orb_{\fcat T/^\rep X}} \fcat A_* \simeq \rlaxlim_{\Orb_{\fcat T/^\rep Y}} \Ran_{f^\flat}^{\rlax} \fcat A_*.$$
Then the functor $f^{\Glo,*}$ is induced by the morphism of diagrams
$$f^{*}_-\colon \fcat A_* \tto \Ran_{f^\flat}^{\rlax} \fcat A_*.$$
By \cref{lemma: right adjoint to oplax limit_lax is strict}, the left adjoint $f^{\Glo}_{\#}$ is equivalent to the functor induced on the right lax limits by the natural transformation
$$f_{\#, -} \colon \Ran_{f^\flat}^{\rlax} \fcat A_* \tto \fcat A_*$$
which is the left mate of the natural transformation $f^*_-$, i.e. $f_{\#,-}$ is obtained from $f^{*}_-$ by passage to the pointwise left adjoints. 

Note that $f^\flat$ is a Cartesian fibration and the strict fiber over $y\colon BG \to Y$ is the groupoid $\Hom_{/Y}(BG, X)$. Therefore, by \Cref{on_KanExtensions_of_GrothFib}, the value of $\Ran_{f^\flat}^{\rlax} \fcat A_*$ at $y$ is equivalent to the functor category $\fcat A(BG)^{\Hom_{/Y}(BG, X)}$ and the natural map of diagrams
$$f^{*}_-\colon \fcat A_* \tto \Ran_{f^\flat}^{\rlax} \fcat A_*$$
evaluated at $y$ is given by the diagonal map
$$f^{*}_y \colon \fcat A(BG) \tto \fcat A(BG)^{\Hom_{/Y}(BG, X)}.$$
Hence, $f_{\#, y}$ is equivalent to the colimit functor.

Now, assuming that the coefficient system $\fcat A$ is $E_1$-monoidal, the left adjoint $f^{\Glo}_{\#}$ of the monoidal functor $f^{\Glo,*}$ is a \emph{left lax} $\LocSys_{\fcat T}^\Glo(Y, \fcat A)$-linear functor. We show that $f^{\Glo}_{\#}$ is \emph{strictly} linear. Because of the pointwise description above for functors $f^{\Glo,*}$ and $f^{\Glo}_{\#}$, it suffices to check that the colimit functor
$$\colim_{\fcat I} \colon \Fun(\fcat I, \fcat C) \tto \fcat C$$
is $\fcat C$-linear for any small diagram $\fcat I$ and cocomplete monoidal category $\fcat{C}$ such that the tensor product preserves small colimits in each variable. The latter is clear.
\end{proof}
\end{prop}

\begin{lem}\label{lemma: base change orbispaces}
Let
\[\xymatrix{
X^\prime \ar[r]^-q \ar[d]^g & X \ar[d]^f \\
Y^\prime \ar[r]^-p & Y
}\]
be a fibered diagram of orbispaces and faithful morphisms. Then the base change maps
$$q^{\Glo}_{\#} \circ g^{\Glo,*} \tto f^{\Glo,*} \circ p^{\Glo}_{\#} , \qquad p^{\Glo,*} \circ f^{\Glo}_* \tto g^{\Glo}_* \circ q^{\Glo,*}$$
are equivalences.

\begin{proof}
The second equivalence follows from the first one by passage to the right adjoints, so it is enough to prove only the first equivalence. By the proof of \Cref{prop_LSGlo_sharp_pushforward}, all functors in question can be computed pointwise in this case. More precisely, it suffices to show that the commutative diagram
\[\xymatrix{
\fcat A(BG)^{\Hom_{/Y}(BG, X^\prime)} & \ar[l] A(BG)^{\Hom_{/Y}(BG, X)}\\
\fcat A(BG)^{\Hom_{/Y}(BG, Y^\prime)} \ar[u] & \ar[l] \ar[u] \fcat A(BG)
}\]
is vertically left adjointable for each $y\colon BG \to Y$. Since
$$\Hom_{/Y}(BG, X^\prime) \simeq \Hom_{/Y}(BG, Y^\prime) \times \Hom_{/Y}(BG, X),$$
this is straightforward.
\end{proof}
\end{lem}

\begin{rem}\label{remark: parametrized global systems}
Suppose that the coefficient system $\fcat{A}$ is symmetric monoidal and $\fcat{A}$ is of Beck--Chevalley type. Then \Cref{lemma: base change orbispaces} and \Cref{prop_LSGlo_sharp_pushforward} imply together that the functor 
$$\LocSys^{\Glo}(-,\fcat{A})\colon (\Type^\Orb)^\op \tto \CAlg(\Prs^\LL)$$
is a presentably symmetric monoidal $\Type^\Orb$-category, see~\cite[Definition~2.6.2.10]{MW25} or~\cite[Definition~2.15]{Cnossen23}.
\end{rem}

\begin{ex}\label{example: orbifamily and pushforward}
Let $X\in \Type^{\Glo}$ be an orbispace. Consider the canonical projection 
$$\pi\colon X_{\fcat{T}}=X\times_{\mathcal{N}}\mathcal{N}_{\fcat T} \tto X,$$
see \cref{definiton: classfier subgroups_family}. As in \cref{remark: family_classfier is left kan}, we have $X_{\fcat T} \times_X X_{\fcat{T}} \simeq X$. Therefore, by \cref{lemma: base change orbispaces}, we have $$\pi^{\Glo}_{\#} \circ \pi^{\Glo,*} \xrightarrow{\simeq} \Id \;\; \text{and} \;\;  \Id \xrightarrow{\simeq} \pi^{\Glo}_{*} \circ \pi^{\Glo,*}.$$
In particular, the functors 
$$\pi^{\Glo}_{\#}\colon \LocSys^{\Glo}(X_{\fcat T},\fcat{A}) \tto \LocSys^{\Glo}(X,\fcat{A}),$$ 
$$\pi^{\Glo}_{*}\colon \LocSys^{\Glo}(X_{\fcat T},\fcat{A}) \tto \LocSys^{\Glo}(X,\fcat{A}),$$ 
are fully faithful. Moreover, by the formula in \cref{prop_LSGlo_sharp_pushforward} and \cref{remark: faithful points family}, a local system $\mathcal{L} \in \LocSys^{\Glo}(X,\fcat{A})$ belongs to the essential image of $\pi^{\Glo}_{\#}$ if and only if $\mathcal{L}^x\simeq 0$ for all faithful points $x\colon BG \to X$ and $G\notin \fcat{T}$. Similarly,
$$(\pi^{\Glo,*}\mathcal{L})^x\simeq
\begin{cases}
\mathcal{L}^{\pi \circ x}, & \mbox{if $G \in \fcat{T}$,}\\
0, & \mbox{otherwise,}
\end{cases}
$$
where $\mathcal{L}\in \LocSys^{\Glo}(X,\fcat{A})$ is globally equivariant local system and $x\colon BG \to X_{\fcat{T}}$ is a faithful morphism. As usual, the right adjoint $\pi^{\Glo}_{*}$ is harder to describe, see e.g. \cref{example: right adjoint proper family}.
\end{ex}

\begin{cor}\label{corollary: finite G-space}
Let $\fcat{A} \colon \Orb^{\op}_{\fcat T} \to \Prs^{\LL}$ be a strongly continuous coefficient system and let $f\colon X \to BG$ be a finite $G$-space. Then the pushforward
$$f^{\Glo}_*\colon \LocSys^{\Glo}_{\fcat T}(X,\fcat{A}) \tto \LocSys^{\Glo}_{\fcat T}(BG,\fcat{A}) $$
preserves small colimits.
\end{cor}

\begin{proof}
Note that the pullback $X\times_{BG}BH$ is a finite $H$-space for any closed subgroup $H\subset G$, $H\in \fcat{T}$. Therefore, by \cref{lemma: base change orbispaces} and \cref{LS_glo_colimits_pointwise}, it suffices to show that the composite $$\ev_{BG}\circ f^{\Glo}_*\colon \LocSys^{\Glo}_{\fcat T}(X,\fcat{A}) \tto \fcat{A}(BG)$$
preserves small colimits for $G \in \fcat T$. This follows by \cref{prop_LSglo_preserves_limits} and \cref{example: star-pullback fixed points}.
\end{proof}

\begin{ex}\label{example: star-pullback fixed points}
Let $G$ be a compact Lie group from $\fcat T$. By the virtue of \Cref{cor_limits_in_LS_Glo_pointwise_on_orbisp}, the evaluation functor
$$\ev_{BG}\colon \LocSys^\Glo(BG,\fcat{A}) \tto \fcat{A}(BG) $$
preserves limits, and so, $\ev_{BG}$ admits a left adjoint
$$\triv^G\colon \fcat{A}(BG) \simeq \lim_{G/H \in \Orb_G^\op} \fcat A(BH) \xymatrix{\ar@{^(->}[r] &} \LocSys^\Glo(BG,\fcat{A})$$
which is the fully faithful embedding of the strict limit into the lax limit. In particular, we have $$(\triv^G(\mathcal{F}))^{y} \simeq y_{\fcat{A}}^{*}(\mathcal{F}) \in \fcat{A}(BK),$$ 
where $\mathcal{F} \in \fcat{A}(BG)$ and $y\colon BK \to BG$ is a faithful point. Note that the natural morphism
$$\triv^H \circ f^*_\fcat{A} \tto f^{\Glo,*} \triv^G$$
is an equivalence for any faithful morphism $f\colon BH \to BG$. Therefore, by taking right adjoints, we obtain the following equivalence
$$(f^{\Glo}_*\mathcal{L})^{BG} \simeq \fcat{A}_*(f)\mathcal{L}^{BH}, \;\; \mathcal{L} \in \LocSys^{\Glo}(BG,\fcat{A}). $$
\end{ex}

\begin{cor}\label{corollary: internal Hom colimits}
Let $f\colon X\to BG$ be a finite $G$-space. Suppose that the coefficient system $\fcat{A}$ is $E_1$-monoidal and strongly continuous. Then the endofunctor
$$[f^{\Glo}_{\#}\mathbbl{1}_X,-]\colon \LocSys^{\Glo}_{\fcat T}(BG,\fcat{A}) \tto \LocSys^{\Glo}_{\fcat T}(BG,\fcat{A}) $$
preserves small colimits.
\end{cor}

\begin{proof}
By \cref{prop_LSGlo_sharp_pushforward}, we have $[f^{\Glo}_{\#}\mathbbl{1}_X,-]\simeq f^{\Glo}_*f^{\Glo,*}(-)$. Finally, \cref{corollary: finite G-space} implies the assertion.
\end{proof}

\subsection{Functoriality along faithful morphisms for \texorpdfstring{$\LocSys^\Glo$}{LSGlo}}\label{subsection: functoriality along faithful morphisms}
In this section, based on the result of the previous section, we establish additional functoriality of $\LocSys^\Glo_{\fcat T}(-, \fcat A)$ for faithful morphism between general global spaces. Throughout this section, we assume that $\fcat{T}$ is a multiplicative global family. 

\begin{lem}\label{corollary: base change simple case glo}
Let
\[\xymatrix{
BH' \ar[r]^-q \ar[d]^g & BH \ar[d]^f \\
BG' \ar[r]^-p & BG
}\]
be a fibered diagram in $\Orb$ such that $f$ is faithful and $p$ is full. Suppose that the coefficient system $\fcat{A}$ is of Beck--Chevalley type. Then the base change maps
$$f^{\Glo,*} \circ p^{\Glo}_* \tto q^{\Glo}_* \circ g^{\Glo,*}, \qquad g^{\Glo}_{\#}\circ q^{\Glo,*} \tto p^{\Glo,*} \circ f^{\Glo}_{\#} $$
are equivalences.
\end{lem}

\begin{proof}
%By presenting the map $BG^\prime \to BG$ as a composite of a full morphism followed by a faithful morphism and by using \cref{lemma: base change orbispaces}, we can assume that the map $BG^\prime \to BG$ is full.
Let $y\colon BK \to BH$ be a faithful point of $BH$, $K\in \fcat{T}$ and set $K' = K\times_G G'$. Then, by \cite[Proposition~6.1.1]{Rezk_GlobalHomotopy}, the following commutative diagram
\[\xymatrix{
BK' \ar[r]^-{r} \ar[d]^x & BK \ar[d]^y \\
BH' \ar[r]^-q \ar[d]^g & BH \ar[d]^f \\
BG' \ar[r]^-p & BG
}\]
consists of fibered squares. Then, by applying \cref{corollary: pushforward for orbits} twice, we have
\begin{align*}
(f^{\Glo,*} p^{\Glo}_*\mathcal{L})(y) \simeq (p^{\Glo}_*\mathcal{L})(f\circ y) \simeq \fcat{A}_*(r)(\mathcal{L}(x\circ g)) \simeq \fcat{A}_*(r)(g^{\Glo,*}\mathcal{L}(x)) \simeq (q^{\Glo}_*g^{\Glo,*}\mathcal{L})(y)
\end{align*}
for every $\mathcal{L} \in \LocSys^\Glo(BG', \fcat A)$.
\end{proof}

%We finish this section by showing that the assumptions for the base change from the previous lemma can be considerably weakened.
\begin{prop}\label{LSGlo_basechange_faithful_pbs_ra}
Let
\[\xymatrix{
X^\prime \ar[r]^-q \ar[d]^g & X \ar[d]^f \\
Y^\prime \ar[r]^-p & Y
}\]
be a fibered diagram of global spaces such that $f$ and $g$ are faithful. Suppose that the coefficient system $\fcat{A}$ is of Beck--Chevalley type. Then the base change map
$$f^{\Glo,*} \circ p^{\Glo}_* \tto q^{\Glo}_* \circ g^{\Glo,*}$$ %, \qquad g^{\Glo}_{\#}\circ q^{\Glo,*} \tto p^{\Glo,*} \circ f^{\Glo}_{\#} $$
is an equivalence.
\begin{proof}
Consider the following commutative square
\begin{equation}\label{equation: LSGlo_basechange eq1}
\xymatrix{\LocSys^\Glo_{\fcat T}(X^\prime, \fcat A) & \ar[l]_-{q^{\Glo,*}} \LocSys^\Glo_{\fcat T}(X, \fcat A) \\
\LocSys^\Glo_{\fcat T}(Y^\prime, \fcat A) \ar[u]^{g^{\Glo,*}} & \LocSys^\Glo_{\fcat T}(Y, \fcat A). \ar[l]_-{p^{\Glo,*}} \ar[u]^{f^{\Glo,*}}
}\end{equation}
We will show that~\eqref{equation: LSGlo_basechange eq1} is horizontally right adjointable. %Then, by \cref{corollary: left adjoint faithful global},~\eqref{equation: LSGlo_basechange eq1} is also left adjointable with respect to vertical arrows. 
By \cite[Corollary~4.7.4.18]{Lur_HA}, a small limit of adjointable squares is adjointable. So, by presenting the global spaces $X, Y$ and $Y^\prime$ as colimits of representable global spaces and by \Cref{prop_LSglo_preserves_limits}, we can assume that $X = BH$, $Y = BG$ and $Y^\prime = BG'$, where $H \subseteq G$ is some subgroup. Moreover, by presenting the map $p\colon BG^\prime \to BG$ as a composite of a full morphism followed by a faithful morphism and by using \cref{lemma: base change orbispaces}, we can assume that the map $p\colon BG^\prime \to BG$ is full. Finally, the assertion follows by \cref{corollary: base change simple case glo}.
\end{proof}
\end{prop}

\begin{cor}\label{corollary: left adjoint faithful global}
Let $f\colon X \to Y$ be a faithful morphism of global spaces. Suppose that the coefficient system $\fcat{A}$ is of Beck--Chevalley type. Then the pullback functor
$$f^{\Glo,*}\colon \LocSys^\Glo_{\fcat T}(Y, \fcat A) \tto \LocSys^\Glo_{\fcat T}(X, \fcat A)$$
admits a left adjoint $f^{\Glo}_{\#}$.

\begin{proof}
Let us present $Y\in \Type^\Glo$ as a colimit $Y\simeq \colim_{\alpha \in \fcat{I}} T_{\alpha}$ of representable global spaces. Set $X_\alpha \simeq X\times_Y T_{\alpha}$ and $f_\alpha\colon X_\alpha \to T_\alpha$ be the faithful projection map, $\alpha \in \fcat{I}$. Then, by the virtue of \cref{cor_LSGlo_LKan_extended} and \cref{prop_LSGlo_sharp_pushforward}, the pullback functor $f^{\Glo,*}$ is a limit $\lim_{\alpha \in \fcat{I}^\op} f^{\Glo,*}_\alpha$. Moreover, by \cref{LSGlo_basechange_faithful_pbs_ra}, the diagram
$$f^{\Glo,*}_\bullet\colon \fcat{I}^\op \tto \Fun(\Delta^1,\Cat), \;\; \alpha \mapsto f^{\Glo,*}_{\alpha}$$
factors through the subcategory $\Fun^{\mathrm{LAd}}(\Delta^1,\Cat)$ of right adjoints and left adjointable natural transformations, see~\cite[Definition~4.7.4.16]{Lur_HA}. Finally, by \cite[Corollary~4.7.4.18]{Lur_HA}, a limit of right adjoint functors along left adjointable natural transformations is a right adjoint, i.e. $f^{\Glo,*}\simeq \lim_{\alpha} f^{\Glo,*}_\alpha$ is a right adjoint.
\end{proof}
\end{cor}

\begin{cor}\label{LSGlo_basechange_faithful_pbs}
Let
\[\xymatrix{
X^\prime \ar[r]^-q \ar[d]^g & X \ar[d]^f \\
Y^\prime \ar[r]^-p & Y
}\]
be a fibered diagram of global spaces such that $f$ and $g$ are faithful. Suppose that the coefficient system $\fcat{A}$ is of Beck--Chevalley type. Then the base change map
%$$f^{\Glo,*} \circ p^{\Glo}_* \tto q^{\Glo}_* \circ g^{\Glo,*}$$ 
$$g^{\Glo}_{\#}\circ q^{\Glo,*} \tto p^{\Glo,*} \circ f^{\Glo}_{\#} $$
is an equivalence. \qed
\end{cor}

\begin{rem}
The result above is the best we can get without additional assumptions on the coefficient system $\fcat A$ or the fibered square of global spaces. Even the assumption that $f$ is faithful will not help, since the projection $\coprod X \to X$ is faithful, but according to \Cref{LSGlo_no_bs_in_general} the base change property does not hold in this case.
\end{rem}

\begin{cor}\label{corollary: left adjoint faithful global linear}
Let $f\colon X \to Y$ be a faithful morphism of global spaces. Suppose that the coefficient system $\fcat{A}$ is $E_1$-monoidal and is of Beck--Chevalley type. Then the $\#$-pushforward
$$f^{\Glo}_{\#}\colon \LocSys^\Glo_{\fcat T}(X, \fcat A) \tto \LocSys^\Glo_{\fcat T}(Y, \fcat A)$$
is strict $\LocSys^\Glo_{\fcat T}(Y, \fcat A)$-linear.
\end{cor}

\begin{proof}
By construction, the functor $f^{\Glo}_{\#}$ is left lax $\LocSys^\Glo_{\fcat T}(Y, \fcat A)$-linear, i.e.\@ there exists a canonical morphism 
$$f^{\Glo}_{\#}(f^{\Glo,*}\mathcal{L} \otimes \mathcal{M}) \tto \mathcal{L} \otimes f^{\Glo}_{\#}(\mathcal{M}). $$
To check that this morphism is an equivalence by \cref{cor_LSGlo_LKan_extended} and \cref{LSGlo_basechange_faithful_pbs} we can assume that $X$ and $Y$ are global orbits. In this case the result is already established in \cref{prop_LSGlo_sharp_pushforward}.
\end{proof}

\begin{cor}\label{corollary: projection formula global Hom}
Let $f\colon X\to Y$ be a faithful morphism of global spaces. Suppose that the coefficient system $\fcat{A}$ is $E_1$-monoidal and is of Beck--Chevalley type. Then the natural transformation
$$[f_{\#}^{\Glo}(\mathcal{L}),\mathcal{M}] \tto f^{\Glo}_*[\mathcal{L},f^{\Glo,*}(\mathcal{M})]$$
is an equivalence for all $\mathcal{L} \in \LocSys^{\Glo}_{\fcat T}(X,\fcat{A})$ and $\mathcal{M} \in \LocSys^{\Glo}_{\fcat T}(Y,\fcat{A})$. \qed
\end{cor}

\begin{cor}[Relative K\"unneth formula]\label{cor_LSGlo_Kunneth}
Let $p\colon X \to B$ and $q\colon Y \to B$ be a pair of faithful morphisms of global spaces and let $r\colon X \times_B Y \to B$ be the fiber product. Then there is a natural equivalence
$$r^{\Glo}_{\#}(\mathcal L \boxtimes \mathcal M) \simeq p^{\Glo}_{\#}(\mathcal L) \otimes q^{\Glo}_{\#}(\mathcal M)$$ 
for every $\mathcal L \in \LocSys^\Glo_{\fcat T}(X, \fcat A)$, $\mathcal M \in \LocSys^\Glo_{\fcat T}(Y, \fcat A)$. 
\qed
\end{cor}

\begin{rem}
Let $H,G$ be a pair of non-trivial compact Lie groups and $\fcat A = \underline{\Sp}$ be a constant coefficient system. Then, by \Cref{LS_glo_on_BG}, we have
$$\LocSys^\Glo(BH, \fcat A) \simeq \Fun(\Orb_H^\op, \Sp), \qquad \LocSys^\Glo(BG, \fcat A) \simeq \Fun(\Orb_G^\op, \Sp)$$
so
$$\LocSys^\Glo(BH, \fcat A) \otimes \LocSys^\Glo(BG, \fcat A) \simeq \Fun(\Orb_H^\op \times \Orb_G^\op, \Sp).$$
Since $\Orb_H \times \Orb_G \not\simeq \Orb_{H\times G}$, we deduce that the categorical K\"{u}nneth formula fails for globally equivariant local systems, i.e.
$$\LocSys^\Glo(BH, \fcat A) \otimes \LocSys^\Glo(BG, \fcat A) \not\simeq \LocSys^\Glo(BH\times BG, \fcat A).$$
\end{rem}

\subsection{Enhanced global sections functor}\label{section: global section}
For a global homotopy type $X$, it seems natural to define the global sections functor
$$\LocSys^\Glo_{\fcat T}(X, \fcat A) \tto \LocSys^\Glo_{\fcat T}(*, \fcat A) \simeq \fcat A(*)$$
as the $*$-pushforward $p^{\Glo}_*$ along the structure map $p\colon X \to *$. However, we will show that $p^{\Glo}_*$ factors naturally through an intermediate category $\fcat A(X)$ (see~\Cref{LSglo_enhacned_GS}) which has more information, especially, if the coefficient system $\fcat A$ is non-affine.
\begin{construction}\label{LSglo_enhacned_GS}
Let $\fcat{A}\colon \Orb^{\op}_{\fcat{T}} \to \Prs^\LL$ be a coefficient system. Abusing notation, we define the functor
$$\fcat{A} \colon  (\Type^{\Glo})^{\op} \tto \Prs^\LL $$
as the right Kan extension of $\fcat{A}$ along the Yoneda embedding $\Orb_{\fcat{T}} \hookrightarrow \Type^{\Glo}$. In particular,
$$\fcat A(X) \simeq \lim_{T \in (\Orb_{\fcat T/X})^\op} \fcat A(T) \in \Prs^\LL$$
for a global space $X \in \Type^{\Glo}$.
\end{construction}

\begin{rem}\label{remark: global sections}
In other words, if the coefficient system $\fcat{A}\colon \Orb^{\op}_{\fcat{T}} \to \CAlg(\Prs^\LL)$ is symmetric monoidal, then the extension $\fcat{A}\colon  (\Type^{\Glo})^{\op} \tto \Prs^\LL$ is a global 2-ring associated with the coefficient system $\fcat{A}$, see \cref{remark: geometric and pregenuine}.
\end{rem}

\begin{defn}\label{LSglo_enhached_GS2}
We denote by
$$\mathbbl 1_{(-)}^\Glo\otimes - \colon \fcat A(-) \tto \LocSys^\Glo_{\fcat T}(-, \fcat A)$$
the natural transformation between the functors $\fcat{A}(-), \LocSys^{\Glo}_{\fcat T}(-,\fcat{A})\colon (\Type^{\Glo})^{\op} \to \Prs^{\LL}$ induced by the fully faithful embedding of the strict limit into the partial left lax limit, see \cref{definition: global local systems}. If the coefficient system $\fcat{A}$ is $E_n$-monoidal for some $0\leq n \leq \infty$, then the natural transformation $\mathbbl 1_{(-)}^\Glo\otimes -$ is (strict) $E_n$-monoidal as well.
\end{defn}
\begin{defn}\label{LSglo_enhached_GS3}
Let $X \in \Type^{\Glo}$ be a global space. We define the \emdef{enhanced global sections functor}
$$\Gamma_{\fcat A}(X, -) \colon \LocSys^\Glo_{\fcat T}(X, \fcat A) \tto \fcat A(X)$$
as the right adjoint to $\mathbbl 1_X^\Glo \otimes -$, see \cref{LSglo_enhached_GS2}. Again, if the coefficient system $\fcat{A}$ is $E_n$-monoidal for some $0 \leq n \leq \infty$, then $\Gamma_{\fcat A}(X, -)$ is right lax $E_n$-monoidal.
\end{defn}
\begin{ex}\label{example: global_ench_BG}
Let $X = BG$, where $G$ is a compact Lie group from $\fcat T$. By the definition, we have a natural equivalence 
$$\triv^G \simeq \mathbbl 1_{BG}^{\Glo} \otimes - ,$$
see \cref{example: star-pullback fixed points}. Therefore, the enhanced global section functor
$$\Gamma_{\fcat A}(BG, -) \colon \LocSys^\Glo_{\fcat T}(BG, \fcat A) \tto \fcat A(BG)$$
is the natural projection onto $\fcat A(BG)$, i.e. $\Gamma_{\fcat A}(BG, \mathcal L) \simeq \mathcal L^{BG},$
see \cref{LSGlo_on_orbispaces}.
\end{ex}
\begin{rem}
By \cref{LSglo_enhached_GS3}, the enhanced global sections commute with $*$-pushforwards. More precisely, for every morphism $f\colon X \to Y$ of global spaces, the diagram
\[\xymatrix{
\LocSys^\Glo_{\fcat T}(X, \fcat A) \ar[rr]^-{\Gamma_{\fcat A}(X, -)} \ar[d]^{f^{\Glo}_*} && \fcat A(X) \ar[d]^{\fcat A_*(f)} \\
\LocSys^\Glo_{\fcat T}(Y, \fcat A) \ar[rr]^-{\Gamma_{\fcat A}(Y, -)} && \fcat A(Y)
}\]
commutes. In particular, if $Y = *$ is a point, we observe that the pushforward $\Gamma(X, -) := p^{\Glo}_*$, $p\colon X \to *$ factors as a composition of $\Gamma_{\fcat A}(X, -)$ followed by $\fcat A_*(f)$. Here we use the equivalence $\LocSys^\Glo_{\fcat T}(*, \fcat A) \areq \fcat A(*)$.
\end{rem}
\begin{ex}
Let $A$ be a preoriented abelian group non-connective spectral stack and consider the coefficient system $\fcat A(BG) := \QCoh(A[\widehat G])$, see \cref{example: coefficient system from preoriented abelian group}. Then, for a global space $X\in \Type^\Glo$, the category $\fcat A(X)$ identifies with $\QCoh(A(X))$, where $A(X)$ is the spectral stack of $A$-tempered cohomology of $X$, see~\cite[Section~6]{GM}. Hence, if $\mathcal L \in \LocSys^\Glo_\ab(X, A)$ is a globally equivariant local system, then the enhanced global sections $\Gamma_A(X, \mathcal L)$ is a quasi-coherent sheaf on $A(X)$ while $\Gamma(X, \mathcal L)$ is merely a spectrum of the global sections (in the quasi-coherent sense) of $\Gamma_A(X, \mathcal L)$. For instance, if $A = E$ is an oriented spectral elliptic curve and $X=BU(1)$, the naive global section functor $\Gamma(BU(1), -)$ loses some information. As usual, if $A$ is affine (e.g.\@ $A = \mathbb G_{m, KU}$), then the difference is not essential.
\end{ex}

We will discuss \cref{example: global_ench_BG} with more details. Let $G$ be a compact Lie group from $\fcat T$. Then, by \cref{LS_glo_colimits_pointwise}, the enhanced global sections functor
$$\Gamma_{\fcat{A}}(BG,-)\simeq \ev_{BG} \colon \LocSys^{\Glo}_{\fcat T}(BG,\fcat{A}) \tto \fcat{A}(BG) $$
preserves colimits. In particular, $\ev_{BG}$ admits a right adjoint $Q^{BG}$. By e.g. \cref{corollary: right adjoint to evaluation on any object}, for any $\mathcal{F}\in \fcat{A}(BG)$, we obtain
$$(Q^{BG}\mathcal{F})^{BH} \simeq 
\begin{cases}
\mathcal{F}, & \mbox{if $H\simeq G$,}\\
0, & \mbox{otherwise.}
\end{cases}
$$
In particular, the right adjoint $Q^{BG}$ is a fully faithful embedding and the left adjoint $\ev_{BG}$ is a smashing localization. 

We will show that the adjoint pair $\ev_{BG} \dashv Q^{BG}$ is a part of a stable symmetric monoidal recollement (see \cite[Definition~A.8.1]{Lur_HA}). Indeed, let $\fcat{P}=\fcat{P}_G$ be an orbifamily of proper subgroups in $G$, see \cref{example: orbifamily associated to G}, and let 
$$\pi\colon (BG)_{\fcat{P}}=BG \times_{\mathcal{N}}\mathcal{N}_{\fcat{P}} \tto BG$$
be the canonical projection, see \cref{definiton: classfier subgroups_family}. Then, by \cref{example: orbifamily and pushforward}, the functor
$$\pi^{\Glo}_{*}\colon \LocSys^{\Glo}_{\fcat T}((BG)_{\fcat{P}},\fcat{A}) \xymatrix{\ar@{^(->}[r] &} \LocSys^{\Glo}_{\fcat T}(BG,\fcat{A})$$
is fully faithful, $\pi^{\Glo,*}\circ Q^{BG} \simeq 0$, and the functors $\pi^{\Glo,*}$ and $\ev_{BG}$ are jointly conservative. By definition, this implies the following assertion.

\begin{prop}[Isotropy separation for $\LocSys^{\Glo}_{\fcat T}(BG,\fcat{A})$]\label{proposition: global sections semi-orthogonal}
The following diagram
\[\xymatrix{
\LocSys^{\Glo}_{\fcat T}((BG)_{\fcat{P}},\fcat{A}) \ar@<-1ex>@{^(->}[r]_-{\pi^{\Glo}_{*}} & 
\ar@<-1ex>[l]_-{\pi^{\Glo,*}} 
\LocSys^{\Glo}_{\fcat T}(BG,\fcat{A}) \ar@<0.5ex>[r]^-{(-)^{BG}} & \ar@{^(->}@<0.5ex>[l]^-{Q^{BG}} \fcat{A}(BG),
}
\]
presents the category $\LocSys^{\Glo}_{\fcat T}(BG,\fcat{A})$ as a stable symmetric monoidal recollement of $\fcat{A}(BG)$ and $\LocSys^{\Glo}_{\fcat T}((BG)_{\fcat{P}},\fcat{A})$. In particular, there is a fiber sequence
$$\pi^{\Glo}_{\#} \pi^{\Glo,*} \mathcal{L} \tto \mathcal{L} \tto Q^{BG}\mathcal{L}^{BG} $$
for every globally equivariant local system $\mathcal{L} \in \LocSys^\Glo_{\fcat T}(BG,\fcat{A})$. \qed
\end{prop}

\begin{ex}\label{example: constant family}
Let $\fcat{A}=\underline{\Sp}$ be the constant coefficient system. Then the recollement of \cref{proposition: global sections semi-orthogonal} recovers the ordinary recollement
\[\xymatrix{
\Fun(\Orb^{\op}_{G, \fcat{P}}, \Sp) \ar@<-1ex>@{^(->}[r]_-{\pi_{*}} & 
\ar@<-1ex>[l]_-{\pi^{*}} 
\Fun(\Orb^{\op}_G,\Sp) \ar@<0.5ex>[r]^-{(-)^{BG}} & \ar@{^(->}@<0.5ex>[l]^-{Q^{BG}} \Sp,
}
\]
where $\Orb_{G, \fcat{P}}$ is the full subcategory of $\Orb_G$ spanned by the non-trivial orbits $G/H$, $H\neq G$. We recall that this recollement governs the classical isotropy fiber sequence
$$E\fcat{P} \otimes X \tto X \tto \widetilde{E\fcat{P}}\otimes X, $$
where $X\in \Sp^{nG} = \Fun(\Orb^{\op}_G,\Sp)$.
\end{ex}

\begin{ex}\label{example: right adjoint proper family}
By using \cref{right_adjoint_llax_limits}, we obtain that
$$(\pi^{\Glo}_*\mathcal{L})^{BG} \simeq \lim_{\substack{\phi\colon BH \to BG\\ BH\not\simeq BG}} \phi_{\fcat{A},*}(\mathcal{L}^{BH}) \in \fcat{A}(BG) $$
for a globally equivariant local system $\mathcal{L} \in \LocSys^{\Glo}((BG)_{\fcat{P}},\fcat{A})$.
\end{ex}

\begin{ex}\label{example: cyclic group global}
Let $G=C_p$ be a cyclic group of a prime order $p$. Then the orbifamily $\fcat{P}=\fcat{P}_{C_p}$ is trivial. In particular, the global space $(BG)_{\fcat{P}} = \underline{BC_p}$ is constant, see \cref{example: constant global space}. Therefore, 
$$\LocSys^{\Glo}((BG)_{\fcat{P}}, \fcat{A})\simeq \LocSys(BC_p,\fcat{A}(*))$$
is the ordinary category of local systems on $BC_p$. Let $\phi \colon * \to BC_p$ be a point of $BC_p$. Then, by \cref{example: right adjoint proper family}, the gluing functor 
$$(\pi^{\Glo}_*)^{BG}\colon \LocSys(BC_p,\fcat{A}(*)) \tto \fcat{A}(BC_p)$$
is given by the formula
$$(\pi^{\Glo}_*\mathcal{L})^{BG} \simeq (\phi_{\fcat{A},*}\mathcal{L})^{hC_p},\;\; \mathcal{L} \in \LocSys(BC_p,\fcat{A}(*)).$$
Finally, we obtain another description of the category $\LocSys^{\Glo}(BC_p, \fcat{A})$ as a right lax limit
$$\LocSys^{\Glo}(BC_p, \fcat{A}) \simeq \rlaxlim\left(\LocSys(BC_p,\fcat{A}(*)) \xrightarrow{(\phi_{\fcat{A},*}(-))^{hC_p}} \fcat{A}(BC_p)\right).$$
\end{ex}

\section{Genuine equivariant local systems}\label{section: genuine equivariant local system}
The goal of this section is to develop a theory of genuine equivariant local systems $\LocSys^\gen_{\fcat T}(-, \fcat A)$. For a compact Lie group $G$ there is a natural identification
$$\LocSys^\gen(BG, \underline{\Sp}) \simeq \Sp^G,$$
where the right hand side is the genuine equivariant stable category from \cite{LMS86}, which is responsible for our choice of notation. Motivated by the construction of $\Sp^G$ one is tempted to define $\LocSys^\gen_{\fcat T}(X, \fcat A)$ by $\otimes$-inverting Thom local systems of vector bundles on $X$, a notion that we discuss in \Cref{section: thom local systems}. Unfortunately, such a naive construction $\LocSys^\ngen_{\fcat T}(-, \fcat A)$ does not work well, since a general global space may have very few vector 
bundles, see \Cref{ex_LS_ngen_for_N}. Instead in \Cref{section: geuine local systems on global spaces} we define $\LocSys^\gen_{\fcat T}(-, \fcat A)$ as a right Kan extension of $\LocSys^\ngen_{\fcat T}(-, \fcat A)$ from orbits $\Orb_{\fcat T}$ to all global spaces. Recall that for $\LocSys^\Glo_{\fcat T}(-, \fcat A)$ we have a similar presentation by \Cref{cor_LSGlo_LKan_extended}.

For $X$ an orbispace we show that the right adjoint in the adjunction
$$\Sigma^\infty_X \colon \LocSys^\Glo_{\fcat T}(X, \fcat A) \xymatrix{\ar@<0.5ex>[r] & \ar@<0.5ex>[l]} \LocSys^\gen_{\fcat T}(X, \fcat A) \colon \Omega^\infty_X$$
is colimit preserving and conservative. Using this we prove that genuine local systems admit similar functoriality as $\LocSys^\Glo_{\fcat T}(-, \fcat A)$, most notably that the pullback functor along a faithful morphism admits a left adjoint, such that the base change and projection formula hold. One particularly important distinction from the globally equivariant case is that even for a limit preserving coefficient system the genuine pullback functor may not admit a left adjoint for non-faithful morphisms, see \Cref{remark:genuine_sharp_does_not_exists}.

It is an old observation that for a compact Hausdorff $G$-space $X$ any equivariant vector bundle locally embeds as a direct summand into a pullback of a bundle on a point. Using this in \Cref{section: LSgen on G-spaces} we show that in this case
$$\LocSys^\ngen_{\fcat T}(X//G, \fcat A) \simeq \LocSys^\gen_{\fcat T}(X//G, \fcat A).$$
We also show in this section that $\LocSys^\gen_{\fcat T}(BG, \fcat A)$ is rigid and atomically generated over $\fcat A(BG)$ with the set of generators indexed by $G$-orbits.

We finish this section by constructing for a global space $X$ a set of symmetric monoidal continuous jointly conservative \emdef{geometric fixed points functors}
$$\Phi^x_{\fcat A} \colon \LocSys^\gen(X, \fcat A) \tto \fcat A(BH),$$
where $x\colon BH \to X$ ranges over points of $X$. They are essentially determined by the property that there is a natural equivalence
$$\Phi_{\fcat A}^x(\Sigma^\infty_X \mathcal L) \simeq \mathcal L^x$$
for $\mathcal L \in \LocSys^\Glo(X, \fcat A)$. We use these functors to establish a semi-orthogonal decomposition for $\LocSys^\gen_{\fcat T}(BG, \fcat A)$ which generalizes classical isotropy separation for $\Sp^G$. This result is interesting in its own right, but will also be our main tool to prove a comparison between genuine and tempered local systems in \Cref{section: comparison of tempered and genuine local systems}.

\subsection{Thom local systems}\label{section: thom local systems}
An important aspect of equivariant homotopy theory is that the naive stabilization $\Sp^{nG} := \Fun(\Orb_G^\op, \Sp)$ of the category of $G$-spaces $\Type^G$ is \emph{not} suitable to work with the (generalized) equivariant cohomology (e.g.\@ because the suspension spectra of finite $G$-CW complexes are not dualizable in $\Sp^{nG}$, hence one cannot hope for Spanier--Whitehead or Poincar\'e/Atiyah dualities in $\Sp^{nG}$). Instead, the category of genuine $G$-spectra is defined as
$$\Sp^G := \Sp^{nG}[\{(\mathbb S^V)^{-1}\}_V],$$
where $V$ ranges over all finite dimensional real representations of $G$, and $\mathbb S^V$ denotes the corresponding representation sphere (a one-point compactification of $V$ with the natural $G$-action). Now note that by \Cref{LS_glo_on_BG}
$$\Sp^{nG} \simeq \LocSys^\Glo(BG, \underline{\Sp}).$$
The goal of this section is to introduce a construction of the so-called globally equivariant Thom local systems in $\LocSys^\Glo_{\fcat T}(X, \fcat A)$, which extends a construction of representation spheres $\mathbb S^V$ to the case of arbitrary global spaces $X$ and not necessarily constant coefficient systems $\fcat A$.
\begin{conv}
In this section, we assume that $\fcat A$ is a symmetric monoidal coefficient system on $\Orb_{\fcat T}$ of Beck--Chevalley type, see \cref{assumpt_limit_prerserving_cs}.
\end{conv}
\begin{defn}
Let $\Vect_\Glo$ denote the global space
$$\Vect_\Glo := \coprod_{n\ge 0} BO(n).$$
For a global space $X$ we define the \emdef{groupoid of vector bundles over $X$} as
$$\Vect_X := \Hom_{\Type^\Glo}(X, \Vect_\Glo).$$
\end{defn}
\begin{rem}
The block-diagonal inclusion $O(n) \times O(m) \inj O(n+m)$ induces a commutative monoid structure on $\Vect_\Glo$ and hence on $\Vect_X$ for any $X \in \Type^\Glo$. For vector bundles $E, E^\prime \in \Vect_X$, we denote the corresponding operation by $E\oplus E^\prime$ and we call it the direct sum of vector bundles.
\end{rem}
\begin{ex}\label{ex_VB_on_BG}
Let $G$ be a compact Lie group. Then a map $BG \to \Vect_\Glo$ factors through $BO(n)$ for some $n$, which by construction corresponds to some continuous group homomorphism $G \to O(n)$ up to conjugation. It follows that the set of objects of the groupoid $\Vect_{BG}$ is in natural bijection with the set of isomorphism classes of finite dimensional continuous real representations of $G$.
\end{ex}

\begin{construction}
Let $V$ be a finite dimensional real vector space. The space $V\setminus \{0\}$ is acted on by $O(V)$ and hence produces a faithful morphism
$$(V\setminus \{0\})\gitq O(V) \tto BO(V)$$
of global spaces. Taking coproduct over $V = \mathbb R^n$ for $n\ge 0$, we obtain a faithful map
$$E^\circ_{\mathrm{univ}} \tto \Vect_\Glo.$$
Let $E$ be a vector bundle over a global space $X$ which is classified by a map $f\colon X \to \Vect_\Glo$. We define $E^\circ \to X$ as a pullback of $E^\circ_{\mathrm{univ}}$ along $f$ and sometimes call it \emdef{$E$ without the zero section}. A fiberwise one-point compactification $\overline E_X$ of $E$ is defined similarly. Note that $\overline E_X$ is equivalent to the fiberwise suspension of $E^\circ$.
\end{construction}

\begin{construction}\label{def_ThomLS_glo}
Let $E$ be a vector bundle over a global space $X$ and let $p\colon E^\circ \to X$ denote the structure morphism. We define the \emdef{Thom local system $\mathbbl 1_X^E$ associated with $E$} as the cofiber of the counit of adjunction
$$\mathbbl 1_X^E := \cofib\left(p^{\Glo}_{\#} \mathbbl 1_{E^\circ} \tto \mathbbl 1_X \right) \in \LocSys^\Glo_{\fcat T}(X, \fcat A),$$
where
$$p^{\Glo}_{\#} \colon \LocSys^\Glo_{\fcat T}(E^\circ, \fcat A) \tto \LocSys^\Glo_{\fcat T}(X, \fcat A)$$
is the left adjoint of the pullback functor $p^{\Glo,*}$ which exists by \Cref{corollary: left adjoint faithful global}. Equivalently, by excision (see \Cref{notation_escision} below), if $\overline p\colon \overline E_X \to X$ is the structure morphism for the fiberwise one-point compactification of $E$, then there is a natural equivalence
$$\mathbbl 1_X^E \simeq \cofib\left(\mathbbl 1_X \tto \overline p_{\#} \mathbbl 1_{\overline E_X}\right),$$
where the map is induced by the infinity section $X \to \overline E_X$.
\end{construction}
\begin{lem}\label{lem_excision_for_global_homology}
The functor
$$C_{*}^{B, \Glo}(-, \fcat A) \colon (\Type^\Glo_{/B})^{\rep} \tto \LocSys^\Glo_{\fcat T}(B, \fcat A), \qquad \quad (p\colon X \to B) \mapsto p_{\#}(\mathbbl 1_X)$$
preserves small colimits.

\begin{proof}
Let $X \simeq \colim_I X_i \in (\Type^\Glo_{/ B})^\rep$. Then by \Cref{prop_LSglo_preserves_limits}
$$\LocSys^\Glo_{\fcat T}(X, \fcat A) \simeq \lim_I \LocSys^\Glo_{\fcat T}(X_i, \fcat A).$$
Let $p \colon X \to B$, $p_{i}\colon X_i \to B$ denote the structure maps. Then, by \Cref{lem_adj_functors_to_limit}, we have a natural equivalence
$$p^{\Glo}_{\#}(\mathcal L) \simeq \colim_{I^\op} p^{\Glo}_{i,\#}(\mathcal L_i)$$
for $\mathcal L \in \LocSys^\Glo_{\fcat T}(X, \fcat A)$, where $\mathcal L_i$ denotes the restriction of $\mathcal L$ to  $X_i$. The assertion follows by setting $\mathcal L = \mathbbl 1_{X}$ .
\end{proof}
\end{lem}
\begin{notation}\label{notation_escision}
In what follows, we will call this property \emdef{excision}.
\end{notation}
\begin{ex}
Continuing \Cref{ex_VB_on_BG}, for a finite dimensional real representation $V$ of $G$ considered as a vector bundle on $BG$, by construction the corresponding Thom local system in
$$\LocSys^\Glo(BG, \underline{\Sp}) \simeq \Sp^{nG}$$
is equivalent to the suspension spectrum of the corresponding representation sphere $\mathbb S^V$.
\end{ex}

This construction has the following familiar properties.
\begin{prop}\label{Thom_LS_basics_Glo}
Let $E$ be a vector bundle over a global space $X$. Then:
\begin{enumerate}
\item If $f\colon X^\prime \to X$ is a map of global spaces then there is a natural equivalence
$$f^*(\mathbbl 1_X^E) \simeq \mathbbl 1_{X^\prime}^{f^*E}.$$

\item For a pair of vector bundles $E^\prime, E^{\prime\prime}$ there is a natural equivalence
$$\mathbbl 1_X^{E^\prime \oplus E^{\prime\prime}} \simeq \mathbbl 1_X^{E^\prime} \otimes \mathbbl 1_X^{E^{\prime\prime}}.$$

\item Let $E$ be a trivial bundle of rank $r$. Then $\mathbbl 1_X^E \simeq \mathbbl 1_{X}[r]$.
\end{enumerate}

\begin{proof}
The first part follows from the fact that $(f^*E)^\circ \simeq f^*(E^\circ)$ and the base change for $\#$-pushforwards (\Cref{LSGlo_basechange_faithful_pbs}).

If $E = E^\prime \oplus E^{\prime\prime}$, we have
$$\overline{(E^\prime \oplus E^{\prime\prime})}_X \simeq \overline{E^\prime}_X \times \overline{E^{\prime\prime}}_X/(\overline{E^\prime}_X \vee \overline{E^{\prime\prime}}_X)$$
so the result follows from the relative K\"unneth formula \Cref{cor_LSGlo_Kunneth} and excision.

If $\rk E = 0$, then $\overline E_X \simeq X\coprod X$ and hence by definition $\mathbbl 1_X^E \simeq \mathbbl 1_X$. The computation of the Thom local system for the trivial bundle $E$ of rank $r>0$ follows by excision since $\overline E_X$ in this case is the $k$-times suspension of the one-point compactification of the trivial bundle. 
\end{proof}
\end{prop}

\subsection{Genuine local systems on global spaces}\label{section: geuine local systems on global spaces}
In this section we construct the category of genuine equivariant local systems and, using the results of \Cref{section: globally equivariant local systems}, we establish similar functoriality of this construction. Throughout this section, we fix a strongly continuous symmetric monoidal coefficient system 
$$\fcat{A}\colon \Orb^{\op}_{\fcat{T}}\to \CAlg(\Prs^\LL)$$
of Beck--Chevalley type, see \cref{assumpt_limit_prerserving_cs}.

\begin{defn}\label{def_genuine_LS}
For a compact Lie group $G$ from $\fcat T$, we define the category of \emdef{genuine local systems on $BG$} as
$$\LocSys^\gen_{\fcat T}(BG, \fcat A) = \LocSys^\Glo_{\fcat T}(BG, \fcat A)[\{\left(\mathbbl 1^V\right)^{-1}\}_{V\in \Vect_{BG}}],$$
where the right hand side denotes the formal $\otimes$-inversion of Thom local systems $\mathbbl 1^V$, see \Cref{definition: object inversion}.
\end{defn}

\begin{ex}
Let $G$ be a compact Lie group. Suppose that the coefficient system $\fcat{A}=\underline{\Sp}$ is constant. Then, by construction, we have
$$\LocSys^\gen(BG, \underline{\Sp}) \simeq \Sp^G,$$
where $\Sp^G$ denotes the genuine stable $G$-equivariant homotopy category.
\end{ex}

\begin{ex}\label{example: geom type gen BG}
Let $G$ be a compact Lie group from $\fcat T$. Suppose that the coefficient system $\fcat{A}$ is of geometric type and let $\fcat{A}_G \in \CAlg(\LocSys^{\Glo}_{\fcat T}(BG,\underline{\fcat{A}(BG)}))$ be a commutative algebra as in \cref{example: locsysglo geometric type}. Then, by construction, we have
$$\LocSys^\gen_{\fcat T}(BG, \fcat{A}) \simeq \Mod_{\fcat{A}_G}\left(\LocSys^{\gen}_{\fcat T}(BG,\underline{\fcat{A}(BG)})\right) \simeq \Mod_{\fcat{A}_G}\left(\Sp^G \otimes \fcat{A}(BG)\right).$$
\end{ex}

\begin{lem}\label{lemma: thom_cons_and_cont_bg}
Let $G$ be a compact Lie group from $\fcat T$ and let $V\in \Vect_{BG}$ be a real representation of~$G$. Then the endofunctor
$$[\mathbbl 1_{BG}^V, -] \colon \LocSys^\Glo_{\fcat T}(BG, \fcat A) \tto \LocSys^\Glo_{\fcat T}(BG, \fcat A)$$
is conservative and preserves small colimits.
\end{lem}

\begin{proof}
We recall that the coefficient system $\fcat{A}$ is strongly continuous. We will prove the first statement by induction on the dimension and the number of connected components of $G$. If $G$ is trivial, then $\mathbbl 1^V_{BG} \simeq \mathbbl 1_{BG}[{\dim V}]$ and the endofunctor $[\mathbbl 1_{BG}^V, -]$ is conservative because the category $\fcat A(*)$ is stable. Let $\mathcal L \in \LocSys^\Glo_{\fcat T}(BG, \fcat A)$ be a globally equivariant local system. If there exists a proper closed subgroup $H \subset G$ such that $\mathcal L|_{BH} \not\simeq 0$, then we have
$$[\mathbbl 1_{BG}^V, \mathcal L]|_{BH} \simeq [\mathbbl 1_{BH}^{V|_{BH}}, \mathcal L|_{BH}] \not\simeq 0$$
by using the inductive assumption and \cref{corollary: left adjoint faithful global linear}. Hence, $[\mathbbl 1_{BG}^V, \mathcal L] \not\simeq 0$. If $\mathcal L|_{BH} \simeq 0$ for all proper subgroups $H$, then
$$[\mathbbl 1_{BG}^V, \mathcal L]^G \simeq [\mathbbl 1[{\dim V^G}], \mathcal L^G]_{\fcat{A}(BG)} \simeq \Sigma^{- \dim V^G}\mathcal{L}^G.$$
We conclude that $[\mathbbl 1_{BG}^V, \mathcal L]$ vanishes if and only if $\mathcal L \simeq 0$.

Finally, the endofunctor $[\mathbbl 1_{BG}^V, -]$ preserves colimits by \cref{corollary: internal Hom colimits} and the fact that $\mathbb S^V$ is a finite $G$-space.
\end{proof}

\begin{rem}\label{remark: genlocal_susp}
By the definition, there is a natural symmetric monoidal functor
$$\Sigma^{\infty}_{BG} \colon \LocSys^{\Glo}_{\fcat T}(BG,\fcat{A}) \tto \LocSys^{\gen}_{\fcat T}(BG,\fcat{A}). $$
Note that $\Sigma^{\infty}_{BG}$ preserves colimits, so $\Sigma^{\infty}_{BG}$ admits a right adjoint
$$\Omega^{\infty}_{BG}\colon \LocSys^{\gen}_{\fcat T}(BG,\fcat{A}) \tto \LocSys^{\Glo}_{\fcat T}(BG,\fcat{A}).$$
\end{rem}

\begin{cor}\label{cor_OmegaInftyNgen}
Let $G$ be a compact Lie group from $\fcat T$. Then the forgetful functor
$$\Omega^\infty_{BG} \colon \LocSys^{\gen}_{\fcat T}(BG, \fcat A) \tto \LocSys^\Glo_{\fcat T}(BG, \fcat A)$$
is conservative and preserves small colimits.

\begin{proof}
We recall that the coefficient system $\fcat{A}$ is strongly continuous. By the limit presentation of $\LocSys^{\gen}_{\fcat T}(BG, \fcat A)$ from \cref{proposition: inversion symmetric},  it is enough to prove that for any real representation $V$ on $BG$ the induced functor
$$[\mathbbl 1_{BG}^V, -] \colon \LocSys^\Glo_{\fcat T}(BG, \fcat A) \tto \LocSys^\Glo_{\fcat T}(BG, \fcat A)$$
is conservative and preserves small colimits, which is proven in \cref{lemma: thom_cons_and_cont_bg}.
\end{proof}
\end{cor}

\begin{notation}\label{notation: categorical fixed points}
Given a subgroup $H \subset G$, we will write 
$$(-)^H\colon \LocSys^{\gen}_{\fcat T}(BG, \fcat A) \tto \fcat A(BH)$$ for the composite 
$$(-)^H\simeq (\Omega^\infty_{BG}(-))^H \simeq \ev_{BH}\circ \Omega^{\infty}_{BG}.$$
This abuse of notation is standard in genuine equivariant homotopy theory.
\end{notation}

Next, following \cite[Section~2.3]{Cnossen23}, we will construct a functor
$$\LocSys^{\gen}_{\fcat T}(-,\fcat{A}) \colon \Orb^{\op}_{\fcat T} \tto \CAlg(\Prs^\LL)$$
which extends the assignment
$$BG \xymatrix{\ar@{|->}[r] &} \LocSys^\gen_{\fcat T}(BG, \fcat A). $$
%extends to a functor from the category of orbits $\Orb^{\op}$ to $\CAlg(\Prs^\LL)$.

\begin{defn}[Definitions~2.43 and 2.44 in \cite{Cnossen23}]\label{definition: augmented categories}
We define $\Cat_{\mathrm{aug}}$ as the full subcategory of $\Ar(\Cat)=\Fun(\Delta^1,\Cat)$ spanned by the morphisms in $\Cat$ corresponding to the fully faithful embeddings $\iota\colon \fcat{I} \hookrightarrow \fcat{C}$, where $\fcat{I}$ is a small category. Similarly, we define the category $\CAlg(\Prs^\LL)_{\mathrm{aug}}$ as the pullback
$$\xymatrix{
\CAlg(\Prs^\LL)_{\mathrm{aug}}\ar[d] \ar[r] & \Cat_{\mathrm{aug}} \ar[d] \\
\CAlg(\Prs^\LL) \ar[r] & \Cat.
}$$
\end{defn}

Note that the forgetful functor $\CAlg(\Prs^\LL)_{\mathrm{aug}} \to \CAlg(\Prs^\LL)$ admits a section
$$(-)_{\mathrm{inv}}\colon \CAlg(\Prs^\LL) \tto \CAlg(\Prs^\LL)_{\mathrm{aug}} $$
which equips a symmetric monoidal category $\fcat{C}$ with its collection of $\otimes$-invertible objects. We recall that the functor $(-)_{\mathrm{inv}}$ admits a left adjoint.

\begin{prop}\label{proposition: functorial inversion}
The functor $(-)_{\mathrm{inv}}$ admits a left adjoint
$$\mathrm{Quot}\colon \CAlg(\Prs^\LL)_{\mathrm{aug}} \tto \CAlg(\Prs^\LL) $$
given on objects by sending a pair $(\fcat{C},\fcat{I})$ to the formal inversion $\fcat{C}[\{X^{-1}\}_{X\in \fcat{I}}]$.
\end{prop}

\begin{proof}
See \cite[Lemma~2.45]{Cnossen23}.
\end{proof}

\begin{construction}\label{construction: lsgen functor}
By \cref{Thom_LS_basics_Glo}, there is a functor
$$\mathrm{Thom}\colon \Orb^{\op}_{\fcat T} \tto \CAlg(\Prs^\LL)_{\mathrm{aug}} $$
which sends an object $T \in \Orb_{\fcat T}$ to the pair $(\LocSys^{\Glo}_{\fcat T}(T,\fcat{A}),\mathrm{Thom}(T))$, where $\mathrm{Thom}(T)$ is a full subcategory of $\LocSys^{\Glo}_{\fcat T}(T,\fcat{A})$ spanned by all Thom local systems. Finally, we define 
$$\LocSys^{\gen}_{\fcat T}(-,\fcat{A})\colon \Orb^{\op}_{\fcat T} \tto \CAlg(\Prs^\LL)$$
as the composite
$$\LocSys^{\gen}_{\fcat T}(-,\fcat{A})\colon\Orb^{\op}_{\fcat T} \xrightarrow{\mathrm{Thom}} \CAlg(\Prs^\LL)_{\mathrm{aug}} \xrightarrow{\mathrm{Quot}} \CAlg(\Prs^\LL).$$
\end{construction}

By construction, a morphism $f\colon T' \to T$ in $\Orb_{\fcat T}$ induces a colimit-preserving symmetric monoidal functor
$$f^{\gen,*}\colon \LocSys^{\gen}_{\fcat T}(T,\fcat{A}) \tto \LocSys^{\gen}_{\fcat T}(T',\fcat{A}) $$
such that the following diagram
\begin{equation}\label{equation: pullback for gen orbits}
\xymatrix{
\LocSys^\Glo_{\fcat T}(T', \fcat A) \ar[r]^-{\Sigma^\infty_{T'}} & \LocSys^\gen_{\fcat T}(T', \fcat A)  \\
\LocSys^\Glo_{\fcat T}(T, \fcat A) \ar[r]^-{\Sigma^\infty_{T}}\ar[u]^{f^{\Glo,*}} & \LocSys^{\gen}_{\fcat T}(T, \fcat A) \ar[u]^{f^{\gen,*}}}
\end{equation}
commutes.

\begin{lem}\label{lemma: sharp for gen orbits}
Let $f\colon T' \to T$ be a faithful morphism in $\Orb_{\fcat T}$. Then the diagram~\eqref{equation: pullback for gen orbits} is horizontally right adjointable, that is, the diagram
$$
\xymatrix{
\LocSys^\Glo_{\fcat T}(T', \fcat A) & \ar[l]_-{\Omega^\infty_{T'}} \LocSys^{\gen}_{\fcat T}(T', \fcat A) \\
\LocSys^\Glo_{\fcat T}(T, \fcat A) \ar[u]_{f^{\Glo,*}} & \ar[l]_-{\Omega^\infty_{T}} \LocSys^{\gen}_{\fcat T}(T, \fcat A) \ar[u]_{f^{\gen,*}}
}$$
commutes.
\end{lem}

\begin{proof}
Without loss of generality, we can assume that $T\simeq BG$ and $T'\simeq BH$, where $H$ is a subgroup of $G$. By construction, there is a natural transformation
$$\theta\colon f^{\Glo,*}\Omega^{\infty}_{BG} \tto \Omega^{\infty}_{BH}f^{\gen,*}.$$
By \cref{cor_OmegaInftyNgen}, the category $\LocSys^{\gen}_{\fcat T}(BG,\fcat{A})$ is generated under colimits by the image of $\Sigma^\infty_{BG}$. Again, by \cref{cor_OmegaInftyNgen}, both sides of $\theta$ commute with colimits. Therefore, it suffices to show that
$$\theta\colon f^{\Glo,*}\Omega^{\infty}_{BG}\Sigma^{\infty}_{BG}\mathcal{L} \tto \Omega^{\infty}_{BH}f^{\gen,*}\Sigma^{\infty}_{BG}\mathcal{L} \simeq \Omega^{\infty}_{BH}\Sigma^{\infty}_{BH}f^{\Glo,*}\mathcal{L} $$
is an equivalence for all $\mathcal{L}\in \LocSys^{\Glo}_{\fcat T}(BG,\fcat{A}).$ However, by \cref{proposition: right adjoint to inversion} and \cref{Thom_LS_basics_Glo}, we have
$$\Omega^{\infty}_{BG}\Sigma^{\infty}_{BG}\mathcal{L} \simeq \colim_{V\in \Vect_{BG}} [\mathbbl{1}^V_{BG},\mathbbl{1}^V_{BG}\otimes \mathcal{L}] $$
and $$\Omega^{\infty}_{BH}\Sigma^{\infty}_{BH}f^{\Glo,*}\mathcal{L} \simeq \colim_{W\in \Vect_{BH}} [\mathbbl{1}^W_{BH},\mathbbl{1}^W_{BH}\otimes f^{\Glo,*}\mathcal{L}].$$
Note that any real $H$-representation $W$ is a direct summand of some $G$-representation $V$. Therefore,
$$\Omega^{\infty}_{BH}\Sigma^{\infty}_{BH}f^{\Glo,*}\mathcal{L} \simeq \colim_{V\in \Vect_{BG}} [f^{\Glo,*}\mathbbl{1}^V_{BH},f^{\Glo,*}\mathbbl{1}^V_{BH}\otimes f^{\Glo,*}\mathcal{L}].$$
Finally, the assertion follows by \cref{corollary: left adjoint faithful global linear}.
\end{proof}

\begin{prop}\label{proposition: existence of left adjoint gen faithful orb}
Let $f\colon T' \to T$ be a faithful morphism in $\Orb_{\fcat T}$. The pullback functor
$$f^{\gen,*}\colon \LocSys^{\gen}_{\fcat T}(T, \fcat A) \tto \LocSys^{\gen}_{\fcat T}(T', \fcat A)$$
preserves small limits, is closed symmetric monoidal, and admits a $\LocSys^{\gen}_{\fcat T}(T, \fcat A)$-linear left adjoint $f^{\gen}_{\#}$. In particular, for all $\mathcal L \in \LocSys^{\gen}_{\fcat T}(T', \fcat A)$ and $\mathcal M, \mathcal N \in \LocSys^{\gen}_{\fcat T}(T, \fcat A)$, we have
$$f^{\gen}_{\#}(\mathcal L \otimes f^{\gen,*}(\mathcal M)) \xymatrix{\ar[r]^-\sim &} f^{\gen}_{\#}(\mathcal L)\otimes \mathcal M, \qquad f^{\gen,*}[\mathcal M, \mathcal N] \xymatrix{\ar[r]^-\simeq &} [f^{\gen,*}\mathcal M, f^{\gen,*} \mathcal N].$$
\end{prop}

\begin{proof}
We will show first that $f^{\gen,*}$ commutes with limits. Indeed, by \cref{cor_OmegaInftyNgen}, it suffices to show that $\Omega^{\infty}_{T'}f^{\gen,*}$ is limit-preserving. However, by \cref{lemma: sharp for gen orbits}, we have $\Omega^{\infty}_{T'}f^{\gen,*} \simeq f^{\Glo,*}\Omega^{\infty}_T$ and the assertion follows by \cref{corollary: left adjoint faithful global}. By the adjoint functor theorem, the pullback $f^{\gen,*}$ admits a left adjoint
$$f^{\gen}_{\#}\colon \LocSys^{\gen}_{\fcat T}(T',\fcat{A}) \tto \LocSys^{\gen}_{\fcat T}(T,\fcat{A}).$$ Moreover, by \cref{lemma: sharp for gen orbits}, the following diagram
\begin{equation}\label{equation: sharp and susp gen orbit}
\xymatrix{
\LocSys^\Glo_{\fcat T}(T', \fcat A) \ar[r]^-{\Sigma^\infty_{T'}} \ar[d]^{f^{\Glo}_{\#}} & \LocSys^{\gen}_{\fcat T}(T', \fcat A) \ar[d]^{f^{\gen}_{\#}} \\
\LocSys^\Glo_{\fcat T}(T, \fcat A) \ar[r]^-{\Sigma^\infty_{T}} & \LocSys^{\gen}_{\fcat T}(T, \fcat A)
}
\end{equation}
commutes.

Next, it suffices to show that the natural transformation
$$f^{\gen}_{\#}(\mathcal L \otimes f^{\gen,*}(\mathcal M)) \tto f^{\gen}_{\#}(\mathcal L)\otimes \mathcal M $$
is an equivalence for all $\mathcal L \in \LocSys^{\gen}_{\fcat T}(T', \fcat A)$ and $\mathcal M \in \LocSys^{\gen}_{\fcat T}(T, \fcat A)$. By \cref{cor_OmegaInftyNgen}, we can assume that $\mathcal{L}\simeq \Sigma^{\infty}_{T'}\mathcal{L}_0$ and $\mathcal{M}\simeq \Sigma^{\infty}_{T}\mathcal{M}_0$ for some $\mathcal{L}_0 \in \LocSys^{\Glo}_{\fcat T}(T',\fcat{A})$ and $\mathcal{M}_0 \in \LocSys^{\Glo}_{\fcat T}(T,\fcat{A})$. Finally, the assertion follows by commutativity of the diagram~\eqref{equation: sharp and susp gen orbit} and \cref{corollary: left adjoint faithful global linear}.
\end{proof}

\begin{warn}\label{remark:genuine_sharp_does_not_exists}
Unfortunately, even if the pullback functor
$$f^{\Glo,*}\colon \LocSys^\Glo_{\fcat T}(T, \fcat A) \tto \LocSys^\Glo_{\fcat T}(T', \fcat A)$$
admits a left adjoint for full morphisms $f\colon T' \to T$, the genuine pullback $f^{\gen,*}$ may \emph{not} admit a left adjoint as the next example shows.

Indeed, let $G$ be a nontrivial compact Lie group and let $p\colon BG \to *$ be the natural map and let $\fcat A = \underline{\Sp}$ be the constant coefficient system. We claim that the pullback functor $p^*\colon \Sp \to \Sp^G$ does not admit a left adjoint. Indeed, let us assume that the left adjoint $p_\#$ exists. Then, since the right adjoint $p^*$ is continuous, $p_{\#}$ preserves compact objects. Moreover, since compact objects in $\Sp^G$ coincide with the dualizable objects, $p_{\#}$ preserves dualizable objects as well. In particular, $p_\#(\mathbb S) \in \Sp$ is dualizable, and so is $p_\#(\mathbb S)^\vee$. However,
$$p_\#(\mathbb S)^\vee \simeq \Hom_{\Sp^G}(\mathbb S, \mathbb S) \simeq \mathbb S^G$$
which is never dualizable for a nontrivial group $G$, since $\mathbb{S}^G$ contains $\Sigma^\infty_+ {BG}$ as a direct summand by the tom Dieck splitting.
\end{warn}

\begin{cor}\label{corollary: base change simple case}
Let
\[\xymatrix{
T'_0 \ar[r]^-q \ar[d]^g & T' \ar[d]^f \\
T_0 \ar[r]^-p & T
}\]
be a fibered diagram in $\Orb_{\fcat T}$ such that $f$ is faithful and $p$ is full. Then the base change maps
$$f^{\gen,*} \circ p^{\gen}_* \tto q^{\gen}_* \circ g^{\gen,*}, \qquad g^{\gen}_{\#}\circ q^{\gen,*} \tto p^{\gen,*} \circ f^{\gen}_{\#} $$
are equivalences.
\end{cor}

\begin{proof}
Recall that the coefficient system $\fcat{A}$ is of Beck--Chevalley type. We will show that the canonical natural transformation
$$\theta \colon f^{\gen,*} \circ p^{\gen}_* \tto q^{\gen}_* \circ g^{\gen,*} $$
is an equivalence. By \cref{cor_OmegaInftyNgen}, it suffices to show that $\Omega^{\infty}_{T'}\theta$ is an equivalence. Finally, the assertion follows by \cref{lemma: sharp for gen orbits} and \cref{LSGlo_basechange_faithful_pbs}.
\end{proof}

\begin{defn}\label{construction: genuine local system}
We define the category of \emdef{genuine equivariant local systems} 
$$\LocSys^\gen_{\fcat T}(-, \fcat A) \colon \Type^{\Glo, \op} \tto \CAlg(\Prs^\LL)$$
as the right Kan extension of the functor 
$$\LocSys^\gen_{\fcat T}(-,\fcat A) \colon \Orb^{\op}_{\fcat T} \tto \CAlg(\Prs^\LL) $$
from \cref{construction: lsgen functor} along the Yoneda embedding $\yoneda \colon \Orb^\op_{\fcat T} \inj \Type^{\Glo, \op}$. More concretely, for a global space $X\in \Type^\Glo$, we have
$$\LocSys^\gen_{\fcat T}(X, \fcat A) \simeq \lim_{T \in (\Orb_{\fcat T /X})^\op} \LocSys^\gen_{\fcat T}(T, \fcat A).$$
As in \Cref{notation_for_LSGlo_allfamily}, for $\fcat T$ being a family consisting of all compact Lie groups, we denote $\LocSys^\gen_{\fcat T}(X, \fcat A)$ simply by $\LocSys^\gen(X, \fcat A)$.
\end{defn}
\begin{rem}\label{rem_LSGen_T_depends_only_on_Type_T}
By construction the category $\LocSys^\gen_{\fcat T}(X, \fcat A)$ depends only on the restriction of $X$ to $\Type^\Glo_{\fcat T}$. This behavior is consistent with \Cref{rem_LSGlo_T_depends_on_SGlo_T}.
\end{rem}
%\begin{rem}
%By \cref{construction: genuine local system}, there is a natural colimit preserving symmetric monoidal functor
%$$\LocSys^\ngen(X, \fcat A) \tto \LocSys^\gen(X, \fcat A).$$
%\end{rem}
%\begin{ex}
%If we consider the constant coefficient system $\fcat A(BG) := \Sp$ then for a compact Lie group $G$ by construction
%$$\LocSys^\gen(BG, \mathcal A) \simeq \Sp^G,$$
%where $\Sp^G$ denotes the genuine stable $G$-equivariant homotopy category.
%\end{ex}

\begin{notation}
Given a global space $X$, we write
$$\Sigma^\infty_{X} \colon \LocSys^\Glo_{\fcat T}(X, \fcat A) \tto \LocSys^{\gen}_{\fcat T}(X, \fcat A)$$
for the canonical symmetric monoidal functor, i.e. $\Sigma^{\infty}_X \simeq \lim_{T\to X} \Sigma^{\infty}_T$, see \cref{remark: genlocal_susp}.  By the adjoint functor theorem, $\Sigma^\infty_{X}$ admits a right adjoint which we will denote by $\Omega^\infty_{X}$.
\end{notation}

\begin{rem}\label{remark: parametrized genuine formal inversion}
\cref{proposition: existence of left adjoint gen faithful orb} and \cref{corollary: base change simple case} imply together that the functor
$$\LocSys^{\gen}(-,\fcat{A})\colon (\Type^{\Orb})^{\op} \tto \CAlg(\Prs^\LL)$$
is a presentably symmetric monoidal $\Type^\Orb$-category, cf. \cref{remark: parametrized global systems}. Moreover, one can check that the natural transformation 
$$\Sigma^{\infty}\colon \LocSys^{\Glo}(-,\fcat{A}) \tto \LocSys^{\gen}(-,\fcat{A})$$
exhibits $\LocSys^{\gen}(-,\fcat{A})$ as the \emph{formal inversion} of $\LocSys^{\Glo}(-,\fcat{A})$ in the full subcategory given by Thom local systems, see~\cite[Section~2.3]{Cnossen23}.
\end{rem}

\begin{rem}
By \cref{construction: genuine local system}, a morphism of global spaces $f\colon X \to Y$ induces a colimit-preserving symmetric monoidal \emdef{pullback functor}
$$f^{\gen,*}\colon \LocSys^\gen_{\fcat T}(Y, \mathcal A) \tto \LocSys^\gen_{\fcat T}(X, \mathcal A).$$
By the adjoint functor theorem, $f^{\gen,*}$ admits the (right lax monoidal) right adjoint
$$f^{\gen}_*\colon \LocSys^\gen_{\fcat T}(X, \mathcal A) \tto \LocSys^\gen_{\fcat T}(Y, \mathcal A).$$
\end{rem}

The construction from \Cref{construction: genuine local system} has the following functoriality properties.

\begin{prop}\label{proposition: genuine sharp easy}
Let $f\colon X \to Y$ be a faithful morphism of global spaces. Then
\begin{enumerate}[label=(\arabic*)]
%\item The pullback functor
%$$f^{\gen,*}\colon \LocSys^{\gen}(Y, \fcat A) \tto \LocSys^{\gen}(X, \fcat A)$$
%preserves small limits, $f^{\gen,*}$ is closed symmetric monoidal, and $f^{\gen,*}$ admits a left adjoint 
%$$f^{\gen}_{\#}\colon \LocSys^{\gen}(X, \fcat A) \tto \LocSys^{\gen}(Y, \fcat A);$$ %In particular, there are natural equivalences 
%$$f^{\gen}_{\#}(\mathcal L \otimes f^{\gen,*}(\mathcal M)) \xymatrix{\ar[r]^-\sim &} f^{\gen}_{\#}(\mathcal L)\otimes \mathcal M,$$
%$$f^{\gen,*}[\mathcal M, \mathcal N] \xymatrix{\ar[r]^-\simeq &} [f^{\gen,*}\mathcal M, f^{\gen,*} \mathcal N]$$
%for all $\mathcal L \in \LocSys^{\gen}(X, \fcat A)$ and $\mathcal M, \mathcal N \in \LocSys^{\gen}(Y, \fcat A)$;

\item The diagram
$$ \xymatrix{
\LocSys^\Glo(X, \fcat A) & \ar[l]_-{\Omega^\infty_{X}} \LocSys^{\gen}(X, \fcat A) \\
\LocSys^\Glo(Y, \fcat A) \ar[u]_{f^{\Glo,*}} & \ar[l]_-{\Omega^\infty_{Y}} \LocSys^{\gen}(Y, \fcat A) \ar[u]_{f^{\gen,*}}
}$$
is commutative;

\item Consider the fiber square
$$\xymatrix{
X^\prime \ar[d]^{g} \ar[r]^-q & X \ar[d]^f \\
Y^\prime \ar[r]^-p & Y,
}$$
where the morphisms $f$ and $g$ are faithful. Then the diagram of categories
$$\xymatrix{
\LocSys^\gen(X^\prime, \fcat A) & \ar[l]_-{q^{\gen,*}} \LocSys^{\gen}(X, \fcat A) \\
\LocSys^\gen(Y^\prime, \fcat A) \ar[u]_{g^{\gen,*}} & \ar[l]_-{p^{\gen,*}} \LocSys^{\gen}(Y, \fcat A) \ar[u]_{f^{\gen,*}}
}$$
is horizontally right adjointable. That is the natural transformation
$$f^{\gen,*} \circ p^{\gen}_* \tto q^{\gen}_* \circ g^{\gen,*}$$
is an equivalence.
\end{enumerate}
\end{prop}

\begin{proof}
First we show that the commutative diagram
\begin{equation}\label{equation: star_global_susp_eq1}
\xymatrix{
\LocSys^\Glo_{\fcat T}(X, \fcat A) \ar[r]^-{\Sigma^\infty_{X}} & \LocSys^\gen_{\fcat T}(X, \fcat A)  \\
\LocSys^\Glo_{\fcat T}(Y, \fcat A) \ar[r]^-{\Sigma^\infty_{Y}}\ar[u]^{f^{\Glo,*}} & \LocSys^{\gen}_{\fcat T}(Y, \fcat A) \ar[u]^{f^{\gen,*}}}
\end{equation}
is horizontally right adjointable. By construction, the diagram~\eqref{equation: star_global_susp_eq1} is the limit of commutative diagrams
\[\xymatrix{
\LocSys^\Glo_{\fcat T}(T'_{\alpha}, \fcat A) \ar[r]^-{\Sigma^\infty_{T'_{\alpha}}} & \LocSys^\gen_{\fcat T}(T'_{\alpha}, \fcat A)  \\
\LocSys^\Glo_{\fcat T}(T_{\alpha}, \fcat A) \ar[r]^-{\Sigma^\infty_{T_{\alpha}}}\ar[u]^{f^{\Glo,*}_{\alpha}} & \LocSys^{\gen}_{\fcat T}(T_{\alpha}, \fcat A), \ar[u]^{f^{\gen,*}_{\alpha}}}
\]
where $T_\alpha \in \Orb_{\fcat T}$. These diagrams are right adjointable by \cref{lemma: sharp for gen orbits}. Since, any limit of right adjointable squares is right adjointable (see~\cite[Corollary~4.7.4.18]{Lur_HA}), so is the diagram~\eqref{equation: star_global_susp_eq1}.

Next, we will prove the second part. Again, by \cite[Corollary~4.7.4.18]{Lur_HA}, a small limit of adjointable squares is adjointable. Note that by \Cref{rem_LSGen_T_depends_only_on_Type_T}, we can assume that all global spaces lie in $\Type^\Glo_{\fcat T}$. So, by presenting $X, Y$ and $Y^\prime$ as colimits of global $\fcat T$-orbits and by \Cref{def_genuine_LS}, we can assume that $X = T'$, $Y = T$ and $Y^\prime = T_0$, $f\colon T' \to T$ is faithful. Then, the map $p\colon T_0 \to T$ is a composite of a full morphism followed by a faithful morphism. We will consider both cases separately. If $p$ is full, then the assertion is the content of \cref{corollary: base change simple case}. If $p$ is faithful, then it suffices to show that 
$$\theta \colon f^{\gen,*} \circ p^{\gen}_* \tto q^{\gen}_* \circ g^{\gen,*} $$
is an equivalence. By \cref{cor_OmegaInftyNgen}, it suffices to show that $\Omega^{\infty}_{T'}\theta$ is an equivalence. Finally, the assertion follows by \cref{lemma: sharp for gen orbits} and \cref{LSGlo_basechange_faithful_pbs_ra}.
\end{proof}

\begin{thm}\label{proposition: genuine sharp}
Let $f\colon X \to Y$ be a faithful morphism of global spaces. Then
\begin{enumerate}[label=(\arabic*)]
\item The pullback functor
$$f^{\gen,*}\colon \LocSys^{\gen}_{\fcat T}(Y, \fcat A) \tto \LocSys^{\gen}_{\fcat T}(X, \fcat A)$$
preserves small limits, $f^{\gen,*}$ is closed symmetric monoidal, and $f^{\gen,*}$ admits a left adjoint 
$$f^{\gen}_{\#}\colon \LocSys^{\gen}_{\fcat T}(X, \fcat A) \tto \LocSys^{\gen}_{\fcat T}(Y, \fcat A);$$ %In particular, there are natural equivalences 
%$$f^{\gen}_{\#}(\mathcal L \otimes f^{\gen,*}(\mathcal M)) \xymatrix{\ar[r]^-\sim &} f^{\gen}_{\#}(\mathcal L)\otimes \mathcal M,$$
%$$f^{\gen,*}[\mathcal M, \mathcal N] \xymatrix{\ar[r]^-\simeq &} [f^{\gen,*}\mathcal M, f^{\gen,*} \mathcal N]$$
%for all $\mathcal L \in \LocSys^{\gen}(X, \fcat A)$ and $\mathcal M, \mathcal N \in \LocSys^{\gen}(Y, \fcat A)$;

\item The diagram
\[\xymatrix{
\LocSys^\Glo_{\fcat T}(X, \fcat A) \ar[r]^-{\Sigma^\infty_{X}} \ar[d]^{f^{\Glo}_{\#}} & \LocSys^\gen_{\fcat T}(X, \fcat A) \ar[d]^{f^{\gen}_{\#}} \\
\LocSys^\Glo_{\fcat T}(Y, \fcat A) \ar[r]^-{\Sigma^\infty_{Y}} & \LocSys^{\gen}_{\fcat T}(Y, \fcat A),
}\]
is commutative;

\item Consider the fiber square
$$\xymatrix{
X^\prime \ar[d]^{g} \ar[r]^-q & X \ar[d]^f \\
Y^\prime \ar[r]^-p & Y,
}$$
where the morphisms $f$ and $g$ are faithful. Then the diagram of categories
\[\xymatrix{
\LocSys^\gen_{\fcat T}(X^\prime, \fcat A) & \ar[l]_-{q^{\gen,*}} \LocSys^{\gen}_{\fcat T}(X, \fcat A) \\
\LocSys^\gen_{\fcat T}(Y^\prime, \fcat A) \ar[u]_{g^{\gen,*}} & \ar[l]_-{p^{\gen,*}} \LocSys^{\gen}_{\fcat T}(Y, \fcat A) \ar[u]_{f^{\gen,*}}
}\]
is vertically left adjointable. That is the natural transformation
$$g^{\gen}_\# \circ q^{\gen,*} \tto p^{\gen,*} \circ f^{\gen}_{\#}$$
is an equivalence.
\end{enumerate}
\end{thm}

\begin{proof}
We will only prove the existence of the left adjoint $f^{\gen}_\#$; two other parts follow formally from \cref{proposition: genuine sharp easy}. The proof is similar to the proof of \cref{corollary: left adjoint faithful global}.

By \Cref{rem_LSGen_T_depends_only_on_Type_T} we can assume that $Y\in \Type^\Glo_{\fcat T}$. Let us present it as a colimit $Y\simeq \colim_{\alpha \in \fcat{I}} T_{\alpha}$ of representable global $\fcat T$-spaces. Set $X_\alpha \simeq X\times_Y T_{\alpha}$ and $f_\alpha\colon X_\alpha \to T_\alpha$ be the faithful projection map, $\alpha \in \fcat{I}$. Then, by the virtue of \cref{construction: genuine local system} and \cref{proposition: existence of left adjoint gen faithful orb}, the pullback functor $f^{\gen,*}$ is a limit $\lim_{\alpha \in \fcat{I}^\op} f^{\gen,*}_\alpha$. Moreover, by \cref{proposition: genuine sharp easy}, the diagram
$$f^{\gen,*}_\bullet\colon \fcat{I}^\op \tto \Fun(\Delta^1,\Cat), \;\; \alpha \mapsto f^{\gen,*}_{\alpha}$$
factors through the subcategory $\Fun^{\mathrm{LAd}}(\Delta^1,\Cat)$ of right adjoints and left adjointable natural transformations, see~\cite[Definition~4.7.4.16]{Lur_HA}. Finally, by \cite[Corollary~4.7.4.18]{Lur_HA}, a limit of right adjoint functors along left adjointable natural transformations is a right adjoint, i.e. $f^{\gen,*}\simeq \lim_{\alpha} f^{\gen,*}_\alpha$ is a right adjoint.
\end{proof}

\begin{cor}\label{cor_OmegaInftygen}
Let $X\in \Type^{\Glo}$ be an orbispace. Then the forgetful functor
$$\Omega^\infty_{X} \colon \LocSys^{\gen}_{\fcat T}(X, \fcat A) \tto \LocSys^\Glo_{\fcat T}(X, \fcat A)$$
is conservative and preserves small colimits.

\begin{proof}
By the second part of \cref{proposition: genuine sharp} and \cref{cor_OmegaInftyNgen}, the right adjoint $\Omega^\infty_X$ is the limit of colimit-preserving conservative functors $\Omega^{\infty}_{BG}$. Therefore, by~\cite[Proposition~5.5.3.13]{Lur_HTT}, $\Omega^{\infty}_X$ is conservative and preserves small colimits.
\end{proof}
\end{cor}

\begin{cor}\label{corollary: projection formula for genuine sharp}
Let $f\colon X\to Y$ be a faithful morphism of global spaces. Then the left adjoint
$$f^{\gen}_{\#}\colon \LocSys^{\gen}_{\fcat T}(X, \fcat A) \tto \LocSys^{\gen}_{\fcat T}(Y, \fcat A)$$ 
is a $\LocSys^{\gen}_{\fcat T}(Y, \fcat A)$-linear. That is there are natural equivalences 
$$f^{\gen}_{\#}(\mathcal L \otimes f^{\gen,*}(\mathcal M)) \xymatrix{\ar[r]^-\sim &} f^{\gen}_{\#}(\mathcal L)\otimes \mathcal M,$$
$$f^{\gen,*}[\mathcal M, \mathcal N] \xymatrix{\ar[r]^-\sim &} [f^{\gen,*}\mathcal M, f^{\gen,*} \mathcal N]$$
for all $\mathcal L \in \LocSys^{\gen}_{\fcat T}(X, \fcat A)$ and $\mathcal M, \mathcal N \in \LocSys^{\gen}_{\fcat T}(Y, \fcat A)$.
\end{cor}

\begin{proof}
We will show that the canonical morphism
$$\theta\colon f^{\gen}_{\#}(\mathcal L \otimes f^{\gen,*}(\mathcal M)) \tto f^{\gen}_{\#}(\mathcal L)\otimes \mathcal M$$
is an equivalence for all $\mathcal L \in \LocSys^{\gen}_{\fcat T}(X, \fcat A)$ and $\mathcal M \in \LocSys^{\gen}_{\fcat T}(Y, \fcat A)$. First suppose that $Y$ be an orbispace. Then $X$ is an orbispaces as well and we can assume by \cref{cor_OmegaInftygen} that $\mathcal{L}\simeq \Sigma^{\infty}_{X}\mathcal{L}_0$ and $\mathcal{M}\simeq \Sigma^{\infty}_{Y}\mathcal{M}_0$ for some $\mathcal{L}_0 \in \LocSys^{\Glo}_{\fcat T}(X,\fcat{A})$ and $\mathcal{M}_0 \in \LocSys^{\Glo}_{\fcat T}(Y,\fcat{A})$. Hence, the assertion follows by the second part of \cref{proposition: genuine sharp} and \cref{corollary: left adjoint faithful global linear}.

Next, if $Y\in \Type^\Glo$ is arbitrary, then $Y\simeq \colim Y_\alpha$, where $Y_{\alpha}$ is a representable global space. Let $p_{\alpha} \colon Y_{\alpha} \to Y$ and $q_{\alpha}\colon X_{\alpha}=Y_{\alpha} \times_Y X \to X$ be the projections. Then it is enough to show that $p^{*}_{\alpha}(\theta)$ is an equivalence for any $\alpha$. However, by the previous paragraph and \cref{proposition: genuine sharp}, we observe
\begin{align*}
p^{\gen,*}_{\alpha}f^{\gen}_{\#}(\mathcal L \otimes f^*(\mathcal M)) &\simeq f^{\gen}_{\alpha,\#} q_{\alpha}^{\gen,*}(\mathcal L \otimes f^{\gen,*}(\mathcal M)) \\
&\simeq f^{\gen}_{\alpha,\#} (q_{\alpha}^{\gen,*}(\mathcal L) \otimes q_{\alpha}^{\gen,*}f^{\gen,*}(\mathcal M)) \\
&\simeq f^{\gen}_{\alpha,\#} (q_{\alpha}^{\gen,*}(\mathcal L) \otimes f^{\gen,*}_{\alpha}p^{\gen,*}_{\alpha}(\mathcal M)) \\
&\simeq f^{\gen}_{\alpha,\#} q_{\alpha}^{\gen,*}(\mathcal L) \otimes p^{\gen,*}_{\alpha}(\mathcal M)\\
&\simeq p^{\gen,*}_{\alpha}(f^{\gen}_{\#}(\mathcal L) \otimes \mathcal M).\qedhere
\end{align*}
\end{proof}

\begin{ex}\label{example: omegainfty is not conservative}
We note that the statement of~\cref{cor_OmegaInftygen} is false if $X$ is \emph{not} an orbispace. Indeed, let $G=C_p$ be a cyclic group of prime order $p$, let $X=\ast \coprod_{BG} \ast \in \Type^{\Glo}$ be the suspension of $BG$, and let $\Mod_{\mathbb{Q}}\colon \Orb^{\op} \to \Prs^\LL$ denote the constant coefficient system with the value $\Mod_{\mathbb{Q}}$. Then, by \cref{prop_LSglo_preserves_limits} and \cref{construction: genuine local system}, the following commutative squares
\[\xymatrix{
\LocSys^\Glo(X, \Mod_{\mathbb{Q}}) \ar[r] \ar[d] & \Mod_{\mathbb{Q}} \ar[d]^-{\pi^{\Glo,*}} \\
\Mod_{\mathbb{Q}} \ar[r]^-{\pi^{\Glo,*}} & \LocSys^\Glo(BG, \Mod_{\mathbb{Q}}),
}\qquad
\xymatrix{
\LocSys^\gen(X, \Mod_{\mathbb{Q}}) \ar[r] \ar[d] & \Mod_{\mathbb{Q}} \ar[d]^-{\pi^{\gen,*}} \\
\Mod_{\mathbb{Q}} \ar[r]^-{\pi^{\gen,*}} & \LocSys^\gen(BG, \Mod_{\mathbb{Q}}),
}\]
are pullbacks, where $\pi\colon BG \to \ast$ is the projection. By \cref{corollary: pushforward for orbits}, the functor $\pi^{\Glo,*}$ is a fully faithful, so $$\LocSys^\Glo(X, \Mod_{\mathbb{Q}}) \simeq \Mod_{\mathbb{Q}}.$$ However, by \cref{example: gen gluing functor rational} below, we have $\LocSys^\gen(BG, \Mod_{\mathbb{Q}}) \simeq \Mod_{\mathbb{Q}} \times \Fun(BC_p,\Mod_{\mathbb{Q}})$ and $\pi^{\gen,*}\simeq \Id\times \pi^*$. Therefore, we identify
$$\LocSys^\gen(X, \Mod_{\mathbb{Q}}) \simeq \Fun(\Sigma(BG_+),\Mod_{\mathbb{Q}}). $$
Moreover, under these identifications, the functor $\Sigma^{\infty}_{X}$ is equivalent to the pullback
$$q^*\colon \Mod_{\mathbb{Q}} \tto \Fun(\Sigma(BG_{+}),\Mod_{\mathbb{Q}})$$
along the projection $q\colon \Sigma(BG_{+}) \to \ast$. Therefore, $\Omega^{\infty}_{X} \simeq q_*$, which is not conservative.

Similarly, one can construct an example of a global space $Y$ such that $\Omega^{\infty}_Y$ is \emph{not} continuous. For instance, let $Y=\bigvee_{i=1}^{\infty}X=\bigvee_{i=1}^{\infty}\Sigma BC_p$ be the wedge of infinite copies of $X$. Then, we have
$$\LocSys^\Glo(Y, \Mod_{\mathbb{Q}}) \simeq \Mod_{\mathbb{Q}}, \qquad  \LocSys^\gen(Y, \Mod_{\mathbb{Q}}) \simeq \Fun\left(\bigvee_{i=1}^{\infty} \Sigma(BG_+),\Mod_{\mathbb{Q}}\right).$$
Again, under these identifications, the functor $\Omega^{\infty}_Y$ is equivalent to the pushforward $\tilde{q}_*$ along the projection $\tilde{q}\colon \bigvee_{i=1}^{\infty} \Sigma(BG_+) \to \ast$. Since $\bigvee_{i=1}^{\infty} \Sigma(BG_+)$ is not (rationally) of finite type, the pushforward $\tilde{q}_*$ does not preserve all colimits.
\end{ex}

\subsection{Genuine local systems on \texorpdfstring{$G$}{G}-spaces}\label{section: LSgen on G-spaces}

%\todo{RELATIVE VERSION! LEAVE ONLY G-SPACES, TWO DESCRIPTIONS}

Our definition of genuine local systems (\cref{construction: genuine local system}) may look ad hoc, so we discuss in this section alternative ways to define $\LocSys^\gen_{\fcat T}(X,\fcat A)$. Motivated by the definition of $\Sp^G$, one is temped to give the following construction.

\begin{defn}\label{definition: naive genuine}
Let $X$ be a global space and let $\fcat A$ be a strongly continuous symmetric monoidal coefficient system of Beck--Chevalley type. Define a \emdef{naive genuine equivariant local systems} on $X$ as
$$\LocSys^{\ngen}_{\fcat T}(X, \fcat A) := \LocSys^\Glo_{\fcat T}(X, \fcat A)[\{\left(\mathbbl 1^E\right)^{-1}\}_{E\in \Vect_X}].$$
\end{defn}

Unfortunately, this construction does not produce satisfactory answer in general, e.g. because of the lack of enough of vector bundles on a general global space, see \Cref{ex_LS_ngen_for_N}. However, it gives the good category in the case of global quotients, see \Cref{prop_LS_gen_vs_ngen_global_quotients}.

%The following example shows that the category obtain from globally equivariant local systems by inverting Thom local systems of all vector bundles does not produce interesting enough theory in general.
\begin{ex}\label{ex_LS_ngen_for_N}
Let $\mathcal N \in \Type^{\Glo}$ be the normal subgroup classifier from \cref{example: normal subgroup classifier}. By  construction, any map $\mathcal N \to BG$ factors through a point for every compact Lie group $G$, see also \cite[Exercise~4.1.1]{Rezk_GlobalHomotopy}. In particular, all vector bundles on $\mathcal N$ are trivial. It follows that
$$\LocSys^\Glo(\mathcal N, \fcat A) \xymatrix{\ar[r]^-\sim &} \LocSys^\ngen(\mathcal N, \fcat A).$$
On the other hand,  $\mathcal{N} \simeq \colim_{\Orb^{\rep}}BG$ by \cref{example: normal subgroup classifier as a colimit}. Therefore,
$$\LocSys^{\gen}(\mathcal N, \fcat A) \simeq \lim_{BG\in\Orb^{\rep,\op}} \LocSys^{\gen}(BG,\fcat A). $$
In particular, if $\fcat A = \underline{\Sp}$ is a constant coefficient system, then
$$\LocSys^{\gen}(\mathcal N, \underline{\Sp}) \simeq \lim_{BG\in\Orb^{\rep,\op}} \Sp^G, $$
cf. \cref{ex_LSglo_for_N}. Once again, a similar category was studied in \cite{LNP25}.
\end{ex}

Another way to define genuine local systems is to make the definition relative to a certain base.
\begin{defn}\label{def_ngenB_LS}
Let $p\colon X \to B$ be a morphism of global homotopy types. We define
$$\LocSys^{\ngen/B}_{\fcat T}(X, \fcat A) := \LocSys^\Glo_{\fcat T}(X, \fcat A)\otimes_{\LocSys^\Glo_{\fcat T}(B, \fcat A)} \LocSys^\ngen_{\fcat T}(B, \fcat A).$$
Equivalently,
$$\LocSys^{\ngen/B}_{\fcat T}(X, \fcat A) := \LocSys^\Glo_{\fcat T}(X, \fcat A)\left[\left\{\left(\mathbbl 1^{p^*E}\right)^{-1}\right\}_{E\in \Vect_B}\right].$$
\end{defn}
\begin{ex}
If $B = *$ is a point, then $\LocSys^{\ngen/*}_{\fcat T}(X, \fcat A) \simeq \LocSys^\Glo_{\fcat T}(X, \fcat A)$. In the other extreme case, $B = X$, we have $\LocSys^{\ngen/X}_{\fcat T}(X, \fcat A) \simeq \LocSys^\ngen_{\fcat T}(X, \fcat A)$.
\end{ex}
\begin{ex}\label{example: relative genuine local systems over BG}
Assume that any vector bundle on $X$ embeds as a direct summand into a pullback of a vector bundle on $B$. Then, for every $E \in \Vect_X$, the Thom local system $\mathbbl 1_X^E$ is already $\otimes$-invertible in $\LocSys^{\ngen/B}_{\fcat T}(X, \fcat A)$, see \cref{Thom_LS_basics_Glo}. Hence the natural map
$$\LocSys^{\ngen/B}_{\fcat T}(X, \fcat A) \tto \LocSys^{\ngen}_{\fcat T}(X, \fcat A)$$
is an equivalence.

For instance, this is the case if $B = BG$, where $G$ is a compact Lie group $G$, and $X = Y\gitq G$, where $Y$ is the underlying $G$-homotopy type of a compact $G$-topological space, see e.g.\@ \cite[Proposition~2.4]{Segal68} and \cref{notation: gitq}.
\end{ex}

\begin{lem}\label{lem_ngenB_colimit_of_types}
Let $X_\bullet\colon \fcat{I} \to \Type^\Glo_{/^\rep B}$ be a small diagram of global spaces faithful over $B \in \Type^\Glo$. Then the functor induced by $*$-pullbacks
$$\LocSys^{\ngen/B}_{\fcat T}\left(\colim_{i\in\fcat I} X_i, \fcat A\right) \tto \lim_{i\in \fcat{I}^{\op}} \LocSys^{\ngen/B}_{\fcat T}(X_i, \fcat A)$$
is an equivalence.

\begin{proof}
By the 2-out-of-3 property for faithful morphisms (see \cref{proposition: basic properties of faithful}), all transition morphisms in $X_{\bullet}$ are faithful. Therefore, by \cref{corollary: left adjoint faithful global},
$$\LocSys^{\Glo}_{\fcat T}\left(\colim_{i\in \fcat{I}}X_i,\fcat A\right) \simeq \lim_{i\in \fcat{I}^{\op}} \LocSys^{\Glo}_{\fcat T}(X_i, \fcat A) \simeq \colim_{i\in \fcat{I}} \LocSys^{\Glo}_{\fcat T}(X_i, \fcat A) \in \Prs^\LL,$$
where the last colimit is taken along the $\#$-pushforwards. Moreover, by \cref{corollary: left adjoint faithful global linear}, the transition functors in the colimit are $\LocSys^{\Glo}_{\fcat T}(B,\fcat{A})$-linear, so 
$$\LocSys^{\Glo}_{\fcat T}\left(\colim_{i\in \fcat{I}}X_i,\fcat A\right) \simeq \colim_{i\in \fcat{I}} \LocSys^{\Glo}_{\fcat T}(X_i, \fcat A) \in \Mod_{\LocSys^{\Glo}_{\fcat T}(B,\fcat{A})}(\Prs^\LL).$$
The assertion follows by tensoring the last equivalence with $\LocSys^{\ngen}_{\fcat T}(B,\fcat{A})$ over $\LocSys^{\Glo}_{\fcat T}(B,\fcat{A})$.
%By the second part of \Cref{lem_basics_ngenB} we have a commutative diagram
%$$\xymatrix{
%\LocSys^{\ngen/B}(\colim_I X_i, \fcat A) \ar[r]\ar[d] & \lim_I \LocSys^{\ngen/B}(X_i, \fcat A) \ar[d] \\ 
%\LocSys^\Glo(\colim_I X_i, \fcat A) \ar[r]^-\sim & \lim_I \LocSys^\Glo(X_i, \fcat A)
%}$$
%where the vertical functors are conservative by \Cref{cor_OmegaInftyNgen} and the lower horizontal functor is an equivalence by \Cref{prop_LSglo_preserves_limits}. It follows that the upper vertical functor is conservative. Moreover, since by the last part of \Cref{lem_basics_ngenB} the Beck--Chevalley condition holds, the result follows from \cite[Theorem 4.7.5.2]{Lur_HA}.
\end{proof}
\end{lem}

\begin{prop}\label{prop_LS_gen_vs_ngen_global_quotients}
Let $G$ be a compact Lie group and let $X \in \Type^G $ be a $G$-homotopy type. Let $X\gitq G \in \Type^\Glo$ denote the corresponding global space. Then the natural functors
$$\LocSys^{\ngen/BG}_{\fcat T}(X \gitq G, \fcat A) \tto \LocSys^\ngen_{\fcat T}(X \gitq G, \fcat A) \tto \LocSys^\gen_{\fcat T}(X \gitq G, \fcat A) $$
are equivalences. %Moreover, both parts are equivalent to
%$$\LocSys^{\ngen/BG}(X \gitq G, \fcat A) = \LocSys^\Glo(X \gitq G, \mathcal A) \otimes_{\LocSys^\Glo(BG, \fcat A)} \LocSys^\gen(BG, \fcat A).$$

\begin{proof}
First, we will show that the natural functor
$$\LocSys^{\ngen/BG}_{\fcat T}(X\gitq G, \fcat A) \tto \LocSys^\ngen_{\fcat T}(X \gitq G, \fcat A)$$
is an equivalence. By definition, it suffices to show that the Thom local system $\mathbbl 1^E_{X\gitq G}$ is already $\otimes$-invertible in $\LocSys^{\ngen/BG}_{\fcat T}(X\gitq G, \fcat A)$ for every $E \in \Vect_{X\gitq G}$. By \cref{lem_ngenB_colimit_of_types}, we have
$$\LocSys^{\ngen/BG}_{\fcat T}(X\gitq G, \fcat A) \simeq \lim_{G/H \in \Orb_{G, \fcat T/X}^\op} \LocSys^{\ngen/BG}_{\fcat T}(BH, \fcat A).$$
Hence, it is enough to show that the fiber $x^{\Glo,*}(\mathbbl 1^E_{X\gitq G})$ is $\otimes$-invertible for every $G$-equivariant map $x\colon G/H \to X$, where $H \in \fcat T$. By \cref{Thom_LS_basics_Glo},
$$x^{\Glo,*}(\mathbbl 1^E_{X\gitq G}) \simeq \mathbbl 1^{x^* E}_{BH} \in \LocSys^{\ngen/BG}(BH, \fcat A).$$
Since
$$\LocSys^{\ngen/BG}_{\fcat T}(BH, \fcat A) \simeq \LocSys^{\ngen}_{\fcat T}(BH, \fcat A)$$
by \cref{example: relative genuine local systems over BG}, the Thom local system $\mathbbl 1^{x^* E}_{BH}$ is $\otimes$-invertible in $\LocSys^{\ngen/BG}_{\fcat T}(BH, \fcat A)$ according to the definition of $\LocSys^{\ngen}_{\fcat T}(BH, \fcat A)$.
%and in the later category this holds by construction.

Finally, we will show that the composite
$$\LocSys^{\ngen/BG}_{\fcat T}(X\gitq G, \fcat A) \tto \LocSys^\gen_{\fcat T}(X \gitq G, \fcat A)$$
is an equivalence. By \cref{def_ngenB_LS} and \cref{def_genuine_LS}, the natural transformation of functors
$$\LocSys^{\ngen/BG}_{\fcat T}(-, \fcat A)|_{\Type^G}, \LocSys^\gen_{\fcat T}(-, \fcat A)|_{\Type^G} \colon \xymatrix{\Type^{G, \op} \ar[r]^-{-\gitq G} & \Type^{\Glo, \op} \ar[r] & \Prs^\LL}$$
is an equivalence when restricted to the full subcategory $\Orb_G \subset \Type^G$ spanned by $G$-orbits. Note that both functors  $\LocSys^{\ngen/BG}_{\fcat T}(-, \fcat A)|_{\Type^G}$ and $\LocSys^\gen_{\fcat T}(-, \fcat A)|_{\Type^G}$ preserve small limits by \cref{lem_ngenB_colimit_of_types} and \cref{construction: genuine local system}, respectively. The assertion follows since any $G$-space is a colimit of $G$-orbits. % Since any $X \in \Type^G$ is a small colimit of $G$-orbits it is enough to show that both $\LocSys^{\ngen}(-, \mathcal A)_{|\Type^G}$ and $\LocSys^\gen(-, \mathcal A)_{|\Type^G}$ preserve small limits. For $\LocSys^\gen(-, \mathcal A)_{|\Type^G}$ this follows by construction and the fact that the functor $\Type^G \to \Type^\Glo$ is colimit preserving. For $\LocSys^{\ngen}(-, \mathcal A)_{|\Type^G}$ we will see this by comparing it with $\LocSys^{\ngen/BG}(X//G, \mathcal A)$ and using \Cref{lem_ngenB_colimit_of_types} in the later case.
\end{proof}
\end{prop}

\begin{ex}\label{example: G-spaces geom type gen}
Let $G$ be a compact Lie group from $\fcat T$ and let $X \in \Type^G$ be a $G$-space. Suppose that the coefficient system $\fcat{A}$ is of geometric type. Then, by \cref{example: locsysglo geom type general}, \cref{example: geom type gen BG}, and \cref{prop_LS_gen_vs_ngen_global_quotients}, we have
$$\LocSys^{\gen}_{\fcat T}(X\gitq G,\fcat{A}) \simeq \Mod_{\fcat{A}_G}\left(\LocSys^{\gen}_{\fcat T}(X\gitq G,\underline{\fcat{A}(BG)})\right), $$
where $\fcat{A}_G \in \CAlg(\LocSys^{\Glo}_{\fcat T}(BG,\underline{\fcat{A}(BG)}))$ is a commutative algebra as in \cref{example: locsysglo geometric type}.
\end{ex}

\begin{cor}\label{corollary: genuine local systems}
Let $X \in \Type^{\Glo}$ be an orbispace. Then there are a natural equivalences
\begin{align*}
\LocSys^{\gen}_{\fcat T}(X,\fcat{A}) &\simeq \lim_{BG\in \Orb^{\rep,\op}}\LocSys^{\gen}_{\fcat T}(X\times_{\mathcal{N}} BG,\fcat{A}) \\
&\simeq \lim_{BG\in \Orb^{\rep,\op}} \LocSys^{\Glo}_{\fcat T}(X\times_{\mathcal{N}} BG,\fcat{A})\otimes_{\LocSys^{\Glo}_{\fcat T}(BG,\fcat{A})}\LocSys^{\gen}_{\fcat T}(BG,\fcat{A}).
\end{align*}
\end{cor}

\begin{proof}
The first equivalence follows by \cref{example: any orbispace is a colimit of G-spaces} and \cref{construction: genuine local system}. The second equivalence follows by \cref{proposition: G-spaces and orbispace} and \cref{prop_LS_gen_vs_ngen_global_quotients}.
\end{proof}

In the rest of the section we show that the category $\LocSys^{\gen}_{\fcat T}(BG,\fcat{A})$ admits a nice set of generators if the coefficient system $\fcat{A}$ is of geometric type. We recall the following classical result for the constant coefficient system first.

\begin{thm}\label{theorem: genuine BG rigid}
Let $G$ be a compact Lie group. Then the category $\LocSys^{\gen}(BG,\underline{\Sp}) \simeq \Sp^G$ is a rigid compactly generated presentably symmetric monoidal category. More precisely, the objects $p^\gen_{\#}\mathbbl 1_{BH}$, where $p\colon BH \to BG$ ranges over the set of faithful morphisms, constitute a set of compact (and dualizable) generators of $\Sp^G$.
\end{thm}

\begin{proof}
Note that $\Hom_{\Sp^G}(p^\gen_{\#}\mathbbl 1_{BH},-) \simeq (\Omega^{\infty}(-))^H$. Therefore, by \cref{cor_OmegaInftyNgen} and \cref{LSGlo_on_orbispaces}, the set $\{p^\gen_{\#}\mathbbl 1_{BH}\}_{p\in \Orb^{\rep}_{/BG}}$ is a set of compact generators. Finally, each generator $p^\gen_{\#}\mathbbl 1_{BH} \simeq \Sigma^{\infty}_{BG} (G/H)_+$ is dualizable by the equivariant Atiyah duality, see~\cite[Theorem~III.5.1]{LMS86}
\end{proof}

By using the theory of \cite{Ramzi26} and \cref{example: geom type gen BG}, we can generalize \cref{theorem: genuine BG rigid} to an arbitrary coefficient system of geometric type. Let $\fcat{V}\in \CAlg(\Prs^\LL)$ be a presentably symmetric monoidal category and let $\fcat{W} \in \CAlg_{\fcat{V}}(\Prs^\LL)$ be a $\fcat{V}$-algebra. Recall that $\fcat{W}$ is called \emdef{locally rigid over $\fcat{V}$} if $\fcat{W}$ is dualizable over $\fcat{V}$ and the multiplication map $\fcat{W}\otimes_{\fcat{V}} \fcat{W} \to \fcat{W}$ admits a continuous strict $\fcat{V}$-linear right adjoint (see~\cite[Definition~4.5]{Ramzi26}). An object $x\in \fcat{W}$ is called \emdef{$\fcat{V}$-atomic} if the internal Hom-functor $$[x,-]_{\fcat{V}} \colon \fcat{W} \tto \fcat{V}$$ commutes with arbitrary colimits and strict $\fcat{V}$-linear (see~\cite[Definition~1.22]{Ramzi24_dualizable}). The $\fcat{V}$-algebra $\fcat{W}$ is called \emdef{rigid over $\fcat{V}$} if $\fcat{W}$ is locally rigid over $\fcat{V}$ and the monoidal unit $\mathbbl{1}_{\fcat{W}} \in \fcat{W}$ is $\fcat{V}$-atomic (see~\cite[Definition~4.36]{Ramzi26}).

Recall that the category $\fcat{W}$ is \emdef{atomically generated} if $\fcat{W}$ is the smallest full $\fcat{V}$-submodule of $\fcat{W}$ containing all atomic objects and closed under colimits (see~\cite[Definition~1.27]{Ramzi24_dualizable}). Finally, we recall from~\cite[Example~4.37]{Ramzi26} that if $\fcat{W}$ is atomically generated, the monoidal unit $\mathbbl{1}_{\fcat{W}} \in \fcat{W}$ is $\fcat{V}$-atomic, and every $\fcat{V}$-atomic object is dualizable, then $\fcat{W}$ is rigid over $\fcat{V}$.

\begin{cor}\label{corollary: genuine BG rigid over ABG}
Let $G$ be a compact Lie group from $\fcat T$. Suppose that the coefficient system $\fcat{A}$ is of geometric type. Then the category $\LocSys^{\gen}_{\fcat T}(BG,\fcat{A}) \in \CAlg_{\fcat{A}(BG)}(\Prs^\LL)$ is a rigid $\fcat{A}(BG)$-atomically generated $\fcat{A}(BG)$-algebra. More precisely, the objects $p^\gen_{\#}\mathbbl 1_{BH}$, where $p\colon BH \to BG$ ranges over the set of faithful morphisms, constitute a set of $\fcat{A}(BG)$-atomic (and dualizable) generators of a $\fcat{A}(BG)$-module $\LocSys^{\gen}_{\fcat T}(BG,\fcat{A})$.
\end{cor}

\begin{proof}
By \cref{theorem: genuine BG rigid}, it is clear that the free $\fcat{A}(BG)$-linear category 
$$\LocSys^{\gen}_{\fcat T}(BG,\underline{\fcat{A}(BG)})\simeq \Sp^G \otimes \fcat{A}(BG)$$
satisfies all the desired properties. By \cref{example: geom type gen BG}, we have
$$\LocSys^{\gen}_{\fcat T}(BG,\fcat{A}) \simeq \Mod_{\fcat{A}_G}\left(\LocSys^{\gen}_{\fcat T}(BG,\underline{\fcat{A}(BG)})\right) $$
where $\fcat{A}_G \in \CAlg(\LocSys^{\Glo}_{\fcat T}(BG,\underline{\fcat{A}(BG)}))$ is the commutative algebra from \cref{example: locsysglo geometric type}. By \cite[Example~4.38]{Ramzi26}, the symmetric monoidal category $\LocSys^{\gen}_{\fcat T}(BG,\fcat{A})$ is rigid over the symmetric monoidal category $\LocSys^{\gen}_{\fcat T}(BG,\underline{\fcat{A}(BG)})$ generated by $\fcat{A}_G$ as a $\LocSys^{\gen}_{\fcat T}(BG,\underline{\fcat{A}(BG)})$-module. Finally, the assertion follows by \cite[Example~4.45]{Ramzi26}.
\end{proof}

\begin{rem}\label{remark: genuine BG norm map}
Let $p\colon BH \to BG$ be a faithful morphism. We recall that the equivariant Atiyah duality used in \cref{theorem: genuine BG rigid} is closely related to the Wirthm\"{u}ller isomorphism relating the \emph{left} adjoint $p^\gen_{\#}$ with the \emph{right} adjoint $p^{\gen}_*$, see \cite[\S II.6]{LMS86} or \cite[\S 4.2]{Cnossen23} for a detailed account. We recall the construction of the \emph{norm map} here. Consider the commutative diagram \begin{equation}\label{eq_star_shriek_BKBHBG}
\xymatrix{
BH \ar[r]^-\Delta & BH \times_{BG} BH \ar[d]^s \ar[r]^-r & BH \ar[d]^p \\
 & BH \ar[r]^-p & BG.
}\end{equation}
Let $D_p\in \LocSys^{\gen}(BH,\fcat{A})$ denote the genuine local system $s^{\gen}_*\Delta^{\gen}_{\#} \mathbbl{1}_{BH}$. Then there is a natural transformation
$$\varepsilon\colon s^{\gen,*} D_p \tto \Delta^{\gen}_{\#} \mathbbl{1}_{BH}.$$
Then by the base change and the projection formulas, we obtain the natural transformation
\begin{align*}
p^{\gen,*}p^{\gen}_{\#}(-\otimes D_p) \simeq r^{\gen}_{\#}s^{\gen,*}(-\otimes D_p) &\simeq r^{\gen}_{\#}(s^{\gen,*}(-)\otimes s^{\gen,*}(D_p)) \\
&\xrightarrow{\varepsilon} r^{\gen}_{\#}(s^{\gen,*}(-)\otimes \Delta^{\gen}_{\#} \mathbbl{1}_{BH}) \\
&\simeq r^{\gen}_{\#}\Delta^{\gen}_{\#}(\Delta^{\gen,*}s^{\gen,*}(-))\\
&\simeq (-).
\end{align*}
Using the adjoint pair $p^{\gen,*}\dashv p^{\gen}_{*}$, we obtain the norm map
\begin{equation}\label{equation: norm map}
\mathrm{Nm}^{\fcat{A}}_p \colon p^{\gen}_{\#}(-\otimes D_p) \tto p^{\gen}_*(-).
\end{equation}
If $\fcat{A}=\underline{\Sp}$ is the constant coefficient system, then the norm map $\mathrm{Nm}^{\Sp}_p$ is an equivalence by~\cite[\S II.6]{LMS86} or~\cite[Theorem~18.6.5]{MS06}. Moreover, the genuine local system $D_p$ is dualizable and $D_p \simeq \mathbb{S}^{-\mathfrak{g}/\mathfrak{h}}$, where $\mathfrak{g}/\mathfrak{h}$ is the tangent space of the $H$-manifold $G/H$ at the point~$eH$. By the same argument as in \cref{corollary: genuine BG rigid over ABG}, one can show that $\Nm_p^\fcat{A}$ is an equivalence for any coefficient system $\fcat{A}$ of geometric type and $D_p \simeq \mathbbl{1}^{-\mathfrak{g}/\mathfrak{h}}_{BH}$.
\end{rem}

\begin{rem}\label{remark: generating set not geom type}
At the moment of writing, we are not aware if the analog of \cref{corollary: genuine BG rigid over ABG} is true for a coefficient system $\fcat{A}$ which is \emph{not} of geometric type. Similarly, we are not aware if the norm map~\eqref{ex_LS_ngen_for_N} is an equivalence for a coefficient system which is \emph{not} of geometric type.
\end{rem}

\subsection{Geometric fixed points}\label{section: geometric fixed points}
In this section, we construct the functor geometric fixed points from the category of genuine local systems on $BG$. Using this functor, we will obtain a semi-orthogonal decomposition  for the category $\LocSys^{\gen}_\fcat{T}(BG,\fcat{A})$ (\cref{corollary: isotropy separation for gen}), which generalizes the classical isotropy separation fiber sequence for $\Sp^G$.

Throughout this section, $\fcat{A}\colon \Orb^{\op}_{\fcat{T}} \to \Prs^{\LL}$ is a stable $E_\infty$-monoidal strongly continuous coefficient system of Beck--Chevalley type, see~\cref{definition: coefficient system}.

\begin{construction}\label{construction: geometric fixed points}
Let $G\in \fcat{T}$ be a compact Lie group and $V$ be a real $G$-representation. Then, by \cref{prop_LSGlo_sharp_pushforward}, there is an equivalence $(\mathbbl{1}^V_{BG})^G \simeq \Sigma^{\dim V^{G}} \mathbbl{1} \in \fcat{A}(BG)$. Therefore, the symmetric monoidal evaluation functor
$$(-)^G\colon \LocSys^{\Glo}_{\fcat{T}}(BG,\fcat{A}) \tto \fcat{A}(BG) $$
maps any Thom local system on $BG$ to an invertible object. By \cref{definition: object inversion} and \cref{def_genuine_LS}, the functor $(-)^G$ factors via the functor
$$\Phi^G_{\fcat{A}} \colon \LocSys^{\gen}_{\fcat{T}}(BG,\fcat{A}) \tto \fcat{A}(BG)$$
which we call the \emdef{geometric fixed points} functor.
\end{construction}

\begin{rem}\label{remark: gf for constant}
Let $\fcat{A}=\underline{\Sp}$ be the constant coefficient system. Then $\LocSys^{\gen}(BG,\underline{\Sp})\simeq \Sp^G$ and the functor $\Phi^G_{\underline{\Sp}}$ coincides with the usual functor of geometric fixed points $\Phi^G$.
\end{rem}

\begin{prop}\label{proposition: basic prop of gfp}
Let $G\in \fcat{T}$ be a compact Lie group. Then the functor $\Phi^{G}_{\fcat{A}}$ satisfies the following properties
\begin{enumerate}
\item $\Phi^G_{\fcat{A}}$ is symmetric monoidal and colimit-preserving;
\item $\Phi^G_{\fcat{A}}\circ \Sigma^{\infty}_{BG} \simeq (-)^G$;
\item let $f\colon BH \to BG$ be induced by an inclusion of a closed proper subgroup $H \subset G$, then $\Phi^G_{\fcat{A}} \circ f^{\gen}_{\#} \simeq 0$;
\item  there is a natural equivalence
$$\Phi^{G}_{\fcat{A}}(\mathcal{L}) \simeq \colim_{V \in \Vect_{BG}} \Sigma^{-\dim V^G}(\mathbbl{1}^V_{BG} \otimes^{\gen} \mathcal{L})^G$$
for a genuine local system $\mathcal{L} \in \LocSys^{\gen}_{\fcat{T}}(BG,\fcat{A})$.
\end{enumerate}
\end{prop}

\begin{proof}
The first two parts follows immediately from \cref{construction: geometric fixed points}. We prove now the third claim. Indeed, by \cref{cor_OmegaInftyNgen}, it suffices to show that $\Phi^G_{\fcat{A}} \circ f^{\gen}_{\#} \circ \Sigma^{\infty}_{BH} \simeq 0. $
However, by \cref{lemma: sharp for gen orbits}, we have 
$$\Phi^G_{\fcat{A}} \circ f^{\gen}_{\#} \circ \Sigma^{\infty}_{BH} \simeq \Phi^G_{\fcat{A}}  \circ \Sigma^{\infty}_{BG} \circ f^{\Glo}_{\#} \simeq (f^{\Glo}_{\#}(-))^G.$$
Finally, the assertion follows by \cref{prop_LSGlo_sharp_pushforward}.

By \cref{lemma: X-linearization of sigma-omega}, there is a natural equivalence
$$\mathcal{L} \simeq \colim_{V\in \Vect_{BG}}\left( (\mathbbl{1}^V_{BG})^{-1} \otimes \Sigma^{\infty}_{BG} \Omega^{\infty}_{BG} (\mathbbl{1}^V_{BG} \otimes \mathcal{L})\right). $$
This implies the last part.
\end{proof}

\begin{cor}\label{corollary: gfp is localization}
Let $G \in \fcat{T}$ be a compact Lie group. Then the functor $\Phi^G_{\fcat{A}}$ is a localization such that the right adjoint
$$\Xi^{G}_{\fcat{A}} \colon \fcat{A}(BG) \tto \LocSys^{\gen}_{\fcat{T}}(BG,\fcat{A}) $$
commutes with small colimits.
\end{cor}

\begin{proof}
By \cref{proposition: basic prop of gfp}(1) and the adjoint functor theorem, the functor $\Phi^G_{\fcat{T}}$ admits a right adjoint
$$\Xi^{G}_{\fcat{A}} \colon \fcat{A}(BG) \tto \LocSys^{\gen}_{\fcat{T}}(BG,\fcat{A}). $$
By passing to right adjoints in \cref{proposition: basic prop of gfp}(2) and \cref{lemma: sharp for gen orbits}, we obtain
$$(\Omega^{\infty}_{BG}\Xi^G_{\fcat{A}}(\mathcal F))^H \simeq 
\begin{cases}
\mathcal F \in \fcat{A}(BG) & \mbox{if $H = G$,}\\
0 & \mbox{if $H$ is a proper subgroup}
\end{cases}
$$
for all $\mathcal F \in \fcat{A}(BG)$. This implies that the functor $\Xi^G_{\fcat{A}}$ commutes with small colimits, see \cref{cor_OmegaInftyNgen}.

Since the Thom local systems $\mathbbl{1}^V_{BG} \in \LocSys^{\gen}_{\fcat{T}}(BG,\fcat{A})$ are $\otimes$-invertible, the canonical morphism
$$\mathbbl{1}^V_{BG} \otimes \Xi^G_{\fcat{A}}(\mathcal F) \tto \Xi^G_{\fcat{A}}(\Phi^G_{\fcat{A}}(\mathbbl{1}^V_{BG}) \otimes \mathcal F) \simeq \Xi^G_{\fcat{A}}(\Sigma^{\dim V^G} \mathcal F)$$
is an equivalence for all $\mathcal F \in \fcat{A}(BG)$ and $V \in \Vect_{BG}$. By the last part of \cref{proposition: basic prop of gfp}, we observe that the counit map $ \Phi^G_{\fcat{A}} \Xi^G_{\fcat{A}} \to \Id$ is an equivalence. Indeed, for any $\mathcal F \in \fcat{A}(BG)$, we have
\begin{align*}
\Phi^G_{\fcat{A}} \Xi^G_{\fcat{A}}(\mathcal F) &\simeq \colim_{V \in \Vect_{BG}} \Sigma^{-\dim V^G}(\mathbbl{1}^V_{BG} \otimes \Xi^{G}_{\fcat{A}}(\mathcal F))^G \simeq \\
&\simeq \colim_{V \in \Vect_{BG}} \Sigma^{-\dim V^G}(\Xi^{G}_{\fcat{A}}(\Sigma^{\dim V^G}\mathcal F))^G \simeq (\Xi^{G}_{\fcat{A}}(\mathcal F))^G \simeq \mathcal F.\qedhere
\end{align*}
\end{proof}

By \cref{proposition: basic prop of gfp}, the counit $\Sigma^{\infty}_{BG}\Omega^{\infty}_{BG} \to \Id$ induces the canonical morphism
$$\theta_{\mathcal{L}} \colon \mathcal{L}^G \tto \Phi^{G}_{\fcat{A}}(\mathcal{L}) $$
for a genuine local system $\mathcal{L} \in \LocSys^{\gen}_{\fcat{T}}(BG,\fcat{A})$.
\begin{thm}\label{theorem: gfp_is_fp}
Let $G\in \fcat{T}$ be a compact Lie group and let $\mathcal{L} \in \LocSys^{\gen}_{\fcat{T}}(BG,\fcat{A})$ be a genuine local system on $BG$. Suppose that $\mathcal{L}^H=0$ for all proper subgroups $H\subset G$. Then the canonical morphism $\theta_{\mathcal{L}}$ is an equivalence.
\end{thm}

\begin{proof}
Let $\mathcal{C} \subset \LocSys^{\gen}_{\fcat{T}}(BG,\fcat{A})$ be a full subcategory spanned by genuine local systems $\mathcal{L}$ such that $\mathcal{L}^H=0$ for all proper subgroups $H\subset G$. By \cref{cor_OmegaInftyNgen}, the subcategory $\mathcal{C}$ is closed under small colimits. Moreover, $\mathcal{C}$ is presentable as the intersection of kernels of colimit-preserving functors between presentable categories.

Similarly, let $\mathcal{C}' \subset \LocSys^{\Glo}_{\fcat{T}}(BG,\fcat{A})$ be a full subcategory spanned by globally equivariant local systems $\mathcal{M}$ such that $\mathcal{M}^H=0$ for all proper subgroups $H\subset G$. By the diagram~\eqref{equation: pullback for gen orbits}, the essential image of $\mathcal{C}'$ under the functor $\Sigma^{\infty}_{BG}$ is contained in the full subcategory $\mathcal{C}$. First, we will show that $\mathcal{C}$ is generated under colimits by $\Sigma^{\infty}_{BG}\mathcal{C'}$. Indeed, let $\mathcal{M} \in \mathcal{C}'$ and $\mathcal{L} \in \mathcal{C}$, then 
$$\Hom_{\LocSys^\gen(BG,\fcat{A})}(\Sigma^{\infty}_{BG}\mathcal{M},\mathcal{L}) \simeq \Hom_{\LocSys^\Glo(BG,\fcat{A})}(\mathcal{M},\Omega^{\infty}_{BG}\mathcal{L}) \simeq \Hom_{\fcat{A}(BG)}(\mathcal{M}^G,\mathcal{L}^G), $$
where the last equivalence follows by \cref{proposition: global sections semi-orthogonal}. Let $\mathcal{M} \in \LocSys^\Glo(BG,\fcat{A})$ be given by the formula
$$
\mathcal{M}^H=
\begin{cases}
\mathcal{L}^G, & \mbox{if $H= G$,}\\
0, & \mbox{if $H$ is proper.}
\end{cases}
$$
Then $\mathcal{M} \in \mathcal{C'}$ and  $$\Hom_{\LocSys^\gen(BG,\fcat{A})}(\Sigma^{\infty}_{BG}\mathcal{M},\mathcal{L}) \not\simeq 0 $$
if $\mathcal{L}\not\simeq 0$. Therefore, $\Sigma^{\infty}_{BG}\mathcal{C'}$ detects equivalence and since $\mathcal{C}$ is presentable and stable, $\Sigma^{\infty}_{BG}\mathcal{C'}$ generates $\mathcal{C}$ under colimits.
%
%This shows that, for every $\mathcal{L} \in \mathcal{C}$, one can find $\mathcal{M} \in \mathcal{C'}$ such that 
%$$\Hom_{\LocSys^\gen(BG,\fcat{A})}(\Sigma^{\infty}_{BG}\mathcal{M},\mathcal{L}) \not\simeq 0. $$

Therefore, by \cref{cor_OmegaInftyNgen}, it suffices to show that $\theta_{\mathcal{L}}$ is an equivalence only for $\mathcal{L}\simeq \Sigma^{\infty}_{BG}\mathcal{M}$, $\mathcal{M}\in \mathcal{C}'$. However, by \cref{corollary: right adjoint to inversion many objects}, we have
$$(\Sigma^{\infty}_{BG}\mathcal{M})^G \simeq (\Omega^{\infty}_{BG}\Sigma^{\infty}_{BG}\mathcal{M})^G \simeq \colim_{V\in \Vect_{BG}} [\mathbbl{1}^V_{BG}, \mathbbl{1}^V_{BG} \otimes \mathcal{M}]^G_{\Glo}. $$
Since $\mathcal{M}\in \mathcal{C}'$ and the tensor product in $\LocSys^{\Glo}(BG,\fcat{A})$ is componentwise, we have
$$\mathbbl{1}^V_{BG} \otimes \mathcal{M} \simeq \Sigma^{\dim V^G} \mathcal{M} \in \mathcal{C}'.$$
Moreover, by the adjunctions from \cref{proposition: global sections semi-orthogonal}, we have
$$[\mathbbl{1}^V_{BG}, \mathbbl{1}^V_{BG} \otimes \mathcal{M}]^G_{\Glo} \simeq [(\mathbbl{1}^V_{BG})^G, \Sigma^{\dim V^G} \mathcal{M}^G]_{\fcat{A}(BG)} \simeq \mathcal{M}^G$$
for any $G$-representation $V$. Therefore, $\theta_{\Sigma^{\infty}_{BG}\mathcal{M}}$ is an equivalence for any globally equivariant local system $\mathcal{M} \in \mathcal{C}'$.
\end{proof}

\begin{defn}\label{definition: gfp orbispace}
Let $X \in \Type^{\Orb}$ be an orbispace and $x\colon BG \to X$ be a point, where $G \in \fcat T$. We define the functor $\Phi^{x}_{\fcat{A}}$ of \emdef{geometric fixed points at $x$} as the composite
$$\Phi^{x}_{\fcat{A}}\colon \LocSys^{\gen}_{\fcat{T}}(X,\fcat{A}) \xrightarrow{x^{\gen,*}} \LocSys^{\gen}_{\fcat{T}}(BG,\fcat{A}) \xrightarrow{\Phi^{G}_{\fcat{A}}} \fcat{A}(BG). $$
\end{defn}

\begin{rem}\label{remark: family of gf constant}
Let $X=BG$ and $\fcat{A}=\underline{\Sp}$ be the constant coefficient system. Then a faithful point $x\colon BH \to BG$ corresponds to a closed subgroup $H$ of $G$ and the functor $\Phi^{x}_{\underline{\Sp}}\colon \Sp^G \to \Sp$ coincides with the usual functor of geometric fixed points $\Phi^H$, cf. \cref{remark: gf for constant}.
\end{rem}

\begin{cor}\label{corollary: gfp are jointly conservative}
Let $X \in \Type^{\Glo}$ be an orbispace. Then the functors $\Phi^x_{\fcat{A}}$, where $x$ ranges over the set of faithful $\fcat T$-point of $X$, are jointly conservative.
\end{cor}

\begin{proof}
By \cref{construction: genuine local system}, we can assume that $X=BG$ is a representable orbispace, where $G \in \fcat T$. In this case, we will prove the assertion by induction on the size (dimension and the number of connected components) of a compact Lie group $G$. If $G$ is trivial, then the statement is trivial as well. 

Let $\mathcal{L} \in \LocSys^{\gen}_{\fcat{T}}(BG,\fcat{A})$ such that $\Phi^{H}_{\fcat{A}}(\mathcal{L})=0$ for all closed subgroups $H \subseteq G$. By induction, $\mathcal{L}^H = 0$ for all proper closed subgroups $H\subset G$. By \cref{theorem: gfp_is_fp}, we obtain $\mathcal{L}^G \simeq \Phi^G_{\fcat{A}}\mathcal{L} \simeq 0$. So, $\mathcal{L}=0$ by \cref{cor_OmegaInftyNgen}.
\end{proof}

\begin{cor}\label{corollary: gfp are jointly conservative global}
Let $X \in \Type^{\Glo}$ be a global space. Then the functors $\Phi^x_{\fcat{A}}$, where $x$ ranges over the set of $\fcat T$-point of $X$, are jointly conservative.
\end{cor}

\begin{proof}
Again, by \cref{construction: genuine local system}, we can assume that $X=BG$ is a representable global space. Therefore, the assertion follows by \cref{corollary: gfp are jointly conservative}.
\end{proof}

\begin{cor}\label{corollary: ep tilde}
Let $G \in \fcat T$ be a compact Lie group. Then the commutative square
$$
\xymatrix{
\LocSys^\Glo_{\fcat T}(BG, \fcat{A}) \ar[r]^-{(-)^G} \ar[d]^{\Sigma^{\infty}_{BG}} & \fcat{A}(BG) \ar[d]^{\simeq}\\
\LocSys^\gen_{\fcat T}(BG, \fcat{A}) \ar[r]^-{\Phi^G_{\fcat{A}}}  & \fcat{A}(BG)
}
$$
is horizontally right adjointable. In other words, the natural transformation
$$\theta\colon \Sigma^{\infty}_{BG} \circ Q^{BG} \tto \Xi^G_{\fcat{A}}$$
is an equivalence. In particular, $\Xi^G_{\fcat{A}}(\mathbbl{1})\simeq \Sigma^\infty_{BG}Q^{BG}(\mathbbl{1})$.
\end{cor}

\begin{proof}
By \cref{corollary: gfp are jointly conservative}, it suffices to show that $\Phi^H_{\fcat{A}}(\theta)$ is an equivalence for all closed subgroups $H\subset G$. By the second part of \cref{proposition: basic prop of gfp} and \cref{corollary: gfp is localization}, $\Phi^H_{\fcat{A}}(\theta)$ is an equivalence if $H=G$. If $H$ is a proper subgroup, then the right hand side of $\Phi^H_{\fcat{A}}(\theta)$ is trivial by the third part of \cref{proposition: basic prop of gfp}. Finally, if $H$ is a proper subgroup, then the left hand side $$\Phi^H_{\fcat{A}}(\Sigma^{\infty}_{BG}Q^{BG}(-)) \simeq (Q^{BG}(-))^H$$ 
of $\Phi^H_{\fcat{A}}(\theta)$ is also trivial by \cref{right_adjoint_llax_limits}.
\end{proof}

\begin{cor}\label{corollary: gfp is smashing localization}
Let $G \in \fcat{T}$ be a compact Lie group. Then the functor $\Phi^G_{\fcat{A}}$ is a smashing localization.
\end{cor}

\begin{proof}
By \cref{corollary: gfp is localization}, it suffices to show that the canonical morphism
$$\theta_{\mathcal{L},\mathcal{M}}\colon \mathcal{L} \otimes \Xi^G_{\fcat{A}}(\mathcal{M}) \tto \Xi^G_{\fcat{A}}(\Phi^G_{\fcat{A}}(\mathcal{L}) \otimes \mathcal{M}) $$
is an equivalence for all $\mathcal{L}\in \LocSys^{\gen}_{\fcat{T}}(BG,\fcat{A})$ and $\mathcal{M}\in \fcat{A}(BG)$. By~\cref{corollary: gfp are jointly conservative}, it is enough to show that $\Phi^H_{\fcat{A}}(\theta_{\mathcal{L},\mathcal{M}})$ is an equivalence for all subgroups $H \subseteq G$. By \cref{corollary: gfp is localization},  $\Phi^G_{\fcat{A}}(\theta_{\mathcal{L},\mathcal{M}})$ is an equivalence. Finally, by \cref{proposition: basic prop of gfp}, $\Phi^H_{\fcat{A}}(\theta_{\mathcal{L},\mathcal{M}})$ is an equivalence for all proper subgroups $H\subset G$ as well, since both sides of $\Phi^H_{\fcat{A}}(\theta_{\mathcal{L},\mathcal{M}})$ are trivial.
\end{proof}

In the rest of section, we will show the category of genuine local systems $\LocSys^{\gen}_{\fcat T}(BG,\fcat{A})$ satisfies the isotropy separation property generalizing the one for the genuine stable category $\Sp^G$.

As in the end of \cref{section: global section}, let $G\in \fcat{T}$ be a compact Lie group, let $\fcat{P}=\fcat{P}_G$ be a family of compact Lie groups which are isomorphic to a proper subgroup of $G$ and let $$\pi\colon (BG)_{\fcat{P}} = BG \times_{\mathcal{N}}\mathcal{N}_{\fcat{P}}\tto BG$$ be the canonical projection. 

\begin{cor}[Isotropy separation for $\LocSys^{\gen}_{\fcat T}(BG,\fcat{A})$]\label{corollary: isotropy separation for gen}
Let $G \in \fcat T$ be a compact Lie group. The following diagram
\[\xymatrix{
\LocSys^{\gen}_{\fcat T}((BG)_{\fcat{P}},\fcat{A}) \ar@<-1ex>@{^(->}[r]_-{\pi^{\gen}_{*}} & 
\ar@<-1ex>[l]_-{\pi^{\gen,*}} 
\LocSys^{\gen}_{\fcat T}(BG,\fcat{A}) \ar@<0.5ex>[r]^-{\Phi^G} & \ar@{^(->}@<0.5ex>[l]^-{\Xi^{G}} \fcat{A}(BG),
}
\]
presents the category $\LocSys^{\gen}_{\fcat T}(BG,\fcat{A})$ as a stable symmetric monoidal recollement of $\fcat{A}(BG)$ and $\LocSys^{\gen}_{\fcat T}((BG)_{\fcat{P}},\fcat{A})$.
In particular, there is a fiber sequence
$$\pi^{\gen}_{\#} \pi^{\gen,*} \mathcal{L} \tto \mathcal{L} \tto \Xi^G_{\fcat{A}}\Phi^G_{\fcat{A}}  \mathcal{L} $$
for any genuine local system $\mathcal{L} \in \LocSys^\gen_{\fcat{T}}(BG,\fcat{A})$.
\end{cor}

\begin{proof}
By \cref{proposition: basic prop of gfp}, the composite $\Phi^G \circ \pi^{\gen}_{\#}$ is trivial, and by \cref{corollary: gfp are jointly conservative}, the functors $\pi^{\gen,*}$ and $\Phi^G_{\fcat{A}}$ are jointly conservative, so the assertion follows. 
\end{proof}

\begin{ex}\label{example: genuine separation constant}
Compare with \cref{example: constant family}. Let $\fcat{A}=\underline{\Sp}$ be the constant coefficient system. Then the recollement of \cref{corollary: isotropy separation for gen} recovers the ordinary recollement
\[\xymatrix{
\Sp^G_{\fcat{P}} \ar@<-1ex>@{^(->}[r] & 
\ar@<-1ex>[l]
\Sp^G \ar@<0.5ex>[r]^-{\Phi^G} & \ar@{^(->}@<0.5ex>[l]^-{\Xi^G} \Sp,
}
\]
for the genuine stable equivariant category $\Sp^G$. Here $\Sp^G_{\fcat{P}}$ is the full subcategory of $\Sp^G$ spanned by $X\in \Sp^G$ such that $\Phi^G(X)\simeq 0$. We recall that this recollement governs the classical isotropy fiber sequence
$$E\fcat{P} \otimes X \tto X \tto \widetilde{E\fcat{P}}\otimes X, $$
where $X\in \Sp^{G} \simeq \LocSys^{\gen}(BG,\underline{\Sp})$.
\end{ex}

\begin{ex}\label{example: gen gluing functor constant}
Compare with \cref{example: right adjoint proper family} and \cref{example: cyclic group global}. Let $G=C_p$ be a cyclic group of a prime order $p$ and suppose $\fcat{A}=\underline{\fcat{C}}$ is a constant coefficient system, $\fcat{C}\in \Prs^{\LL}_{\mathrm{st}}$. We will compute the gluing functor
$$\Phi^G(\pi^{\gen}_*(-))\colon \LocSys(BC_p,\fcat{C}) \tto \fcat{C}. $$
By applying \cref{corollary: isotropy separation for gen} to $\mathcal L = \pi_*^\gen(-)$ and taking $G$-invariants, there is a fiber sequence
$$(\pi_{\#}^{\gen}(-))^G \tto (\pi_{*}^{\gen}(-))^G \tto \Phi^G(\pi^{\gen}_*(-))$$
of functors. We note that
\begin{enumerate}
\item the fiber $(\pi_{\#}^{\gen}(-))^G$ commutes with arbitrary colimits by \cref{cor_OmegaInftyNgen};
\item the middle term $(\pi_{*}^{\gen}(-))^G$ is equivalent to $(-)^{hC_p}$ by \cref{example: cyclic group global};
\item the composite $\Phi^G(\pi^{\gen}_* \circ f_{\#}(-))$ is trivial where $f\colon *\to BC_p$ is a point. Indeed, by \cref{proposition: basic prop of gfp}, we have
\begin{align*}
\Phi^G(\pi^{\gen}_*f_{\#}\mathcal{L}) &\simeq \colim_{\substack{V\in \Vect_{BC_p},\\ V^{C_p}=0}} (\mathbbl 1^V_{BC_p} \otimes^{\gen} \pi^{\gen}_* f_{\#}\mathcal{L})^{C_p} \\
&\simeq \colim_{\substack{V\in \Vect_{BC_p},\\ V^{C_p}=0}} (\pi^{\gen}_* f_{\#}\Sigma^{\dim V}\mathcal{L})^{C_p}\\ 
&\simeq \colim_{\substack{V\in \Vect_{BC_p},\\ V^{C_p}=0}} \Sigma^{\dim V}(\pi^{\gen}_* f_{\#}\mathcal{L})^{C_p}
\end{align*}
for any $\mathcal{L} \in \fcat{C}$. Since the induced map $S^{\dim V} \to S^{\dim W}$ is trivial for an inclusion $V\hookrightarrow W$, $V\neq W$, we obtain $\Phi^G(\pi^{\gen}_*f_{\#}\mathcal{L})\simeq 0$.
\end{enumerate}
By using \cite[Theorem~1.8.7]{Raksit26}, we obtain 
$$\Phi^G(\pi^{\gen}_*(-)) \simeq (-)^{tC_p} $$
as expected from the classical case, see~\cite[Example~2.23]{AMR17} or \cite{GM95}. Finally, we obtain another description of the category $\LocSys^{\gen}(BC_p, \underline{\fcat{C}})$ as a right lax limit
$$\LocSys^{\gen}(BC_p, \underline{\fcat{C}}) \simeq \rlaxlim\left(\LocSys(BC_p,\fcat{C}) \xrightarrow{(-)^{tC_p}} \fcat{C}\right).$$
\end{ex}

\begin{ex}\label{example: gen gluing functor geom type}
As in the previous example, let $G=C_p$, but assume additionally that the coefficient system $\fcat{A}$ is of geometric type. Then, by \cref{example: geom type gen BG}, we have an equivalence 
$$\LocSys^{\gen}(BC_p,\fcat{A}) \simeq \Mod_{\fcat{A}_{C_p}}(\LocSys^{\gen}(BC_p, \underline{\fcat{A}(BC_p)}))$$
for the commutative algebra $\fcat{A}_{C_p}\in \CAlg\left(\LocSys^{\Glo}(BC_p, \underline{\fcat{A}(BC_p)})\right)$ from \cref{example: locsysglo geometric type}. Since the last equivalence in \cref{example: gen gluing functor constant} is symmetric monoidal, we obtain
$$\LocSys^{\gen}(BC_p,\fcat{A}) \simeq \rlaxlim\left(\LocSys(BC_p,\fcat{A}(*)) \xrightarrow{(\phi_{\fcat{A},*}(-))^{tC_p}} \fcat{A}(BC_p)\right),$$
where $\phi\colon *\to BC_p$ is a point of $BC_p$. 
\end{ex}

\begin{ex}\label{example: gen gluing functor rational}
As before, let $G=C_p$, but assume additionally that the coefficient system $\fcat{A}\colon \Orb^{\op} \to \Prs^\LL_{\mathbb{Q}}$ is $\mathbb{Q}$-linear and is of geometric type. Then the gluing functor $(\phi_{\fcat{A},*}(-))^{tC_p} \simeq 0$ is trivial and we obtain the splitting
$$(\pi^*, \Phi^{C_p})\colon \LocSys^{\gen}(BC_p,\fcat{A}) \xrightarrow{\simeq} \LocSys(BC_p,\fcat{A}(*)) \times \fcat{A}(BC_p).$$

At the moment of writing, we are not aware if the assumption on geometric type is necessary here.
\end{ex}

\section{Complex periodic coefficients}\label{section: complex periodic coefficients}
Suppose that the coefficient system $\fcat{A}$ is associated with a preoriented abelian group object $A \in \PreAb(\SpDM_S^{\nc,\flat,\qcs})$ over $S$, see \cref{ssect_coefs_from_PreAbStk}. If both $S$ and $A$ are nc-affine and $G$ is a compact Lie group, then the category $\LocSys^{\gen}(BG,\fcat{A})$ of genuine local systems identifies with $\Mod_{\Sigma^{\infty}_{BG}\fcat{A}_G}(\Sp^G)$, see \cref{example: G-spaces geom type gen}. In other words, the category of genuine local systems recovers \emph{only} genuine equivariant cohomology theories which are suspensions of (semi-)naive equivariant cohomology theories. In particular, if $A=\mathbb{G}_{m,KU}$ (see \cref{example: multiplicative group}), the category $\LocSys^{\gen}(BG,\fcat{A})$ does \emph{not} recover the genuine equivariant $K$-theory. 

The main goal of this section is to introduce a fix to the problem above. Namely, we will define the category of \emdef{tempered local systems} $\LocSys^{\temp}(X,A)$ on a global space $X\in \Type^\Glo$. % and establish some fundamental properties of the construction in \cref{section: tempered locsys}.
Our construction is inspired by~\cite{Lur_Ell3}, but deals with compact abelian Lie groups of \emph{positive dimension} as well. As in loc.\@ cit., the category $\LocSys^{\temp}(X,A)$ is well-behaved only if $A$ is \emdef{oriented} and we will introduce oriented abelian group objects in spectral stacks over $S$ in \cref{section: oriented abelian group objects}. Informally, a preoriented abelian group object $A$ is oriented if the formal completion $\widehat{A}$ of $A$ at the unit identifies via the canonical map with the cospectrum $\cSpec_S(C_*(BU(1),\mathcal{O}_S))$. We will recall all the relevant notions in \cref{section: formal completions}. Also, in the same section, we show that the formal completion of spectral stacks behaves as expected if $A$ is \emdef{fiber-smooth} over $S$. The assumption that $A$ is oriented forces the coalgebra $C_*(BU(1),\mathcal{O}_S)$ to be \emdef{smooth}, which is automatic if $S$ is \emdef{locally complex periodic}, see \cref{section: complex oriented coefficient systems} for the definition and certain special properties.

Starting from \cref{section: as comparison}, we deal only with a coefficient system $\fcat{A}$ associated with an oriented abelian group object $A$ over $S$; we call coefficient systems arising this way \emdef{complex oriented}. In \cref{section: as comparison}, we establish the crucial property of complex oriented coefficient systems which resembles the \emdef{Atiyah--Segal completion theorem}, see \cref{lem_TempAS_structure_sheaf}. In \cref{section: tempered locsys}, we define the category $\LocSys^\temp(X,A)$ of tempered local systems as a full subcategory of $\LocSys^\Glo(X,A)$ and we show that this subcategory is closed under pullbacks and the inclusion admits a left adjoint if $X$ is an orbispace. In \cref{section: monoidal structure}, we extend this result to arbitrary global spaces by computing the left orthogonal to $\LocSys^\temp(X,A)$ in $\LocSys^\Glo(X,A)$. The left orthogonal is also a $\otimes$-ideal, so the (componentwise) monoidal structure on  $\LocSys^\Glo(X,A)$ descends to the monoidal structure on $\LocSys^\temp(X,A)$.

Finally, in \cref{section: tempered local of tori}, given a compact torus $T\cong U(1)^{\times r}$, we identify the category of tempered local systems $\LocSys^{\temp}(BT,A)$ with the category $\QCoh(A^{\times r})$ of quasi-coherent sheaves on $A^{\times r}$.

\subsection{Formal completions in spectral algebraic geometry}\label{section: formal completions}
In this section we review several constructions and results about formal completions in spectral algebraic geometry which are necessary to work with oriented abelian group objects over a general spectral base stack.
\begin{construction}\label{constr_complement}
Let $X$ be a non-connective spectral stack and let $Y \inj X$ be a monomorphism. Define a \emdef{complement $X\setminus Y$ of $Y$ in $X$} to be a subprestack of $X$ such that for a test non-connective affine spectral scheme $T$ the groupoid $(X\setminus Y)(T)$ is the union of connected components of $X(T)$ for which the fibered product $T \times_X Y$ is empty.
\end{construction}
\begin{rem}\label{rem_complements_pullback}
Let $X^\prime \to X$ be a map and let $Y^\prime := X^\prime \times_X Y$. Then it follows directly from the definition that the complement $X^\prime\setminus Y^\prime$ is equivalent to the pullback $X^\prime\times_X (X\setminus Y)$.
\end{rem}
\begin{ex}\label{ex_compl_in_aff}
Let $X$ be a spectral scheme and let $Z \inj X$ be a closed embedding. We claim that the complement $X\setminus Z$ of $Z$ in $X$ in the sense of \Cref{constr_complement} is represented by the open subscheme $U \subseteq X$ complementary to $Z$. Indeed, since the map $U \inj X$ is a monomorphism in spectral stacks, a map from a test affine spectral scheme $T \to X$ factors through $U$ if and only if the induced projection $j\colon U_T := T\times_X U \to T$ is an equivalence. The latter map is an equivalence if and only if the fiber of 
$$\mathcal O_T \tto j_* \mathcal O_{U_T}$$
vanishes. This last condition holds if and only if the complementary closed subscheme $Z_T = T\times_X Z$ is empty, since the fiber of the map above is set theoretically concentrated on $Z_T$. It follows that $U$ represents precisely the functor which we denoted $X\setminus Z$.
\end{ex}
\begin{rem}
Let $Z \inj X$ be a representable closed embedding of spectral stacks. It follows from \Cref{rem_complements_pullback} and \Cref{ex_compl_in_aff} that the inclusion map $X\setminus Z \inj X$ is a representable open embedding.
\end{rem}

%In that follows the following simple assertion about locality on $X$ of the complement construction will be quite useful.
We observe that the construction of the complement is local on $X$, which will be used a lot in the sequel.
\begin{lem}
Let $X \simeq \colim X_{\alpha}$ be a non-connective spectral stack and let $Y \to X$ be a map. Denote $Y_\alpha := X_\alpha \times_X Y$. Then
$$X\setminus Y \simeq \colim_\alpha X_\alpha \setminus Y_\alpha.$$

\begin{proof}
Let $X^\prime \to X$ be a map and let $X_\alpha^\prime := X^\prime \times_X X_\alpha$. %Then since $\SpStk$ is a topos, $X^\prime \simeq \colim_\alpha X_\alpha^\prime$.
Then since colimits in $\SpStk^\nc$ are universal, $X^\prime \simeq \colim_\alpha X_\alpha^\prime$. 
Apply this observation with $X^\prime = X\setminus Y$ and use that by \Cref{rem_complements_pullback}
\[(X\setminus Y) \times_X X_\alpha \simeq X_\alpha \setminus Y_\alpha.\qedhere\]
\end{proof}
\end{lem}

In this work we will be most interested in the following special case of the complement construction.
\begin{construction}\label{constr_compl_and_comp}
Let $Z \to X$ be a morphism of non-connective spectral stacks such that the induced map of connective stacks $Z_{\ge 0} \to X_{\ge 0}$ is a closed embedding and the natural map
$$Z \tto Z_{\ge 0} \times_{X_{\ge 0}} X$$
is an equivalence. We will call such maps \emdef{non-connective closed embeddings}. We define an \emdef{open complement $X\setminus Z$} as in \Cref{constr_complement} and \emdef{formal completion $X^\wedge_Z$} as a complement of the open embedding $(X\setminus Z) \inj X$.

Note that by \Cref{rem_complements_pullback} these are the pullback of the corresponding constructions for the inclusion of connective stacks $Z_{\ge 0} \inj X_{\ge 0}$, namely
$$X\setminus Z \simeq (X_{\ge 0}\setminus Z_{\ge 0}) \times_{X_{\ge 0}} X, \qquad X^\wedge_Z \simeq ( (X_{\ge 0})^\wedge_{Z_{\ge 0}}) \times_{X_{\ge 0}} X.$$
\end{construction}
\begin{ex}
Let $X$ be a spectrum of a connective $E_\infty$-ring $R$ and let $Z \inj X$ be a closed embedding. Unwinding the definitions, we observe that the data of a map 
$$\Spec A \tto X^\wedge_Z$$
from the test connective $E_\infty$-ring $A$ is equivalent to the data of a map $\alpha \colon R \to A$ of $E_\infty$-rings such that the induced map on $\pi_0$ lands the underlying topological space of $\Spec \pi_0(A)$ to the underlying topological space of $Z_0 \subseteq X_0$.

If we assume additionally that the ideal $I$ corresponding to the closed embedding $Z_0 \inj X_0$ is finitely generated, then the above condition holds if and only if $\pi_0(\alpha)$ maps some power of $I$ to $0$. In other words, $\pi_0(\alpha)$ must be continuous with respect to the $I$-adic topology on the source and with the discrete topology on the target. It follows that
$$X^\wedge_Z \simeq \Spf R,$$
where the right hand side denotes the formal spectrum of $R$ with $I$-adic topology as defined in \cite[Construction 8.1.1.10]{Lur_SAG}.
\end{ex}

\begin{defn}\label{def_loc_tors_compl}
Let $Z \to X$ be a map of non-connective spectral stacks as in \Cref{constr_compl_and_comp} and let $j\colon X\setminus Z \to X$ be the corresponding open embedding. Assume additionally that $j_{\ge 0}$ is quasi-compact. It follows that $j$ is quasi-compact and (automatically) separated, hence by \Cref{lem_push_is_nice_for_rep_qcqs_alg_sp_morphs} and \Cref{rem_scalloped_in_connective_case} we deduce that the pushforward functor
$$j_* \colon \QCoh(X\setminus Z) \tto \QCoh(X)$$
is continuous, satisfies base change, and the projection formula for $j_*$ holds. Moreover, since $j$ is an embedding, it follows that $j_*$ is fully faithful by the base change argument.

We define the category $\QCoh(X)^{Z\mdef\loc}$ of \emdef{$Z$-local quasi-coherent sheaves on $X$} as the essential image of the functor $j_*$. Equivalently,
$$\QCoh(X)^{Z\mdef\loc} \simeq \Mod_{j_* \mathcal O_{X\setminus Z}} \QCoh(X).$$
We define the categories $\QCoh(X)^{Z\mdef\tors}$ of \emdef{$Z$-torsion} and $\QCoh(X)^{Z\mdef\complete}$ of \emdef{$Z$-complete} sheaves as the left and right orthogonal of $\QCoh(X)^{Z\mdef\loc}$, respectively.
\end{defn}
\begin{rem}\label{rem_tors_and_comp_as_kernels}
Since the inclusion of $Z$-local sheaves admits a left adjoint $\mathcal F \mapsto \mathcal F \otimes_{\mathcal O_X} j_* \mathcal O_{X\setminus Z}$, the category of $Z$-torsion sheaves identifies with the kernel of this functor. Similarly, the category of $Z$-complete sheaves identifies with the kernel of the functor $\HHom_{\mathcal O_X}(j_*\mathcal O_{X\setminus Z}, -)$.
\end{rem}
\begin{rem}\label{rem_locality_of_loc_and_tors}
Let $f\colon Y \to X$ be a map of non-connective spectral stacks and let $Z \inj X$ be as in \Cref{constr_compl_and_comp}. Then the pullback of $Z$-local (resp.\@ pushforward of $Z_Y:= Y \times_X Z$-local) object via $f^*$ (resp.\@ $f_*$) is $Z_Y$-local (resp.\@ $Z$-local). It follows formally that $f^*$-pullback (resp.\@ $f_*$-pushforward) of $Z$-torsion (resp.\@ $Z_Y$-complete) sheaf is $Z_Y$-torsion (resp.\@ $Z$-complete). In particular, the properties of being $Z$-local and $Z$-torsion can be checked locally on $X$.
\end{rem}
\begin{prop}\label{prop_semi_ort_tors_loc_cpl}
Let $Z \to X$ be a non-connective closed embedding with the quasi-compact open complement (see \cref{def_loc_tors_compl}). Then
\begin{enumerate}
\item the full subcategory $\QCoh(X)^{Z\mdef\tors}$ is presentable and it is a $\otimes$-ideal in $\QCoh(X)$. Moreover, the inclusion $\QCoh(X)^{Z\mdef\tors} \inj \QCoh(X)$ admits a continuous $\QCoh(X)$-linear right adjoint
$$(-)_Z^\approx \colon \QCoh(X) \tto \QCoh(X)^{Z\mdef\tors};$$

\item the full subcategory $\QCoh(X)^{Z\mdef\complete}$ is presentable and the inclusion
$\QCoh(X)^{Z\mdef\complete} \inj \QCoh(X)$
admits a left adjoint \emdef{$Z$-completion functor}
$$(-)^\wedge_Z \colon \QCoh(X) \tto \QCoh(X)^{Z\mdef\complete};$$

\item the pairs of categories
\begin{equation}\label{eq_tors_loc_comp_semi_ort}
(\QCoh(X)^{Z\mdef\tors}, \QCoh(X)^{Z\mdef\loc}) \quad \text{and} \quad (\QCoh(X)^{Z\mdef\loc}, \QCoh(X)^{Z\mdef\complete})
\end{equation}
form semi-orthogonal decompositions of $\QCoh(X)$. In particular, the composite functors
\begin{equation}\label{eq_tors_comp_mutual_equiv}
\xymatrix{\QCoh(X)^{Z\mdef\tors} \ar@{^(->}@<0.5ex>[r] & \ar@<0.5ex>[l]^-{(-)_Z^\approx} \QCoh(X) \ar@<0.5ex>[r]^-{(-)_Z^\wedge} & \ar@{^(->}@<0.5ex>[l]\QCoh(X)^{Z\mdef\complete}}
\end{equation}
are mutually inverse equivalences.
\end{enumerate}

\begin{proof}
By \Cref{rem_tors_and_comp_as_kernels}, we have
$$\QCoh(X)^{Z\mdef\tors} \simeq \ker\left(-\otimes j_*\mathcal O_{X\setminus Z} \colon \QCoh(X) \tto \QCoh(X) \right).$$
Since $\Prs^\LL$ is closed under limits, it follows that $\QCoh(X)^{Z\mdef\tors}$ is presentable. Similarly, since $\Prs^\RR$ is closed under limits, $\QCoh(X)^{Z\mdef\complete}$ is presentable as well. In particular, since $\QCoh(X)^{Z\mdef\tors}$ (resp.\@ $\QCoh(X)^{Z\mdef\complete}$) is closed in $\QCoh(X)$ under colimits (resp.\@ under limits) the inclusion admits a right (resp.\@ left) adjoint by the adjoint functor theorem.

We will show that $Z$-torsion sheaves form a $\otimes$-ideal. Let $\mathcal F$ be a $Z$-torsion sheaf and let $\mathcal G$ be an arbitrary sheaf on $X$. We claim that $\mathcal F \otimes \mathcal G \in \QCoh(X)^{Z\mdef\tors}$. This is formally equivalent to the fact that the internal $\HHom_X(\mathcal G, \mathcal H)$ is $Z$-local for any $Z$-local $\mathcal H$. But this is indeed the case since $Z$-local objects are by definition $\mathcal O_{X\setminus Z}$-modules, and $\HHom_X(\mathcal G, \mathcal H)$ has a natural structure of $\mathcal O_{X\setminus Z}$-module co-induced from the one on $\mathcal H$.

Next, we will show that the right adjoint functor $(-)^\approx_Z$ is $\QCoh(X)$-linear. Indeed, by tensoring the counit map
$$\mathcal O_X^\approx:= (\mathcal O_X)^\approx_Z \tto \mathcal O_X$$
with $\mathcal F \in \QCoh(X)$, we obtain a cofiber sequence
$$\xymatrix{\mathcal O_X^\approx\otimes \mathcal F \ar[r] & \mathcal F \ar[r] & C_{\mathcal F},}$$
where the left hand side is torsion, since $Z$-torsion modules form a $\otimes$-ideal, and the right hand side is $Z$-local, since $Z$-local sheaves also form a $\otimes$-ideal. It follows that
$$\mathcal F\otimes \mathcal O_X^\approx \areq \mathcal F^\approx_Z,$$
hence the functor $(-)_Z^\approx$ is continuous and $\QCoh(X)$-linear.

The last assertion is a general result about semi-orthogonal decompositions, see e.g.\@ \cite[Proposition 7.2.1.10]{Lur_SAG}.
\end{proof}
\end{prop}
\begin{rem}
In the notation of the previous proposition it is easy to see that $\mathcal O_X^\approx$ is the fiber of the natural map $\mathcal O_X \to j_*\mathcal O_{X\setminus Z}$. More generally, for $\mathcal F \in \QCoh(X)$ there is a fiber sequence
$$\xymatrix{\mathcal F \otimes \mathcal O_X^\approx \simeq \mathcal F^\approx_Z \ar[r] & \mathcal F \ar[r] & j_*j^* \mathcal F \simeq \mathcal F \otimes j_*\mathcal O_{X\setminus Z}.}$$
\end{rem}
\begin{rem}
For $\mathcal F \in \QCoh(X)$ there is a natural equivalence
$$\mathcal F_Z^\wedge \simeq \HHom_X(\mathcal O_X^\approx, \mathcal F).$$
\end{rem}

\begin{ex}
Let $Z \to X$ be a non-connective closed embedding with the quasi-compact open complement and let $X = \Spec R$ for some $E_\infty$-ring $R$. By assumption, the map $Z_0 \to X_0$ is a closed embedding cut out by a finitely generated ideal $I = (f_1, f_2, \ldots, f_n) \subseteq \pi_0 R$. Therefore, by \Cref{ex_compl_in_aff}, the open complement $X\setminus Z$ (defined as a functor of points) is the base change of the open embedding $U_{\ge 0} \inj \Spec R_{\ge 0}$ along $\Spec R \to \Spec R_{\ge 0}$, where $U_{\ge 0}$ is the complement of $Z_{\ge 0}$ in the sense of spectral schemes.

In particular, by \Cref{rem_tors_and_comp_as_kernels} and the Mayer--Vietoris argument, a module $M \in \Mod_R$ is $Z$-torsion if and only if
$$M\otimes_R R[f_i^{-1}] \simeq 0$$
for all $i$. Since $f_i$ generate $I$, the last condition is equivalent for $M$ to be $I$-nilpotent in the sense of \cite[Definition 7.1.1.6]{Lur_SAG}. By passing to the left orthogonal (resp.\@ right orthogonal), we obtain that the category of $Z$-local (resp.\@ $Z$-complete) sheaves coincides with the category of $I$-local (resp.\@ $I$-complete) modules from loc.\@ cit. 
\end{ex}
%\begin{prop}
%Let $Z \to X$ be as in \Cref{constr_compl_and_comp} and let $f\colon X^\prime \to X$ be a morphism. Denote $Z^\prime := X^\prime \times_X Z$. Then the natural map
%$$f^*(\mathcal F_Z^\approx) \tto f^*(\mathcal F)_{Z^\prime}^\approx$$
%is an equivalence for any $\mathcal F \in \QCoh(X)$.
%
%\begin{proof}
%By applying $f^*$ to the cofiber sequence
%$$\xymatrix{\mathcal F^\approx_Z \ar[r] & \mathcal F \ar[r] & C_{\mathcal F},}$$
%and using that both torsion and local sheaves are stable under pullbacks, we see that
%$$f^*(\mathcal F^\approx_Z) \simeq f^*(\mathcal F)_{Z^\prime}^\approx.\qedhere$$ 
%\end{proof}
%\end{prop}

The following proposition shows that the properties of being $Z$-local and $Z$-torsion are local on $X$.
\begin{prop}\label{prop_loc_and_tors_on_colim}
Let $Z \to X$ be a non-connective closed embedding with the quasi-compact open complement and let $X \simeq \colim_\alpha X_\alpha$. Then, under the equivalence
$$\QCoh(X) \simeq \lim_\alpha \QCoh(X_\alpha),$$
the inclusions
$$\QCoh(X)^{Z\mdef\loc} \inj \QCoh(X) \quad\text{and}\quad \QCoh(X)^{Z\mdef\tors} \inj \QCoh(X)$$
correspond to the limits of the inclusions
$$\QCoh(X_\alpha)^{Z_\alpha\mdef\loc} \inj \QCoh(X_{\alpha}) \quad\text{and}\quad \QCoh(X_\alpha)^{Z_\alpha\mdef\tors} \inj \QCoh(X_\alpha),$$
where $Z_\alpha := X_\alpha \times_X Z$.

\begin{proof}
The assertion about $Z$-local objects is true by the base change for $j_*$ and the fact that the pullback of $X\setminus Z$ along the map $X_\alpha \to X$ is $X_\alpha\setminus Z_\alpha$. The result about $Z$-torsion sheaves follows formally from the previous assertion and the stability of torsion sheaves under pullbacks.
\end{proof}
\end{prop}

The situation with complete sheaves is a bit more complicated, since they are not preserved by pullbacks in general.
\begin{defn}
Let $Z \to X$ be a non-connective closed embedding with the quasi-compact open complement and let $f\colon X^\prime \to X$ be a morphism. Denote $Z^\prime := X^\prime \times_X Z$. Define the \emdef{completed pullback $f^{\widehat *}$} as
$$f^{\widehat *} \colon \QCoh(X)^{Z\mdef\complete} \tto \QCoh(X^\prime)^{Z^\prime\mdef\complete}, \qquad \mathcal F \mapsto f^*(\mathcal F)^\wedge_{Z^\prime}.$$
\end{defn}
\begin{cor}\label{cor_comp_shv_on_colim}
Let $Z \to X$ be a non-connective closed embedding with the quasi-compact open complement and let $X \simeq \colim_\alpha X_\alpha$. Denote $Z_\alpha := X_\alpha \times_X Z$. Then
$$\QCoh(X)^{Z\mdef\complete} \simeq \lim_\alpha \QCoh(X_\alpha)^{Z_\alpha\mdef\complete},$$
where the transition functors in the limit diagram are given by completed pullbacks.

\begin{proof}
Let $f\colon X^\prime \to X$ be a map and denote $Z^\prime := X^\prime \times_X Z$. Since the cofiber of the natural map
$$\mathcal F^\approx_Z \tto \mathcal F, \; \mathcal F \in \QCoh(X)$$
is $Z$-local, and since local sheaves are stable under pullbacks, we deduce that the natural map
$$f^*(\mathcal F^\approx_Z)^\wedge_{Z^\prime} \tto f^*(\mathcal F)^\wedge_{Z^\prime} = f^{\widehat *}(\mathcal F)$$
is an equivalence.  This observation and the equivalence between the full subcategories of complete and torsion objects from \cref{prop_semi_ort_tors_loc_cpl} show that the diagram $\QCoh(X_\alpha)^{Z_\alpha\mdef\complete}$ with completed pullbacks as the transition functors is equivalent to the diagram of $\QCoh(X_\alpha)^{Z_\alpha\mdef\tors}$ with the usual pullbacks as the transition functors. Finally, the result follows from \Cref{prop_loc_and_tors_on_colim}.
\end{proof}
\end{cor}

Our next goal is to relate the geometric notion of formal completion with the categorical one. Let $Z \to X$ be a non-connective closed embedding with the quasi-compact open complement. First, we note that the pullback
$$(X\setminus Z) \times_X X^\wedge_Z$$
is empty. In particular, all sheaves on $X^\wedge_Z$ are complete and, by \Cref{rem_locality_of_loc_and_tors}, the image of 
$$\widehat i_* \colon \QCoh(X_Z^\wedge) \tto \QCoh(X),$$
where $\widehat i\colon X_Z^\wedge \to X$ is the natural inclusion, lies in the category of $Z$-complete sheaves. In the setting of derived algebraic geometry, one shows that this functor induces an equivalence between quasi-coherent sheaves on the formal completion $X^\wedge_Z$ and complete sheaves on $X$, see e.g. \cite[Proposition~7.1.3]{GaitsRoz_dgindschemes}. However, in the spectral setting, we do not know if this functor is fully faithful or essentially surjective even for connective affine $X$, see \cite[Remark 8.3.4.5]{Lur_SAG}. One can deduce formally from \cite[Theorem 8.3.4.4]{Lur_SAG} and \Cref{cor_comp_shv_on_colim} that for a connective spectral stack $\widehat i_*$ restricts to an equivalence of connective subcategories. But unfortunately, this is not enough for the present paper.

Hence, in the sequel, we will restrict our attention to a class of smooth (in an appropriate sense) non-connective spectral stacks, where essentially the classical argument works.

\begin{construction}\label{constr_coSpec}
Let $X$ be a spectral stack and let $C$ be a cocommutative coalgebra in $\QCoh(X)$. We define the \emdef{cospectrum $\cSpec_X(C)$ of $C$} to be a spectral prestack over $X$ such that for a connective $E_\infty$-ring $R$ and a map $f\colon \Spec R \to X$ the space of lifts of $f$ to a map to $\cSpec_X(C)$ is the space of grouplike elements in $C_R := f^* C$, i.e.\@ a space of cocommutative coalgebra maps $\Hom_{\cCAlg_R}(R, C_R)$.

For a non-connective spectral stack $X$ and a cocommutative coalgebra $C$ in $\QCoh(X)$ such that the natural map
$$C_{\ge 0} \otimes_{\mathcal O_{X_{\ge 0}}} \mathcal O_X \tto C$$
is an equivalence, we define
$$\cSpec_X(C) := \cSpec_{X_{\ge 0}}(C_{\ge 0}) \times_{X_{\ge 0}} X.$$
\end{construction}
\begin{ex}
Let $X$ be a non-connective spectral stack and let $C$ be a cocommutative coalgebra in $\QCoh(X)$ as in \Cref{constr_coSpec} which is dualizable as a sheaf. Then there is a natural equivalence
$$\cSpec_X(C) \simeq \Spec_X(C^\vee).$$
\end{ex}
\begin{defn}
Let $X$ be a non-connective spectral stack. A cocommutative coalgebra $C$ in $\QCoh(X)$ is called \emdef{smooth} if for each $E_\infty$-ring $R$ and a map $f\colon \Spec R \to X$ the pullback $C_R:=f^* C$ is smooth in the sense of \cite[Definition 1.2.4]{Lur_Ell2}; i.e. $C_R$ is flat as an $R$-module and $\pi_0(C_R)$ is isomorphic to a divided power coalgebra on a finite rank projective $\pi_0(R)$-module. We denote the full subcategory of smooth coalgebras over $X$ by $\cCAlg^\sm_X$.
\end{defn}
\begin{defn}
A \emdef{formal hyperplane over a non-connective spectral stack $X$} is a prestack of the form $\cSpec_X(C)$ for some smooth coalgebra $C$. We write $\Hyp_X$ for the full subcategory of $\SpPStk_{/X}$ spanned by formal hyperplanes.
\end{defn}
\begin{rem}\label{rem_cCAlgSm_vs_FHyp}
Let $X \simeq \colim X_\alpha$. By construction, we have
$$\cCAlg_X^\sm \simeq \lim_\alpha \cCAlg_{X_\alpha}^\sm \quad \text{and} \qquad \Hyp_X \simeq \lim_\alpha \Hyp_{X_\alpha}.$$
Recall from \cite[Proposition 1.5.9]{Lur_Ell2} that the cospectrum functor induces an equivalence
$$\cCAlg_X^\sm \areq \Hyp_{X}$$
if $X$ is affine. By descent, the same is true for a general non-connective spectral stack $X$. We will write
$$\Dist_X \colon \Hyp_X \tto \cCAlg_X^\sm$$
the functor inverse to $\cSpec_X$.
\end{rem}
\begin{defn}\label{definition: fiber-smooth}
A morphism $p\colon X \to S$ of non-connective spectral stacks is called \emdef{fiber-smooth} if it is representable in non-connective Deligne--Mumford stacks, flat, and for each non-connective $E_\infty$-ring $R$ and a map $\Spec(R)\to S$, the connective cover $X_{R, \ge 0}=(\Spec(R)\times_{S} X)_{\ge 0}$ is fiber-smooth over $\Spec(R_{\ge 0})$ in the sense of~\cite[Definition~11.2.5.5]{Lur_SAG}; i.e. $X_{R, \ge 0}$ is flat over $R_{\ge 0}$, locally finitely presented, and its base change to any algebraically closed field (which is discrete by flatness) is regular in the sense of the classical algebraic geometry.
\end{defn}
\begin{rem}\label{rem_fcomp_fsmooth}
Let $X \to S$ be a fiber-smooth separated morphism of non-connective spectral stacks equipped with a section $s\colon S \to X$. Then, combining \Cref{rem_cCAlgSm_vs_FHyp} with \cite[Proposition 1.5.15]{Lur_Ell2}, we find that the formal completion of $X$ along $s$ is a formal hyperplane over~$S$.
\end{rem}

\begin{lem}\label{lem_fsmooth_dualizable}
Let $p\colon X \to S$ be a fiber-smooth non-connective spectral stack over a non-connective spectral stack $S$ equipped with a section $s\colon S \to X$. Then the sheaf $s_* \mathcal O_S \in \QCoh(X)$ is dualizable. Moreover, if $p$
is of relative dimension $1$, then the fiber of $\mathcal O_X \to s_* \mathcal O_S$ is flat and $\otimes$-invertible.

\begin{proof}
By flatness, $s_* \mathcal O_S$ is the pullback of $s_{\ge 0 *} (\mathcal O_{S_{\ge 0}})$ along the natural map $X \to X_{\ge 0}$. Hence it is enough to prove the assertion when $S$ and $X$ are connective. Moreover, the assertion is local on $S$, so we can assume that $S$ is affine.

Furthermore, we claim that it is enough to prove that the restriction of $s_{\ge 0 *} (\mathcal O_{S_{\ge 0}})$ along the embedding $i\colon X_{0} \inj X$ is dualizable. In fact, even more generally, let $\mathcal F \in \QCoh(X)_{\ge 0}$ be a connective sheaf. Then we claim that $\mathcal F$ is dualizable if and only if $i^* \mathcal F$ is dualizable. Indeed, since by assumption $X$ is represented by a spectral Deligne--Mumford stack, $X$ admits an \'etale atlas $U_\alpha$, where each $U_\alpha$ is affine. By flatness of the map $U_\alpha \to X$, we observe that the square
\[\xymatrix{
U_{\alpha, 0} \ar[r]\ar[d] & U_\alpha \ar[d] \\
X_0 \ar[r] & X
}\]
is fibered. Since the property of being dualizable is local, it follows that we can assume that $X = \Spec R$ for some connective $E_\infty$-ring $R$. In this case the result is standard: let $M$ be connective $R$-module such that $M\otimes_R \pi_0(R)$ is perfect as $\pi_0(R)$-module. We claim that $M$ is perfect. Indeed, let $n$ be the Tor amplitude of $M$. Note that $n$ is finite, since $M\otimes_R \pi_0 R$ is perfect. If $n=0$, then $M$ is flat and so $M\otimes_R \pi_0 R$ is flat and compact, hence projective of finite rank. By \cite[Corollary 7.2.2.19]{Lur_HA}, so is $M$. For $n>0$, we argue as in \cite[Proposition 7.2.4.23(4)]{Lur_HA} and we observe that the fiber of any map
$$R^{\oplus k} \tto M,$$
which induces surjection on $\pi_0$ (such map exists since $\pi_0 M \simeq \pi_0(M\otimes_R \pi_0 R)$), has Tor-amplitude at most $n-1$, hence the fiber is perfect by induction. It follows that $M$ is a finite extension of perfect modules, so $M$ is perfect as well.

By flatness and base change, the restriction of $s_{\ge 0 *} (\mathcal O_{S_{\ge 0}})$ to $X_0$ is equivalent to $s_{0*}(\mathcal O_{S_0})$. By smoothness of $X_0$ over $S_0$, the embedding $S_0 \inj X_0$ is a regular, so $s_{0*}(\mathcal O_{S_0})$ is dualizable. Moreover, if $p$ is of relative dimension $1$, then the ideal sheaf $\mathcal I_0$ of $S_0$ in $X_0$ is a line bundle. Since by the previous argument the fiber $\mathcal I$ of the map $\mathcal O_X \to s_* \mathcal O_S$ is dualizable and since the pullback functor
$$i^*\colon \QCoh(X) \tto \QCoh(X_0)$$
is symmetric monoidal and conservative when restricted to eventually connective sheaves, the sheaf $\mathcal I$ is $\otimes$-invertible. It is also flat, since by definition $\mathcal{I}$ is flat if and only if the restriction $i^*(\mathcal I)$ is flat.
\end{proof}
\end{lem}

\begin{lem}\label{lem_QCoh_vs_complSh_fsmooth}
Let $p\colon X \to S$ be a separated fiber-smooth morphism of non-connective spectral stacks equipped with a section $s\colon S \to X$. Let $\widehat i \colon \widehat X \to X$ denotes the inclusion of the formal completion of $X$ along $S$. Then the pushforward functor
$$\widehat i_*\colon \QCoh(\widehat X) \tto \QCoh(X)$$
is fully faithful with essential image being the category of $S$-complete sheaves.

\begin{proof}
By \cite[Proposition 1.5.15]{Lur_Ell2}, the formal completion $\widehat X$ is a formal hyperplane over $S$. Let $\mathcal P$ denote the vector bundle on $S_0$ dual to $\Prim \pi_0 \Dist(\widehat X)$. For every $n \in \mathbb Z_{\ge 0}$, we write $\mathcal A_n$ for the unique flat $\mathcal O_S$-algebra such that
$$\pi_0(\mathcal A_n) \simeq \Sym^{\le n}(\mathcal P).$$
Note that each $\mathcal A_n$ is dualizable as a sheaf. Moreover, by construction, there is a natural cocommutative coalgebra map $\mathcal A_n^\vee \to \Dist(\widehat X)$. By passing to cospectra, we obtain a compatible family of maps
$$\Spec_S \mathcal A_n \simeq \cSpec_S(\mathcal A_n^\vee) \tto \cSpec( \Dist_S(\widehat X)) \simeq \widehat X.$$
Note that
\begin{equation}\label{eq_formal_hyp_as_colim}
\indlim \Spec_S \mathcal A_n \areq \widehat X
\end{equation}
is an equivalence. Indeed, by universality of colimits in topoi, this claim can be checked locally on $S$, and in the affine case the claim is \cite[Proposition 1.5.8]{Lur_Ell2}.

By composing with the inclusion $\widehat X \to X$, we obtain a family of maps
$$i_n \colon \Spec_S(\mathcal A_n) \tto X,$$
which are non-connective affine since $X$ is separated over $S$. We denote
$$\mathcal O_S^{(n)} := i_{n*} \mathcal O_{\Spec_S \mathcal A_n}.$$
For example, $\mathcal O_S^{(0)} \simeq s_* \mathcal O_S$. Arguing as in \Cref{lem_fsmooth_dualizable}, we observe that the sheaves $\mathcal O_S^{(n)}$ are dualizable, and by construction  the map
$$\mathcal O_S^{(n+1)} \tto \mathcal O_S^{(n)}$$
is a square-zero extension for every $n\geq 0$. By \eqref{eq_formal_hyp_as_colim}, we have
$$\QCoh(\widehat X) \simeq \prolim \Mod_{\mathcal A_i} \QCoh(S) \simeq \Mod_{\mathcal O_S^{(n)}} \QCoh(X).$$
Under this identification, the pushforward $\widehat i_*$ corresponds to the functor
$$\{\mathcal F_n\}_{n\ge 0} \mapsto \prolim \mathcal F_n,$$
and its left adjoint $\widehat i^*$ is equivalent to
$$\mathcal F \mapsto \{\mathcal F\otimes_{\mathcal O_X} \mathcal O_S^{(n)}\}_{n\ge 0}.$$

First, we will show that $\widehat i_*$ is fully faithful. Let $\{\mathcal F_n\}_{n\ge 0}$ be a compatible family of $\mathcal O_S^{(n)}$-modules. We will show that the natural map
\begin{equation}\label{eq_counit_QCoh_hyp}
\prolim(\mathcal F_m) \otimes_{\mathcal O_X} \mathcal O_S^{(n)} \tto \mathcal F_n
\end{equation}
is an equivalence for every $n\ge 0$. Using dualizability of $\mathcal O_S^{(n)}$, we rewrite the left hand side as
$$\prolim(\mathcal F_m \otimes_{\mathcal O_X} \mathcal O_S^{(n)}) \simeq \prolim\left(\mathcal F_m \otimes_{\mathcal O_S^{(m)}} \mathcal O_S^{(m)} \otimes_{\mathcal O_X} \mathcal O_S^{(n)}\right).$$
Under this identification, \eqref{eq_counit_QCoh_hyp} identifies with the map induced on limits by the following map of pro-objects
$$\mathcal F_m \otimes_{\mathcal O_S^{(m)}} \mathcal O_S^{(m)} \otimes_{\mathcal O_X} \mathcal O_S^{(n)} \tto \mathcal F_m \otimes_{\mathcal O_S^{(m)}} \mathcal O_S^{(n)} \simeq \mathcal F_n.$$
The fiber of this map is pro-zero by nilpotency of the maps $\mathcal O_S^{(m)} \to \mathcal O_S^{(n)}$, hence \eqref{eq_counit_QCoh_hyp} is an equivalence.

Since the pullback $S \times_X \widehat X$ is empty, any sheaf in $\QCoh(\widehat X)$ is automatically complete, and hence $\widehat i_*$ lands in the category of $S$-complete sheaves. We will show that $\widehat i_*$ is essentially surjective. Since $\widehat i_*$ is fully faithful, it is enough to show that the left adjoint $\widehat i^*$ is conservative. In fact, we claim that if $\mathcal F$ is an $S$-complete sheaf in $\QCoh(X)$ such that
\begin{equation}\label{eq_QCohFCompFSmooth_1}
\mathcal F \otimes_{\mathcal O_X} s_* \mathcal O_S \simeq 0,
\end{equation}
then $\mathcal F \simeq 0$. First, we will reduce the claim to a non-connective affine base. Let $S^\prime \to S$ be a map and let
$$f\colon X^\prime := S^\prime \times_S X \to X$$
be the induced map. By definition, the fiber of the natural map
$$f^*(\mathcal F) \tto f^{\widehat *}(\mathcal F) = (f^*(\mathcal F))^\wedge_{S^\prime}$$
is $S^\prime$-local. Since $s^\prime_* \mathcal O_{S^\prime}$ is $S^\prime$-torsion, the fiber of
$$f^*(\mathcal F) \otimes_{\mathcal O_{X^\prime}} s^\prime_* \mathcal O_{S^\prime} \tto f^{ \widehat*}(\mathcal F) \otimes_{\mathcal O_{X^\prime}} s^\prime_* \mathcal O_{S^\prime}$$
is both local and torsion, hence vanishes. It follows that if \eqref{eq_QCohFCompFSmooth_1} holds, then
$$f^{\widehat *}(\mathcal F) \otimes_{\mathcal O_{X^\prime}} s^\prime_* \mathcal O_{S^\prime} \simeq f^*(\mathcal F) \otimes_{\mathcal O_{X^\prime}} s^\prime_* \mathcal O_{S^\prime} \simeq f^*(\mathcal F\otimes_{\mathcal O_X} s_*\mathcal O_S) \simeq 0,$$
i.e.\@ the property \eqref{eq_QCohFCompFSmooth_1} is stable under completed pullbacks. It follows that the assertion is local on $S$, so we can assume that $S$ is non-connective affine.

Then, $X$ is representable by a non-connective spectral algebraic space and hence admits an \'etale covering $\{j_\alpha \colon U_\alpha \to X\}$ by non-connective affine schemes $U_\alpha$. Let us denote by $i_\alpha\colon Z_\alpha \inj U_\alpha$ the fiber product $S \times_X U_\alpha$. Again, since $i_{\alpha*} \mathcal O_{Z_\alpha}$ is torsion in $\QCoh(U_\alpha)$ there is a natural equivalence
$$0 \simeq j_\alpha^*(\mathcal F \otimes_{\mathcal O_X} s_* \mathcal O_S) \simeq j_\alpha^*(\mathcal F)\otimes_{\mathcal O_{U_\alpha}} i_{\alpha*} \mathcal O_{Z_\alpha} \areq j_\alpha^{\widehat *}(\mathcal F) \otimes_{\mathcal O_{U_\alpha}} i_{\alpha*} \mathcal O_{Z_\alpha}.$$
Moreover, by refining if necessary, we can assume that the ideal of the closed embedding of the classical affine schemes $(Z_\alpha)_0 \inj (U_\alpha)_0$ is generated by a regular sequence.

Hence we reduced the claim to the following assertion. Let $R$ be an $E_\infty$-algebra and let $A \to B$ be a map of flat $R$-algebras such that the induced ring morphism $\pi_0(A) \to \pi_0(B)$ is surjective with ideal $I$ generated by a regular sequence $f_1, \ldots, f_d$. Let $M \in \Mod_A$ be an $I$-complete module such that $M\otimes_A B \simeq 0$. Then $M \simeq 0$. To see this note that by flatness
$$B \simeq A/f_1 \otimes_A \ldots \otimes_A A/f_d,$$
where $A/f_i$ denotes the cofiber of the multiplication by $f_i$. Denote $M_0 := M$ and set
$$M_i := M_{i-1} \otimes_A A/f_i$$
for  $0< i \le d$. We will prove by descending induction on $i$ that $M_i \simeq 0$. The base of induction $i=d$ is known by assumption. Now assume that $M_i \simeq 0$. 
Since $M_{i-1}$ is obtained from $M$ as a finite colimit, it is $I$-complete and in particular $f_i$-complete, which means that the limit
$$\prolim(\xymatrix{\ldots \ar[r]^-{\cdot f_i} & M_{i-1} \ar[r]^-{\cdot f_i} & M_{i-1}})$$
vanishes. But since $M_{i-1}/f_i = M_i \simeq 0$, $f_i$ acts invertibly on $M_{i-1}$, hence $M_{i-1} \simeq 0$.
\end{proof}
\end{lem}
\begin{rem}\label{rem_locality_in_smooth_case}
It follows from the proof of \Cref{lem_QCoh_vs_complSh_fsmooth} that a sheaf $\mathcal F \in \QCoh(X)$ is $S$-local if and only if $s^*\mathcal F \simeq 0$. Indeed, the only if direction follows from the base change and the fact that $S \times_X (X\setminus S) \simeq \eset$. Conversely, note that $s^*\mathcal F \simeq 0$ if and only if
$$s_*s^*\mathcal F \simeq \mathcal F \otimes_{\mathcal O_X} s_*\mathcal O_S$$
vanishes. We will show that this implies $\mathcal F^\wedge_S \simeq 0$, i.e. $\mathcal F$ is $S$-local. By the proof of \Cref{lem_QCoh_vs_complSh_fsmooth}, there is an equivalence
$$\mathcal F^\wedge_S \simeq \prolim \mathcal F \otimes_{\mathcal O_X} \mathcal O_S^{(n)}.$$
However, by construction, every $\mathcal O_S^{(n)}$ admits a finite filtration whose quotients are $s_* \mathcal O_S$-modules. This implies that each term in the limit above vanishes. So, $\mathcal F^\wedge_S \simeq 0$.
\end{rem}
\begin{cor}\label{cor_QCoh_formal_hyp_conservative}
Let $\widehat p\colon H \to S$ be a formal hyperplane over a non-connective spectral stack~$S$. Then the pushforward functor
$$\widehat p_*\colon \QCoh(H) \tto \QCoh(S)$$
factors through an equivalence
$$\QCoh(H) \areq \left(\Mod_{\widehat p_*\mathcal O_H} \QCoh(S)\right)^{S\mdef\complete}.$$
In particular, $\widehat p_*$ is conservative.

\begin{proof}
Let us denote by $\mathcal A$ the unique flat $E_\infty$-algebra in $\QCoh(S)$ such that
$$\pi_0(\mathcal A) \simeq \Sym ((\Prim \pi_0 \Dist_S(H))^\vee),$$
where $\Prim \pi_0 \Dist_S(H)$ is the vector bundle of primitive elements in $\pi_0 \Dist_S(H)$. Note that the morphism
$$X := \Spec_S \mathcal A \xymatrix{\ar[r]^-p &} S$$
is fiber-smooth and admits a section. Moreover, the formal completion of $X$ identifies naturally with $H$. It follows from \Cref{lem_QCoh_vs_complSh_fsmooth} that the pushforward $\widehat p_*$ factors through a fully faithful embedding into $S$-complete modules in $\QCoh(X)$ followed by $p_*$. Moreover, each sheaf in the image of $\widehat p_*$ is naturally a module over the completion of $\mathcal A$, which is equivalent to $\widehat p_* \mathcal O_H$ by construction.
\end{proof}
\end{cor}

\subsection{Complex oriented coefficient systems}\label{section: complex oriented coefficient systems}

We recall first some basic notions of complex orientable and complex periodic geometry here. We refer the reader to \cite[Section~4]{BDL25} for a detailed account. 

\begin{defn}[Definition~4.1.1 in \cite{Lur_Ell2}]\label{definition: complex orientable}
An $E_\infty$-ring $R$ is \emdef{complex orientable} if the unit map
$e\colon \mathbb{S} \to R$ factors as a composite
$$e\colon \mathbb{S} \simeq \Sigma^{\infty-2}\mathbb{CP}^1 \tto \Sigma^{\infty-2}\mathbb{CP}^\infty \xrightarrow{u} R. $$
In this case, we will say that $u$ is a \emdef{complex orientation} of $R$.
\end{defn}
In what follows, we will need a globalization of \cref{definition: complex orientable} to the non-affine case. To motivate it, we recall that the complex orientations of $R$ correspond one-to-one to the trivializations of the $R$-linear Thom local system $R^{\mathbb{C}(1)}$ of the universal line bundle $\mathbb{C}(1)$ on $BU(1) \simeq \mathbb{CP}^\infty$, see e.g.~\cite{ABGHR14}. Essentially, by definition, $R^{\mathbb{C}(1)}$ as an object of 
$$\LocSys(BU(1), \Mod_R) \simeq \Mod_R^{BU(1)}$$
is equivalent to the ``augmentation ideal'' $\Sigma^\infty U(1) \otimes R$ of the $R$-linear group algebra $R[U(1)] := \Sigma^\infty_+ U(1) \otimes R$. Using this interpretation, we can extend \cref{definition: complex orientable} to an arbitrary non-connective spectral stack.
\begin{defn}\label{def_complex_or_ss}
A non-connective spectral stack $S$ is called \emdef{complex orientable} if the $U(1)$-representation $\Sigma^\infty U(1) \otimes \mathcal O_S$ in $\QCoh(S)^{BU(1)}$ is trivial. If $S$ is complex orientable, a choice of trivialization produces the following cofiber sequence
$$\xymatrix{\Sigma \mathcal O_S \simeq \Sigma^\infty U(1)\otimes \mathcal O_S \ar[r] & \mathcal O_S[U(1)] \ar[r] & \mathcal O_S}$$
in $\QCoh(S)^{BU(1)}$. The connecting morphism in this sequence corresponds to a class in $u\in H^2(S, \mathcal O_S^{hU(1)})$. As before, we will say that $u$ is a \emdef{complex orientation of $S$}.
\end{defn}
\begin{rem}\label{pb_complex_or}
Let $f \colon S^\prime \to S$ be a map of non-connective spectral stacks, and assume that $S$ is complex orientable. Since $\Sigma^\infty U(1) \otimes \mathcal O_{S^\prime}$ is equivalent to the pullback of the corresponding representation from $S$, the stack $S^\prime$ is also complex orientable. Moreover, the pullback along $f$ induces a map
$$H^2(S, \mathcal O_S^{hU(1)}) \tto H^2(S^\prime, \mathcal O_{S^\prime}^{hU(1)}),$$
which maps the complex orientations of $S$ to the complex orientations of $S^\prime$.
\end{rem}

\begin{defn}[Definition~4.1.5 in \cite{Lur_Ell2}]\label{definition: weakly 2-periodic}
An $E_\infty$-ring $R$ is \emdef{weakly $2$-periodic} if $R$ satisfies two conditions
\begin{itemize}
\item $\pi_2(R)$ is a projective module of rank $1$ over the commutative ring $\pi_0(R)$;
\item for every integer $n$, the canonical map $\pi_2(R)\otimes_{\pi_0(R)}\pi_n(R)\to \pi_{n+2}(R)$ is an isomorphism.
\end{itemize}
\end{defn}

\begin{defn}[Definition~4.1.8 in \cite{Lur_Ell2}]\label{definition: complex periodic}
An $E_\infty$-ring $R$ is \emdef{complex periodic} if $R$ is complex orientable and weakly $2$-periodic. 
\end{defn}

\begin{lem}\label{lemma: descent of complex periodic}
Let $f\colon R \to R'$ be a faithfully flat map of $E_\infty$-ring spectra. Then $R$ is complex periodic (resp.\@ complex orientable) if and only if $R'$ is.

\begin{proof}
Compare with~\cite[Lemma~4.1.2.2]{BDL25}. By~\cite[Remark 4.1.10]{Lur_Ell2}, if $R$ is complex periodic, then so is $R'$. Vice versa, assume that $R'$ is weakly $2$-periodic. By flatness $\pi_2(R) \otimes_{\pi_0(R)} \pi_0(R') \simeq \pi_2(R')$. Since $R'$ is weakly $2$-periodic, $\pi_2(R')$ is invertible as a $\pi_0(R')$-module. Since the map $\pi_0(R) \to \pi_0(R')$ is faithfully flat, $\pi_2(R)$ is invertible as a $\pi_0(R)$-module, i.e. $R$ is weakly $2$-periodic.

Recall that the homotopy ring spectrum $R$ is complex orientable if and only if the Atiyah--Hirzebruch spectral sequence computing $R$-cohomology of $\mathbb{CP}^\infty$ degenerates at the second page. Since $f\colon R \to R'$ is faithfully flat, the induced homomorphism $\pi_*(f)$ on the homotopy groups is injective. Therefore, if the Atiyah--Hirzebruch spectral sequence for $R'$-cohomology of $\mathbb{CP}^\infty$ degenerates at the second page, then so does the spectral sequence for $R$-cohomology of $\mathbb{CP}^\infty$, see~\cite[Lemma~4.1.2.3]{BDL25}.
\end{proof}
\end{lem}

\begin{defn}[Definition~4.0.0.1 in~\cite{BDL25}]\label{definition: complex periodic stack}
A non-connective spectral stack $S$ is called \emdef{locally complex orientable} (resp. \emdef{locally complex periodic}) if there exists a jointly surjective family of maps $\{\Spec R_i \to S\}_{i\in I}$ such that each $R_i$ is complex orientable (resp. complex periodic).
\end{defn}
%\begin{ex}
%For any complex orientable $E_\infty$-ring spectrum $A$ the corresponding spectral scheme $\Spec A$ is obviously locally complex orientable.
%\end{ex}
\begin{ex}
Let $\Mell^{\ori}$ denote the non-connective spectral stack which classifies oriented elliptic curves, see~\cite[Definition~7.2.9]{Lur_Ell2}. Then $\Mell^{\ori}$ is locally complex periodic. Indeed, by \cite[Remark 7.2.8]{Lur_Ell2}, a map $\Spec R \to \Mell^{\ori}$ exists only if the $E_\infty$-ring $R$ is complex periodic. So any cover of $\Mell^{\ori}$ by affine spectral schemes satisfies the conditions of \cref{definition: complex periodic stack}.
\end{ex}
\begin{rem}\label{loc_comp_or_over}
It follows formally by \cref{lemma: descent of complex periodic} that if $S$ a locally complex orientable (resp. periodic) non-connective spectral stack then any nc-spectral stack $X$ over $S$ is also locally complex orientable (resp. periodic).
\end{rem}

\iffalse
\begin{notation}
For the rest of this section we fix an oriented abelian group non-connective spectral stack $A$ over a locally complex periodic base $S$ \todo{ref}.
\end{notation}
\fi

We will recall a crucial fact that the functor of homotopy orbits is conservative in a locally complex orientable case.

\begin{lem}\label{lem_coinvariants_conservativity_com_or}
Let $S$ be a locally complex orientable non-connective spectral stack and let $\mathcal F$ be a $G$-equivariant quasi-coherent sheaf on $S$, where $G$ is a connected compact abelian Lie group. Then $\mathcal F_{hG} \simeq 0$ if and only if $\mathcal F \simeq 0$.

\begin{proof}
See~\cite[Theorem~7.43]{MNN17}, we repeat the argument for reader's convenience. The assertion is local on $S$, so we can assume without loss of generality that $S = \Spec R$ for some (complex orientable) $E_\infty$-ring $R$ admitting a ring homomorphism $MU\to R$. Also, we note that $G\cong U(1)^{\times n}$. We have a cofiber sequence of $U(1)$-equivariant $R$-modules
$$\xymatrix{\Sigma^\infty U(1) \otimes R \ar[r] & R[U(1)] \ar[r] & R,}$$
where $R[U(1)] = \Sigma^\infty_+ U(1) \otimes R$ is the $R$-linear group algebra of $U(1)$. However, as a local system on $BU(1)$, the $U(1)$-representation $\Sigma^\infty U(1)$ is equivalent to the unit sphere bundle for the tautological line bundle over $BU(1)$. Then a complex orientation of $R$ induces an $U(1)$-equivariant trivialization $\Sigma^\infty U(1)\otimes R \simeq \Sigma R$. Hence, $R[U(1)]$ is an extension of $R$-modules with a trivial $U(1)$-action. By taking $n$-th tensor power of $R[U(1)]$, we deduce that $R[G]$ has a finite filtration whose associated graded pieces are shifts of $R$ equipped with a trivial $G$-action. It follows that $\mathcal{F}$ has a finite filtration with associated graded pieces being shifts of $\mathcal{F}_{hG} \simeq \mathcal{F}\otimes_{R[G]} R$, hence it vanishes if and only if $\mathcal{F}_{hG} \simeq 0$.
\end{proof}
\end{lem}

\subsection{Oriented abelian group objects}\label{section: oriented abelian group objects}
The basic input necessary to construct tempered equivariant cohomology and tempered local system is the notion of oriented abelian group object in spectral stacks.
\begin{defn}[Definition~5.17 in \cite{GM}]\label{definition: complex oriented abelian group object}
Let $p\colon A \to S$ be a preoriented (see \Cref{definition: preab nc stack}) fiber-smooth abelian group non-connective spectral stack over a non-connective spectral stack $S$\footnote{Recall that by abuse of notation we implicitly assume that $A$ is a relative flat quasi-compact separated non-connective spectral Deligne--Mumford stack over $S$, see \Cref{definition: ab nc stack}.}. By \Cref{rem_preor_induces_U_1_eq_lift} a preorientation of $A$ induces a lift of the unit section $e\colon S \to A$
to a $U(1)$-equivariant map, which by definition factors through the formal completion $\widehat A$ of $A$ at $e$. By \Cref{rem_fcomp_fsmooth}, $\widehat A$ is a formal hyperplane over $S$, hence this datum is equivalent to a datum of a $U(1)$-equivariant map of smooth cocommutative coalgebras
$$\mathcal O_{S} \tto \Dist_{S}(\widehat A).$$
Since the forgetful functor $\cCAlg \QCoh(S) \to \QCoh(S)$ preserves colimits, this map corresponds to a non-equivariant map of cocommutative coalgebras
$$\alpha \colon (\mathcal O_{S})_{hU(1)} \tto \Dist_{S}(\widehat A).$$
We say that the preorientation is an \emdef{orientation} if $\alpha$ is an equivalence. In this case we call $A$ an \emdef{oriented} abelian group object.
\end{defn}
\begin{rem}\label{rem_locality_of_orient}
All constructions involved in the definition of oriented abelian group object $A$ over $S$ are compatible with base changes in $S$. In particular, oriented abelian group objects are stable under pullbacks and the condition of preorientation to be an orientation can be checked locally (e.g.\@ after a base change to any non-connective affine spectral stack mapping to $S$).
\end{rem}

\begin{rem}\label{remark: orientation vs lurie}
Let $S = \Spec R$ be a non-connective affine base spectral stack. Then a preorientation of $A$ induces a preorientation of the pointed hyperplane $\widehat A$ in the sense of \cite[Definition 4.3.1]{Lur_Ell2}. Additionally, if we suppose that $R$ is complex periodic, then %Moreover, orientability of $A$ implies that the coalgebra $C_*(\mathbb{CP}^\infty, R)$ is smooth, and hence \todo{why?} that $R$ is complex periodic. 
$A$ is orientable in the sense of \cref{definition: complex oriented abelian group object} if and only if the pointed formal hyperplane $\widehat A$ is orientable in the sense of \cite[Definition~4.3.9]{Lur_Ell2} by virtue of Proposition~4.3.23 in loc.\@cit. %It also follows that an oriented abelian group object can exist only over a locally complex periodic base (\Cref{definition: complex periodic stack}).
\end{rem}

\begin{rem}\label{remark: smooth vs complex periodic}
Let $S = \Spec R$ be a non-connective affine base spectral stack. Then orientability of $A$ in the sense of \cref{definition: complex oriented abelian group object} implies that the coalgebra $C_*(\mathbb{CP}^\infty, R)$ is smooth. At the moment of writing, we are not aware whether or not the smoothness of  $C_*(\mathbb{CP}^\infty, R)$ implies that $R$ is complex periodic.
\end{rem}

\begin{defn}\label{definition: complex oriented coefficient systems}
%Let $p\colon A \to S$ be a preoriented abelian group object in non-connective spectral stacks over $S$. Suppose that $p$ is \todo{separated?} and flat. 
In the setting of~\cref{example: coefficient system from preoriented abelian group}, we will say that the coefficient system $$\fcat{A} \colon \Orb^{\op}_{\ab} \to \CAlg(\Prs^{\LL}_{\mathrm{st}}),\;\; T \mapsto \QCoh(A[\widehat{T}])$$ is \emdef{complex oriented} if the preoriented abelian group object $p\colon A \to S$ is an \emph{oriented} in the sense of \cref{definition: complex oriented abelian group object}.
\end{defn}

\begin{rem}\label{remark: complex oriented good properties}
By \cref{corollary: preab is geom type}, any complex oriented coefficient system is of geometric type. 
\end{rem}

\begin{ex}\label{example: additive group is complex oriented}
We will use the notation from \cref{example: additive group}. Let $R=\mathbb{Q}[\beta,\beta^{-1}]$, $|\beta|=2$, $S=\Spec(R)$ and let $\mathbb{G}_{a,R}= \Spec(R[t])$ be a (strict) affine line over $S$. By~\cref{example: additive group}, $p\colon \mathbb{G}_{a,R} \to S$ is a preoriented fiber-smooth separated abelian group object in non-connective spectral stacks over $S$. We will show that $\mathbb{G}_{a,R}$ is an \emph{oriented} abelian group object; in particular, $\mathbb{G}_{a,R}$ defines a complex oriented coefficient system.

%Since $R$ is a $\mathbb{Q}$-algebra, we can identify all $E_\infty$-algebras with their homotopy groups \todo{what does this mean?}. 
Since $R$ is complex periodic, the homotopy quotient $R_{hU(1)}$ is a smooth coalgebra and we have 
$$\cSpec(R_{hU(1)})\simeq \Spf(R^{hU(1)}),$$
where $R^{hU(1)}\simeq \mathbb{Q}[\beta^{\pm}][[u]]$ and $|u|=-2$ is the complex orientation of $R$. By \cref{remark: orientation vs lurie}, it suffices to show that the map 
$$\alpha \colon \Spf(R^{hU(1)}) \tto \widehat{\mathbb{G}}_{a,R}\simeq \Spf(\mathbb{Q}[\beta^\pm][[t]])$$
given by the preorientation is an equivalence. Since both sides are affine, it is enough to show that the map
$$\alpha^*\colon \mathbb{Q}[\beta^\pm][[t]] \tto \mathbb{Q}[\beta^{\pm}][[u]] $$
is an equivalence. We claim that $\alpha^*(t)=\beta u$, which will finish the proof. Indeed, since $\alpha^*$ is a map of Hopf algebras, $\alpha^*(t)$ is a primitive element. Therefore, $\alpha^*(t)$ is linear at $u$. Moreover, by construction, the quotient $\alpha^*(t)/u \in \pi_{2}(R)$ is the image of the coordinate $t$ under the composite~\eqref{equation: preorientation of affine line and coordinate}, i.e. $\alpha^*(t)/u = \beta$ according to our choice of the preorientation. 

%Indeed, by construction, $R[t]$ is a smooth coalgebra in $\Mod_R$ and $\Dist_S(\widehat{A})\simeq R[t]$. Moreover, since $R$ is complex periodic, the homotopy quotient $R_{hU(1)}$ is a smooth coalgebra as well and 
%$$\pi_*(R_{hU(1)})=H_*(\mathbb{CP}^\infty, \mathbb{Q}[\beta^{\pm}])\cong \mathbb{Q}[\beta^{\pm}]\langle b_i \;|\; i\geq 0\rangle, $$
%where $b_0=1$, and $b_i \in H_{2i}(\mathbb{CP}^{\infty}, \mathbb{Q}[\beta^{\pm}]) \cong H_{2i}(\mathbb{CP}^{i}, \mathbb{Q}[\beta^{\pm}])$ maps to $1\in \mathbb{Q}$ under the isomorphisms of $H_{2i}(\mathbb{CP}^{i}/\mathbb{CP}^{i-1}, \mathbb{Q}[\beta^{\pm}])\cong H_{2i}(S^{2i}, \mathbb{Q}[\beta^{\pm}])\cong \mathbb{Q}$. By \cref{definition: complex oriented abelian group object}, it suffices to show that the map
%$$\pi_*(\alpha)\colon \pi_*(R_{hU(1)}) \tto \pi_*(\Dist_S(\widehat{A})) \cong \mathbb{Q}[\beta^{\pm},t] $$
%given by the preorientation from \cref{example: additive group} is an isomorphism. However, $\pi_*(\alpha)$ is $\mathbb{Q}[\beta^{\pm}]$-linear and sends a generator $b_{i}$ to the product $(\beta t)^i$ \todo{elaborate?} for $i\geq 0$. A straightforward dimension count shows that $\pi_*(\alpha)$ is an isomorphism.
\end{ex}

\begin{ex}\label{example: greenlees is complex oriented}
Similarly to \cref{example: additive group is complex oriented}, the preoriented fiber-smooth separated abelian group object $p\colon \mathbb{G}_R \to S=\Spec(\mathbb{Q}[\beta,\beta^{-1}])$ from \cref{example: greenlees} is an oriented abelian group object. In particular, $\mathbb{G}_{R}$ defines a complex oriented coefficient system, which encodes rational equivariant elliptic cohomology.
\end{ex}

\begin{ex}\label{example: multiplicative group is complex oriented}
Consider the strict multiplicative group $\mathbb{G}_{m,KU}$  over the (periodic) complex $K$-theory $KU$, see \cref{example: multiplicative group}. Then $p\colon \mathbb{G}_{m,KU} \to S=\Spec(KU)$ is a preoriented fiber-smooth separated abelian group object in non-connective spectral stacks over $S$ with the preorientation as in \cite[Theorem~4.1]{GM} or, equivalently, as in \cite[Remark~4.3.8]{Lur_Ell2}. We will show that $\mathbb{G}_{m,KU}$ is an \emph{oriented} abelian group object, in particular, $\mathbb{G}_{m,KU}$ defines a complex oriented coefficient system, which encodes the equivariant complex $K$-theory.

Indeed, by~\cite[Proposition~4.3.21]{Lur_Ell2}, the given preorientation corresponds to the element $e\in \pi_2(\Omega^{\infty}\mathbb{G}_{m}(KU))$. Moreover, by~\cite[Section~6.5]{Lur_Ell2}, the element $e$ corresponds further to the complex line bundle $\mathcal{O}(1)$ over $S^2\simeq \mathbb{CP}^1$. Therefore, the induced Bott map $\beta_e\colon \Sigma^2 KU \to KU$ (see~\cite[Construction~4.3.7]{Lur_Ell2}) is given by the multiplication on the first Chern class $c_1(\mathcal{O}(1))\in \pi_2(KU)$. In particular, by the classical Bott periodicity, $\beta_e$ is an equivalence. Finally, by~\cite[Proposition~4.3.23]{Lur_Ell2}, this implies that $\mathbb{G}_{m,KU}$ is an oriented abelian group object.
\end{ex}

\begin{ex}\label{example: deformation is oriented}
Consider the preoriented abelian group scheme $\mathbb{G}_{R}$ over $R=\mathbb{Q}[\varepsilon,\beta^{\pm}]$ from \cref{example: deformation of Gm to Ga}. We will show that $\mathbb{G}_{R}$ is an \emph{oriented} abelian group object; in particular, $\mathbb{G}_{R}$ defines a complex oriented coefficient system, which deforms rational equivariant complex $K$-theory to rational equivariant ordinary cohomology.

Similarly to \cref{example: additive group is complex oriented}, it suffices to show that the map 
$$\alpha \colon \Spf(R^{hU(1)}) \tto \widehat{\mathbb{G}}_R\simeq \widehat{\mathbb{G}}_{a,R}\simeq \Spf(\mathbb{Q}[\varepsilon,\beta^\pm][[t]])$$
given by the preorientation is an equivalence. Since both sides are affine, it is enough to show that the map
$$\alpha^*\colon \mathbb{Q}[\varepsilon,\beta^\pm][[t]] \tto \mathbb{Q}[\varepsilon,\beta^{\pm}][[u]] $$
is an equivalence. By the construction of preorientation, we have $\alpha^*(t)=\beta u$, which implies the claim.
\end{ex}

\begin{ex}\label{example: oriented elliptic curve}
The most interesting examples of complex oriented coefficient systems are provided by \emph{spectral elliptic curves}, see~\cite[Section~5]{GM} and~\cite{Lur_Ell1}. We recall basic notions here.

%We recall the definition of a \emph{variety} over an $E_\infty$-ring $R$, see~\cite[Definition~1.1.1]{Lur_Ell1}. 
A morphism $p\colon X \to \Spec(R)$ is a \emdef{variety over an $E_\infty$-ring $R$} if $p$ is flat and the connective cover $p_{\geq 0}\colon X_{\geq 0} \to \Spec(R_{\geq 0})$ is proper, locally almost of finite presentation, geometrically reduced, and geometrically connected. We say that $X$ is \emdef{of relative dimension~$g$} if, for every field $k$ and for every homomorphism $\pi_0R \to k$, the base change $X_k:=X_{\geq 0}\times_{R_{\geq 0}} k$ has Krull dimension~$g$.

Similarly, $p\colon X \to S$ is a \emdef{variety of relative dimension~$g$ over a non-connective spectral stack~$S$} if the pullback $X \times_S \Spec(R) \to \Spec(R)$ is a variety of relative dimension $g$ over $R$ for every $E_\infty$-ring $R$. A \emdef{(strict) abelian variety} is a (strict) abelian group object in the category of varieties and a strict abelian variety of relative dimension~$1$ is called an \emdef{elliptic curve}. By~\cite[Proposition~1.4.7]{Lur_Ell1}, any strict abelian variety is fiber-smooth.

We say that a preoriented abelian variety $p\colon A \to S$ over a non-connective spectral stack $S$ is \emdef{oriented} if $p$ is an oriented abelian group object in the sense of \cref{definition: complex oriented abelian group object}. By definition, an oriented abelian variety must be an elliptic curve. Moreover, if the base $S$ is affine, then an elliptic curve $A$ over $S$ is oriented if and only if $A$ is oriented in the sense of~\cite[Definition~7.2.7]{Lur_Ell2}. Therefore, any oriented elliptic curve $A\to S$ over a locally complex periodic base $S$ induces a complex oriented coefficient system
$$\fcat{A}\colon \Orb^{\op}_{\ab} \tto \CAlg(\Prs^\LL). $$
The coefficient system $\fcat{A}$ encodes equivariant elliptic cohomology associated with the elliptic curve $A$, see e.g.~\cite[Section~6]{GM}

In particular, the universal oriented elliptic curve 
$$p\colon E^{\mathrm{univ}} \to \mathcal{M}^{\mathrm{or}}_{\mathrm{Ell}}$$
over the moduli stack of oriented elliptic curves $\mathcal{M}^{\mathrm{or}}_{\mathrm{Ell}}$ (see~\cite[Proposition~7.2.10]{Lur_Ell2}) induces the universal example of elliptic coefficient systems
$$\mathcal{TMF}\colon \Orb^{\op}_{\ab} \tto \CAlg(\Prs^\LL), \;\; T \mapsto \QCoh(E^{\mathrm{univ}}[\widehat{T}]).$$
The notation $\mathcal{TMF}$ emphasizes the fact that the associated cohomology theory is the theory of equivariant topological modular forms, see~\cite{GM}.
\end{ex}

\subsection{Atiyah--Segal comparison for tempered cohomology}\label{section: as comparison}
In this section we collect several technical results about the relation between genuine (resp. geometric) and homotopy (resp. Tate) fixed points in the context of Lurie--Gepner--Meier's tempered cohomology constructed from an oriented abelian group object. The main result of this section is the tempered version of the Atiyah--Segal completion theorem \cite[Theorem 2.1]{AtiyahSegal_Completion}. Roughly, the definition of tempered local systems axiomatizes this comparison and some of the auxiliary results of this section are used to establish good properties of the category of tempered local systems. Our strategy is close to the one in \cite[Section 4.6]{Lur_Ell3}, but since we work with non necessarily affine abelian group objects over the base which is only locally complex periodic (e.g. Lurie's moduli stack of oriented spectral elliptic curves), some proofs require different arguments.

We start with the following special case of the inclusion of the trivial group to $U(1)$, which holds essentially by definition of oriented abelian stacks.

\begin{lem}\label{proposition: orientation fixed points}
Let $p\colon A \to S$ be an oriented non-connective spectral abelian group stack over a locally complex periodic base $S$ with the unit section $i\colon S\inj A$. Then the natural map
\begin{equation}\label{equation: orientation}
\mathcal{O}_A \tto (i_*\mathcal{O}_S)^{hU(1)} \; \in \; \QCoh(A)
\end{equation}
induced by the orientation identifies the right hand side with the completion (in the sense of \Cref{prop_semi_ort_tors_loc_cpl}) of the left hand side at $S \inj A$.

\begin{proof}
By \Cref{lem_QCoh_vs_complSh_fsmooth} and \Cref{cor_QCoh_formal_hyp_conservative}, the global section functor
$$\Gamma_S(A, -) := p_* \colon \QCoh(A) \tto \QCoh(S)$$
is conservative when restricted to the full subcategory of $S$-complete sheaves. Hence, it is enough to prove that the map
$$\theta\colon \mathcal{O}_{\widehat A} \simeq (\mathcal{O}_A)^\wedge_S \tto (i_*\mathcal{O}_S)^{hU(1)}$$
becomes an equivalence after applying $\Gamma_S(A, -)$. However, by construction, the functor $\Gamma_S(A, -)$ identifies $\Gamma_S(A, (\mathcal{O}_A)^\wedge_S)$ with the global sections $\Gamma_S(\widehat A,\mathcal{O}_{\widehat A})$ of the structure sheaf of $\widehat A$ and $\Gamma_S(A,(i_*\mathcal{O}_S)^{hU(1)})$ with $U(1)$-homotopy invariants of $p_* i_* \mathcal O_S \simeq \mathcal O_S$, i.e.\@ with the global sections (relative to $S$) of the structure sheaf over $\cSpec( (\mathcal O_S)_{hU(1)})$. Therefore, the map $\Gamma_S(A,\theta)$ is an equivalence, since the map 
$$\cSpec( (\mathcal O_S)_{hU(1)}) \tto \widehat A$$
induced by the orientation is an equivalence by the definition of an oriented abelian group object.
\end{proof}
\end{lem}

To proceed further, we have to prove a sequence of technical results about the functor of homotopy fixed points and the functor of Tate invariants in $\QCoh(A[\widehat G])$.

\begin{prop}\label{Tate_of_coinduced} Let $\fcat C$ be a presentable stable category and let $H \subseteq G$ be an inclusion of compact abelian Lie groups. Then for any $X\in \fcat C^{BH}$ there is a natural equivalence
$$X^{tH} \simeq (\coInd_H^GX)^{tG},$$
where $\coInd_H^G\colon \fcat C^{BH} \to \fcat C^{BG}$ denotes the right adjoint to the forgetful functor $\fcat C^{BG} \to \fcat C^{BH}$.
\end{prop}

\begin{proof}
Let $T\colon \fcat{C}^{BH} \to \fcat{C}$ denote the functor given by $X\mapsto (\coInd_H^GX)^{tG}$. There is a natural transformation
$$\theta \colon X^{hH} \simeq (\coInd_H^GX)^{hG} \tto (\coInd_H^GX)^{tG} \simeq T(X). $$
By~\cite[Theorem~1.8.7]{Raksit26}, it suffices to show that the fiber $\fib(\theta)$ preserves colimits and $T(X)$ vanishes for $H$-induced $X$. By construction, $\fib(\theta)\simeq (\coInd_H^G(-))_{hG}$ preserves colimits since both functors $\coInd_H^G$ and $(-)_{hG}$ preserve colimits. Finally, let $X\simeq \Ind_e^H Y$ for some $Y \in \fcat{C}$. Then, since $\fcat{C}$ is stable, we have $X \simeq \Ind_e^H Y \simeq \Sigma^{\dim H}\coInd_e^H Y$ and 
\[(\coInd_H^GX)^{tG} \simeq \Sigma^{\dim H}(\coInd_e^GY)^{tG} \simeq \Sigma^{\dim H -\dim G} (\Ind_e^GY)^{tG}\simeq 0.\qedhere\]
\end{proof}

\begin{prop}\label{proposition: tate fixed points is a colimit}
Let $\fcat C$ be a presentable stable category. Then for $X \in \fcat C^{BU(1)}$ there is a natural equivalence
$$X^{tU(1)} \simeq \indlim\left(\mathbb S^{\mathbb C(1)^{\oplus n}} \otimes X\right)^{hU(1)},$$
where $\mathbb C(1)$ denotes the tautological $U(1)$-representation and the maps in the diagram are induced by the standard closed embeddings $\mathbb C(1)^{\oplus n} \inj \mathbb C(1)^{\oplus n+1}$, $(x_1, \ldots, x_n) \mapsto (x_1, \ldots, x_n, 0)$.

\begin{proof}
Let $T\colon \fcat{C}^{BU(1)}\to \fcat{C}$ denote the functor given by the formula 
$$T(X) = \indlim\left(\mathbb S^{\mathbb C(1)^{\oplus n}} \otimes X\right)^{hU(1)}$$ 
and let $F(X)$ be the fiber of the natural transformation $X^{hU(1)} \to T(X)$. By~\cite[Theorem~1.8.7]{Raksit26}, it suffices to show that:
\begin{itemize}
\item $T(X)$ vanishes for $U(1)$-induced $X$.
\item $F(-)$ preserves colimits.
\end{itemize}

For the first part, since coinduction and induction differ by a shift, it is enough to show that $T(\coInd^{U(1)}_e(Y)) \simeq 0$ for any $Y \in \fcat C$. By the projection formula, we have
$$\mathbb S^{\mathbb C(1)^{\oplus n}} \otimes \coInd^{U(1)}_e(Y) \simeq \coInd_e^{U(1)}(\Sigma^{2n} Y).$$
Hence,
$$\left(\mathbb S^{\mathbb C(1)^{\oplus n}} \otimes \coInd^{U(1)}_e(Y)\right)^{hU(1)} \simeq \Sigma^{2n} Y.$$
Note that the transition maps in the diagram for $T(\coInd_e^{U(1)}(Y))$ are induced by the embeddings $\mathbb S^{2n} \inj \mathbb S^{2(n+1)}$ and hence are nullhomotopic. This implies that $T(\coInd_e^{U(1)}(Y)) \simeq 0$.

For the last part, let $F_n\colon \fcat{C}^{BU(1)}\to \fcat{C}$ denote the functor given by the formula
$$F_n(X) = \fib\left(X^{hU(1)} \tto \left(\mathbb S^{\mathbb C(1)^{\oplus n}} \otimes X\right)^{hU(1)}\right), \;\; n\geq 0.$$
Note that $F(X) \simeq \colim F_n(X)$, so it is enough to show that $F_n$ preserves colimits for each $n\geq 0$. The functor $F_n$ can be rewritten as
$$F_n(X) \simeq \left(\Sigma^{-1}\cofib(\mathbb S^0 \to \mathbb S^{\mathbb C(1)^{\oplus n}})\otimes X\right)^{hU(1)}.$$
Finally, $\cofib(\mathbb S^0 \to \mathbb S^{\mathbb C(1)^{\oplus n}})$ has a finite $U(1)$-cell decomposition whose cells are equipped with a free action. Therefore, $F_n(X)$ admits a finite filtration whose quotients are shifts of $X$. This implies that $F_n$ preserves colimits. 
\end{proof}
\end{prop}

\begin{cor}\label{cor_TateU1_complex_or}
Let $X$ be a complex orientable non-connective spectral stack (\Cref{def_complex_or_ss}), and let $u \in H^2(X, \mathcal O_X^{hU(1)})$ be a complex orientation of $X$. Then, for every $\mathcal F \in \QCoh(X)^{BU(1)}$, there is a natural equivalence
$$\mathcal F^{tU(1)} \simeq \mathcal F^{hU(1)}[u^{-1}].$$

\begin{proof}
By \cref{proposition: tate fixed points is a colimit}, we have
$$\mathcal F^{tU(1)} \simeq \indlim\left(\mathbb S^{\mathbb C(1)^{\oplus n}} \otimes \mathcal F\right)^{hU(1)}.$$
Since the Thom local system of a direct sum of vector bundles is a tensor product of Thom local system of the summands, the map
$$\mathbb S^{\mathbb C(1)^{\oplus n}} \tto \mathbb S^{\mathbb C(1)^{\oplus n+1}},$$
induced by the embedding of representations $\mathbb C(1)^{\oplus n} \inj \mathbb C(1)^{\oplus n+1}$ as a first $n$ summands, is given by $\Id_{\mathbb S^{\mathbb C(1)^{\oplus n}}}\otimes u_0$, where
$$u_0\colon \mathbb S \simeq \mathbb S^{\mathbb C(1)^{\oplus 0}} \tto \mathbb S^{\mathbb C(1)}$$
is the corresponding map for $n=0$. Since a fiberwise one-point compactification of a vector bundle is a suspension of a fiberwise spherization, we have that 
$$\Sigma \Sigma^\infty U(1) \simeq \mathbb S^{\mathbb C(1)},$$
and under this identification the map $u_0$ corresponds to the connecting morphism in the cofiber sequence 
$$\xymatrix{\Sigma^\infty U(1) \ar[r] & \mathbb S[U(1)] \ar[r] & \mathbb S}$$
of local systems on $BU(1)$. Finally, by \Cref{def_complex_or_ss}, a choice of a complex orientation $u$ identifies $\mathbb S^{\mathbb C(1)^{\oplus n}}\otimes \mathcal O_X$ with $\Sigma^{2n} \mathcal O_X$ and $u_0\otimes \mathcal O_X$ with multiplication by $u$.
\end{proof}
\end{cor}

\begin{lem}\label{LStemp_pullback_Tate_Invs}
Let $A$ be an oriented non-connective spectral abelian group stack over $S$ with the unit section $i\colon S\inj A$ and let $g\colon A' = A\times_{S'} S \to A$ (resp. $i'\colon S' \to A')$ be the pullback of $A$ (resp. of $i$) along a morphism $f\colon S' \to S$. Then the canonical morphism
$$\theta_h \colon g^*((i_*\mathcal{O}_{S})^{hU(1)}) \tto (i'_* \mathcal{O}_{S'})^{hU(1)} \in \QCoh(A')$$
exhibits $(i'_* \mathcal{O}_{S'})^{hU(1)}$ as the $S^\prime$-completion of $g^*((i_*\mathcal{O}_S)^{hU(1)})$. In particular, the fiber of the canonical morphism
$$\theta_t \colon g^*((i_*\mathcal{O}_{S})^{tU(1)}) \tto (i'_* \mathcal{O}_{S'})^{tU(1)} \in \QCoh(A') $$
is $S^\prime$-local.

\begin{proof}
Since $A$ is oriented, its pullback $A'$ is oriented as well, see~\cref{rem_locality_of_orient}. Therefore, by~\cref{proposition: orientation fixed points}, $(i_*\mathcal{O}_{S})^{hU(1)}$ (resp. $(i'_*\mathcal{O}_{S'})^{hU(1)}$) is the $S$-completion (resp.\@ $S'$-completion) of the structure sheaf $\mathcal{O}_S$ (resp.\@ $\mathcal{O}_{S'}$). In particular, $(i'_*\mathcal{O}_{S'})^{hU(1)}$ is $S'$-complete, and it suffices to show that the induced morphism
$$\Hom_{A'}((i'_* \mathcal{O}_{S'})^{hU(1)},\mathcal{F}) \tto \Hom_{A'}(g^*((i_*\mathcal{O}_{S})^{hU(1)}), \mathcal{F})$$
is an equivalence for every $S'$-complete quasi-coherent sheaf $\mathcal{F}$. However,
\begin{align*}
\Hom_{A'}(g^*((i_*\mathcal{O}_{S})^{hU(1)}),\mathcal{F}) &\simeq \Hom_{A}((i_*\mathcal{O}_{S})^{hU(1)}, g_*\mathcal{F}) \\
&\simeq \Hom_{A}(\mathcal{O}_{A},g_*\mathcal{F}) \\
&\simeq \Hom_{A'}(\mathcal{O}_{A'},\mathcal{F}) \\
&\simeq \Hom_{A'}((i'_* \mathcal{O}_{S'})^{hU(1)},\mathcal{F}),
\end{align*}
which finishes the proof of the first part of the lemma. For the second part, we note that $\fib(\theta_h) \in \QCoh(A')^{S\mdef\llocal}$ by the first part and $\fib(\theta_t)\simeq \fib(\theta_h)$.
\end{proof}
\end{lem}

\begin{lem}\label{LStemp_TateInvs_are_local_circle}
Let $p\colon A \to S$ be an oriented non-connective spectral abelian group stack with the unit section $i \colon S \inj A$. Then, for every $\mathcal F \in \QCoh(S)^{BU(1)}$, the Tate construction $(i_*\mathcal F)^{tU(1)}$ is $S$-local.

\begin{proof}
First, assume that $S = \Spec(R)$ is non-connective affine. Let
$$u \in H^{2}(S, \mathcal O_S^{hU(1)})$$
be a complex orientation. By \Cref{pb_complex_or}, the stack $A$ is also complex orientable, and let us denote by $\widetilde u \in H^2(A, \mathcal O_A^{hU(1)})$ the pullback of $u$. Since $(i_*\mathcal F)^{tU(1)}$ is a module over $(i_*\mathcal O_S)^{tU(1)}$, it is enough to prove the assertion for $\mathcal F = \mathcal O_S$. By \Cref{cor_TateU1_complex_or}, we have
$$(i_*\mathcal O_S)^{tU(1)} \simeq (i_* \mathcal O_S)^{hU(1)}[\widetilde u^{-1}].$$
Hence, by \Cref{rem_locality_in_smooth_case}, it is enough to show that $i^*(\widetilde u)$ vanishes. To see this, we note that $\widetilde u$ acts on $(i_* \mathcal O_S)^{hU(1)}$ by multiplication with its image $\mu_*(\widetilde{u})$ under the algebra map 
$$\mu\colon \mathcal O_A^{hU(1)} \tto (i_*\mathcal O_S)^{hU(1)}.$$  
Under the isomorphism
$$H^2(A, (i_*\mathcal O_S)^{hU(1)}) \cong H^{2}(S, \mathcal O_S^{hU(1)}),$$
the element $\mu_*(\widetilde{u})$ corresponds to $u$. Next, by \Cref{cor_QCoh_formal_hyp_conservative}, we can identify $S$-complete sheaves on $A$ with the full subcategory of $p_* \mathcal O_{\widehat A}$-modules in $\QCoh(S)$. Moreover, by orientability of $A$, we have an equivalence $p_* \mathcal O_{\widehat A} \simeq \mathcal O_S^{hU(1)}$. Under these identifications, the restriction functor $i^*$ corresponds to the tensor product functor
$$-\otimes_{\mathcal O_S^{hU(1)}} \mathcal O_S.$$
Hence, we reduced to show that the image of $u$ under the restriction map
$$C^*(\mathbb{CP}^\infty, R) = \mathcal O_S^{hU(1)} \tto \mathcal O_S = C^*(*, R)$$
vanishes. The latter follows easily from the definition of a complex orientation, e.g. the restriction of $u$ to $H^2(\mathbb S^2, R) \simeq \overline H^2(\mathbb S^2, R) \oplus H^2(*, R)$ is a generator of reduced cohomology, hence its further restriction to a point vanishes.

For a general base $S$, let $f\colon S^\prime \to S$ be a map such that $S^\prime \simeq \Spec(R^\prime)$ is affine. Let $g\colon A' = A\times_{S'} S \to A$ be the pullback of $A$ along a morphism $f$. By \Cref{rem_tors_and_comp_as_kernels}, it suffices to show that $g^*((i_*\mathcal F)^{tU(1)})$ is $S^\prime$-local. However, by \cref{LStemp_pullback_Tate_Invs}, the fiber of the canonical morphism
$$g^*((i_*\mathcal F)^{tU(1)}) \tto (i'_* g^*\mathcal F))^{tU(1)} $$
is $S^\prime$-local and $(i'_*(g^*\mathcal F))^{tU(1)}$ is $S^\prime$-local by the previous paragraph. Therefore, $g^*((i_*\mathcal F)^{tU(1)})$ is $S^\prime$-local as well.
\end{proof}
\end{lem}

In the sequel, we will use the following notation.
\begin{notation}\label{nota_tors_loc_comp_BH_BG}
Let $BH \to BG$ be a faithful morphism in $\Orb_\ab$ and let $A$ be a non-connective oriented abelian group spectral stack. We denote by
$$\QCoh(A[\widehat G])^{(H,G)\mdef\tors}, \QCoh(A[\widehat G])^{(H,G)\mdef\loc}, \QCoh(A[\widehat G])^{(H,G)\mdef\complete}$$
the categories of torsion, local, and complete sheaves (in the sense of \Cref{def_loc_tors_compl}) on $A[\widehat G]$ corresponding to the non-connective closed embedding $A[\widehat H] \inj A[\widehat G]$. We call the objects of these categories $(H,G)$-torsion, $(H,G)$-local, and $(H,G)$-complete sheaves respectively. Finally, we will write 
$$(-)^\wedge_{(H,G)}\colon \QCoh(A[\widehat G]) \tto \QCoh(A[\widehat G])^{(H,G)\mdef\complete}$$
for the completion functor at $A[\widehat{H}]$.
\end{notation}

\begin{prop}\label{LStemp_TateInvs_are_local}
Let $H, G$ be a pair of compact abelian Lie groups and let $BH \to BG$ be a faithful morphism. Let $i\colon A[\widehat H] \inj A[\widehat G]$ denote the induced non-connective closed embedding. Then, for every $\mathcal F \in \QCoh(A[\widehat H])^{BG/H}$, the Tate construction $(i_*\mathcal F)^{tG/H}$ is $(H,G)$-local.

\begin{proof}
Note that %for each $\mathcal F \in \QCoh(A[\widehat{H}])^{BG/H}$
the Tate construction $(i_*\mathcal F)^{tG/H}$ is a module over $(i_*\mathcal O_{A[\widehat H]})^{tG/H}$. Hence we can assume that $\mathcal F = \mathcal O_{A[\widehat H]}$ without loss of generality. Next we will explain how to reduce the case of an arbitrary pair $(H,G)$ to the case $(\{e\},U(1))$ considered in \cref{LStemp_TateInvs_are_local_circle}. Indeed, we claim that
\begin{enumerate}
\item \label{LStemp_TateInvs_are_local_1} Let
\[\xymatrix{
H^\prime \ar@{^(->}[r] \ar@{->>}[d] & G^\prime \ar@{->>}[d] \\
H \ar@{^(->}[r] & G
}\]
be a fibered square of compact abelian Lie groups such that the horizontal arrows are embeddings and the vertical arrows are surjections. Then if the locality of the Tate invariants holds for $H \subseteq G$ then it is also holds for $H^\prime \subseteq G^\prime$.

\item \label{LStemp_TateInvs_are_local_2} Let $K \subseteq H \subseteq G$ be compact abelian Lie groups and assume that locality of the Tate invariants holds for $K \subseteq G$. Then the locality of the Tate invariants holds for $K \subseteq H$.

\item \label{LStemp_TateInvs_are_local_3} Let $K \subseteq H \subseteq G$ be compact abelian Lie groups and assume that locality of the Tate invariants holds for $K \subseteq H$ and  $H \subseteq G$. Then the locality of the Tate invariants holds for $K \subseteq G$.
\end{enumerate}

By~\ref{LStemp_TateInvs_are_local_1}, we deduce that the assertion for the pair $H \subseteq G$ follows from the assertion for $\{e\} \subseteq G/H$, hence we can assume $H$ is trivial. By embedding $G$ into some torus and using~\ref{LStemp_TateInvs_are_local_2}, we can further assume that $G$ is a torus. By considering filtration of a torus with associated quotients isomorphic to $U(1)$ and by using~\ref{LStemp_TateInvs_are_local_3} and~\ref{LStemp_TateInvs_are_local_1}, we can assume that $G = U(1)$ is a one-dimensional torus.

Now we explain how to prove the claims. For~\ref{LStemp_TateInvs_are_local_1}, we consider the fiber diagram 
\[\xymatrix{
A[\widehat{H^\prime}] \ar@{^(->}[r]^-{i^\prime} \ar[d] & A[\widehat{G^\prime}] \ar[d]^p \\
A[\widehat H] \ar@{^(->}[r]^-i & A[\widehat G].
}\]
By the base change \Cref{lem_push_is_nice_for_rep_qcqs_alg_sp_morphs}, we have 
$$p^* i_* \mathcal O_{A[\widehat H]} \simeq i^\prime_*\mathcal O_{A[\widehat{H^\prime}]}$$
and hence $(i^\prime_*(\mathcal O_{A[\widehat{H^\prime}]}))^{tG^\prime/H^\prime}$ is a module over $p^*( (i_*\mathcal O_{A[\widehat H]})^{tG/H})$ (here we use that $G^\prime/H^\prime \cong G/H$). Since local objects are stable under pullbacks (see \Cref{rem_locality_of_loc_and_tors}), it follows that $(i^\prime_*(\mathcal O_{A[\widehat{H^\prime}]}))^{tG^\prime/H^\prime}$ is concentrated on the complement
$$p^{-1}(A[\widehat G]\setminus A[\widehat H]) = A[\widehat{G^\prime}]\setminus A[\widehat{H^\prime}],$$
i.e. the module $(i^\prime_*(\mathcal O_{A[\widehat{H^\prime}]}))^{tG^\prime/H^\prime}$ is $(H^\prime, G^\prime)$-local.

For \ref{LStemp_TateInvs_are_local_2} and \ref{LStemp_TateInvs_are_local_3}, we consider the following commutative diagram
\[
\xymatrixcolsep{20mm}
\xymatrix{
A[\widehat{K}] \ar@{^(->}[r]^-{i} \ar@{^(->}[rd]^-{k} & A[\widehat{H}] \ar@{^(->}[d]^-{j} \\
 & A[\widehat G].
}\]
of non-connective closed embeddings. For~\ref{LStemp_TateInvs_are_local_2}, we assume that the assertion holds for the pair $K\subseteq G$. For $\mathcal F \in \QCoh(A[\widehat K])$, we need to prove that $(i_{*}\mathcal F)^{tH/K}$ is concentrated on $A[\widehat H]\setminus A[\widehat K]$. It is enough to prove that $j_{*}((i_{*}\mathcal F)^{tH/K})$ is supported on the complement $A[\widehat G]\setminus A[\widehat K]$. By \Cref{Tate_of_coinduced}, we have
$$(i_{*}\mathcal F)^{tH/K} \simeq (i_{*}\coInd_{H/K}^{G/K}(\mathcal F))^{tG/K}$$
and since $j_{*}$ preserves Tate constructions (as it preserves both limits and colimits), we obtain
$$j_{*}((i_{*}\mathcal F)^{tH/K}) \simeq (k_{*}\coInd_{H/K}^{G/K}(\mathcal F))^{tG/K}.$$
By assumption for the pair $K \subseteq G$, the right hand side is $(K, G)$-local. So, $(i_{*}\mathcal F)^{tH/K}$ is $(K, H)$-local.

We will prove the last claim~\ref{LStemp_TateInvs_are_local_3}. %let
%$$i \colon A[\widehat K] \inj A[\widehat H],\qquad j \colon A[\widehat H] \inj A[\widehat G]$$
%denote the natural inclusions induced by the orientation of $A$. \todo{Extended the argument} 
Since the norm map 
$$\norm_{G/K}\colon \Sigma^{\dim G/K}(j_* i_* \mathcal F)_{hG/K} \tto (j_* i_* \mathcal F)^{hG/K} $$
for the group $G/K$ can be presented as the composite 
%$$\Sigma^{\dim G/H}(\Sigma^{\dim H/K}(j_* i_* \mathcal F)_{hH/K})_{hG/H} \xrightarrow{N_{H/K}} \Sigma^{\dim G/H}((j_* i_* \mathcal F)^{hH/K})_{hG/H} \xrightarrow{N_{G/H}} ((j_* i_* \mathcal F)^{hH/K})^{hG/H}$$
\[
\xymatrixcolsep{20mm}
\xymatrix{
\Sigma^{\dim G/H}(\Sigma^{\dim H/K}(j_* i_* \mathcal F)_{hH/K})_{hG/H} \ar[r]^-{\norm_{H/K}} \ar[rd]_-{\norm_{G/K}}& \Sigma^{\dim G/H}((j_* i_* \mathcal F)^{hH/K})_{hG/H} \ar[d]^-{\norm_{G/H}} \\
& ((j_* i_* \mathcal F)^{hH/K})^{hG/H}
}\]
of the norm maps for the subgroup $H/K$ and the quotient group $G/H$, we obtain the following cofiber sequence
$$\xymatrix{\Sigma^{\dim G/H} (j_*(i_*\mathcal F)^{tH/K})_{hG/H} \ar[r] & (j_* i_* \mathcal F)^{tG/K} \ar[r] & (j_* (i_*\mathcal F)^{hH/K})^{tG/H}}.$$
By assumption for the pair $H \subseteq G$, $(j_* (i_*\mathcal F)^{hH/K})^{tG/H}$ is $(H,G)$-local and hence is $(K,G)$-local. By assumption for the pair $K \subseteq H$, the sheaf $(i_*\mathcal F)^{tH/K}$ is $(K, H)$-local. Hence, the pushforward $j_*(i_*\mathcal F)^{tH/K}$ is $(K, G)$-local. Since local objects are closed under colimits, both $\Sigma^{\dim G/H} (j_*(i_*\mathcal F)^{tH/K})_{hG/H}$ and $(j_* i_* \mathcal F)^{tG/K}$ are $(K,G)$-local.
\end{proof}
\end{prop}

\begin{thm}\label{lem_TempAS_structure_sheaf}
Let $i\colon BH \inj BG$ be a faithful morphism, where $H, G$ are compact abelian Lie groups. Then the canonical map
$$\theta_{H,G}\colon \mathcal{O}_{A[\widehat{G}]} \tto \left(i_*\mathcal{O}_{A[\widehat{H}]}\right)^{hG/H}$$
exhibits the right hand side as the $(H, G)$-completion of the left hand side.
\end{thm}

\begin{proof}
We proceed as in \cref{LStemp_TateInvs_are_local}. Note that the assertion is true for the pair $(\{e\},U(1))$ by \cref{proposition: orientation fixed points}. Next, we explain how to reduce the general case to the case $(\{e\},U(1))$. Indeed, we claim that
\begin{enumerate}
\item \label{claim:str_sheaf_temp_1} Let
\[\xymatrix{
H^\prime \ar@{^(->}[r] \ar@{->>}[d] & G^\prime \ar@{->>}[d] \\
H \ar@{^(->}[r] & G
}\]
be a fibered square of compact abelian Lie groups such that the horizontal arrows are embeddings and the vertical arrows are surjections. Suppose that $\theta_{H,G}$ is the $(H,G)$-completion. Then $\theta_{H',G'}$ is the $(H',G')$-completion.

\item \label{claim:str_sheaf_temp_2} Let $K \subseteq H \subseteq G$ be compact abelian Lie groups and assume that $\theta_{K,H}$ is the $(K, H)$-completion and $\theta_{H,G}$ is the $(H,G)$-completion. Then $\theta_{K,G}$ is the $(K,G)$-completion.

\item \label{claim:str_sheaf_temp_3} Let $K \subseteq H \subseteq G$ be compact abelian Lie groups, $G$ is connected, and assume that $\theta_{K,G}$ is the $(K,G)$-completion and $\theta_{H,G}$ is the $(H,G)$-completion. Then $\theta_{K,H}$ is the $(K,H)$-completion.
\end{enumerate}
By \ref{claim:str_sheaf_temp_1}, the assertion for a pair $H \subseteq G$ follows from the assertion for the pair $* \subseteq G/H$. Hence, we can assume that $H=\{e\}$. By embedding $G$ into a torus and by \ref{claim:str_sheaf_temp_3}, we can assume that $G$ is a torus. By considering filtration of a torus with associated quotients isomorphic to $U(1)$ and by \ref{claim:str_sheaf_temp_1} and \ref{claim:str_sheaf_temp_2}, we can further assume that $G = U(1)$.

We prove now the reduction steps. Note that the target of $\theta_{H,G}$ is always $(H,G)$-complete, so it suffices to check that the fiber $\fib(\theta_{H,G})$ is $(H,G)$-local. For \ref{claim:str_sheaf_temp_1}, we consider the fiber diagram 
\[\xymatrix{
A[\widehat{H^\prime}] \ar@{^(->}[r]^-{i^\prime} \ar[d] & A[\widehat{G^\prime}] \ar[d]^p \\
A[\widehat H] \ar@{^(->}[r]^-i & A[\widehat G].
}\]
By the base change, we obtain the following commutative diagram
\[
\xymatrixcolsep{20mm}
\xymatrix{
\mathcal{O}_{A[\widehat{G'}]}\simeq p^*(\mathcal{O}_{A[\widehat{G}]}) \ar[r]^-{p^*(\theta_{H,G})} \ar[rd]_-{\theta_{H',G'}}& p^*\left(\left(i_*\mathcal{O}_{A[\widehat{H}]}\right)^{hG/H} \right) \ar[d]^-{\zeta^h} \\
& \left(i'_*\mathcal{O}_{A[\widehat{H'}]}\right)^{hG'/H'}.
}\]
Here we use that $G/H\cong G'/H'$. By assumption, the fiber of the horizontal arrow $p^*(\theta_{H,G})$ is $(H^\prime,G^\prime)$-local. Furthermore, the vertical map $\zeta^h$ fits into the following commutative diagram
\[
%\xymatrixcolsep{20mm}
\xymatrix{
\Sigma^{\dim{G/H}}p^*\left(\left(i_*\mathcal{O}_{A[\widehat{H}]}\right)_{hG/H} \right) \ar[r] \ar[d]^-{\zeta_h} & p^*\left(\left(i_*\mathcal{O}_{A[\widehat{H}]}\right)^{hG/H} \right) \ar[d]^-{\zeta^h} \ar[r] & p^*\left(\left(i_*\mathcal{O}_{A[\widehat{H}]}\right)^{tG/H} \right) \ar[d]^-{\zeta^t}\\
\Sigma^{\dim{G'/H'}}\left(i'_*\mathcal{O}_{A[\widehat{H'}]}\right)_{hG'/H'} \ar[r] & \left(i'_*\mathcal{O}_{A[\widehat{H'}]}\right)^{hG'/H'} \ar[r] & \left(i'_*\mathcal{O}_{A[\widehat{H'}]}\right)^{tG'/H'}.
}
\]
Since $\zeta_h$ is an equivalence, we have $\fib(\zeta^h) \simeq \fib(\zeta^t)$. Finally, the fiber $\fib(\zeta^t)$ is $(H',G')$-local as the fiber of a map between two $(H',G')$-local objects, see \cref{LStemp_TateInvs_are_local}. This implies the first part.
%Furthermore, by \cref{LStemp_TateInvs_are_local}, the fiber of the vertical arrow is also $(H^\prime,G^\prime)$-local \todo{loc.cit, is about Tate invariants}. This implies the first part.

For \ref{claim:str_sheaf_temp_2} and \ref{claim:str_sheaf_temp_3}, we consider the commutative diagram
\[
\xymatrixcolsep{20mm}
\xymatrix{
A[\widehat{K}] \ar@{^(->}[r]^-{i} \ar@{^(->}[rd]^-{k} & A[\widehat{H}] \ar@{^(->}[d]^-{j} \\
 & A[\widehat G].
}\]
By the base change, we obtain the following commutative diagram
\[
\xymatrixcolsep{20mm}
\xymatrix{
\mathcal{O}_{A[\widehat{G}]} \ar[r]^-{\theta_{K,G}} \ar[d]^-{\theta_{H,G}} & \left(k_*\mathcal{O}_{A[\widehat{K}]}\right)^{h{G/K}} \ar[d]^-{\simeq} \\
\left(j_*\mathcal{O}_{A[\widehat{H}]}\right)^{h{G/H}} \ar[r]^-{(j_*\theta_{K,H})^{h{G/H}}} & \left(j_*\left(i_*\mathcal{O}_{A[\widehat{K}]}\right)^{h{H/K}}\right)^{hG/H}.
}
\]
Therefore, we obtain the fiber sequence
\begin{equation}\label{equation:fiber_seq_compl}
\fib(\theta_{H,G}) \tto \fib(\theta_{K,G}) \tto \left(j_*\fib(\theta_{K,H})\right)^{hG/H}
\end{equation}
in the category $\QCoh(A[\widehat{G}])$.

We prove \ref{claim:str_sheaf_temp_2}. By assumption, $\fib(\theta_{H,G})$ is $(H,G)$-local and $\fib(\theta_{K,H})$ is $(K,H)$-local. Since pushforwards and limits of local modules are local, $\left(j_*\fib(\theta_{K,H})\right)^{hG/H}$ is $(K,G)$-local. Therefore, by the fiber sequence~\eqref{equation:fiber_seq_compl}, $\fib(\theta_{K,G})$ is $(K,G)$-local.

Now we prove \ref{claim:str_sheaf_temp_3}. By assumption, $\fib(\theta_{H,G})$ is $(H,G)$-local and $\fib(\theta_{K,G})$ is $(K,G)$-local. Therefore, by the fiber sequence~\eqref{equation:fiber_seq_compl}, $\left(j_*\fib(\theta_{K,H})\right)^{hG/H}$ is $(K,G)$-local. We will show that $\fib(\theta_{K,H})$ is $(K,H)$-local. Indeed, let $M:=\fib(\theta_{K,H})^\wedge_{(K,H)}$ denote the $(K,H)$-completion of $\fib(\theta_{K,H})$. It suffices to show $M\simeq 0$. However, since $\left(j_*\fib(\theta_{K,H})\right)^{hG/H}$ is $(K,G)$-local, we have $(j_*M)^{hG/H}\simeq 0$. Therefore, we obtain
$$\Sigma^{\dim G/H+1}(j_*M)_{hG/H}\simeq (j_*M)^{tG/H}.$$
However, the left hand side is always an $(H,G)$-torsion and the right hand side is $(H,G)$-local by \cref{LStemp_TateInvs_are_local}. So, $(j_*M)_{hG/H}\simeq 0$ as well. Now, by \cref{lem_coinvariants_conservativity_com_or} and the connectedness assumption for $G$, we obtain $j_*M\simeq 0$. Finally, since the pushforward $j_*$ is conservative, $M\simeq 0$ as well. This finishes the proof of the last part.
%Finally, by \cref{lem_coinvariants_conservativity_com_or}, this implies that $\fib(\theta_{K,H})$ is $(K,H)$-local \todo{not clear}.
\end{proof}

\begin{cor}\label{lem_TempAS_uber}
Let $i\colon BH \inj BG$ be a faithful morphism, where $H, G$ are compact abelian Lie groups. Then, for every $\mathcal F \in \QCoh(A(BG))$, the natural map
$$\theta^{\mathcal{F}}_{H,G}\colon \mathcal F \tto \left(i_*i^*\mathcal F\right)^{hG/H}$$
exhibits the right hand side as the $A[\widehat H]$-completion of the left hand side.

\begin{proof}
We will mimic the proof of~\cite[Lemma~4.6.12]{Lur_Ell3}. The construction
$$\mathcal{F} \mapsto (i_*i^*\mathcal{F})^{hG/H} $$
determines a functor $F \colon \QCoh(A[\widehat{G}]) \to \QCoh(A[\widehat{G}])^{(H,G)\mdef\complete}$. 
Let $\mathcal{C} \subseteq \QCoh(A[\widehat{G}])$ be the full subcategory spanned by those objects $\mathcal{F}$ for which the map $\theta_{H,G}^{\mathcal{F}}\colon \mathcal{F} \to F(\mathcal{F})$ exhibits $F(\mathcal{F})$ as an $(H,G)$-completion. By \cref{lem_TempAS_structure_sheaf}, the structure sheaf $\mathcal{O}_{A[\widehat{G}]} \in \mathcal{C}$. Consequently, to show that $\mathcal{C}=\QCoh(A[\widehat{G}])$, it suffices to show that $\mathcal{C}$ is a $\otimes$-ideal. For this, it will suffice to show that the functor $F$ preserves colimits and is strictly $\QCoh(A[\widehat{G}])$-linear.

By \cref{LStemp_TateInvs_are_local}, we can factor the functor $F$ as the composition of continuous functors
$$\QCoh(A[\widehat{G}]) \xrightarrow{F'} \QCoh(A[\widehat{G}]) \xrightarrow{(-)^{\wedge}_{(H,G)}} \QCoh(A[\widehat{G}])^{(H,G)\mdef\complete} $$
where $F'$ is given as follows
$$\mathcal{F} \mapsto \Sigma^{\dim(G/H)}(i_*i^*\mathcal{F})_{h{G/H}}. $$
It is left to note that $F'$ is $\QCoh(A[\widehat{G}])$-linear by the projection formulas, and the completion functor $(-)^\wedge_{(H,G)}$ is symmetric monoidal, hence in particular $\QCoh(A[\widehat G])$-linear.
\end{proof}
\end{cor}

\subsection{Tempered local systems}\label{section: tempered locsys}
The goal of this section is to introduce a full subcategory $\LocSys^\temp(X, A) \subseteq \LocSys^\Glo_{\ab}(X, A)$ of \emdef{$A$-tempered local systems} and establish its various nice properties for a global space $X$ and an oriented non-connective spectral abelian group stack $A$ over $S$. This is a version of Lurie's theory of tempered local systems from \cite[Section 5]{Lur_Ell3}, which works for infinite compact groups, but requires as an input an honest abelian spectral group object $A$ instead of only a spectral $p$-divisible group.

For the rest of the section we fix an oriented non-connective spectral abelian group stack $A$ over a non-connective spectral stack $S$, see \Cref{definition: complex oriented abelian group object}. Recall that using $A$ as an input, in \Cref{ssect_coefs_from_PreAbStk} we constructed a coefficient system of geometric type indexed by compact abelian Lie groups. We denote the corresponding category of globally equivariant local systems by $\LocSys^\Glo_\ab(-, A)$.

\begin{construction}\label{def_LS_temp}
Let $H, G$ be a pair of compact abelian Lie groups and let $i\colon BH \to BG$ be a faithful morphism. By abuse of notation, let us denote the induced map $A[\widehat H] \to A[\widehat G]$ also by $i$. By \cref{lem_grp_embeding_induces_aff_morph}, $i$ is a non-connective affine morphism.

Let $X \in \Type^\Glo$ be a global homotopy type and let $\mathcal L \in \LocSys^\Glo_\ab(X, A)$ be a globally equivariant local system over $X$. By definition, the pushforward $i_*\mathcal L(BH)$ is acted on by the quotient group $G/H$ and the structure map in $\QCoh(A[\widehat G])$
$$\mathcal L^{BG} \tto i_*\mathcal L^{BH}$$
is $G/H$-equivariant for the trivial $G/H$-action on the source. It follows that there is a natural comparison map
\begin{equation}\label{eq_def_LS_temp_non_complete}
\mathcal L^{BG} \tto (i_*\mathcal L^{BH})^{hG/H}.
\end{equation}
Note that since the full subcategory of $\QCoh(A[\widehat G])$ spanned by $(H,G)$-complete objects (see \Cref{nota_tors_loc_comp_BH_BG}) is closed under limits the right hand side of the map above is $(H,G)$-complete, so there is an induced comparison map
\begin{equation}\label{eq_def_LS_temp}
(\mathcal L^{BG})^\wedge_{(H,G)} \tto (i_*\mathcal L^{BH})^{hG/H}.
\end{equation}

We define the category $\LocSys^\temp(X, A)$ of \emdef{$A$-tempered local systems on $X$} as a full subcategory of $\LocSys^\Glo_\ab(X, A)$ spanned by the objects $\mathcal L$ such that the comparison map \eqref{eq_def_LS_temp}
is an equivalence for every faithful morphism $i\colon BH \to BG$ over $X$.
\end{construction}

\begin{rem}\label{LSTemp_cond_transitivity}
Let $i\colon BK \to BH$ and $j\colon BH \to BG$ be a pair of faithful morphisms over $X$. Let $\mathcal L \in \LocSys^\Glo_\ab(X, A)$ be a globally equivariant local system such that the comparison map~\eqref{eq_def_LS_temp} is an equivalence for both $i$ and $j$. Then the comparison map~\eqref{eq_def_LS_temp} for the composite $j\circ i$ is an equivalence as well. Indeed, the comparison map~\eqref{eq_def_LS_temp_non_complete} for $j\circ i$ factors as the following composite
\[\xymatrix{
\mathcal L^{BG} \ar[r] & (j_*\mathcal L^{BH})^{hG/H} \ar[r] & (j_* i_* \mathcal L^{BK})^{hG/K}.
}\]
By assumption, the fiber of the left map is $(H,G)$-local and hence is $(K,G)$-local. Also the fiber of the right map is equivalent to the $(j_*\mathcal F)^{hG/H}$, where $\mathcal F$ is the fiber of the comparison map $\mathcal L^{BH} \to (i_*\mathcal L^{BK})^{hH/K}$ for $i$. By assumption $\mathcal F$ is $(K,H)$-local, hence $j_* \mathcal F$ is $(K,G)$-local. Since local object are closed under limits, we deduce that $(j_*\mathcal F)^{hG/H}$ is $(K,G)$-local. It follows that the comparison map \eqref{eq_def_LS_temp} is an equivalence for $j\circ i$.
\end{rem}

\begin{ex}
Since we assume that $A$ is oriented, the monoidal unit $\mathbbl{1}^{\Glo}_X \in \LocSys^\Glo_\ab(X, A)$ is tempered for any global space $X \in \Type^{\Glo}$, see \cref{lem_TempAS_structure_sheaf}.
\end{ex}

First we record some formal properties of \cref{def_LS_temp}.
\begin{prop}\label{pb_preserves_LStemp}
Let $f\colon X \to Y$ be a morphism of global spaces. Then the pullback functor $f^{\Glo,*}\colon \LocSys^\Glo_\ab(Y, A) \to \LocSys^\Glo_\ab(X, A)$ preserves the full subcategory spanned by tempered local systems.

\begin{proof}
Let  $BH \to BG$  be a faithful morphism over $X$ and let $\mathcal L \in \LocSys^\Glo_\ab(Y, A)$ be a tempered local system. Then, by construction, the comparison map \eqref{eq_def_LS_temp} for $f^{\Glo,*} \mathcal L$ can be identified with the comparison map for $\mathcal L$, where we consider $BH \to BG$ over $Y$ via the postcomposition with $f\colon X \to Y$. If follows that if $\mathcal L$ is tempered, then so is $f^{\Glo,*} \mathcal L$.
\end{proof}
\end{prop}

By \cref{pb_preserves_LStemp} and since the assignment $X \mapsto \LocSys^\Glo_\ab(X, A)$ is functorial, we deduce that the assignment $X \mapsto \LocSys^\temp(X, A)$ is also functorial. Similarly to $\LocSys^\Glo_\ab(-, A)$, we have the following descent property for tempered local systems, see \cref{prop_LSglo_preserves_limits}.
\begin{prop}\label{proposition: temp descent}
The functor
$$\LocSys^\temp(-, A) \colon \Type^{\Glo, \op} \tto \Prs^\LL$$
preserves small limits.

\begin{proof}
Let $X_\alpha$ be a small diagram of global spaces with colimit $X$. By \Cref{prop_LSglo_preserves_limits}, we have
$$\LocSys^\Glo_\ab(X, A) \simeq \lim \LocSys^\Glo_\ab(X_\alpha, A).$$
Hence, it is enough to prove that $\mathcal L \in \LocSys^\Glo_\ab(X, A)$ is tempered if and only if its restriction to each $X_\alpha$ is tempered. The only if direction follows by \cref{pb_preserves_LStemp}. Conversely, assume that the restriction of $\mathcal L$ to each $X_\alpha$ is tempered. Note that any faithful morphism $i\colon BH \to BG$ over $X$ factors through some $X_\alpha$, so the comparison map for $i$ is an equivalence since it depends only on the values of $\mathcal L$ on $BH$ and $BG$ and the map between them.
\end{proof}
\end{prop}

As in the case of $\LocSys^\Glo$, the construction simplifies for $X$ an orbispace, cf. \cref{LSGlo_on_orbispaces}.
\begin{prop}\label{proposition: tempered for orbi}
Let $X$ be an orbispace. Then $\mathcal L \in \LocSys^\Glo_\ab(X, A)$ is tempered if and only if the comparison map \eqref{eq_def_LS_temp} is an equivalence for $BH \to BG$ \emph{faithful} over $X$.

\begin{proof}
Let $\mathcal L \in \LocSys^\Glo_\ab(X, A)$ be an object such that the comparison map \eqref{eq_def_LS_temp} is an equivalence for $BH \to BG$ \emph{faithful} over $X$. We will show that~\eqref{eq_def_LS_temp} is an equivalence for an arbitrary faithful map $i\colon BH \to BG$ over $X$. Since $X$ is an orbispace, there is the following commutative square over $X$
\[\xymatrix{
BH \ar@{^(->}[r]^-i\ar@{->>}[d]^q & BG \ar@{->>}[d]^p \\
BH^\prime \ar@{^(->}[r]^-{i^\prime} & BG^\prime
}\]
such that the horizontal arrows are faithful, the vertical arrows are full, and $BH^\prime$, $BG^\prime$ are faithful over $X$, see \cite[Proposition 4.3.1]{Rezk_GlobalHomotopy}. Let $H^{\prime\prime} = H^\prime \times_{G^\prime} G$, then $i$ factors through~$BH^{\prime\prime}$ and we can refine the square above to the following commutative diagram
\[\xymatrix{
BH \ar@{^(->}[r]^-j \ar[rd]^-q & BH^{\prime\prime} \ar@{^(->}[r]^-{i^{\prime\prime}} \ar@{->>}[d]^r & BG \ar@{->>}[d]^p \\
& BH^\prime \ar@{^(->}[r]^-{i^\prime} & BG^\prime
}\]
such that the right square is fibered, in particular $G^\prime/H^\prime \cong G/H^{\prime\prime}$. By \Cref{LSTemp_cond_transitivity}, it suffices to show that the comparison map \eqref{eq_def_LS_temp} is an equivalence separately for $j$ and $i^{\prime\prime}$.

Consider the case $i^{\prime\prime}$ first. By the definition of $\LocSys^\Glo_\ab(X, A)$, there are natural equivalences
$$\mathcal L^{BH^{\prime\prime}} \simeq r^*\mathcal L^{BH^\prime}, \qquad \mathcal L^{BG} \simeq p^*\mathcal L^{BG^\prime}.$$
Then, by the base change, we have
$$(i^{\prime\prime}_* \mathcal L^{BH^{\prime\prime}})^{hG/H^{\prime\prime}} \simeq (i^{\prime\prime}_* r^* \mathcal L^{BH^{\prime}})^{hG/H^{\prime\prime}} \simeq (p^* i^{\prime}_* \mathcal L^{BH^{\prime}})^{hG/H^{\prime\prime}}.$$
Using the locality of the Tate invariants (see \Cref{LStemp_TateInvs_are_local}) and the fact that $p^*$ preserves colimits, the right hand side of the last identity identifies with the $(H^{\prime\prime}, G)$-completion of $p^*((i^{\prime}_* \mathcal L^{BH^{\prime}})^{hG^{\prime}/H^{\prime}})$. Moreover, since the comparison map  \eqref{eq_def_LS_temp} is an equivalence for orbits faithful over $X$, we have the equivalence
$$p^*((i^{\prime}_* \mathcal L^{BH^{\prime}})^{hG^{\prime}/H^{\prime}}) \simeq p^*\left((\mathcal L^{BG^\prime})^\wedge_{(H^\prime, G^\prime)}\right).$$
It is left to note that 
$$(\mathcal L^{BG})^\wedge_{(H^{\prime\prime}, G)} \simeq p^*(\mathcal L^{BG^\prime})^\wedge_{(H^{\prime\prime}, G)} \xymatrix{\ar[r]^-\sim &} \left(p^*\left((\mathcal L^{BG^\prime})^\wedge_{(H^\prime, G^\prime)}\right)\right)^\wedge_{(H^{\prime\prime}, G)},$$
where the first equivalence follows from the definition of $\LocSys^\Glo_\ab(X, A)$ and the second equivalence is the general fact (which can be deduced e.g. from the fact that $p^*$ sends $(H^\prime, G^\prime)$-local objects to $(H^{\prime\prime}, G)$-local objects).

Finally, we will show that the comparison map \eqref{eq_def_LS_temp} is an equivalence for $j$. We note that  there are natural equivalences
$$\mathcal L^{BH} \simeq q^* \mathcal L^{BH^\prime} \simeq j^* r^* \mathcal L^{BH^\prime}, \qquad \mathcal L^{BH^{\prime\prime}} \simeq r^* \mathcal L^{BH^\prime}$$
by the definition of $\LocSys^\Glo_\ab(X, A)$. Under these identifications, the comparison map
$$(r^* \mathcal L^{BH^\prime})^\wedge_{(H, H^{\prime\prime})} \simeq (\mathcal L^{BH^{\prime\prime}})^\wedge_{(H, H^{\prime\prime})} \tto (j_* \mathcal L^{BH})^{hH^{\prime\prime}/H} \simeq (j_* j^* r^* \mathcal L^{BH^\prime})^{hH^{\prime\prime}/H}$$
is induced by the unit of adjunction $\mathcal F \to j_* j^* \mathcal F$. So, it is an equivalence by \Cref{lem_TempAS_uber}.
\end{proof}
\end{prop}
\begin{cor}
Let $X$ be an orbispace and denote by $\Orb_{\ab/^\rep X}$ the subcategory of $\Orb_{\ab/X}$ spanned by the orbits with a faithful map to $X$. Then the restriction along the embedding
$$i\colon \Orb_{\ab/^\rep X} \xymatrix{\ar@{^(->}[r] &} \Orb_{\ab/X}$$ induces a fully faithful embedding
$$\LocSys^\temp(X, A) \xymatrix{\ar@{^(->}[r] &} \llaxlim_{BG \in (\Orb_{\ab/^\rep X})^\op} \QCoh(A(BG)),$$
whose essential image consists of the objects such that the comparison map \eqref{eq_def_LS_temp} is an equivalence.

\begin{proof}
This follows from \Cref{LSGlo_on_orbispaces} and \cref{proposition: tempered for orbi}.
\end{proof}
\end{cor}

\begin{prop}\label{LStemp_presentable}
Let $X$ be a global space. Then the category $\LocSys^\temp(X, A)$ is presentable and it is closed under colimits in $\LocSys^\Glo_\ab(X, A)$. Moreover, if $X$ is an orbispace, then the embedding functor $\LocSys^\temp(X, A) \inj \LocSys^\Glo_\ab(X, A)$ admits a left adjoint.

\begin{proof}
Let $i\colon BH \to BG$ be a faithful morphism over $X$. Consider the functor
$$F_i\colon \LocSys^\Glo_\ab(X, A) \tto \QCoh(A(BG))^{(H,G)\mdef\complete},$$
$$F_i(\mathcal L) = \fib\left(\mathcal L^{BG} \tto (i_*\mathcal L^{BH})^{hG/H}\right)^\wedge_{(H,G)}.$$
We will show that the functor $F_i$ preserves colimits by showing that both sides preserve colimits. Indeed, by \Cref{LS_glo_colimits_pointwise}, the evaluation functors
$$\ev_{BH}\colon \LocSys^\Glo_\ab(X, A) \tto \QCoh(A(BH)), \qquad \ev_{BG}\colon \LocSys^\Glo_\ab(X, A) \tto \QCoh(A(BG))$$
preserve colimits as well as the pushforward $i_*$ along the closed embedding. Next, we recall that the norm map
$$\left(\Sigma^{\dim G/H} (i_*\mathcal L^{BH})_{hG/H}\right)_{(H,G)}^\wedge \tto (i_* \mathcal L^{BH})^{hG/H}$$
is an equivalence as the Tate invariants are $(H,G)$-local, see \Cref{LStemp_TateInvs_are_local}. This implies that the functor $F_i$ is colimit preserving. It follows that the full subcategory 
$$\ker F_i \subseteq \LocSys^\Glo_\ab(X, A)$$
is presentable (as a limit of a diagram of presentable categories and colimit preserving functors) and it is closed under colimits in $\LocSys^\Glo_\ab(X, A)$. Since $\LocSys^\temp(X, A)$ is the intersection of $\ker F_i$ over the set of all faithful morphisms $i\colon BH\to BG$ over $X$, $\LocSys^\temp(X, A)$ is presentable and it is closed under colimits in $\LocSys^\Glo_\ab(X, A)$.

Finally, suppose $X$ is an orbispace, we will show that the full subcategory $\LocSys^\temp(X, A)$ of $\LocSys^\Glo_\ab(X, A)$ is closed under limits. Indeed, by \Cref{cor_limits_in_LS_Glo_pointwise_on_orbisp}, the condition~\eqref{eq_def_LS_temp} is closed under limits if $BG \to X$ is faithful. Since, by \cref{proposition: tempered for orbi}, it is sufficient to check the condition~\eqref{eq_def_LS_temp} only for faithful points of $X$, we deduce that $\LocSys^\temp(X, A)$ is closed in $\LocSys^\Glo_\ab(X, A)$ under limits and hence the inclusion functor admits a left adjoint by the adjoint functor theorem.
\end{proof}
\end{prop}

\begin{prop}\label{from_glob_sects_are_tempered}
Let $X\in \Type^{\Glo}$ be a global space. Then, for every $\mathcal F \in \QCoh(A(X))$, the globally equivariant local system $\mathbbl{1}^{\Glo}_X \otimes \mathcal F \in \LocSys^\Glo_\ab(X, A)$ is tempered, see \cref{LSglo_enhached_GS2}. %In particular, the monoidal unit of $\LocSys^\Glo_\ab(X, A)$ is tempered.
\end{prop}

\begin{proof}
By \cref{proposition: temp descent}, we can assume that $X=BG$ is a representable global space, where $G$ is a compact abelian Lie group. In this case, the assertion follows by \cref{lem_TempAS_uber}.
\end{proof}

\subsection{Monoidal structure}\label{section: monoidal structure}
In this section, given a global space $X$, we construct a symmetric monoidal structure on $\LocSys^\temp(X, A)$. For this end, we will prove that the inclusion 
$$\LocSys^\temp(X, A) \inj \LocSys^\Glo_\ab(X, A)$$
admits a left adjoint (compare with \cref{LStemp_presentable}) and describe its kernel.
\begin{defn}\label{definition: proper torsion}
Let $G$ be a compact abelian Lie group and let $A$ be an oriented non-connective spectral abelian group stack over a locally complex orientable non-connective spectral stack $S$. Define $\QCoh(A(BG))^\tors$ to be the smallest full stable subcategory of $\QCoh(A(BG))$ closed under small colimits and containing $\QCoh(A(BG))^{(H,G)\mdef\tors}$ for all proper closed subgroup $H \subsetneq G$.
\end{defn}
\begin{rem}\label{rem_torsion_stable_under_pb}
Let $f\colon BH \to BG$ be a full map. Then the pullback functor
$$f^*\colon \QCoh(A(BG)) \tto \QCoh(A(BH))$$
maps $\QCoh(A(BG))^\tors$ into $\QCoh(A(BH))^\tors$.
\end{rem}
\begin{defn}
Let $X$ be a global space. We say that $\mathcal L \in \LocSys^\Glo_\ab(X, A)$ is \emdef{null} if the evaluation $\mathcal L^x \in \QCoh(A(BG))$ lies in the full subcategory $\QCoh(A(BG))^\tors$ for each $x\colon BG \to X$, where $G$ is a compact abelian Lie group. We denote by $\LocSys^\nnull(X, A)$ the full subcategory of $\LocSys^\Glo_\ab(X, A)$ spanned by the null local systems.
\end{defn}
\begin{rem}\label{rem_functoriality_and_descent_LSnull}
Let $f\colon X \to Y$ be a morphism of global spaces. Since the condition for a local system to be null is pointwise, the pullback functor
$$f^{\Glo,*}\colon \LocSys^\Glo_\ab(Y, A) \tto \LocSys^\Glo_\ab(X, A)$$
maps $\LocSys^\nnull(Y, A)$ into $\LocSys^\nnull(X, A)$.

Moreover, let $\{X_\alpha\}$ be a small diagram of global spaces. Combining the previous observation with \Cref{proposition: temp descent}, we deduce that $*$-pullbacks induce an equivalence
$$\LocSys^\nnull(\colim_\alpha X_\alpha, A) \areq \lim_\alpha \LocSys^\nnull(X_\alpha, A).$$
\end{rem}
\begin{rem}\label{rem_nul_on_orbisp}
Let $X$ be an orbispace. By \cite[Proposition 4.3.1]{Rezk_GlobalHomotopy}, a morphism $BG \to X$ factors through a faithful map $B\overline G \to X$. It follows from \Cref{rem_torsion_stable_under_pb} that $\mathcal L \in \LocSys^\Glo_\ab(X, A)$ is null if and only if $\mathcal L^x$ is torsion for each \emph{faithful} point $x\colon BG \to X$.
\end{rem}

Let $X$ be an orbispace. By \Cref{LStemp_presentable}, the inclusion of the full subcategory $$\LocSys^\temp(X, A) \inj \LocSys^\Glo_\ab(X, A)$$ admits a left adjoint $L_X$. In particular, there is a semi-orthogonal decomposition
$$\left\langle \ker L_X, \LocSys^\temp(X, A)\right\rangle = \LocSys^\Glo_\ab(X, A).$$
We will describe the kernel of $L_X$ in \cref{thm_Ltemp_kernel}. Its proof requires two auxiliary observations.

\begin{lem}\label{lem_sharp_pf_preserves_null}
Let $f\colon X \to Y$ be a faithful morphism of global spaces. Then the $\#$-pushforward functor (see \Cref{prop_LSGlo_sharp_pushforward})
$$f^{\Glo}_{\#}\colon \LocSys^\Glo_\ab(X, A) \tto \LocSys^\Glo_\ab(Y, A)$$
preserves the full subcategories of null local systems.

\begin{proof}
First assume that $Y$ is an orbispace. Let $\mathcal L \in \LocSys^\nnull(X, A)$. For a faithful point $y \colon BG \to Y$ with abelian $G$, we have
$$f^{\Glo}_{\#}(\mathcal L)^y \simeq \colim_{x \in \Hom_{/Y}(BG, X)} \mathcal L^x$$
by the description of $f^{\Glo}_{\#}$ from \Cref{prop_LSGlo_sharp_pushforward}. Since the full subcategory $\QCoh(A(BG))^\tors$ is closed under colimits, we observe that $f^{\Glo}_{\#}(\mathcal L)^y$ is torsion. By \Cref{rem_nul_on_orbisp}, this implies that $f^{\Glo}_{\#}(\mathcal L) \in \LocSys^\nnull(Y, A)$. 

If $Y$ is an arbitrary global space, then $Y$ is a colimit of orbispaces and the lemma follows by \cref{rem_functoriality_and_descent_LSnull}. 
\end{proof}
\end{lem}

\begin{lem}\label{lem_Ltemp_test_null}
Let $G$ be a compact abelian Lie group. For a faithful morphism $i\colon BH \to BG$, consider the functor
$$F_i\colon \LocSys^\Glo_\ab(BG, A) \tto \QCoh(A(BG)), \qquad F_i(\mathcal L) := \fib\left(\mathcal L^{BG} \tto (i_*\mathcal L^{BH})^{hG/H}\right). $$
Then the functor $F_i$ admits a left adjoint $S_i$ such that
$$S_i(\mathcal{F})^{y} \simeq
\begin{cases}
y^*\mathcal{F} \in \QCoh(A(BK)) & \text{if $K\not \subseteq H$},\\
0 & \text{if $K\subseteq H$}
\end{cases}
$$
for a faithful morphism $y\colon BK \to BG$. In particular, if $\mathcal{F} \in \QCoh(A(BG))^\tors$, then $S_i(\mathcal{F}) \in \LocSys^{\nnull}(BG,A)$.

\begin{proof}
By \Cref{LStemp_presentable} and \Cref{cor_limits_in_LS_Glo_pointwise_on_orbisp}, the functor $F_i$ preserves arbitrary limits. Hence, by the adjoint functor theorem, $F_i$ admits a left adjoint $S_i$. We will compute $S_i(\mathcal{F})^{y}$ in both cases.

Suppose that $K \not \subseteq H$. Then, by \cref{right_adjoint_llax_limits}, the evaluation functor
$$(-)^y\colon \LocSys^\Glo_\ab(BG,A) \simeq \rlaxlim_{BG'\in\Orb^{\rep}_{/BG}} \QCoh(A(BG')) \tto \QCoh(A(BK)) $$
admits a right adjoint $Q^y$ such that $(Q^y(\mathcal{F}))^{BG}\simeq y_*\mathcal{F}$ and $(Q_y(\mathcal{F}))^{BH}\simeq 0$. In particular, $F_i \circ Q^y \simeq y_*$. By taking left adjoints, we obtain the first case.

Suppose that $K \subseteq H$. Then it suffices to show that $i^{\Glo,*}\circ S_i\simeq 0$, or equivalently, $F_i \circ i^{\Glo}_* \simeq 0$. However, by \cref{example: star-pullback fixed points}, we have $(i^{\Glo}_*\mathcal{L})^{BG}\simeq i_*(\mathcal{L}^{BH})$ and $(i^{\Glo}_* \mathcal{L})^{BH} \simeq (i^{\Glo,*}i^{\Glo}_* \mathcal{L})^{BH} \simeq \Coind_e^{G/H}(\mathcal{L}^{BH})$. Moreover, the canonical morphism
$$i_*\mathcal{L}^{BH} \simeq (i^{\Glo}_* \mathcal{L})^{BG} \tto (i_*(i^{\Glo}_* \mathcal{L})^{BH})^{hG/H} \simeq i_*\mathcal{L}^{BH} $$
is an equivalence for any $\mathcal{L} \in \LocSys^\Glo_\ab(BH,A)$.
\iffalse
Consider the full subcategory $\Orb_{G,\leq H} \subset \Orb_G \simeq \Orb^{\rep}_{/BG}$ spanned by the orbits $G/G'$ such that $G' \subseteq H$. Again, by \cref{right_adjoint_llax_limits}, the projection
$$\LocSys^\Glo_\ab(BG,A) \simeq \lim^{\rlax}_{BG'\in\Orb^{\rep}_{/BG}} \QCoh(A(BG')) \tto \lim^{\rlax}_{G/G'\in\Orb_{G, \leq H}} \QCoh(A(BG')) $$
admits a right adjoint $Q^H$ such that $(Q^H \mathcal{L})^{BG}\simeq i_*\mathcal{L}^{BH}$, $(Q^H \mathcal{L})^{BH} \simeq \Coind_e^{G/H}(\mathcal{L}^{BH})$, and the canonical morphism
$$i_*\mathcal{L}^{BH} \simeq (Q^H \mathcal{L})^{BG} \tto (i_*(Q^H \mathcal{L})^{BH})^{hG/H} \simeq i_*\mathcal{L}^{BH} $$
is an equivalence for any $\mathcal{L} \in \lim^{\rlax}_{G/G'\in\Orb_{G, \leq H}} \QCoh(A(BG'))$. In particular, $F_i \circ Q^H \simeq 0$. Finally, by taking left adjoints, we obtain the second case.
\fi
\end{proof}
\end{lem}

\begin{thm}\label{thm_Ltemp_kernel}
Let $X$ be an orbispace. Then the kernel of the localization functor
$$L_X\colon \LocSys^\Glo_\ab(X, A) \tto \LocSys^\temp(X, A)$$
is $\LocSys^\nnull(X, A)$.

\begin{proof}
First, we claim that the intersection of $\LocSys^\temp(X, A)$ and $\LocSys^\nnull(X, A)$ consists only of the trivial local system $0$. To see this, let $\mathcal L \in \LocSys^\Glo_\ab(X, A)$ be both tempered and null. Assume that $\mathcal L \not\simeq 0$ and let $G$ be the smallest (with respect to the dimension and the number of connected components) compact abelian Lie group such that there exists a map $x\colon BG \to X$ and $\mathcal L^{BG} \not\simeq 0$. Note that $G$ must be non-trivial since $\QCoh(A(*))^\tors \simeq \{0\}$. Since $\mathcal L$ is tempered, we know that for every faithful $i\colon BH \inj BG$, $H\not\cong G$, the fiber of the comparison map
$$\mathcal L^{BG} \tto (i_*\mathcal L^{BH})^{hG/H}$$
is $(H,G)$-local. However, by assumption, the right hand side above vanishes. Hence $\mathcal L^{BG}$ is $(H,G)$-local for all proper subgroups $H \subsetneq G$. Since the intersection of $\QCoh(A(BG))^{(H,G)\mdef\mathrm{loc}}$ over all proper subgroups $H \subsetneq G$ is precisely the right orthogonal of $\QCoh(A(BG))^\tors$, we have $\mathcal L^{BG} \simeq 0$. So, this is a contradiction.

It is left to show that $\ker L_X \subseteq \LocSys^\nnull(X, A)$. Let $x\colon BG \to X$ be a faithful point of $X$ with abelian $G$, and let $i\colon BH \to BG$ be a faithful morphism over $X$. Then the functor
$$F_i\colon \LocSys^\Glo_\ab(X, A) \tto \QCoh(BG),$$ 
$$F_i(\mathcal L) = \fib\left(\mathcal L^{BG} \tto (i_*\mathcal L^{BH})^{hG/H}\right)$$
is limit preserving by \Cref{LStemp_presentable} and \Cref{cor_limits_in_LS_Glo_pointwise_on_orbisp}. By the adjoint functor theorem, $F_i$ admits a left adjoint $S_i$, see \cref{lem_Ltemp_test_null}. We note that $\mathcal L \in \LocSys^\Glo_\ab(X, A)$ is tempered if and only if  $F_i(\mathcal L)$ is $(H,G)$-local for all $i$ as above, which is equivalent to
$$\Hom_{\LocSys^\Glo_\ab(X, A)}(S_i(\mathcal F), \mathcal L) \simeq 0$$
for all $\mathcal F \in \QCoh(A(BG))^{(H,G)\mdef\tors}$. Therefore, the kernel $\ker L_X$ is generated under colimits by the objects $S_i(\mathcal F)$, where $\mathcal F \in \QCoh(A(BG))^{(H,G)\mdef\tors}$. Since the category $\LocSys^\nnull(X, A)$ is closed under colimits, it is enough to show that $S_i(\mathcal F) \in \LocSys^\nnull(X, A)$. Note that the functor $F_i$ factors through the restriction
$$x^{\Glo,*}\colon \LocSys^\Glo_\ab(X, A) \tto \LocSys^\Glo_\ab(BG, A),$$
hence the left adjoint $S_i$ factors through the $\#$-pushforward
$$x^{\Glo}_{\#} \colon \LocSys^\Glo_\ab(BG, A) \tto \LocSys^\Glo_\ab(X, A).$$
By \Cref{lem_sharp_pf_preserves_null}, $\#$-pushforward preserves null local systems, hence we can assume without loss of generality that $X=BG$. In this case, the functor $S_i$ admits a concrete description (see \Cref{lem_Ltemp_test_null}) from which it is clear that it preserves null local systems.
\end{proof}
\end{thm}

As the first corollary, we deduce the following result on functoriality of tempered local systems for orbispaces.
\begin{lem}\label{lem_pb_commutes_with_tempereization_orbisps}
Let $f\colon X \to Y$ be a (not necessarily faithful) morphism of orbispaces. Then the commutative diagram (see \Cref{pb_preserves_LStemp})
\[\xymatrix{
\LocSys^\temp(Y, A) \ar[r]^-{f^{\temp,*}}\ar@{^(->}[d]^{i_Y} & \LocSys^\temp(X, A) \ar@{^(->}[d]^{i_X} \\
\LocSys^\Glo_\ab(Y, A) \ar[r]^-{f^{\Glo,*}} & \LocSys^\Glo_\ab(X, A)
}\]
is vertically left adjointable and horizontally right adjointable. In other words, there are natural equivalences
$$f^{\temp,*} \circ L_Y \simeq L_X \circ f^{\Glo,*}, \qquad f^{\Glo}_{*} \circ i_X \simeq i_Y \circ f^{\temp}_{*}.$$

\begin{proof}
By passing to adjoints, these two assertions are equivalent to each other, so we prove only the latter. Let $\mathcal L$ be a tempered local system on $X$. We will show that $f^\Glo_{*} \mathcal L$ is tempered. By \Cref{thm_Ltemp_kernel} (applied to $Y$), it is enough to show that
$$\Hom_{\LocSys^\Glo_\ab(Y, A)}(\mathcal M, f^{\Glo}_{*} \mathcal L) \simeq \Hom_{\LocSys^\Glo_\ab(X, A)}(f^{\Glo,*} \mathcal M, \mathcal L)$$
vanishes for all $\mathcal M \in \LocSys^\nnull(Y, A)$. Since the null local systems are preserved under pullbacks by \Cref{rem_functoriality_and_descent_LSnull}, the result follows by \Cref{thm_Ltemp_kernel} again (now applied to $X$).
\end{proof}
\end{lem}

Our next goal is to extend these results to arbitrary global spaces.
\begin{cor}\label{LStemp_admits_left_adj_global}
$\quad$
\begin{enumerate}
\item Let $X$ be a global space. Then the inclusion $i_X\colon \LocSys^\temp(X, A) \inj \LocSys^\Glo_\ab(X, A)$ admits a left adjoint $L_X$. The kernel of $L_X$ is $\LocSys^\nnull(X, A)$.

\item Let $f\colon X \to Y$ be a morphism of global spaces. Then the commutative diagram
\[\xymatrix{
\LocSys^\temp(Y, A) \ar[r]^-{f^{\temp,*}}\ar@{^(->}[d]^{i_Y} & \LocSys^\temp(X, A) \ar@{^(->}[d]^{i_X} \\
\LocSys^\Glo_\ab(Y, A) \ar[r]^-{f^{\Glo,*}} & \LocSys^\Glo_\ab(X, A)
}\]
is vertically left adjointable and horizontally right adjointable. That is there are natural equivalences
$$f^{\temp,*} \circ L_Y \simeq L_X \circ f^{\Glo,*}, \qquad f^{\Glo}_{*} \circ i_X \simeq i_Y \circ f^{\temp}_{*}.$$
\end{enumerate}

\begin{proof}
One can present $X$ as a colimit of a small diagram of representable global spaces, $X \simeq \colim BG_\alpha$. Then, by \Cref{prop_LSglo_preserves_limits}, \Cref{proposition: temp descent}, and \Cref{pb_preserves_LStemp}, we have
$$i_X \simeq \lim_\alpha i_{BG_\alpha}.$$
Moreover, since the left adjoints $L_{BG_\alpha}$ are compatible with pullbacks by \Cref{lem_pb_commutes_with_tempereization_orbisps}, the limit
$$L_X \simeq \lim_\alpha L_{BG_\alpha}$$
is well defined and it is a left adjoint to $i_X$. Moreover,
$$\ker L_X \simeq \lim_\alpha \ker L_{BG_\alpha}.$$
By \Cref{thm_Ltemp_kernel} and \Cref{rem_functoriality_and_descent_LSnull}, we deduce that
$$\ker L_X \simeq \lim_\alpha \LocSys^\nnull(BG_\alpha, A) \simeq \LocSys^\nnull(X, A).$$
The second part can be proved now in the same way as the second part of \Cref{lem_pb_commutes_with_tempereization_orbisps}.
\end{proof}
\end{cor}

Finally, by using the semi-orthogonal decomposition above, we can construct the monoidal structure on $\LocSys^\temp$.
\begin{thm}\label{thm_CMon_str_on_LStemp}
Let $X$ be a global space. Then there exists an essentially unique symmetric monoidal structure on $\LocSys^\temp(X, A)$ such that the localization functor
$$L_X \colon \LocSys^\Glo_\ab(X, A) \tto \LocSys^\temp(X, A)$$
is symmetric monoidal. For a morphism of global spaces $f\colon X \to Y$ the pullback functor
$$f^{\temp,*}\colon \LocSys^\temp(Y, A) \tto \LocSys^\temp(X, A)$$
is naturally symmetric monoidal. More precisely, the functor
$$\LocSys^\temp(-, A) \colon \Type^{\Glo, \op} \tto \Prs^\LL$$
lifts naturally to a functor to $\CAlg(\Prs^\LL)$ compatibly with the corresponding lift for $\LocSys^\Glo_\ab(-, A)$.

\begin{proof}
By \cite[Proposition 2.2.1.9]{Lur_HA}, the first assertion follows from the fact that $\LocSys^\nnull(X, A)$ is a $\otimes$-ideal in $\LocSys^\Glo_\ab(X, A)$. This holds because the monoidal structure on $\LocSys^\Glo_\ab(X, A)$ is pointwise and the torsion subcategory $\QCoh(A(BG))^\tors$ is a $\otimes$-ideal in $\QCoh(A(BG))$ for each compact abelian Lie group $G$.

The second assertion follows formally from the second part of \Cref{LStemp_admits_left_adj_global} and the fact that the functor
$$\LocSys^\Glo_\ab(-, A) \colon \Type^{\Glo, \op} \tto \Prs^\LL$$
factors through the forgetful functor $\CAlg(\Prs^\LL) \to \Prs^\LL$.
\end{proof}
\end{thm}
\begin{cor}\label{cor_sharp_pf_LStemp}
Let $f\colon X \to Y$ be a faithful morphism of global spaces. Then the pullback functor
$$f^{\temp,*}\colon \LocSys^\temp(Y, A) \tto \LocSys^\temp(X, A)$$
admits a $\LocSys^\temp(Y, A)$-linear left adjoint $f^{\temp}_{\#}$. Moreover, for a fiber square of global spaces
\[\xymatrix{
X^\prime \ar[r]^-q \ar[d]^g & X \ar[d]^f \\
Y^\prime \ar[r]^-p & Y,
}\]
the natural base change map
$$g^{\temp}_{\#} \circ q^{\temp,*} \tto p^{\temp,*} \circ f^{\temp}_{\#}$$
is an equivalence.

\begin{proof}
By \Cref{prop_LSGlo_sharp_pushforward}, the pullback functor
$$f^{\Glo,*} \colon \LocSys^\Glo_\ab(Y, A) \tto \LocSys^\Glo_\ab(X, A)$$
admits a $\LocSys^\Glo_\ab(Y, A)$-linear left adjoint $f^{\Glo}_{\#}$. It follows formally from \Cref{pb_preserves_LStemp} and \Cref{LStemp_admits_left_adj_global} that the functor
$$f^{\temp}_{\#} = L_Y \circ f^{\Glo}_{\#} \circ i_X$$
is a left adjoint to $f^{\temp,*}$. First, we will show that the functor $f^{\temp}_{\#}$ is $\LocSys^{\temp}(Y,A)$-linear. More precisely, we show that the natural map
$$f^{\temp}_{\#}(\mathcal L \otimes^\temp f^{\temp,*}\mathcal M) \tto f^{\temp}_{\#}(\mathcal L) \otimes^\temp \mathcal M$$
is an equivalence for all $\mathcal L \in \LocSys^\temp(X, A)$ and $\mathcal M \in \LocSys^\temp(Y, A)$. The left hand side is equivalent to
$$f^{\temp}_{\#}(\mathcal L) \otimes^\temp \mathcal M \simeq L_Y(L_Y f^{\Glo}_{\#} \mathcal L \otimes^\Glo \mathcal M) \simeq L_Y(f^{\Glo}_{\#} \mathcal L \otimes^\Glo \mathcal M) \simeq L_Y f^{\Glo}_{\#}(\mathcal L \otimes^\Glo f^{\Glo,*}\mathcal M),$$
where the first equivalence is the construction of tempered monoidal product, the second equivalence follows from the fact that $\ker L_Y$ is a $\otimes$-ideal, and the last equivalence is the projection formula in the $\LocSys^\Glo$-case. Finally, by \Cref{lem_sharp_pf_preserves_null}, the natural map
$$L_Y f^{\Glo}_{\#}(\mathcal L \otimes^\Glo f^{\Glo,*}\mathcal M) \tto L_Y f^{\Glo}_{\#}(L_X(\mathcal L \otimes^\Glo f^{\Glo,*}\mathcal M)) \simeq f^{\temp}_{\#}(\mathcal L \otimes^\temp f^{\Glo,*} \mathcal M)$$
is an equivalence.

The result about the base change follows formally from the corresponding assertion in the global case (\Cref{LSGlo_basechange_faithful_pbs}) and the second part of \Cref{LStemp_admits_left_adj_global}.
\end{proof}
\end{cor}

\subsection{Tempered local systems on classifying stacks of tori}\label{section: tempered local of tori}
In this section,  we work out a concrete description of the category $\LocSys^\temp(BT, A)$ for a connected compact abelian Lie group $T$ in algebro-geometric terms of  the stack $A(BT)$.

Recall from \Cref{LSglo_enhached_GS3} that there is a natural enhancement of the global section functor
$$\Gamma_A(X, -) \colon \LocSys^\Glo_\ab(X, A) \tto \QCoh(A(X))$$
for a global space $X$. By abuse of notation, we will denote the composite of $\Gamma_A(X, -)$ with the inclusion
$$i_X \colon \LocSys^\temp(X, A) \inj \LocSys^\Glo_\ab(X, A)$$
also by $\Gamma_A(X, -)$.
\begin{thm}\label{LStemp_on_tori}
Let $A$ be an oriented non-connective spectral abelian group stack over a locally complex periodic base $S$. Then, for a connected compact abelian Lie group $T$, the global section functor
$$\Gamma_A(BT, -) \colon \LocSys^\temp(BT, A) \tto \QCoh(A(BT))$$
is an equivalence.

\begin{proof}
By construction and \Cref{thm_CMon_str_on_LStemp}, the functor $\Gamma_A(BT, -)$ admits a symmetric mo\-noi\-dal left adjoint 
$$\mathbbl 1_{BT}^\temp\otimes - \colon \QCoh(A(BT)) \to \LocSys^\temp(BT, A)$$ which is given by the composition
$$\mathbbl 1_{BT}^\temp\otimes - \simeq L_{BT}(\mathbbl 1_{BT}^\Glo \otimes -).$$
First, we claim that $\mathbbl 1_{BT}^\temp\otimes -$ is fully faithful. Indeed, by \Cref{from_glob_sects_are_tempered}, the object $\mathbbl 1_{BT}^\Glo \otimes \mathcal F$ is tempered for every $\mathcal F \in \QCoh(A[\widehat T])$, hence
$$\Gamma_A(BT, \mathbbl 1_{BT}^\temp \otimes \mathcal F) \simeq \Gamma_A(BT, \mathbbl 1_{BT}^\Glo \otimes \mathcal F).$$
So, the assertion follows since the functor $\mathbbl 1_{BT}^\Glo \otimes -$ is fully faithful, see \Cref{LSglo_enhacned_GS}.

To finish the proof, it is enough to show that $\Gamma_A(BT, -)$ is conservative. Indeed, let $\mathcal{L}\in \LocSys^\temp(BT,A)$ be a tempered local system such that $\mathcal{L}^{BT} \simeq 0$. We will show that $\mathcal{L}^{BH}\simeq 0$ for all proper closed subgroups $H\subsetneq T$ as well. %For this end, assuming that $\mathcal L^{BT} \simeq 0$ we need to prove that $\mathcal L^{BH} \simeq 0$ for all proper subgroups $H$ of $T$.
Let $i\colon A[\widehat H] \inj A[\widehat T]$ be the map induced by a faithful morphism $BH \inj BT$. Since $\mathcal{L}$ is tempered, we have
$$0 \simeq (\mathcal L^{BT})_{(H,T)}^\wedge \xymatrix{\ar[r]^-\sim & } (i_{*}\mathcal L^{BH})^{hT/H}.$$
Moreover, by the norm cofiber sequence, we observe
$$\xymatrix{\Sigma^{\dim{T/H}} (i_* \mathcal L^{BH})_{hT/H} \ar[r] & (i_* \mathcal L^{BH})^{hT/H} \ar[r] & (i_* \mathcal L^{BH})^{tT/H}}.$$
So, the vanishing of the middle term implies that
$$(i_* \mathcal L^{BH})_{hT/H} \simeq \Sigma^{-\dim T/H-1} (i_* \mathcal L^{BH})^{tT/H}.$$
Since, by \Cref{LStemp_TateInvs_are_local}, the right hand side is $(H,T)$-local and the left hand side is $(H,T)$-torsion both must be trivial. By \Cref{loc_comp_or_over}, $A[\widehat T]$ is locally complex orientable. So, $i_* \mathcal L^{BH} \simeq 0$ by the virtue of \Cref{lem_coinvariants_conservativity_com_or}. Since $i_*$ is conservative, $\mathcal L^{BH} \simeq 0$.
\end{proof}
\end{thm}

\begin{rem}\label{remark: non-conservative in general}
We warn the reader that the functor 
$$\Gamma_A(BG,-)\colon \LocSys^{\temp}(BG,A) \tto \QCoh(A[\widehat{G}]) $$
is \emph{not} necessary conservative for a non-connected compact abelian Lie group $G$. Indeed, let $G=C_p$ and $A=\mathbb{G}_{a,R}$, where $R=\Spec(\mathbb{Q}[\beta,\beta^{-1}])$, see \cref{example: additive group}. Then, by e.g. \cref{example: cyclic group global}, a tempered local system $\mathcal{L}\in \LocSys^{\temp}(BC_p,A)$ is a triple 
$$\mathcal{L}=(\mathcal{L}^{C_p}\in \Mod_R,\mathcal{L}^{e}\in \LocSys(BC_p,\Mod_R), \alpha\colon \mathcal{L}^{C_p} \to (\mathcal{L}^e)^{hC_p})$$
such that $\alpha$ is an equivalence. So, as before, $\Gamma_A(BG,\mathcal{L})\simeq 0$ implies $(\mathcal{L}^e)^{hC_p} \simeq 0$. However, there are many non-trivial local systems $\mathcal{L}^{e}\in \LocSys(BC_p,\Mod_R)$ such that $(\mathcal{L}^{e})^{hC_p}\simeq 0$, e.g.~$\mathcal{L}^e=V \otimes_{\mathbb{Q}} R$ such that $V$ is a non-trivial rational $C_p$-representation and $V^{C_p} = 0$.
\end{rem}

\section{Comparison of tempered and genuine local systems}\label{section: comparison of tempered and genuine local systems}
The main goal of this section is to compare the category $\LocSys^\temp(X, A)$ of tempered local systems with the category $\LocSys^\gen_\ab(X, A)$ of genuine local systems. Namely, we show in \cref{theorem: temp are modules in gen} that $\LocSys^\temp(X, A)$ is a \emph{smashing localization} of $\LocSys^\gen_\ab(X, A)$ if $X$ is an orbispace. As a corollary, if $X=BG$ for a compact Lie group $G$ and $A$ is nc-affine, then there are identifications 
$$\LocSys^\temp(BG, A)\simeq \Mod_{R_G}(\Sp^G) \;\; \text{and} \;\; \LocSys^\gen_\ab(BG, A) \simeq \Mod_{\Sigma^\infty_{BG} \Omega^{\infty}_{BG}R_G}(\Sp^G)$$
for a genuine equivariant cohomology theory $R_G \in \CAlg(\Sp^G)$ and $R_G$ is an \emph{idempotent algebra} over the suspension object $\Sigma^\infty_{BG} \Omega^{\infty}_{BG}R_G$. So, although the category $\LocSys^\gen_\ab(X, A)$ does not recover genuine equivariant cohomology theories in general, a certain smashing localization of it does, at least in the complex periodic case.

In \cref{section: thom is tensor-invertible}, we compute the tempered Thom local systems using \cref{LStemp_on_tori} and we show that those are $\otimes$-invertible in $\LocSys^\temp(X, A)$. Therefore, the localization 
$$L_X\colon \LocSys^{\Glo}_\ab(X,A) \to \LocSys^\temp(X, A)$$
factors through the functor 
$$L^\gen_X\colon \LocSys^{\gen}_\ab(X,A) \to \LocSys^\temp(X, A).$$
In \cref{section: smashing localization}, we show that $L^\gen_X$ admits a right adjoint $\beta_X$ and we compute the essential image of $\beta_X$ in terms of geometric fixed points in \cref{proposition: geometric fixed points of temp}. Using this computation, we deduce that $\beta_X$ is strictly $\LocSys^{\gen}_\ab(X,A)$-linear in \cref{proposition: gen-to-temp is smashing} which implies the equivalence 
$$\LocSys^\temp(X, A) \simeq \Mod_{\beta_X(\mathbbl{1})}\LocSys^{\gen}_\ab(X,A)$$ of the main \cref{theorem: temp are modules in gen}. By the previously obtained description of the essential image for $\beta_X$, we also show that $\beta_X(\mathbbl{1})$ is an idempotent algebra.

In \cref{section: dexterity}, we study the category $\LocSys^{\temp}(X,A)$ as a module in $\Prs^\LL$ over the symmetric monoidal category $\LocSys^{\Glo}_\ab(X,A)\in \CAlg(\Prs^\LL)$. Namely, in \cref{theorem: dextery faithful tempered},  we show that the diagram
\[
\xymatrix{
\LocSys^\Glo_\ab(X, A) \ar[r]^-{L_{X}}  & \LocSys^\temp(X, A) \\
\LocSys^\Glo_\ab(Y, A) \ar[r]^-{L_{Y}} \ar[u]^{f^{\Glo,*}} & \LocSys^{\temp}(Y, A) \ar[u]_{f^{\temp,*}}
}\]
is a pushout square in $\CAlg(\Prs^\LL)$ if $f\colon X\to Y$ is a faithful morphism of orbispaces. This observation simplifies some constructions in our future work.
%In this section we show that $\LocSys^\temp(X, A)$ is a smashing localization of $\LocSys^\gen_\ab(X, A)$.
\subsection{\texorpdfstring{$\otimes$}{Tensor}-invertibility of tempered Thom local systems}\label{section: thom is tensor-invertible}
\begin{notation}
Let $G$ be a compact Lie group. For a $G$-homotopy type $X$ with the structure (faithful) map $p\colon X \gitq G \to BG$ let us denote
$$C_*^{BG}(X, A) := p^{\temp}_{\#}(\mathbbl 1_{X\gitq G}),\ C^*_{BG}(X, A) := p^{\temp}_*(\mathbbl 1_{X \gitq G}) \quad \in \quad \LocSys^\temp(BG, A).$$
Let $\overline{C}_*^{BG}(X, A)$ denote the fiber of the natural map $C_*^{BG}(X, A) \to C_*^{BG}(*, A) \simeq \mathbbl 1_{BG}$. The object $\overline C^*_{BG}(X, A)$ is defined similarly. We note that 
$$C^*_{BG}(X, A) \simeq [C_*^{BG}(X, A), \mathbbl 1_{BG}^\temp]$$
by the projection formula for tempered $\#$-pushforward (\Cref{cor_sharp_pf_LStemp}), where $[-,-]$ is the internal $\Hom$ in $\LocSys^\temp(BG, A)$.
\end{notation}
\begin{rem}
Let $G$ be a compact Lie group and let $V$ be a finite dimensional real representation of $G$. Then the Thom local system $\mathbbl 1_{BG}^{V, \temp}$ is equivalent to $\overline C^{BG}_*(\mathbb S^V, A)$, where $\mathbb S^V$ is the corresponding representation sphere.
\end{rem}
\begin{rem}
Let $G$ be a compact abelian Lie group and let $X$ be a $G$-space. Then
$$\Gamma_A(BG, C^*_{BG}(X, A)) \simeq p_* \mathcal O_{A(X \gitq G)},$$
where $p\colon A(X \gitq G) \to A(BG)$ denotes the induced map on $A$-cohomology, see \cref{section: global section}. In particular, if $G = T$ is a torus and $H \subseteq T$ is a closed subgroup, then the object $C^*_{BT}(T/H, A) \in \LocSys^\temp(BT, A)$ corresponds under the identification of \Cref{LStemp_on_tori}
$$\Gamma_A(BT, -) \colon \LocSys^\temp(BT, A) \xymatrix{\ar[r]^-\sim & } \QCoh(A(BT))$$
to the structure sheaf of the closed embedding $A(BH) \inj A(BT)$.
\end{rem}

The proof of the following assertion is the same as for \Cref{lem_excision_for_global_homology}.
\begin{lem}\label{lem_excision_for_tempered_homology}
The functor
$$C_*^{BG}(-, A) \colon \Type^G \tto \LocSys^\temp(BG, A)$$
preserves small colimits. \qedhere
\end{lem}

\begin{notation}
Let $X$ be a global space and let $E$ be a vector bundle on $X$. We denote $$\mathbbl 1_X^{E, \temp} := L_X \mathbbl 1_X^{E,\Glo},$$
the image of the Thom local system from \Cref{def_ThomLS_glo} under the localization functor
$$L_{X}\colon \LocSys^\Glo_\ab(X, A) \tto \LocSys^\temp(X, A)$$
of \Cref{LStemp_admits_left_adj_global}. Equivalently,
$$\mathbbl 1_X^{E, \temp} \simeq \cofib\left(p_{\#}^{\temp}\mathbbl 1_{E^\circ} \tto \mathbbl 1_X\right),$$
where $p\colon E^\circ \to X$ denotes the structure map.
\end{notation}
\begin{rem}\label{Thom_LS_basics_temp}
Let $f\colon X \to Y$ be a map of global spaces and let $E$ be a vector bundle on $Y$. Then there is a natural equivalence
$$\mathbbl 1_X^{f^*E,\temp} \simeq f^*\mathbbl 1_Y^{E,\temp}.$$
For $E_1$, $E_2$ a pair of vector bundles on $X$ there is a natural equivalence
$$\mathbbl 1_X^{E_1\oplus E_2, \temp} \simeq \mathbbl 1_X^{E_1, \temp} \otimes 1_X^{E_2, \temp}.$$ 
One can formally deduce these assertions from the $\LocSys^\Glo$-case \Cref{Thom_LS_basics_Glo}.
\end{rem}
We will need the following basic computation.
\begin{prop}\label{Thom_for_characters}
Let $T$ be a torus and let $\chi\colon T \to \mathbb C^\times$ be a character. Then, under the identification $\LocSys^\temp(BT, A) \simeq \QCoh(A(BT))$ of \Cref{LStemp_on_tori}, we have
$$\mathbbl 1_{BT}^{\chi,\temp} \simeq \begin{cases}
\mathcal O_{A(BT)}(\chi), & \chi \text{ is non-trivial},\\
\mathcal O_{A(BT)}[2], & \chi \text{ is trivial}.
\end{cases}$$
Here, $\mathcal O_{A(BT)}(\chi)$ is the line bundle corresponding to the divisor $A(B\ker \chi) \inj A(BT)$ induced by the group embedding $\ker \chi \inj T$.

\begin{proof}
We will show first that $\mathbbl 1_{BT}^{\chi,\temp}$ is dualizable. Indeed, since any $T$-character is a restriction of the tautological representation of $U(1)$, we can assume without loss of generality that $T = U(1)$ and $\chi = \mathbb C(1)$. The representation sphere $\mathbb S^{\mathbb C(1)}$ admits a $U(1)$-equivariant cell decomposition with one cell being $\mathbb S^0$ and another being $U(1)$ with the tautological action. By \cref{lem_excision_for_tempered_homology}, we have $\mathbbl 1_{BU(1)}^{\mathbb C(1),\temp}$ is an extension of $\mathbbl 1_{BU(1)}^\temp$ and $i^{\temp}_{\#} \mathbbl 1^\temp_*$, where
$$i\colon * \tto BU(1)$$
is a basepoint. Since dualizable objects are closed under extensions, it suffices to show that $i^{\temp}_{\#} \mathbbl 1^\temp_*$ is dualizable. However, under the identifications
$$\LocSys^\temp(*, A) \simeq \QCoh(S) \quad\text{and}\quad \LocSys^\temp(BU(1), A) \simeq \QCoh(A),$$
the $\#$-pushforward corresponds to the left adjoint $i^{*, L}$ of the $*$-pullback functor
$$i^* \colon \QCoh(A) \tto \QCoh(S).$$
The map
$$i\colon S \tto A$$
is non-connective affine, and by \Cref{lem_fsmooth_dualizable} the sheaf $i_* \mathcal O_S$ is dualizable. It follows that 
$$i_{\#}(\mathbbl 1_S) \simeq i^{*, L}(\mathcal O_S) \simeq i_*(\mathcal O_S)^\vee,$$
and is dualizable.

Now, we finish the computation. If the character $\chi$ is trivial, then $\mathbb S^{\chi} \simeq \mathbb S^2$ as $T$-spaces. Therefore, by \Cref{lem_excision_for_tempered_homology}, we have
$$\overline C_*^{BT}(\mathbb S^2, A) \simeq \mathcal O_{A(BT)}[2].$$
If $\chi$ is non-trivial, then there is a pushout square of $T$-homotopy types
\[\xymatrix{
T/\ker \chi \ar[r]\ar[d] & {*} \ar[d] \\
{*} \ar[r] & \mathbb S^\chi.
}\]
Applying $C^*_{BT}(-, A)$ to this square we obtain a fiber sequence
\[\xymatrix{\overline C^*_{BT}(\mathbb S^\chi, A) \ar[r] & \mathcal O_{A(BT)} \ar[r] & i_* \mathcal O_{A(B\ker \chi)},}\]
where $i\colon A(B\ker \chi) \inj A(BT)$ is the natural embedding. It follows that 
$$\overline C^*_{BT}(\mathbb S^\chi, A) \simeq \mathcal O_{A(BT)}(-\chi).$$
Finally, since $\overline C_*^{BT}(\mathbb S^\chi, A)$ is dualizable, we have 
\[\overline C_*^{BT}(\mathbb S^\chi, A) \simeq \left(\overline C^*_{BT}(\mathbb S^\chi, A)\right)^\vee \simeq\mathcal O_{A(BT)}(\chi).\qedhere\]
\end{proof}
\end{prop}

\begin{lem}\label{lemma: temp_thom is invertible}
Let $X$ be a global space and let $E$ be a vector bundle over $X$. Then $\mathbbl 1_X^{E, \temp}$ is $\otimes$-invertible in $\LocSys^\temp(X, A)$.

\begin{proof}
Since $\LocSys^\temp(X, A)$ is a limit of $\LocSys^\temp(BG, A)$ for various orbits mapping to $X$ it is enough to show that $x^* \mathbbl 1_X^{E,\temp}$ is $\otimes$-invertible for all compact abelian Lie groups $G$ and maps $x\colon BG \to X$. Since
$$x^* \mathbbl 1_X^{E,\temp} \simeq \mathbbl 1_{BG}^{x^* E, \temp},$$
we can assume that $X = BG$. Moreover, we can embed $G$ into some torus $T$ and extend $E$ to a representation of $T$. Since a finite dimensional representation of $T$ splits as a direct sum of characters and by the multiplicativity of the Thom local systems (see \Cref{Thom_LS_basics_temp}), we can assume that $E$ is rank $1$. The result then follows from the computation of \Cref{Thom_for_characters}.
\end{proof}
\end{lem}

\begin{cor}\label{corollary: existence of gen-to-temp local}
Let $X$ be a global space. Then the localization functor
$$L_X \colon \LocSys^\Glo_\ab(X, A) \tto \LocSys^\temp(X, A)$$
factors through a symmetric monoidal colimit-preserving functor $L_X^{\gen}$:
\[\xymatrix{
\LocSys^\Glo_\ab(X, A) \ar[d] \ar[r]^-{L_X} & \LocSys^\temp(X, A) \\
\LocSys^\gen_\ab(X, A). \ar@{-->}[ru]_-{L_X^{\gen}}
}\]
Moreover, the functor $L_{-}^\gen$ commutes with $*$-pullbacks.

\begin{proof}
By \cref{construction: genuine local system} and by \Cref{proposition: temp descent}, we can assume that $X = BG$ for some compact abelian Lie group $G$. In this case the existence $L_X^{\gen}$ is provided by the universal property of \cref{definition: object inversion}. Indeed, for any vector bundle $E$ on $BG$, the corresponding Thom local system $\mathbbl 1_{BG}^{E,\temp} = L_{BG} \mathbbl 1_{BG}^{E,\Glo}$ is $\otimes$-invertible in $\LocSys^\temp(BG, A)$ by \cref{lemma: temp_thom is invertible}.
\end{proof}
\end{cor}

\begin{cor}\label{corollary: gen-to-temp local is functor along sharp}
Let $f\colon X\to Y$ be a faithful morphism of global spaces. Then the commutative square
\[
\xymatrix{
\LocSys^\gen_\ab(X, A) \ar[r]^-{L^{\gen}_{X}}  & \LocSys^\temp(X, A) \\
\LocSys^\gen_\ab(Y, A) \ar[r]^-{L^{\gen}_{Y}} \ar[u]^{f^{\gen,*}} & \LocSys^{\temp}(Y, A) \ar[u]^{f^{\temp,*}}
}\]
is vertically left adjointable. In other words, the natural transformation
$$f^{\temp}_{\#} \circ L^{\gen}_X \tto L^{\gen}_Y \circ f^{\gen}_{\#} $$
is an equivalence.
\end{cor}

\begin{proof}
By \cref{proposition: temp descent}, \cref{cor_sharp_pf_LStemp}, and the third part of \cref{proposition: genuine sharp}, we can assume that $Y=BG$ for some compact abelian Lie group $G$. Then $X$ is an orbispace. Hence, by \cref{cor_OmegaInftygen}, it suffices to show that 
$$f^{\temp}_{\#} (L^{\gen}_X \mathcal{M}) \tto L^{\gen}_Y(f^{\gen}_{\#} \mathcal{M}) $$
is an equivalence for $\mathcal{M} \simeq \Sigma^{\infty}_{X} \mathcal{M}_0$, $\mathcal{M}_0 \in \LocSys^\Glo_\ab(X,A)$. Finally, the assertion follows by the second part of \cref{proposition: genuine sharp} and \cref{pb_preserves_LStemp}.
\end{proof}

\subsection{\texorpdfstring{$L^\gen_X$}{Lgen} is a smashing localization}\label{section: smashing localization}
After the preparation from the previous subsection, we prove here the main result of this section \Cref{theorem: temp are modules in gen}.
\begin{construction}\label{construction: right adjoint to temp localization}
Let $X$ be a global space. By \cref{corollary: existence of gen-to-temp local}, the functor
$$L_X^{\gen}\colon \LocSys^\gen_\ab(X,A) \tto \LocSys^{\temp}(X,A) $$
admits a right adjoint
$$\beta_X \colon \LocSys^{\temp}(X,A) \tto \LocSys^\gen_\ab(X,A). $$
By \cref{corollary: gen-to-temp local is functor along sharp}, the functor $\beta_{-}$ commutes with $*$-pullbacks along \emph{faithful} morphisms.
\end{construction}

\begin{lem}\label{lemma: beta is cons and cont}
Let $X$ be an orbispace. Then the functor $\beta_X$ is conservative and preserves small colimits.
\end{lem}

\begin{proof}
By \cref{corollary: existence of gen-to-temp local}, the composite $\Omega^{\infty}_X \circ \beta_X$ is equivalent to the canonical fully faithful embedding
$$\LocSys^{\temp}(X,A) \tto \LocSys^\Glo_\ab(X,A). $$
Hence, the lemma follows by \cref{cor_OmegaInftygen} and \cref{LStemp_presentable}.
\end{proof}

\begin{defn}\label{definition: proper localization}
Let $G$ be a compact Lie group. We say that a morphism $\theta\colon \mathcal{F} \to \mathcal{F}'$ in $\QCoh(A(BG))$ \emph{exhibits $\mathcal{F}'$ as the $(\fcat{P}, G)$-localization of $\mathcal{F}$} if
\begin{enumerate}
\item $\mathcal{F}'$ is $(H,G)$-local for all proper closed subgroups $H\subset G$;
\item the fiber of $\theta$ lies in $\QCoh(A(BG))^{\tors}$, see \cref{definition: proper torsion}.
\end{enumerate}
\end{defn}

Recall the geometric fixed points functor construction \Cref{construction: geometric fixed points}.
\begin{lem}\label{lemma: fiber of cat-to-geom is tors}
Let $X=BG$ and $\mathcal{L} \in \LocSys^\gen_\ab(BG,A)$ be a genuine local system. Then the fiber of the canonical morphism
$$\theta \colon \mathcal{L}^G \tto \Phi^G\mathcal{L} \in \QCoh(A(BG))$$
lies in $\QCoh(A(BG))^{\tors}$.
\end{lem}

\begin{proof}
By \cref{corollary: isotropy separation for gen}, there is an equivalence 
$$\fib(\theta)\simeq (\pi^{\gen}_{\#}\pi^{\gen,*}\mathcal{L})^G,$$
where $$\pi\colon (BG)_{\fcat{P}}=BG\times_{\mathcal{N}} \mathcal{N}_{\fcat{P}}\tto BG$$
is the projection. Since the full subcategory $\QCoh(A(BG))^{\tors}$ is closed under colimits in $\QCoh(A(BG))$, it suffices to show that 
$$(f^{\gen}_{\#} \mathcal{L})^G \in \QCoh(A(BG))^{\tors}$$ for all $\mathcal{L} \in \LocSys^\gen_\ab(BH,A)$ and $f\colon BH \to BG$, where $H \subset G$ is a proper closed subgroup. By \cref{cor_OmegaInftygen}, the category $\LocSys^\gen_\ab(BH,A)$ is generated under colimits by the suspension objects. Hence, we can further assume that $\mathcal{L} \simeq \Sigma^{\infty}_{BH}\mathcal{L}_0$ for some globally equivariant local system $\mathcal{L}_0 \in \LocSys^\Glo_\ab(BH,A)$.

By \cref{lemma: sharp for gen orbits} and \cref{corollary: right adjoint to inversion many objects}, we have
\begin{align*}
(f^{\gen}_{\#} \Sigma^{\infty}_{BH} \mathcal{L}_0)^G &\simeq (\Sigma^{\infty}_{BG}f^{\Glo}_{\#}  \mathcal{L}_0)^G \\
&\simeq \colim_{V \in \Rep(G)} [\mathbbl 1^V_{BG},\mathbbl 1^V_{BG}\otimes f^{\Glo}_\# \mathcal{L}_0]^G \\
&\simeq \colim_{V \in \Rep(G)} [\mathbbl 1^V_{BG},f^{\Glo}_\#(\mathbbl 1^{f^*V}_{BH}\otimes \mathcal{L}_0)]^G.
\end{align*}
Therefore, it is enough to show that 
$$[\mathbbl 1^V_{BG},f^{\Glo}_\#\mathcal{L}_0]^G \in \QCoh(A(BG))^{\tors}$$
for all globally equivariant local systems $\mathcal{L}_0 \in \LocSys^\Glo_\ab(BH,A)$ and for all $G$-representations $V$. However, since the one-point compactification $\mathbb S^V$ is a finite $G$-space, it suffices to show that
$$[p^{\Glo}_{\#}\mathbbl{1},f^{\Glo}_{\#}\mathcal{L}_0]^G \in \QCoh(A(BG))^{\tors} $$
for all globally equivariant local systems $\mathcal{L}_0 \in \LocSys^\Glo_\ab(BH,A)$, all faithful morphisms $p\colon BK \to BG$, and $f\colon BH \to BG$, where $H \subset G$ is a proper closed subgroup. Finally, we observe
$$[p^{\Glo}_{\#}\mathbbl{1},f^{\Glo}_{\#}\mathcal{L}_0]^G \simeq (p^{\Glo}_{*}p^{\Glo,*}f^{\Glo}_{\#}\mathcal{L}_0)^G \simeq 
\begin{cases}
p_*\left((p^{\Glo,*}f^{\Glo}_{\#}\mathcal{L}_0)^K\right), & \mbox{if $K$ is a proper subgroup,}\\
0, & \mbox{if $K=G$}
\end{cases}
$$
by using the projection formula for $p^{\Glo}_{\#}$ (see \cref{corollary: projection formula global Hom}) together with \cref{example: star-pullback fixed points} and  \cref{prop_LSGlo_sharp_pushforward} for final computations.
\end{proof}

\begin{prop}\label{proposition: geometric fixed points of temp}
Let $X=BG$ and $\mathcal{L} \in \LocSys^{\temp}(BG,A)$ be a tempered local system. Set $\beta=\beta_{BG}$. Then the canonical morphism
$$\theta\colon \mathcal{L}^G\simeq (\Omega^{\infty}_{BG}(\beta\mathcal{L}))^G \tto \Phi^G(\beta\mathcal{L}) \in \QCoh(A(BG)) $$
exhibits the geometric fixed points $\Phi^G(\beta\mathcal{L})$ as the $(\fcat{P}, G)$-localization of $\mathcal{L}^G$.
\end{prop}

\begin{proof}
By \cref{lemma: fiber of cat-to-geom is tors}, the fiber $\fib(\theta) \in \QCoh(A(BG))^{\tors}$. We will show that $\Phi^G(\beta\mathcal{L})$ is $(H,G)$-local for every proper closed subgroup $H\subset G$. Let $\chi\colon G \to U(1)$ be a non-trivial character such that $H\subset \ker(f)$. Denote by $f$ the induced map
$$A(B\chi) \colon A(BG) \tto A(U(1)) \simeq A.$$
Consider the fiber sequence 
$$\mathcal{O}_A(-e) \xrightarrow{u^{-1}} \mathcal{O}_A \tto i_*\mathcal{O}_S $$
of quasi-coherent sheaves over $A$ from \cref{lem_fsmooth_dualizable}. Recall that $\mathcal{O}_A(-e)$ is $\otimes$-invertible and let us denote by 
$$u\colon \mathcal{O}_{A} \tto \mathcal{O}_{A}(e)$$
the global section of the line bundle $\mathcal{O}_{A}(e)$ dual to $u^{-1}$. By \cref{rem_locality_in_smooth_case}, $\Phi^G(\beta\mathcal{L})$ is $(H,G)$-local if and only if the tensor product $\Phi^G(\beta\mathcal{L}) \otimes f^*i_*\mathcal{O}_S$ is trivial. Equivalently, $\Phi^G(\beta\mathcal{L})$ is $(H,G)$-local if and only if the map
$$f^*(u)\otimes \Phi^G(\beta\mathcal{L})\colon \Phi^G(\beta\mathcal{L}) \tto f^*\mathcal{O}_{A}(e) \otimes \Phi^G(\beta \mathcal{L})$$
is an equivalence. Moreover, we have 
$$f^*\mathcal{O}_{A}(e) \simeq (\mathbbl{1}^{\Glo} \otimes f^*\mathcal{O}_{A}(e))^G \simeq \Phi^G(\Sigma_{BG}^{\infty} (\mathbbl{1}^{\Glo} \otimes f^*\mathcal{O}_{A}(e))),$$ 
see \cref{example: star-pullback fixed points}, \cref{example: global_ench_BG}, and \cref{proposition: basic prop of gfp}.
Hence, we obtain
$$f^*\mathcal{O}_{A}(e) \otimes \Phi^G(\beta \mathcal{L}) \simeq \Phi^G(\Sigma^{\infty}_{BG}(\mathbbl{1}^{\Glo} \otimes f^*\mathcal{O}_{A}(e)) \otimes \beta \mathcal{L} ). $$
Since $f^*\mathcal{O}_{A}(e) \in \QCoh(A(BG))$ is $\otimes$-invertible and $\beta$ is right lax $\LocSys^\gen_\ab(BG,A)$-linear, we obtain
$$ \Phi^G(\Sigma^{\infty}_{BG}(\mathbbl{1}^{\Glo} \otimes f^*\mathcal{O}_{A}(e)) \otimes \beta \mathcal{L} ) \simeq \Phi^G\beta((\mathbbl{1}^{\temp} \otimes f^*\mathcal{O}_{A}(e))\otimes \mathcal{L}).$$
By \cref{Thom_for_characters}, we have $\mathbbl{1}^{\temp} \otimes f^*\mathcal{O}_{A}(e) \simeq \mathbbl{1}^{\chi,\temp}_{BG}$, and so,
$$f^*(u)\otimes \Phi^G(\beta\mathcal{L}) \simeq \Phi^G\beta(a_\chi \otimes \mathcal{L}),$$
where $a_\chi\colon \mathbbl{1}^{\temp} \to \mathbbl{1}^{\chi,\temp}$ is the Euler class of the representation $\chi$ induced by the $G$-equivariant inclusion $\{0\} \inj \chi$.

However, by \cref{proposition: basic prop of gfp}, we have
$$\Phi^G(\beta\mathcal{L}) \simeq \colim_{V\in \Vect_{BG}}\Sigma^{-\dim V^G}(\mathbbl{1}^{V,\gen}_{BG}\otimes^{\gen} \beta \mathcal{L})^G,$$
where the transition maps are given by the Euler classes $a_V \colon \Sigma^{\dim V^G} \mathbbl{1}^\gen \to \mathbbl{1}^{V,\gen}$, $V\in \Vect_{BG}$. Since the genuine local systems $\mathbbl{1}^{V,\gen}_{BG}$ are $\otimes$-invertible, we can continue
\begin{align*}
\Phi^G(\beta\mathcal{L}) &\simeq \colim_{V\in \Vect_{BG}}\Sigma^{-\dim V^G}(\beta(\mathbbl{1}^{V,\temp}_{BG}\otimes^{\temp} \mathcal{L}))^G \\
&\simeq \colim_{V\in \Vect_{BG}}\Sigma^{-\dim V^G}(\mathbbl{1}^{V,\temp}_{BG}\otimes^{\temp} \mathcal{L})^G.
\end{align*}
Therefore, for every $G$-representation $V$ and every tempered local system $\mathcal{L}\in \LocSys^\temp(BG,A)$, the map $\Phi^G\beta(a_V \otimes \mathcal{L})$ induced by the Euler class $a_V$ is an equivalence. In particular, $f^*(u)\otimes \Phi^G(\beta\mathcal{L}) \simeq \Phi^G\beta(a_\chi \otimes \mathcal{L})$ is an equivalence as well, which finishes the proof.
\end{proof}

\begin{rem}
Suppose that $G$ is a finite abelian group. Then \cref{proposition: geometric fixed points of temp} implies that the functor 
$$\Phi^G\beta \colon \LocSys^{\temp}(BG,A) \tto \QCoh(A(BG))^{(\fcat{P}, G)\mdef\llocal}$$
coincides with the functor of geometric fixed points introduced in~\cite[Definition~15.11]{GLP26}.
\end{rem}

\begin{cor}\label{corollary: geometric fixed points and product}
Let $X=BG$ and $\mathcal{L} \in \LocSys^{\temp}(BG,A)$ be a tempered local system. Set $\beta=\beta_{BG}$. Then for any globally equivariant local system $\mathcal{M}_0 \in \LocSys^\Glo_\ab(X,A)$, the canonical morphism
$$\theta\colon \mathcal{M}^G_0 \otimes \mathcal{L}^G \xrightarrow{\theta_1} (L_{BG}\mathcal{M}_0 \otimes^{\Glo} \mathcal{L})^G \xrightarrow{\theta_2} (L_{BG}\mathcal{M}_0 \otimes^{\temp} \mathcal{L})^G \xrightarrow{\theta_3} \Phi^G\beta(L_{BG}\mathcal{M}_0 \otimes^{\temp} \mathcal{L})$$ 
exhibits $\Phi^G\beta(L_{BG}\mathcal{M}_0 \otimes^{\temp} \mathcal{L})$ as the $(\fcat{P}, G)$-localization of $\mathcal{M}^G_0 \otimes \mathcal{L}^G$.
\end{cor}

\begin{proof}
By \cref{proposition: geometric fixed points of temp}, $\Phi^G\beta(L_{BG}\mathcal{M}_0 \otimes^{\temp} \mathcal{L})$ is $(H,G)$-local for all proper closed subgroups $H\subset G$. Hence, it suffices to check that $\fib(\theta) \in \QCoh(A(BG))^{\tors}$. By \cref{proposition: geometric fixed points of temp}, the fiber of $\theta_3$ is in $\QCoh(A(BG))^{\tors}$. Finally, the fibers $\fib(\theta_1)$ and $\fib(\theta_2)$ are in $\QCoh(A(BG))^{\tors}$ by \cref{thm_CMon_str_on_LStemp} and \cref{thm_Ltemp_kernel}.
\end{proof}

\begin{prop}\label{proposition: gen-to-temp is smashing}
Let $X$ be an orbispace. Then the functor $\beta_X$ is strictly $\LocSys^\gen_\ab(X,A)$-linear. That is the natural morphism
$$\mathcal{M} \otimes^{\gen} \beta_X\mathcal{L} \tto \beta_X(L^{\gen}_X \mathcal{M} \otimes^{\temp} \mathcal{L}) $$
is an equivalence for all $\mathcal{M} \in \LocSys^\gen_\ab(X,A)$ and $\mathcal{L}\in \LocSys^{\temp}(X,A)$.

\begin{proof}
By \cref{construction: genuine local system} and \cref{corollary: gen-to-temp local is functor along sharp}, we can assume that $X=BG$ for some compact abelian Lie group $G$. Then, by \cref{corollary: gfp are jointly conservative}, it suffices to show that
$$\Phi^G(\mathcal{M} \otimes^{\gen} \beta\mathcal{L}) \tto \Phi^G\beta(L^{\gen}_{BG} \mathcal{M} \otimes^{\temp} \mathcal{L}) $$
is an equivalence for all $\mathcal{M} \in \LocSys^\gen_\ab(BG,A)$ and $\mathcal{L}\in \LocSys^{\temp}(BG,A)$. Here, $\beta=\beta_{BG}$. By \cref{cor_OmegaInftygen}, we can assume further that $\mathcal{M} \simeq \Sigma^{\infty}_{BG}\mathcal{M}_0$ for some $\mathcal{M}_0 \in \LocSys^\Glo_\ab(BG,A)$. Therefore, it is enough to show that
$$\mathcal{M}^G_0 \otimes \Phi^G\beta\mathcal{L} \tto \Phi^G\beta(L_{BG} \mathcal{M}_0 \otimes^{\temp} \mathcal{L}) $$
is an equivalence in $\QCoh(A(BG))$ for all $\mathcal{M}_0 \in \LocSys^\Glo_\ab(BG,A)$ and $\mathcal{L}\in \LocSys^{\temp}(BG,A)$. However, by \cref{corollary: geometric fixed points and product}, both sides are the $(\fcat{P}, G)$-localization of $\mathcal{M}_0^G \otimes \mathcal{L}^G$.
\end{proof}
\end{prop}

\begin{thm}\label{theorem: temp are modules in gen}
Let $X$ be an orbispace. Then the adjoint pair
$$L^{\gen}_X \colon \LocSys^\gen_\ab(X,A) \xymatrix{\ar@<0.5ex>[r] & \ar@<0.5ex>[l]} \LocSys^{\temp}(X,A) : \beta_{X}$$
induces a natural equivalence
$$\LocSys^{\temp}(X,A) \simeq \Mod_{\beta_{X}(\mathbbl 1)} \LocSys^\gen_\ab(X,A). $$
Moreover, the commutative algebra $\beta_{X}(\mathbbl 1)$ is idempotent. In particular, the functor $L^{\gen}_X$ is a smashing localization.
\end{thm}

\begin{proof}
By \cref{lemma: beta is cons and cont}, the right adjoint $\beta_X$ satisfies the conditions of the Barr--Beck--Lurie monadicity theorem~\cite[Theorem~4.7.3.5]{Lur_HA}. Therefore, the category $\LocSys^{\temp}(X,A)$ is monadic over $\LocSys^\gen_\ab(X,A)$. By virtue of \cref{proposition: gen-to-temp is smashing},  the monad is given by tensoring with the algebra $\beta_X(\mathbbl{1})$.

Finally, we will prove that the algebra $\beta_X(\mathbbl{1})$ is idempotent, i.e. the multiplication
$$\mu_X\colon \beta_X(\mathbbl {1}) \otimes \beta_X(\mathbbl{1}) \tto \beta_X(\mathbbl{1})$$
is an equivalence. By \cref{corollary: gen-to-temp local is functor along sharp}, we can assume that $X=BG$ for some compact abelian Lie group $G$. Set $\beta=\beta_{BG}$ and $\mu=\mu_{BG}$. Then, by \cref{corollary: gfp are jointly conservative}, it is enough to show that
$$\Phi^G\mu\colon \Phi^G(\beta\mathbbl{1}) \otimes \Phi^G(\beta\mathbbl{1})\tto \Phi^G(\beta\mathbbl{1}) \in \QCoh(A(BG))$$
is an equivalence. By \cref{proposition: geometric fixed points of temp}, $\Phi^G(\beta\mathbbl{1})$ is the $(\fcat{P}, G)$-localization of $\mathbbl{1} \in \QCoh(A(BG))$, which proves the assertion. 
\end{proof}

\begin{cor}\label{corollary: temp in genuine}
Let $X$ be an orbispace. Then the functor
$$\beta_X \colon \LocSys^{\temp}(X,A) \tto \LocSys^\gen_\ab(X,A)$$
is fully faithful and a genuine local system $\mathcal{L} \in \LocSys^\gen_\ab(X,A)$ belongs to the essential image of $\beta_X$ if and only if $\Phi^x\mathcal{L} \in \QCoh(A[\widehat{G}])$ is $(\fcat{P}, G)$-local for all faithful points $x\colon BG \to X$, where $G$ is abelian. \qed
\end{cor}

\begin{rem}\label{remark: smashing finite groups}
For finite abelian groups, \cref{theorem: temp are modules in gen} and \cref{corollary: temp in genuine} were proven in~\cite[Theorem~15.27 and Proposition~15.28]{GLP26}. We point out that the statement in~\cite{GLP26} is even more general in a certain way, since they allow the coefficient system to originate from an oriented divisible group and do not require this divisible group to be induced by any abelian group object. Of course, the price of this approach is that~\cite{GLP26} had to restrict themselves with finite groups.
\end{rem}

\begin{ex}\label{temp_LS_on_BCp}
Recall from \cref{example: gen gluing functor geom type} that there is an equivalence
$$\LocSys^\gen_\ab(BC_p,A) \simeq \rlaxlim\left(\LocSys(BC_p,\QCoh(S)) \xrightarrow{e_*(-)^{tC_p}} \QCoh(A[C_p])\right),$$
where $e\colon S \to A[C_p]$ is the unit section. By \cref{LStemp_TateInvs_are_local}, the functor $e_*(-)^{tC_p}$ factors through the full subcategory $\QCoh(A[C_p])^{(\fcat{P}, C_p)\mdef\llocal}$ of $(e,C_p)=(\fcat{P}, C_p)$-local objects, i.e.
$$e_*(-)^{tC_p}\colon \LocSys(BC_p,\QCoh(S)) \tto \QCoh(A[C_p])^{(\fcat{P}, C_p)\mdef\llocal} \hookrightarrow \QCoh(A[C_p]).$$
Then, by \cref{corollary: temp in genuine}, we have
$$\LocSys^{\temp}(BC_p,A) \simeq \rlaxlim\left(\LocSys(BC_p,\QCoh(S)) \xrightarrow{e_*(-)^{tC_p}} \QCoh(A[C_p])^{(\fcat{P}, C_p)\mdef\llocal}\right).$$
Under this equivalence and the equivalence of \cref{example: cyclic group global}, the functor $$\Omega^\infty_{BC_p}\beta_{BC_p}\colon \LocSys^{\temp}(BC_p,A) \tto \LocSys^\Glo_\ab(BC_p,A)$$ sends an object 
$$\mathcal{L} = (\Phi^{C_p}\mathcal{L}\in \QCoh(A[C_p])^{(\fcat{P}, C_p)\mdef\llocal}, \mathcal{L}^e \in \LocSys(BC_p,\QCoh(S)),\Phi^{C_p}\mathcal{L} \to (e_*\mathcal{L}^e)^{tC_p})$$
of the right lax limit to the triple $$(\Phi^{C_p}\mathcal{L} \times_{(e_*\mathcal{L}^e)^{tC_p}}(e_*\mathcal{L}^e)^{hC_p}, \mathcal{L}^e, \Phi^{C_p}\mathcal{L} \times_{(e_*\mathcal{L}^e)^{tC_p}}(e_*\mathcal{L}^e)^{hC_p} \to (e_*\mathcal{L}^e)^{hC_p}).$$
Note that this triple satisfies the condition~\eqref{def_LS_temp} by \cref{LStemp_TateInvs_are_local} and the assumption on $\Phi^{C_p}\mathcal{L}$. Finally, by arguing as in \cite[Remark~5.0.1]{Lur_Ell3}, one can check independently that the functor $\Omega^\infty_{BC_p}\beta_{BC_p}$ is fully faithful and that its essential image coincides with the full subcategory of tempered local systems defined in \cref{def_LS_temp}.
\end{ex}

\begin{rem}\label{remark: omegainfty is not ff}
Let $X$ be an orbispace and consider the following commutative diagram
\[\xymatrix{
\LocSys^\temp(X, A) \ar[d]^-{\beta_X} \ar[r]^-{i_X} & \LocSys^\Glo_\ab(X, A) \\
\LocSys^\gen_\ab(X, A). \ar[ru]_-{\Omega^{\infty}_X}
}\]
Note that both functors $\beta_X$ and $\iota_X$ are fully faithful embeddings, but the functor $\Omega^{\infty}_X$ is \emph{not} fully faithful, except $X\in \Type$ is a constant global space. Indeed, the functor $\Omega^{\infty}_X$ is fully faithful if and only if the counit map $$\eta\colon \Sigma^{\infty}_X\Omega^{\infty}_X \tto \Id$$
is an equivalence. In particular, if $\Omega^{\infty}_X$ is fully faithful, then the natural map
$\mathcal{L}^x \to \Phi^x\mathcal{L}$
is an equivalence for all genuine local systems $\mathcal{L} \in \LocSys^\gen_\ab(X, A)$ and all faithful points $x\colon BG \to X$ with abelian $G$. Therefore, by taking $\mathcal{L}=\beta_X(\mathbbl{1}^\temp)$ and by \cref{proposition: geometric fixed points of temp}, the $(\fcat{P}, G)$-localization
$$\mathcal{O}_{A(BG)} \tto L_{(\fcat{P},G)}\mathcal{O}_{A(BG)} \in \QCoh(A(BG))$$
is an equivalence for all isotropy groups of $X$. But this is the case only if $G$ is a trivial group. 
\end{rem}

\begin{ex}\label{example: temp is not local}
We note that the statement of~\cref{theorem: temp are modules in gen} is false if $X$ is \emph{not} an orbispace. Indeed, the global space $X$ from \cref{example: omegainfty is not conservative} provides a counterexample. Recall that $X=\ast \coprod_{BG} \ast \in \Type^{\Glo}$ is the suspension of representable global space $BG$ where $G=C_p$ is a cyclic group of prime order $p$. Let $R=\mathbb{Q}[\beta^{\pm 1}] \in \CAlg(\Sp)$ and let $A=\mathbb{G}_{a,R}$ be the oriented (strict) additive group (see~\cref{example: additive group} and \cref{example: additive group is complex oriented}).
\iffalse
Then we have the following identifications
$$A(e)\simeq\Spec(\mathbb{Q}[\beta^{\pm 1}]), \;\; A(BC_p) \simeq \Spec(\mathbb{Q}[\beta^{\pm 1}]), \;\; A(B\mathbb{T}) \simeq \Spec(\mathbb{Q}[\beta^{\pm 1},t]) $$
where $\beta \in \pi_2(R)$ is the Bott class and $t \in \pi_0(R)$ is a generator of the augmentation ideal. \todo{digression on orientations} The $\mathbb{T}$-equivariant map
$$\mathbb{Q}[\beta^{\pm 1},t] \tto \mathbb{Q}[\beta^{\pm 1}] $$
is induced by the preorientation morphism
$$\theta\colon \mathbb{Q}[\beta^{\pm 1},t] \tto (\mathbb{Q}[\beta^{\pm 1}])^{h\mathbb{T}} \simeq \mathbb{Q}[\beta^{\pm 1}][[u]] $$
which maps $t$ to $\beta u$, $|u|=-2$. In particular, $\theta$ is an equivalence after $t$-completion and $\mathbb{G}_{a,R}$ is an oriented strict group.
\fi

Recall from \cref{example: omegainfty is not conservative} that
$$\LocSys^\Glo_\ab(X,A)\simeq \LocSys^\Glo_\ab(*,A) \simeq \Mod_{R}$$ and
$$\LocSys^\gen_\ab(BG,A)\simeq \Mod_{R}\times \Fun(BG,\Mod_{R})\simeq \Fun(BG_+,\Mod_{R}).$$
So, $\LocSys^\gen_\ab(X,A)\simeq \Fun(\Sigma(BG_+),\Mod_R)$ and $\LocSys^\temp(X, A)$ being a full subcategory of
$$\LocSys^\Glo_\ab(X, A) \simeq \LocSys^\Glo_\ab(*, A) \simeq \Mod_R$$
must be equivalent to $\LocSys^\temp(*, A) \simeq \Mod_R$. Under these identifications, the functor $L_X^{\gen}$ is equivalent to the pullback
$$i^*\colon \Fun(\Sigma(BG_+),\Mod_R) \tto \Mod_R, $$
where $i\colon \ast \to \Sigma(BG_+)$ is a point of the connected space $\Sigma(BG_+)$. Therefore, the functor $\beta_X$ is equivalent to the pushforward $i_*$. Finally, the pushforward $i_*$ is neither fully faithful nor continuous.
\end{ex}

\begin{cor}\label{corollary: temp-sharp is gen-sharp global}
Let $f\colon X\to Y$ be a faithful morphism of global spaces. Then the commutative square
\[\xymatrix{
\LocSys^\gen_\ab(X, A) \ar[r]^-{L^{\gen}_{X}} \ar[d]^{f^{\gen}_\#} & \LocSys^\temp(X, A) \ar[d]^{f^{\temp}_\#}\\
\LocSys^\gen_\ab(Y, A) \ar[r]^-{L^{\gen}_{Y}}  & \LocSys^{\temp}(Y, A) 
}\]
is horizontally right adjointable. That is a natural transformation
$$f^{\gen}_{\#} \circ \beta_X \tto \beta_Y \circ f^{\temp}_{\#} $$
is an equivalence.

\begin{proof}
First assume that $X$ and $Y$ are orbispaces. Since $L^\gen_Y$ is a localization, it suffices to show that the genuine $\#$-pushforward $f^{\gen}_\#$ preserves tempered local systems. Let $\mathcal{L} \in \LocSys^{\temp}(X,A)$ be a tempered local system, we will show that $f^{\gen}_\#(\mathcal{L})$ is also tempered. By \cref{theorem: temp are modules in gen}, it suffices to show that
$$f^{\gen}_\#(\mathcal{L}) \tto f^{\gen}_\#(\mathcal{L})\otimes \beta_Y(\mathbbl{1}) $$
is an equivalence. However, by the projection formula, we have $$f^{\gen}_\#(\mathcal{L})\otimes \beta_Y(\mathbbl{1}) \simeq f^{\gen}_\#(\mathcal{L}\otimes f^{\gen,*}\beta_Y(\mathbbl{1}))\simeq f^{\gen}_\#(\mathcal{L}\otimes \beta_X(\mathbbl{1})) \simeq f^{\gen}_\#(\mathcal{L}),$$
which finishes the proof.

For general $X$ and $Y$, the assertion follows by the previous part and the descent properties as in \cref{construction: genuine local system} and \cref{proposition: temp descent} because any limit of horizontally right adjointable squares is right adjointable, see \cite[Corollary~4.7.4.18]{Lur_HA}.
\end{proof}
\end{cor}

We finish this subsection with the following application for tempered local systems similar to \cref{corollary: genuine BG rigid over ABG}.
\begin{cor}\label{cor_rigidity_of_LStemp}
Let $G$ be a compact abelian Lie group. Then the category $\LocSys^{\temp}(BG,A) \in \CAlg_{\QCoh(A(BG))}(\Prs^\LL)$ is a rigid $\QCoh(A(BG))$-atomically generated $\QCoh(A(BG))$-algebra. More precisely, the objects $p^\temp_{\#}\mathbbl 1_{BH}$, where $p\colon BH \to BG$ ranges over the set of faithful morphisms, form a set of $\QCoh(A(BG))$-atomic (and dualizable) generators of a $\QCoh(A(BG))$-module $\LocSys^{\temp}(BG,A)$.
\end{cor}

\begin{proof}
The first assertion follows directly from \cref{corollary: genuine BG rigid over ABG} and \cref{theorem: temp are modules in gen}. The assertion about the set of generators follows by \cref{corollary: temp-sharp is gen-sharp global}.
\end{proof}

\begin{rem}\label{remark: tempered local system comp gen finite}
For a finite abelian group $G$, the assertion of \cref{cor_rigidity_of_LStemp} can be deduced from the main results of~\cite{Lur_Ell3}. Indeed, the objects $p^{\temp}_{\#} \mathbbl 1_{BH}$, $H\subset G$ are compact generators by~\cite[Corollary~5.3.2]{Lur_Ell3}. Moreover, they are dualizable by~\cite[Proposition~7.8.8]{Lur_Ell3}. We point out that~\cite{Lur_Ell3} does not rely on the classical Atiyah duality, whereas we used it in the proof of \cref{corollary: genuine BG rigid over ABG} and \cref{cor_rigidity_of_LStemp}. 
\end{rem}

\begin{rem}\label{remark: ambidexterity for elliptic cohomology}
Let $G$ be a compact abelian Lie group, let $p\colon BG \to *$ be the natural full map, and let $A$ be an oriented abelian group spectral stack over $S$. By~\cite[Section~7]{Lur_Ell3}, if $G$ is finite, then the pullback functor $$p^{\temp,*}\colon \QCoh(S)\simeq \LocSys^{\temp}(*,A)\tto \LocSys^{\temp}(BG,A)$$
admits a left adjoint $p_{\#}^{\temp}$. Moreover, by loc.cit., there is a natural equivalence $p_{\#}^{\temp}\simeq p_{*}^{\temp}$ between the left and right adjoints, cf. \cref{remark: genuine BG norm map}. We point out that this observation is wrong if $G$ is positive-dimensional without additional assumptions on $A$. Indeed, suppose that $G=U(1)$. Then, by \cref{LStemp_on_tori}, we can identify $\LocSys^{\temp}(BG,A)$ with $\QCoh(A)$. Moreover, the tempered pullback functor $p^{\temp,*}$ identifies with the $*$-pullback functor
$$p^*\colon \QCoh(S) \tto \QCoh(A) $$
along the projection $p\colon A \to S$.

Finally, if $A=\mathbb{G}_{a, \mathbb Q[\beta^{\pm 1}]}$ or $A=\mathbb{G}_{m,KU}$ (see \cref{example: additive group is complex oriented} and \cref{example: multiplicative group is complex oriented}), then the pullback $p^*$ does not commute with infinite products, and so $p^*$ does not even admit a left adjoint, see \cref{example: Gm does not have left adjoint}. However, if $A$ is an \emph{oriented spectral elliptic curve} (see \cref{example: oriented elliptic curve}), then one can check that $p^*$ \emph{admits} a left adjoint 
$$p_{\#}\colon \QCoh(A) \tto \QCoh(S).$$
Moreover, by the virtue of Grothendieck--Serre duality, there is a natural equivalence $p_\# \simeq \Sigma p_*$ between the left and (shifted) right adjoints. In other words, we observed the equivalence $p^{\temp}_\# \simeq \Sigma p^{\temp}_*$ between tempered \emph{elliptic} homology and cohomology for $BU(1)$, which is not visible for ordinary or $K$-theoretic (co)homology. We are going to return to this ambidexterity phenomena in the future work.
\end{rem}

\subsection{Dexterity}\label{section: dexterity}
Let $G$ be a compact Lie group and let $f\colon X \to Y$ be a map of orbispaces faithful over $BG$. Then by the results of \Cref{section: LSgen on G-spaces} we know that the comparison morphism
\begin{align*}
\theta_f\colon \LocSys^\Glo(X, \fcat A) \otimes_{\LocSys^\Glo(Y, \fcat A)} \LocSys^\gen(Y, \fcat A) &\tto \LocSys^\gen(X, \fcat A), \\
\mathcal{L} \boxtimes \mathcal{M} &\xymatrix{\ar@{|->}[r] &}  \Sigma^\infty_X\mathcal{L} \otimes f^{\gen,*}\mathcal{M}
\end{align*}
is an equivalence (since any vector bundle on $X$ is locally a summand of a vector bundle pulled back from $Y$ in this case). Unfortunately, we do not know whether this assertion is true for a general faithful morphism of global spaces. The goal of this section is to prove \Cref{theorem: dextery faithful tempered}, which asserts that an analogous comparison map
\begin{align*}
\alpha_f\colon \LocSys^\Glo(X, \fcat A) \otimes_{\LocSys^\Glo(Y, \fcat A)} \LocSys^\temp(Y, \fcat A) &\tto \LocSys^\temp(X, \fcat A),\\
\mathcal{L} \boxtimes \mathcal{M} &\xymatrix{\ar@{|->}[r] &}  L_X\mathcal{L} \otimes f^{\gen,*}\mathcal{M}
\end{align*}
is an \emph{equivalence} if both $X$ and $Y$ are orbispaces. Note that the assertion for $\alpha_f$ would follow if $\theta_f$ is an equivalence by using \cref{theorem: temp are modules in gen}.

One important advantage of dealing with tempered local systems instead of genuine ones is provided by the following observation.

\begin{lem}\label{lemma: relative tensor product is a full subcategory}
Let $f\colon X\to Y$ be a morphism of global spaces. Then the canonical functor
$$\alpha'_f=\Id\otimes L_Y \colon \LocSys^\Glo_\ab(X,A) \tto \LocSys^\Glo_\ab(X,A) \otimes_{\LocSys^\Glo_\ab(Y,A)} \LocSys^{\temp}(Y,A) $$
admits a fully faithful right adjoint $\beta'_f$, i.e. $\alpha_f'$ is a Bousfield localization. % Moreover, $\mathcal{L} \in \LocSys^\Glo_\ab(X,A)$ belongs to the essential image of $i_f$ if and only if 
%$$\Hom(\mathcal{M}_0\otimes f^{\Glo,*}\mathcal{M}_1,\mathcal{L}) \simeq 0$$
%for all $\mathcal{M}_0 \in \LocSys^\Glo_\ab(X,A)$ and $\mathcal{M}_1 \in \LocSys^{\nnull}(Y,A)$.
\end{lem}

\begin{proof}
\iffalse
Note that the functor $\alpha'_f$ is a geometric realization of the functors
$$\alpha'_f\simeq \colim_{n\in \Delta^{\op}}\left(\LocSys^\Glo_\ab(X,A) \otimes \LocSys^\Glo_\ab(Y,A)^{\otimes n} \otimes L_Y\right).$$
By \cref{LStemp_admits_left_adj_global} and \cref{LStemp_presentable}, the functor $L_Y$ is strongly continuous. Therefore, the right adjoint $\beta'_f$ is the totalization
$$\beta'_f\simeq \lim_{n\in \Delta}\left(\LocSys^\Glo_\ab(X,A) \otimes \LocSys^\Glo_\ab(Y,A)^{\otimes n} \otimes i_Y\right). $$
Since a limit of fully faithful functors is fully faithful, the right adjoint $i_f$ is fully faithful.
\fi

By \cref{LStemp_admits_left_adj_global} and \cref{thm_CMon_str_on_LStemp}, we have the following pushout square
$$ 
\xymatrix{
\LocSys^\nnull(Y, A) \ar[r]\ar[d] & \LocSys^\Glo_\ab(Y, A) \ar[d]^{L_Y} \\
0 \ar[r] & \LocSys^\temp(Y, A)
}
$$
in $\Mod_{\LocSys^\Glo_\ab(Y, A)}(\Prs^\LL)$. By tensoring with the $\LocSys^\Glo_\ab(Y, A)$-module $\LocSys^\Glo_\ab(X, A)$, we obtain the pushout square
$$ 
\xymatrix{
\LocSys^\Glo_\ab(X,A) \otimes_{\LocSys^\Glo_\ab(Y,A)}\LocSys^\nnull(Y, A) \ar[r]\ar[d] & \LocSys^\Glo_\ab(X,A) \otimes_{\LocSys^\Glo_\ab(Y,A)}\LocSys^\Glo_\ab(Y, A) \ar[d]^{\Id\otimes L_Y} \\
0 \ar[r] & \LocSys^\Glo_\ab(X,A) \otimes_{\LocSys^\Glo_\ab(Y,A)}\LocSys^\temp(Y, A)
}
$$
in $\Prs^\LL$. Since a pushout in $\Prs^\LL$ of a Bousfield localization is a Bousfield localization, the assertion follows. %This implies that the essential image of $i_f$ is the right orthogonal to the essential image of the functor
%\begin{align*}
%\LocSys^\Glo_\ab(X,A) \otimes_{\LocSys^\Glo_\ab(Y,A)} \LocSys^{\nnull}(Y,A) &\tto \LocSys^\Glo_\ab(X,A) \\
%(\mathcal{M}_0, \mathcal{M}_1) &\mapsto \mathcal{M}_0 \otimes f^{\Glo,*}\mathcal{M}_1.
%\end{align*}
%This implies the last assertion.
\end{proof}

\begin{rem}\label{remark: warning about kernel}
Note that the kernel of $\alpha'_f$ is the essential image of $\LocSys^\Glo_\ab(X,A) \otimes_{\LocSys^\Glo_\ab(Y,A)}\LocSys^\nnull(Y, A)$ in $\LocSys^\Glo_\ab(X,A)$. However, the functor 
$$
\LocSys^\Glo_\ab(X,A) \otimes_{\LocSys^\Glo_\ab(Y,A)}\LocSys^\nnull(Y, A) \tto \LocSys^\Glo_\ab(X,A)
$$
is usually \emph{not} fully faithful.
\end{rem}

We will prove \cref{theorem: dextery faithful tempered} by reducing to a special case. Namely, we note that if $f=g\circ h$ is a composite such that $\alpha_g$ is an equivalence, then $\alpha_f$ is an equivalence if and only if  $\alpha_h$ is. We show first that $\alpha_f$ is an equivalence if $f\colon BG \to \mathcal{N}$ is the unique faithful morphism, where $G$ is a compact abelian Lie group and $\mathcal{N}$ is the normal subgroup classifier from \cref{example: normal subgroup classifier}. Note that the morphism $f$ factors as follows
$$f\colon BG \xrightarrow{g} \mathcal{N}_{\leq G} \xrightarrow{h} \mathcal{N},$$
where $\mathcal{N}_{\leq G}$ is the $(\leq \!\!\!G)$-normal subgroup classifier, see \cref{definiton: classfier subgroups_family}. We will deal with morphisms like $h$ first.

\begin{prop}\label{proposition: dexterity mono}
Let $\pi\colon Y \to X$ be a faithful monomorphism of orbispaces. Then the canonical functor
$$\theta_{\pi} \colon \LocSys^\Glo_\ab(Y, \fcat A) \otimes_{\LocSys^\Glo_\ab(X, \fcat A)} \LocSys^\gen_\ab(X, \fcat A) \tto \LocSys^\gen_\ab(Y, \fcat A)$$
is an equivalence.

\begin{proof}
By \cref{prop_LSGlo_sharp_pushforward}, the left adjoint functor
$$\pi^{\Glo}_{\#}\colon \LocSys^\Glo_\ab(Y, \fcat A) \tto \LocSys^\Glo_\ab(X, \fcat A) $$
is strictly $\LocSys^\Glo_\ab(X, \fcat A)$-linear. Therefore, $\pi^{\Glo}_{\#}$ induces the functor
$$\pi_{\#}^{\Glo}\otimes \Id \colon \LocSys^\Glo_\ab(Y, \fcat A) \otimes_{\LocSys^\Glo_\ab(X, \fcat A)} \LocSys^\gen_\ab(X, \fcat A) \tto \LocSys^\gen_\ab(X, \fcat A). $$
By the second part of \cref{proposition: genuine sharp}, we obtain the following commutative diagram
\begin{equation}\label{equation: tensor product genuine_eq1}
\xymatrix{
\LocSys^\Glo_\ab(Y, \fcat A) \otimes_{\LocSys^\Glo_\ab(X, \fcat A)} \LocSys^\gen_\ab(X, \fcat A) \ar[r]^-{\pi^{\Glo}_{\#} \otimes \Id}\ar[d]^-{\theta_{\pi}} & \LocSys^\gen_\ab(X, \fcat A) \\ 
\LocSys^\gen_\ab(Y, \fcat A) \ar[ru]_-{\pi^{\gen}_{\#}}.&
}
\end{equation}
By assumption, $Y\simeq Y\times_X Y$. So, we have $\pi^{\Glo,*}\circ \pi^{\Glo}_{\#}\simeq \Id$ and $\pi^{\gen,*}\circ\pi^{\gen}_{\#}\simeq \Id$ by \cref{lemma: base change orbispaces} and by the third part of \cref{proposition: genuine sharp} respectively. Therefore, $\pi^{\Glo}_{\#}$ and $\pi^{\gen}_{\#}$ are fully faithful embeddings, and so is $\pi^{\Glo}_{\#}\otimes \Id$. Finally, the commutative diagram~\eqref{equation: tensor product genuine_eq1} implies that the functor $\theta_{\pi}$ is fully faithful.

It is left to show that $\theta_{\pi}$ is essentially surjective. By \cref{cor_OmegaInftygen}, the category $\LocSys^\gen_\ab(Y,\fcat{A})$ is generated under colimits by suspension objects $\Sigma^{\infty} \mathcal{L}$, $\mathcal{L} \in \LocSys^\Glo_\ab(Y,\fcat{A})$ which are in the essential image of $\theta_{\pi}$. Since $\theta_{\pi}$ is fully faithful and commutes with colimits, the assertion follows.
\end{proof}
\end{prop}

\begin{cor}\label{lemma: family_classfier proj_and_cont}
Let $\fcat{T}$ be a family of compact Lie groups and let $X$ be an orbispace. Let
$$\pi\colon X_{\fcat{T}}=X\times_{\mathcal{N}}\mathcal{N}_{\fcat{T}} \tto X$$
be the canonical projection, see \cref{definiton: classfier subgroups_family}. Then the canonical functor
$$\theta_{\pi} \colon \LocSys^\Glo_\ab(X_{\fcat{T}}, \fcat A) \otimes_{\LocSys^\Glo_\ab(X, \fcat A)} \LocSys^\gen_\ab(X, \fcat A) \tto \LocSys^\gen_\ab(X_{\fcat{T}}, \fcat A)$$
is an equivalence.
\end{cor}

\begin{proof}
By \cref{remark: family_classfier is left kan}, $X_{\fcat{T}}\simeq X_{\fcat{T}}\times_X X_{\fcat{T}}$. So, the result follows from~\cref{proposition: dexterity mono}.
\end{proof}

Now let $f\colon BG \to Y=\mathcal{N}_{\leq G}$ be the unique faithful morphism. To show that $\alpha_f$ is an equivalence, we will use the isotropy separation, i.e. the recollements from \cref{proposition: global sections semi-orthogonal} and \cref{corollary: isotropy separation for gen}. Indeed, the category $\LocSys^{\Glo}_{\ab}(BG,\fcat{A})$ splits into two $\LocSys^{\Glo}_{\ab}(Y,\fcat{A})$-modules $\fcat{A}(BG)$ and $\LocSys^{\Glo}_{\ab}((BG)_\fcat{P},\fcat{A})$. We will deal with each of these modules separately.

\begin{lem}\label{lemma: isotropy separation geometric points}
Let $G$ be a compact abelian Lie group and let $f\colon BG \to Y=\mathcal{N}_{\leq G}$ be the unique faithful morphism. We consider $\fcat{A}(BG)\in\Prs^\LL$ as the $\LocSys^\Glo_\ab(Y, \fcat A)$-module via the symmetric monoidal composite
$$\LocSys^\Glo_\ab(Y, \fcat A) \xrightarrow{f^{\Glo,*}} \LocSys^\Glo_\ab(BG, \fcat A) \xrightarrow{(-)^{BG}}	\fcat{A}(BG).$$
Then the functor
\begin{align*}
\fcat A(BG) \otimes_{\LocSys^\Glo_\ab(Y, \fcat A)} \LocSys^\gen_\ab(Y, \fcat A) &\tto \fcat{A}(BG), \\
\mathcal{L} \boxtimes \mathcal{M} &\xymatrix{\ar@{|->}[r] &} \mathcal{L} \otimes \Phi^f\mathcal{M}
\end{align*}
is an equivalence.
\end{lem}

\begin{proof}
Let $\fcat{A}(BG)_{Y} \in \CAlg(\Prs^\LL)$ denote the following limit category
$$\fcat{A}(BG)_{Y} = \rlaxlim_{\substack{T \in (\Orb_{/^\rep Y})^\op \\ T\simeq BG}} \fcat A(T) \simeq \lim_{\substack{T \in (\Orb_{/^\rep Y})^\op \\ T\simeq BG}} \fcat A(T) \in \CAlg(\Prs^\LL). $$
By construction, the evaluation functor $(-)^{BG}\colon \LocSys^\Glo_\ab(Y, \fcat A) \to \fcat{A}(BG)$ factors as follows
$$(-)^{BG}\colon \LocSys^\Glo_\ab(Y, \fcat A) \xrightarrow{\phi^*} \fcat{A}(BG)_Y \xrightarrow{\psi^*} \fcat{A}(BG),$$
where $\phi^*$ is the restriction of the right lax limit indexed by \emph{all} orbits over $Y$ to the right lax limit indexed by the full subcategory of orbits over $Y$ equivalent to $BG$, and $\psi^*$ is the evaluation at the distinguished object $T=BG$ over $Y$. Note that both functors $\phi^*$ and $\psi^*$ are symmetric monoidal. By \cref{right_adjoint_llax_limits}, the functor $\phi^*$ admits a fully faithful right adjoint $\phi_*$ such that
$$(\phi_*(X))^{BH} \simeq 
\begin{cases}
\psi^*(X), & \mbox{if $H\simeq G$,}\\
0, & \mbox{if $H \neq G$.}
\end{cases}
$$
In particular, $\phi_*$ preserves arbitrary colimits and is strictly $\LocSys^\Glo_\ab(Y,\fcat{A})$-linear. Hence, $\phi^*$ is a smashing localization and
$$\fcat{A}(BG)_Y \simeq \Mod_{\phi_*(\mathbbl{1})}(\LocSys^\Glo_\ab(Y,\fcat{A})). $$

Therefore,
\begin{align*}
\fcat{A}(BG)_Y \otimes_{\LocSys^\Glo_\ab(Y, \fcat A)} \LocSys^\gen_\ab(Y, \fcat A) &\simeq \Mod_{\Sigma^{\infty}_Y\phi_*(\mathbbl{1})}(\LocSys^\gen_\ab(Y,\fcat{A}))\\
&\simeq \lim_{g\colon BH\to Y\in \Orb_{/^{\rep} Y}}\left(\Mod_{g^{\gen,*}\Sigma^{\infty}_Y\phi_*(\mathbbl{1})}(\LocSys^\gen_\ab(BH,\fcat{A}))\right).
\end{align*}
By construction and \cref{corollary: ep tilde}, we have
$$g^{\gen,*}\Sigma^{\infty}_Y(\phi_*(\mathbbl{1})) \simeq \Sigma^{\infty}_{BH}g^{\Glo,*}(\phi_*(\mathbbl{1}))
\simeq  \begin{cases}
\Sigma^{\infty}_{BG} Q^{BG}(\mathbbl{1}) \simeq \Xi^G_{\fcat{A}}(\mathbbl{1}), & \mbox{if $H\simeq G$,}\\
0, & \mbox{if $H \neq G$.}
\end{cases}
$$
This implies that
$$
\fcat{A}(BG)_Y \otimes_{\LocSys^\Glo_\ab(Y, \fcat A)}\LocSys^\gen_\ab(Y, \fcat A)
\simeq \lim_{\substack{g\colon T\to Y\in \Orb_{/^{\rep} Y} \\ T\simeq BG}}\left(\Mod_{\Xi^T_{\fcat{A}}(\mathbbl{1})}(\LocSys^\gen_\ab(T,\fcat{A}))\right).
$$
By \cref{corollary: gfp is smashing localization}, we obtain
$$\fcat{A}(BG)_Y \otimes_{\LocSys^\Glo_\ab(Y, \fcat A)}\LocSys^\gen_\ab(Y, \fcat A)
\simeq \lim_{\substack{g\colon T\to Y\in \Orb_{/^{\rep} Y} \\ T\simeq BG}}\fcat{A}(T) \simeq \fcat{A}(BG)_Y.
$$
Finally, we have
\begin{align*}
\fcat A(BG) \otimes_{\LocSys^\Glo_\ab(Y, \fcat A)} \LocSys^\gen_\ab(Y, \fcat A) &\simeq \fcat{A}(BG)\otimes_{\fcat{A}(BG)_Y} \left( \fcat{A}(BG)_Y \otimes_{\LocSys^\Glo_\ab(Y, \fcat A)} \LocSys^\gen_\ab(Y, \fcat A)\right) \\
&\simeq \fcat{A}(BG) \otimes_{\fcat{A}(BG)_Y} \fcat{A}(BG)_Y \simeq \fcat{A}(BG).\qedhere
\end{align*}
\end{proof}

\begin{cor}\label{corollary: isotropy separation geometric points tempered}
Let $G$ be a compact abelian Lie group and let $f\colon BG \to Y=\mathcal{N}_{\leq G}$ be the unique faithful morphism. Then the functor
\begin{align*}
\QCoh(A(BG)) \otimes_{\LocSys^\Glo_\ab(Y, A)} \LocSys^\temp(Y, A) &\tto \QCoh(A(BG))^{(\fcat{P}, G)\mdef\llocal}, \\
\mathcal{L} \boxtimes \mathcal{M} &\xymatrix{\ar@{|->}[r] &} \mathcal{L} \otimes \Phi^f\mathcal{M}
\end{align*}
is an equivalence.

\begin{proof}
Since $\Phi^{f}(\mathbbl{1})\simeq L_{(\fcat{P},G)}\mathcal{O}_{A(BG)}$, the assertion follows by \cref{lemma: isotropy separation geometric points} and \cref{theorem: temp are modules in gen}.
\end{proof}
\end{cor}

\begin{lem}\label{lemma: dexterity isotropy separation tempered}
Let $G$ be a compact abelian Lie group and let $f\colon BG \to Y=\mathcal{N}_{\leq G}$ be the unique faithful morphism. Then the commutative square
\[\xymatrix{
\LocSys^\Glo_\ab(BG,A) \otimes_{\LocSys^\Glo_\ab(Y, A)} \LocSys^\temp(Y, A) \ar[d]^-{(-)^G\otimes \Id} \ar[r]^-{\alpha_f}& \LocSys^\temp(BG,A) \ar[d]^-{\Phi^G} \\
\QCoh(A(BG)) \otimes_{\LocSys^\Glo_\ab(Y, A)} \LocSys^\temp(Y, A) \ar[r]^-{\sim} & \QCoh(A(BG))^{(\fcat{P}, G)\mdef\llocal}
}\]
is horizontally right adjointable.
\end{lem}

\begin{proof}
%\todo{Change notation}
By \cref{corollary: isotropy separation geometric points tempered}, the lower horizontal arrow is an equivalence. Let $\widetilde{\Phi}^G$ denote the functor
$$
\widetilde{\Phi}^G \simeq (-)^G\otimes \Phi^f \colon \LocSys^\Glo_\ab(BG,A) \otimes_{\LocSys^\Glo_\ab(Y, A)} \LocSys^\temp(Y, A) \tto \QCoh(A(BG))^{(\fcat{P}, G)\mdef\llocal}
$$
induced by $(-)^G\otimes \Phi^f$. In other words, $\widetilde{\Phi}^G$ is the composite of the left vertical arrow and the lower horizontal equivalence. Let $\beta_f$ denote the right adjoint to $\alpha_f$. We will show that the natural morphism
$$\theta\colon \widetilde{\Phi}^G\beta_f \tto \Phi^G $$
is an equivalence.

Consider the commutative diagram
\[\xymatrix@C=5em{
\LocSys^\Glo_\ab(BG,A) \ar[r]^-{\alpha_f^\prime=\Id \otimes L_Y} \ar[d]^-{(-)^G} &\LocSys^\Glo_\ab(BG,A) \otimes_{\LocSys^\Glo_\ab(Y, A)} \LocSys^\temp(Y, A) \ar[d]^-{\widetilde{\Phi}^G} \\
\QCoh(A(BG)) \ar[r]^-{L_{(\fcat{P},G)}} &
\QCoh(A(BG))^{(\fcat{P}, G)\mdef\llocal}.
}
\]
Note that $\alpha_f\circ \alpha_f' \simeq L_{BG}$ is the symmetric monoidal localization, see \cref{thm_CMon_str_on_LStemp}. Let $\beta'_f$ be the right adjoint to $\alpha'_f$ and let $\iota_{(\fcat{P},G)}$ denote the right adjoint to the localization $L_{(\fcat{P},G)}$. Note that the fiber $\fib(\theta_1)$ of the natural transformation
$$\theta_1\colon (\beta'_f(-))^G \tto \iota_{(\fcat{P},G)} \widetilde{\Phi}^G$$
is in $\QCoh(A(BG))^{\tors}$. Indeed, it suffices to show that $L_{(\fcat{P},G)}\theta_1$ is an equivalence. However, by the commutativity of the diagram, we have $L_{(\fcat{P},G)}\theta_1 \simeq \widetilde{\Phi}^G (\eta)$, where
$$\eta\colon \alpha'_f \beta'_f \tto \Id $$
is the counit of the adjoint pair $\alpha'_f \dashv \beta'_f$. By \cref{lemma: relative tensor product is a full subcategory}, $\beta'_f$ is fully faithful. In particular, $\eta$ is an equivalence and so is $L_{(\fcat{P},G)}\theta_1$.

We also note that the fully faithful embedding 
$$\iota_{BG}\colon \LocSys^{\temp}(BG,A) \tto \LocSys^\Glo_\ab(BG,A) $$
factors as $\iota_{BG}\simeq \beta'_f \circ \beta_f$. Therefore, the natural transformation
$$\theta_2\colon (\iota_{BG}(-))^G \tto \iota_{(\fcat{P},G)}\Phi^G $$
factors as follows
$$\theta_2 \colon (\iota_{BG}(-))^G  \simeq (\beta'_f\circ \beta_f(-))^G \xrightarrow{\theta_1} \iota_{(\fcat{P},G)}\widetilde{\Phi}^G\beta_f \xrightarrow{\iota_{(\fcat{P},G)}\theta} \iota_{(\fcat{P},G)}\Phi^G. $$
By \cref{lemma: fiber of cat-to-geom is tors}, the fiber $\fib(\theta_2)$ is in $\QCoh(A(BG))^\tors$. Therefore, the fiber $\fib(\iota_{(\fcat{P},G)}\theta)$ is in $\QCoh(A(BG))^{\tors}$ as well. Finally, since $\iota_{(\fcat{P},G)}$ is a fully faithful embedding of $(\fcat{P},G)$-local objects, the fiber $\fib(\iota_{(\fcat{P},G)}\theta)$ is trivial and $\theta$ is an equivalence.
\end{proof}

\begin{rem}\label{remark: failure}
At the moment of writing, we do not know whether the commutative square
\[\xymatrix{
\LocSys^\Glo_\ab(BG,\fcat{A}) \otimes_{\LocSys^\Glo_\ab(Y, \fcat{A})} \LocSys^\gen(Y, \fcat{A}) \ar[d]^-{(-)^G\otimes \Id} \ar[r]^-{\theta_f}& \LocSys^\gen(BG,\fcat{A}) \ar[d]^-{\Phi^G} \\
\fcat{A}(BG) \otimes_{\LocSys^\Glo_\ab(Y, \fcat{A})} \LocSys^\gen(Y, \fcat{A}) \ar[r]^-{\sim} & \fcat{A}(BG)
}\]
is horizontally right adjointable. The lack of this property is the reason why we are not able to prove \cref{theorem: dextery faithful tempered} for genuine local systems.
\end{rem}

The next two lemmas are needed to deal with the second piece of the recollement.

\begin{lem}\label{lemma: dexterity isotropy separation tempered induction}
Let $B \in \Type^\Glo$ be a global space and let $f\colon X \to Y$ be a morphism of global spaces $g\colon X \to B$ and $h\colon Y \to B$ faithful over $B$. Then the commutative square
\[\xymatrix{
\LocSys^\Glo_\ab(Y,\fcat{A}) \otimes_{\LocSys^\Glo_\ab(B, \fcat{A})} \LocSys^\gen_\ab(B, \fcat A) \ar[d]^-{f^{\Glo,*}\otimes \Id} \ar[r]^-{\theta_h}& \LocSys^\gen_\ab(Y,\fcat A) \ar[d]^-{f^{\gen,*}} \\
\LocSys^\Glo_\ab(X, \fcat A) \otimes_{\LocSys^\Glo_\ab(B, \fcat A)} \LocSys^\gen_\ab(B, \fcat A) \ar[r]^-{\theta_{g}} & \LocSys^\gen_\ab(X,\fcat A) 
}\]
is horizontally right adjointable.
\end{lem}

\begin{proof}
It suffices to show that the commutative diagram is vertically left adjointable. By the projection formula for $f^{\Glo}_{\#}$, see \cref{corollary: left adjoint faithful global linear}, the left adjoint for the left vertical map is $f^{\Glo}_{\#}\otimes \Id$. By the projection formula (\cref{proposition: genuine sharp}) for $f^{\gen}_{\#}$, the assertion follows.
\end{proof}

\begin{lem}\label{lemma: dexterity descent}
Let $X_\bullet\colon I \to \Type^\Glo_{/B}$ be a small diagram of global spaces faithful over $B \in \Type^\Glo$. Then the functor induced by $*$-pullbacks
$$\LocSys^\Glo_\ab(\colim_{I }X_\bullet, \fcat{A}) \otimes_{\LocSys^\Glo_\ab(B, \fcat{A})} \LocSys^\gen_\ab(B, \fcat{A}) \tto \lim_{I^{\op}}\left(\LocSys^\Glo_\ab(X_\bullet, \fcat{A}) \otimes_{\LocSys^\Glo_\ab(B, \fcat{A})} \LocSys^\gen_\ab(B, \fcat{A}) \right)$$
is an equivalence.
\end{lem}

\begin{proof}
By the descent property as in \cref{prop_LSglo_preserves_limits} and the projection formula \cref{corollary: left adjoint faithful global linear}, we have
$$\LocSys^\Glo_\ab(\colim_{I }X_\bullet, \fcat{A}) \simeq \lim_{I^{\op}}\LocSys^\Glo_\ab(X_\bullet, \fcat{A}) \simeq \colim_{I }\LocSys^\Glo_\ab(X_\bullet, \fcat{A}) $$
in the category of $\LocSys^\Glo_\ab(B, \fcat{A})$-modules, where the colimit is taken along $\#$-pushforward. This implies the assertion.
\end{proof}

\begin{prop}\label{proposition: dexterity comparison}
Let $G$ be a compact abelian Lie group and let $f\colon BG \to Y=\mathcal{N}_{\leq G}$ be the unique faithful morphism. Then the comparison functor
$$\alpha_f\colon \LocSys^\Glo_\ab(BG, A) \otimes_{\LocSys^\Glo_\ab(Y, A)} \LocSys^\temp(Y, A) \tto \LocSys^{\temp}(BG,A) $$
is an equivalence.
\end{prop}

\begin{proof} 
By the adjoint functor theorem, the functor $\alpha_f$ admits a right adjoint $\beta_f$ and the generators $L_{BG}\mathcal{L} \in \LocSys^{\temp}(BG,A)$, $\mathcal{L} \in \LocSys^\Glo_\ab(BG,A)$ belong to the essential image of $\alpha_f$. Therefore, its right adjoint $\beta_f$ is conservative. We will show that $\alpha_f$ is fully faithful, or equivalently, the unit map 
$$\eta\colon \Id\tto \beta_f\alpha_f $$
is an equivalence.

We will prove that $\eta$ is an equivalence by induction on the size of $G$. If $G=\{e\}$ is trivial, then $\alpha_f$ is an equivalence because $f$ is the identity morphism. Assume now that $\alpha_f$ is an equivalence for all proper subgroups $H\subsetneq G$. Let $\pi\colon (BG)_{\fcat{P}}=BG\times_{\mathcal{N}}\mathcal{N}_{\fcat{P}} \to BG$ be the projection. Then the functors
$$\pi^{\Glo,*}\colon \LocSys^\Glo_\ab(BG,A)\tto \LocSys^\Glo_\ab((BG)_{\fcat{P}},A),$$
$$(-)^G\colon \LocSys^\Glo_\ab(BG,A)\tto \QCoh(A(BG)) $$
are jointly conservative. Let $\fcat{V}$ denote $\LocSys^\Glo_\ab(Y, A)$ and $\fcat{M}$ denote $\LocSys^\temp(Y, A)$. Since the functors $\pi^{\Glo,*}, (-)^G$ are part of the strict $\fcat{V}$-linear recollement (see \Cref{proposition: global sections semi-orthogonal})
\[\xymatrix{
\LocSys^\Glo_\ab((BG)_{\fcat{P}},A) \ar@<0.5ex>@{^(->}[r]^-{\pi^{\Glo}_{\#}} & 
\ar@<0.5ex>[l]^-{\pi^{\Glo.*}} 
\LocSys^\Glo_\ab(BG,A) \ar@<0.5ex>[r]^-{(-)^G} & \ar@{^(->}@<0.5ex>[l] \QCoh(A(BG)),
}
\]
the induced functors 
$$\widetilde{\pi}^{*}=\pi^{\Glo,*}\otimes \Id \colon \LocSys^\Glo_\ab(BG,A)\otimes_{\fcat{V}} \fcat{M}\tto \LocSys^\Glo_\ab((BG)_{\fcat{P}},A)\otimes_{\fcat{V}} \fcat{M}, $$
$$\widetilde{\Phi}^G\simeq (-)^G\otimes \Id \colon \LocSys^\Glo_\ab(BG,A)\otimes_{\fcat{V}} \fcat{M}\tto \QCoh(A(BG))\otimes_{\fcat{V}} \fcat{M} \simeq \QCoh(A(BG))^{(\fcat{P},G)\mdef\llocal} $$
are jointly conservative as well. By \cref{corollary: isotropy separation geometric points tempered} and \cref{lemma: dexterity isotropy separation tempered}, the natural transformation $\widetilde{\Phi}^G(\eta)$ is an equivalence. Therefore, it suffices to show that $\widetilde{\pi}^*(\eta)$ is an equivalence. However, by \cref{lemma: dexterity isotropy separation tempered induction}, we have $\widetilde{\pi}^*(\beta_f\alpha_f)\simeq \beta_{f\circ\pi }\alpha_{f\circ \pi}$. Hence, $\widetilde{\pi}^*(\eta)$ is an equivalence by the inductive assumption and \cref{lemma: dexterity descent}.
\end{proof}

\begin{thm}\label{theorem: dextery faithful tempered}
Let $f\colon X\to Y$ be a faithful morphism of orbispaces. Then the comparison functor
$$\alpha_f\colon \LocSys^\Glo_\ab(X, A) \otimes_{\LocSys^\Glo_\ab(Y, A)} \LocSys^\temp(Y, A) \tto \LocSys^{\temp}(X,A) $$
is an equivalence.
\end{thm}

\begin{proof}
By \Cref{rem_LSGlo_T_depends_on_SGlo_T} and the definition of $\LocSys^\temp$ we can assume without loss of generality that $X$ and $Y$ lie in $\Type^{\Orb, \ab}$. Let $g\colon Y \to Z$ be another faithful morphism of orbispaces. Note that if $\alpha_{g}$ and $\alpha_{g\circ f}$ are equivalences, then $\alpha_f$ is an equivalence as well. Therefore, we can assume without loss of generality that $f\colon X\to \mathcal{N}$ is the unique faithful morphism to the normal subgroup classifier $\mathcal{N}$. Also, by \cref{lemma: dexterity descent}, we can assume that $X=BG$ for some compact abelian Lie group $G$. In this case, the assertion follows by \cref{proposition: dexterity comparison} and \cref{lemma: family_classfier proj_and_cont}. 
\end{proof}

\appendix
\numberwithin{thm}{section}
\numberwithin{equation}{section}
\crefalias{section}{appendix}

\section{Reminder on spectral algebraic geometry}\label{appendix: sag}
We briefly review necessary notation from spectral algebraic geometry. In this work it is more convenient for us to take a functor of points approach to spectral stacks (see e.g.\@ \cite{GaitsRozI} for a similar theory in the context of derived algebraic geometry) instead of Lurie's locally ringed topoi \cite[Part I]{Lur_SAG}. We refer readers to \cite{BDL25} for a more thorough discussion of the subject.
%\begin{notation}
%Let $\CAlg$ denotes the category of $E_\infty$-algebras in spectra and let $\CAlg_{\ge 0}$ denotes the full subcategory spanned $E_\infty$-algebras in connective spectra. We denote by $\tau_{\ge 0}\colon \CAlg \to \CAlg_{\ge 0}$ the right adjoint of the inclusion.
%\end{notation}
\begin{defn}\label{def_tau_SpStk}
Let $\kappa$ be a regular cardinal. Let $\CAlg^\kappa$ (resp.\@ $\CAlg_{\ge 0}^\kappa$) denote the category of $\kappa$-compact $E_\infty$-ring spectra (resp.\@ $\kappa$-compact connective $E_\infty$-ring spectra). Let $\tau$ be a subcanonical Grothendieck topology on $\CAlg$ which restricts to $\CAlg^\kappa$ and $\CAlg^\kappa_{\ge 0}$ for all sufficiently large regular cardinals $\kappa$. Then we define
$$\SpStk_{\tau, \kappa} := \Shv_\tau(\CAlg_{\ge 0}^{\kappa,\op}, \Type), \qquad \SpStk^\nc_{\tau, \kappa} := \Shv_\tau(\CAlg^{\kappa,\op}, \Type).$$
If $\lambda \ge \kappa$, then the inclusion $\CAlg_{\ge 0}^{\kappa} \inj \CAlg_{\ge 0}^{\lambda}$ induces colimit preserving embeddings
$$\SpStk_{\tau, \kappa} \inj \SpStk_{\tau, \lambda}, \qquad \SpStk_{\tau, \kappa}^\nc \inj \SpStk_{\tau, \lambda}^\nc.$$
We define the (large) category $\SpStk_\tau$ of \emdef{spectral $\tau$-stacks} and the category $\SpStk^\nc_{\tau}$ of \emdef{non-connective spectral $\tau$-stacks} as follows
$$\SpStk_\tau := \colim_\kappa \SpStk_{\tau, \kappa}, \qquad \SpStk^\nc_{\tau} := \colim_\kappa \SpStk_{\tau, \kappa}^\nc.$$
For $R \in \CAlg_{\ge 0}$ (resp.\@ $\CAlg$), we will denote by $\Spec R$ the image of $R$ in $\SpStk$ (resp.\@ $\SpStk^\nc$) under the Yoneda embedding.

If $\tau$ is the trivial topology on $\CAlg$, then we denote the corresponding sheaf category by $\SpPStk$ (resp.\@ $\SpPStk^\nc$) and call it (resp.\@ \emdef{non-connective}) \emdef{spectral prestacks}. By construction, for a general $\tau$, there are fully faithful forgetful functors
$$\SpStk_\tau \inj \SpPStk, \qquad \SpStk_\tau^\nc \inj \SpPStk^\nc$$
which admit left exact (i.e.\@ finite limit preserving) left adjoints
$$L_\tau\colon \SpPStk \to \SpStk_\tau, \qquad L_\tau^\nc \colon \SpPStk^\nc \to \SpStk_\tau^\nc.$$
\end{defn}
\begin{rem}
Equivalently, $\SpPStk$ (resp.\@ $\SpPStk^\nc$) is the full subcategory of $\Fun(\CAlg_{\ge 0}, \Type)$ (resp.\@ $\Fun(\CAlg, \Type)$) spanned by the functors which are equivalent to a small colimit of representables.
\end{rem}
\begin{rem}\label{SpPStk_emb_SpPStkNc}
The inclusion $\CAlg_{\ge 0}^{\kappa} \inj \CAlg^{\kappa}$ induces a fully faithful functor
$$\SpStk_{\tau,\kappa} \inj \SpStk_{\tau,\kappa}^\nc.$$
By passing to the colimit in $\kappa$, we obtain a similar embedding
$$\SpStk_\tau \inj \SpStk^\nc_\tau,$$
which admits a left adjoint which we will denote by
$$X \xymatrix{\ar@{|->}[r] &} X_{\ge 0}.$$ 
\end{rem}
\begin{warn}
The categories $\SpStk_\tau$ and $\SpStk_\tau^\nc$ are \emph{not} $\infty$-topoi. Nevertheless, any small $\SpStk_\tau$- (resp.\@ $\SpStk_\tau^\nc$-) valued diagram lands in some $\SpStk_{\tau, \kappa}$ (resp.\@ $\SpStk_{\tau, \kappa}^\nc$) for large enough $\kappa$, so in many respects one can essentially work with these categories as if they were $\infty$-topoi. For instance, we refer the reader to \cite[Proposition 2.1.2.9]{BDL25} for the proofs of several nice properties of $\SpStk_\tau$ and $\SpStk_\tau^\nc$, e.g.\@ the universality of colimits.
\end{warn}
\begin{notation}
We denote the Yoneda embedding
$$\CAlg^\op \tto \SpStk^{\nc}$$
by $\Spec$. A non-connective spectral prestack $S$ is \emdef{non-connective affine} or just \emdef{nc-affine} (resp.\@ \emdef{affine}) if $S$ is equivalent to $\Spec R$ for some (resp.\@ connective) $E_\infty$-ring $R$. A morphism $f\colon X \to Y$ of non-connective spectral prestacks is called \emdef{non-connective affine} or \emdef{nc-affine} if the fibered product $S\times_Y X$ is nc-affine for every nc-affine $S$ and every map $S \to Y$. A non-connective affine map $f\colon X \to Y$ is \emdef{affine} if additionally for any nc-affine $S$ mapping to $Y$ with the fiber product $T:= S\times_Y X$ the commutative square
\[\xymatrix{
T \ar[r]\ar[d] & T_{\ge 0} \ar[d] \\
S \ar[r] & S_{\ge 0}
}\]
is fibered.
\end{notation}

Next we review the theory of quasi-coherent sheaves in this setting.
\begin{defn}
For a regular cardinal $\kappa$ we define a functor
$$\QCoh\colon \SpPStk_\kappa^{\nc, \op} \tto \CAlg(\Prs^\LL)$$
as the right Kan extension of the functor $A \mapsto \Mod_A(\Sp)$ along the Yoneda embedding. If $\lambda \ge \kappa$, then the obtained functors are naturally compatible with each other, so we can glue to a functor
$$\QCoh \colon \SpPStk^{\nc, \op} \tto \CAlg(\Prs^\LL).$$
For a non-connective spectral prestack $X\in \SpPStk^\nc$, we will call $\QCoh(X)$  the category of \emdef{quasi-coherent sheaves on $X$}. The category of quasi-coherent sheaves on a connective spectral stack can be defined similarly or via the inclusion of \Cref{SpPStk_emb_SpPStkNc}.
\end{defn}
\begin{notation}
Let $f\colon X \to Y$ be a morphism of non-connective spectral prestacks. Then, by construction, there is a symmetric monoidal colimit preserving pullback functor
$$f^*\colon \QCoh(Y) \tto \QCoh(X).$$
By the adjoint functor theorem, the pullback $f^*$ admits a right lax monoidal right adjoint
$$f_* \colon \QCoh(X) \tto \QCoh(Y).$$
\end{notation}
\begin{rem}\label{rem_QCoh_fpqc}
Since the functor $R \mapsto \Mod_R$ satisfies fpqc-descent, the functor $\QCoh$ factors through the category of non-connective spectral stacks in the fpqc topology $\SpStk^\nc_\fpqc$. In particular, if $f\colon X \to Y$ is such that the induced map $L_\fpqc X \to L_\fpqc Y$ is an equivalence, then so is
$$f^* \colon \QCoh(Y) \tto \QCoh(X).$$
\end{rem}

Finally, we review the relation between this notion of spectral prestacks and Lurie's (non-connective) spectral Deligne--Mumford stacks $\SpDM$ (resp.\@ $\SpDM^\nc$), see~\cite[Section~1.4.4]{Lur_SAG}. Let $X$ be a non-connective spectral Deligne--Mumford stack. Then we will write
$$h^\nc_X \colon \CAlg \tto \Type$$
for the functor
$$h^\nc_X(R) := \Hom_{\SpDM^\nc}(\Spet R, X).$$
If $X$ is connective, then the spectral prestack $h_X$ is defined similarly.
\begin{prop}
The functors
$$h_- \colon \SpDM \tto \SpPStk, \qquad h^\nc_- \colon \SpDM^\nc \tto \SpPStk^\nc$$
are fully faithful and land in the full subcategories of (non-connective) \'etale spectral stacks. Moreover, if $X$ is a (non-connective) spectral Deligne--Mumford stack, then there is a natural equivalence
$$\QCoh(X) \simeq \QCoh(h^\nc_X),$$
where the category $\QCoh(X)$ on the left is defined in \cite[Definition 2.2.2.1]{Lur_SAG}.
\end{prop}

\begin{proof}
See e.g.~\cite[Proposition 1.6.4.2]{Lur_SAG} or~\cite[Theorem 3.2.1.3]{BDL25} for the proof.
\end{proof}

\begin{conv}
In our work we are more interested in the category of quasi-coherent sheaves on non-connective abelian group spectral stacks rather than in the object itself. By \Cref{rem_QCoh_fpqc}, it follows that in \Cref{def_tau_SpStk} we can pick any Grothendieck topology $\tau$ which is coarser or equal to fpqc and consider $\SpStk_\tau^\nc$ as an ambient geometric category. In fact, the trivial topology would do. For the rest of the paper, $\SpStk$ (resp.\@ $\SpStk^\nc$) will always denote the category of (resp.\@ non-connective) spectral $\tau$-stacks for any one of these choices of $\tau$.
\end{conv}

\section{Cartesian fibrations and lax limits}\label{appendix: cart_lax_limits}
In this section we review the notion of (both left and right) lax limits as well as some properties of lax Kan extensions used in the main text. The main result is \cref{on_KanExtensions_of_GrothFib}.

%In this section we review the notion of partial (op)lax (co)limits in the category of categories and prove some technical facts used in the work.
\begin{defn}[Definition~3.9 in \cite{haugseng2021}]\label{lax_functor_category}
Let $\fcat{C}$ be a small $(\infty,1)$-category and let $\fcat{W}$ be an $(\infty,2)$-category. We write $\Fun(\fcat{C},\fcat{W})^{\llax}$ (resp.\@ $\Fun(\fcat{C},\fcat{W})^{\rlax})$ for the category which has functors as objects and left (resp.\@ right) lax natural transformations as morphisms. Similarly, we write $\Fun(\fcat{C},\fcat{W})^{\mathrm{strict}}$ for the functor category which has strict natural transformations as morphisms.
\end{defn}

\begin{rem}[Remark~3.10 in \cite{haugseng2021}]\label{lax_category_remark}
%\todo{Change the orientation. Remark that we use the opposite to Haugseng.}
Informally, a \emph{left lax} natural transformation $\eta$ between functors $F,G\colon \fcat{C} \to \fcat{W}$ assigns to every morphism $f\colon X \to Y$ in $\fcat{C}$ a left lax commutative square
\begin{equation}\label{lax_category_remark_eq1}
\xymatrix{
F(X) \ar[r]^{\eta_X}\ar[d]_{F(f)} & G(X) \ar[d]^{G(f)} \ar@{=>}[ld]_{\eta_f} \\
F(Y) \ar[r]^{\eta_Y} & G(Y);
}
\end{equation}
while a \emph{right lax} natural transformation assigns a right lax commutative square
\begin{equation}\label{lax_category_remark_eq2}
\xymatrix{
F(X) \ar[r]^{\eta_X}\ar[d]_{F(f)} & G(X) \ar[d]^{G(f)}  \\
F(Y)\ar@{=>}[ru]^{\eta_f} \ar[r]^{\eta_Y} & G(Y).
}
\end{equation}
\end{rem}

For an $(\infty,2)$-category $\fcat{W}$, a functor $p\colon \fcat{C} \to \fcat{D}$ induces the natural functor
$$p^*\colon \Fun(\fcat{D},\fcat{W})^{\llax} \to \Fun(\fcat{C},\fcat{W})^{\llax} 
\;\;\text{(resp. $p^*\colon \Fun(\fcat{D},\fcat{W})^{\rlax} \to \Fun(\fcat{C},\fcat{W})^{\rlax}$)} $$
between the left (resp. right) lax functor categories.

\begin{notation}\label{definition: lax_kan_extension}
We will write 
$$\Ran_p^{\llax}\colon \Fun(\fcat{C},\fcat{W})^{\llax} \to \Fun(\fcat{D},\fcat{W})^{\llax}$$ 
for the right adjoint to the pullback functor $p^*$, if $p^*$ admits a right adjoint. Similarly, we will write
$$\Ran_p^{\rlax}\colon \Fun(\fcat{C},\fcat{W})^{\rlax} \to \Fun(\fcat{D},\fcat{W})^{\rlax}$$ 
for the right adjoint to $p^*$, if $p^*$ admits a right adjoint.
\end{notation}

\begin{rem}\label{remark: existence of right adjoints}
At the moment of writing, we are not aware of a sufficient condition for the existence of $\Ran^{\llax}_p$ (resp.\@ $\Ran^{\rlax}_p$) for a given functor $p$ even in the case $\fcat W = \Cat$, see~\cite{AyalaFrancis20}, \cite[Section~5.6]{AHM26}, or \cite[Section~4.2]{Abellan2023} for a detailed account.
\end{rem}

\begin{defn}\label{lax_kan_extensions_definition}
Let $p\colon \fcat{C} \to \fcat{D}$ and $F\colon \fcat{C} \to \fcat{W}$ be functors. We define the \emdef{left lax right Kan extension 
$$\Ran^{\llax}_p F\colon \fcat{D} \to \fcat{W}$$ of $F$ along $p$} by evaluating the right adjoint functor $\Ran^{\llax}_p$ of~\cref{definition: lax_kan_extension} at $F$. We define the \emdef{right lax right Kan extension 
$$\Ran^{\rlax}_p F\colon \fcat{D} \to \fcat{W}$$ of $F$ along $p$} in a similar way.
\end{defn}

\iffalse
\begin{ex}\label{lax_kan_adjoint}
Let $p\colon \fcat{C} \to \fcat{D}$ be a functor which has a right (resp. left) adjoint $q$. Then the (op)lax left (resp. right) Kan extension $\Lan_p^{\mathrm{(op)lax}}$ (resp. $\Ran_p^{\mathrm{(op)lax}}$) is given by $q^*$. \todo{CHECK}
\end{ex}
\fi

\begin{defn}\label{lax_limit}
Let $F\colon \fcat{C} \to \fcat{W}$ be a functor. We define the \emdef{left lax limit} $$\llaxlim\left(\xymatrix{\fcat{C} \ar[r]^-{F} & \fcat{W}}\right)$$
as the \emph{left} lax right Kan extension $\Ran^{\llax}_{\pi} F$ of $F$ along the canonical projection $\pi\colon \fcat{C} \to \{*\}$. We define the \emdef{right lax limit} $$\rlaxlim\left(\xymatrix{\fcat{C} \ar[r]^-{F} & \fcat{W}}\right)$$
in a similar way as the \emph{right} lax right Kan extension.
%Similarly, we define the \emph{(op)lax colimit}
%$$\colim{}^{\mathrm{(op)lax}}\left(\xymatrix{\fcat{C} \ar[r]^-{F} & \fcat{W}}\right)$$
%as the (op)lax left Kan extension $\Lan^{\mathrm{(op)lax}}_{\pi} F$ of $F$ along the canonical projection $\pi$.
\end{defn}

Let $p\colon \fcat{C} \to \fcat{D}$ be a functor. Consider the induced pullback functor
$$p^*\colon \Cat_{/\fcat{D}} \to \Cat_{/\fcat{C}}, \;\; (\fcat{E} \to \fcat{D}) \mapsto \fcat{E}\times_{\fcat{D}}\fcat{C}$$
between the slice categories $\Cat_{/\fcat{D}}$ and $\Cat_{/\fcat{C}}$. The pullback $p^*$ admits a left adjoint
$$p_!\colon \Cat_{/\fcat{C}} \to \Cat_{/\fcat{D}} $$
which is given by the formula $(\fcat{E} \to \fcat{C})\mapsto (\fcat{E} \to \fcat{C} \xrightarrow{p} \fcat{D})$. Note that the functor $p^*$ admits a right adjoint if and only if $p$ is an \emph{exponential fibration} and $p_*$ is usually rather complicated, see \cite{AyalaFrancis20} for a detailed account.

\begin{ex}\label{example: right adjoint to projection}
Compare with \cite[Example~2.2.2]{AyalaFrancis20}. Let $\pi\colon \fcat{C} \to \{*\}$ be the canonical projection. Then the right adjoint $\pi_*\colon \Cat_{/\fcat{C}} \to \Cat$ is given by the category of sections of the structure morphism, i.e. $$\pi_*(\fcat{E} \to \fcat{C}) \simeq \Fun_{/\fcat{C}}(\fcat{C},\fcat{E}).$$
\end{ex}

\begin{ex}\label{example: lax limits as sections}
Let $\Cat^{\mathrm{cart}}_{/\fcat{C}^{\op}}$ denote the full subcategory of the slice category $\Cat_{/\fcat{C}^{\op}}$ spanned by Cartesian fibrations. By e.g. \cite[Theorem~5.3.1]{haugseng2021lax}, there is the unstraightening equivalence
$$\Unst\colon \Fun(\fcat{C}, \Cat)^{\rlax} \xrightarrow{\simeq} \Cat^{\mathrm{cart}}_{/\fcat{C}^{\op}}.$$
Under the equivalence above, the restriction functor
$$\Cat \to \Fun(\fcat{C},\Cat)^{\rlax}, \;\; G \mapsto G\circ \pi $$
corresponds to the pullback $\pi^*\colon \Cat \to \Cat^{\mathrm{cart}}_{/\fcat{C}^{\op}}\subset \Cat_{/\fcat{C}^{\op}}$
along the canonical projection $\pi\colon \fcat{C}^{\op} \to \{*\}$. Therefore, by \cref{example: right adjoint to projection}, we have
$$\rlaxlim\left(\xymatrix{\fcat{C} \ar[r]^-{F} & \Cat}\right) \simeq \Fun_{/\fcat{C}^\op}(\fcat{C}^\op,\Unst(F))$$
for a functor $F\colon \fcat{C} \to \Cat$. Similarly, let $\Cat^{\mathrm{cocart}}_{/\fcat{C}}$ denote the full subcategory of the slice category $\Cat_{/\fcat{C}}$ spanned by coCartesian fibrations. Then there is the unstraightening equivalence
$$\coUnst\colon \Fun(\fcat{C}, \Cat)^{\llax} \xrightarrow{\simeq} \Cat^{\mathrm{cocart}}_{/\fcat{C}}$$
and we obtain
$$\llaxlim\left(\xymatrix{\fcat{C} \ar[r]^-{F} & \Cat}\right) \simeq \Fun_{/\fcat{C}}(\fcat{C},\coUnst(F))$$
for a functor $F\colon \fcat{C} \to \Cat$.
\end{ex}

\begin{rem}\label{remark: lax and weighted}
We note that any left (resp.\@ right) lax limit of the diagram $F\colon \fcat{C} \to \fcat{W}$ can be expressed as a \emph{weighted} limit in the category $\fcat{W}$, see~\cite[Section~5.2]{AHM26}. In particular, if $G\colon \fcat{W} \to \fcat{W}'$ is a functor of $(\infty,2)$-categories which preserves arbitrary weighted limits, then
$$\llaxlim\left(\xymatrix{\fcat{C} \ar[r]^-{F} & \fcat{W}  \ar[r]^-{G} & \fcat{W}'}\right) \simeq G\left(\llaxlim\left(\xymatrix{\fcat{C} \ar[r]^-{F} & \fcat{W}}\right)\right) $$
and similarly for the right lax limit. 

As a consequence, let $F\colon \fcat{C} \to \Prs^\LL$ (resp. $G\colon \fcat{C} \to \Prs^\RR$) be a diagram of presentable categories and left (resp. right) adjoints and let $\widetilde{F}\colon \fcat{C} \to \Cat$ be the same diagram considered as a diagram of (large) categories and arbitrary functors. Then
\begin{align*}
\llaxlim\left(\xymatrix{\fcat{C} \ar[r]^-{F} & \Prs^\LL}\right) &\simeq \llaxlim\left(\xymatrix{\fcat{C} \ar[r]^-{\widetilde{F}} & \Cat}\right), \\
\llaxlim\left(\xymatrix{\fcat{C} \ar[r]^-{F} & \Prs^\RR}\right) &\simeq \llaxlim\left(\xymatrix{\fcat{C} \ar[r]^-{\widetilde{F}} & \Cat}\right)
\end{align*}
and similarly for the right lax limits.
\end{rem}

\begin{lem}\label{lemma: base change for slice categories}
Let
$$
\xymatrix{
\widetilde{\fcat{C}} \ar[r]^f \ar[d]_{q} & \fcat{C} \ar[d]^p \\
\widetilde{\fcat{D}} \ar[r]^{g} & \fcat{D}.
}
$$
be a fibered square of categories. Then there are natural equivalences $g^*p_! \simeq q_! f^*$ and, if $p$ is exponential, $g^*p_* \simeq q_* f^*$. %of functors between the slice categories $\Cat_{/\fcat{C}}$ and $\Cat_{/\widetilde{\fcat{D}}}$.
\end{lem}

\begin{proof}
The first equivalence follows directly from the definitions and the second equivalence follows from the first one by taking the right adjoints. Note that $q$ is an exponential fibration as a base change of the exponential fibration $p$, see~\cite[Corollary~2.3.1]{AyalaFrancis20}.
\end{proof}

\begin{lem}\label{lemma: right adjoint along cartesian}
$\quad$
\begin{enumerate}
\item Suppose that $p\colon \fcat{C} \to \fcat{D}$ is a Cartesian fibration, then $p_*$ exists and sends coCartesian fibrations over $\fcat{C}$ to coCartesian fibrations over $\fcat{D}$.
\item Suppose that $p\colon \fcat{C} \to \fcat{D}$ is a coCartesian fibration, then $p_*$ exists and sends Cartesian fibrations over $\fcat{C}$ to Cartesian fibrations over $\fcat{D}$.
\end{enumerate}
\end{lem}

\begin{proof}
The first part follows by~\cite[Corollary~3.2.2.12]{Lur_HTT}, see also \cite[Proposition~B.4.10]{Lur_HA}. The second part follows from the first one by applying an autoequivalence
$$\Cat_{/\fcat{C}} \simeq \Cat_{/\fcat{C}^{\op}}, \;\; (\fcat{E} \to \fcat{C}) \mapsto (\fcat{E}^{\op}\to \fcat{C}^{\op}), $$
which switches Cartesian and coCartesian fibrations.
\end{proof}

\begin{cor}\label{on_KanExtensions_of_GrothFib}
Let $p\colon \fcat C \to \fcat D$ be a Cartesian fibration and let $F\colon \fcat C \to \Cat$ be a functor. Then there are natural equivalences
\begin{align*}
(\Ran_p^\rlax F)(Y) & \simeq \rlaxlim\left(\xymatrix{\fcat C_Y \ar[r] & \fcat C \ar[r]^-{F} & \Cat}\right),\\
(\Ran_p^\llax F)(Y) & \simeq \llaxlim\left(\xymatrix{\fcat C_Y \ar[r] & \fcat C \ar[r]^-{F} & \Cat}\right),
\end{align*}
where $\fcat C_Y=\{Y\}\times_{\fcat D} \fcat C$ denotes the strict (i.e.~non-lax) fiber of $p$ over $Y$.
\end{cor}

\begin{proof}
For the first part, let $\Cat^{\mathrm{cart}}_{/\fcat{C}^{\op}}$ denote the full subcategory of the slice category $\Cat_{/\fcat{C}^{\op}}$ spanned by Cartesian fibrations. By e.g.\@ \cite[Theorem~5.3.1]{haugseng2021lax}, there is the unstraightening equivalence
$$\Unst\colon \Cat^{\mathrm{cart}}_{/\fcat{C}^{\op}} \xymatrix{\ar[r]^-\sim &} \Fun(\fcat{C}, \Cat)^{\rlax}.$$
%here the target category has functors as objects and oplax natural transformations as morphisms.
Under the equivalence above, the restriction functor
$$\Fun(\fcat{D},\Cat)^{\rlax} \to \Fun(\fcat{C},\Cat)^{\rlax}, \;\; G \mapsto G\circ p $$
corresponds to the pullback $(p^{\op})^*\colon \Cat^{\mathrm{cart}}_{/\fcat{D}^{\op}} \to \Cat^{\mathrm{cart}}_{/\fcat{C}^{\op}}$
along $p^{\op}\colon \fcat{C}^{\op} \to \fcat{D}^{\op}$. Since $p^{\op}$ is a coCartesian fibration, \Cref{lemma: right adjoint along cartesian} implies that the right adjoint $$(p^{\op})_*\colon \Cat_{/\fcat{C}^{\op}} \to \Cat_{/\fcat{D}^{\op}}$$ preserves Cartesian fibrations, i.e. $(p^{\op})_*$ maps $\Cat^{\mathrm{cart}}_{/\fcat{C}^{\op}}$ onto $\Cat^{\mathrm{cart}}_{/\fcat{D}^{\op}}$. It follows that $\Ran_p^\rlax$ exists and is given by $(p^\op)_*$ under the unstraightening equivalence.

Let $\Unst(F) \in \Cat^{\mathrm{cart}}_{/\fcat{C}^{\op}}$ be the Cartesian fibration corresponding to the functor $F\colon \fcat{C} \to \Cat$. For $Y \in \fcat D$, we then have an equivalence
$$(\Ran_p^{\rlax}F)(Y) \simeq \{Y\}\times_{\fcat{D}^{\op}} p^{\op}_*(\Unst(F)) $$
By applying \Cref{lemma: base change for slice categories} to the fibered square
\[\xymatrix{
(\fcat{C}_Y)^{\op} \ar[r] \ar[d]^{\pi} & \fcat{C}^{\op} \ar[d]^{p^{\op}} \\
\{Y\} \ar[r] & \fcat{D}^{\op},
}\]
we obtain $$\{Y\}\times_{\fcat{D}^{\op}} p^{\op}_*(\Unst(F)) \simeq \pi_*\left(\Unst(F)|_{(\fcat{C}_Y)^{\op}}\right.)$$ Finally, the latter is equivalent to the right lax limit $\rlaxlim\left(\xymatrix{\fcat C_Y \ar[r] & \fcat C \ar[r]^-{F} & \Cat}\right)$. 

The proof of the second assertion is similar. Let $\Cat^{\mathrm{cocart}}_{/\fcat{C}}$ denote the full subcategory of the slice category $\Cat_{/\fcat{C}}$ spanned by coCartesian fibrations. By \cite[Theorem~5.3.1]{haugseng2021lax} again, there is the unstraightening equivalence
$$\coUnst\colon \Cat^{\mathrm{cocart}}_{/\fcat{C}} \xymatrix{\ar[r]^-\sim &} \Fun(\fcat{C}, \Cat)^{\llax}.$$
The rest of the proof follows by \Cref{lemma: base change for slice categories} and \Cref{lemma: right adjoint along cartesian}.
\end{proof}

\begin{ex}\label{example: oplax kan extension slice}
Let $\fcat{C}_{/X}$ be the slice category over an object $X\in \fcat{C}$ and let $F\colon \fcat{C}\to \Cat$ be a functor. We write $p\colon \fcat{C}_{/X}\to \fcat{C}$ for the canonical projection. Let us compute the right lax right Kan extension $$\Ran^{\rlax}_p(F\circ p) \colon \fcat{C} \tto \Cat$$ of the restricted functor $F\circ p$ along the projection $p$. By \Cref{on_KanExtensions_of_GrothFib}, we have an equivalence
\begin{equation}\label{equation: example, slice, eq1}
\left(\Ran^{\rlax}_p(F\circ p)\right)(Y) \simeq \rlaxlim\left(\xymatrix{(\fcat{C}_{/X})_Y \ar[r] & \fcat{C}_{/X} \ar[r]^-{p} & \fcat{C} \ar[r]^-{F} & \Cat}\right),
\end{equation}
where $Y\in \fcat{C}$. We note that the strict fiber $(\fcat{C}_{/X})_Y= \{Y\} \times_{\fcat C} \fcat{C}_{/X}$ is equivalent to the Hom-groupoid $\Hom(Y,X)$, and the diagram in~\eqref{equation: example, slice, eq1} factors as the following composite 
$$\Hom(Y,X) \tto \{Y\} \xrightarrow{Y\mapsto F(Y)} \Cat. $$
%Here the last functor is the constant functor. 
Therefore, we obtain an equivalence
\begin{equation}\label{equation: example, slice, eq2}
\left(\Ran^{\rlax}_p(F\circ p)\right)(Y) \simeq \Fun(\Hom(Y,X)^{\op},F(Y)).
\end{equation}
Moreover, if $f\colon Y\to Z \in \fcat{C}$ is a morphism in $\fcat{C}$, then the induced functor
$$\left(\Ran^{\rlax}_p(F\circ p)\right)(f)\colon \left(\Ran^{\rlax}_p(F\circ p)\right)(Y) \tto \left(\Ran^{\rlax}_p(F\circ p)\right)(Z) $$
corresponds under the equivalence~\eqref{equation: example, slice, eq2} to the composite
$$\Fun(\Hom(Y,X)^{\op},F(Y)) \xrightarrow{F(f)} \Fun(\Hom(Y,X)^{\op},F(Z)) \xrightarrow{\Hom(f,X)^*} \Fun(\Hom(Z,X)^{\op},F(Z)).$$
Finally, we consider the right lax natural transformation $\eta\colon F\to \Ran^{\rlax}_p(F\circ p)$ which is the unit map of the adjunction from \Cref{definition: lax_kan_extension}. One can see that, under the equivalence~\eqref{equation: example, slice, eq2}, $\eta$ is given by the functors of constant diagrams,~i.e.
$$\eta_Y \simeq \const_{F(Y)} \colon F(Y) \tto \Fun(\Hom(Y,X)^{\op}, F(Y)). $$
We note that $\eta$ is a \emph{strict} natural transformation in this example.
\end{ex}

%Let us denote by $\Cat^R$ (resp. $\Cat^L$) the non-full subcategory of $\Cat$ with the same objects and only functors which are right (resp. left) adjoints.

\iffalse
\begin{lem}\label{lemma: (op)lax Kan extension prl vs cat}
Let $F\colon \fcat{C} \to \Prs^\LL$ be a diagram of presentable categories and left adjoints and let $\widetilde{F}\colon \fcat{C} \to \Cat$ be the same diagram considered as a diagram of (large) categories and arbitrary functors. Suppose that $p\colon \fcat{C} \to \fcat{D}$ be a functor. Then there are natural equivalences
$$\Ran_p^{\mathrm{l(r)lax}}F \simeq \Ran_p^{\mathrm{l(r)lax}}\widetilde{F} \in \Fun(\fcat{D}, \Cat). $$
In particular, $\Ran_p^{\mathrm{l(r)lax}}\widetilde{F} \colon \fcat{D} \to \Cat$ is a diagram of presentable categories and left adjoints.
\end{lem}

\begin{proof}
Follows by~\cite[Propositions~5.5.3.6 and~5.5.3.13]{Lur_HTT}.
\end{proof}
\fi

\begin{cor}\label{corollary: right oplax along cocartesian}
Let $p\colon \fcat{C} \to \fcat{D}$ be a coCartesian fibration, let $F\colon \fcat{C} \to \Prs^\LL$ be a functor, and let $F^R\colon \fcat{C}^{\op} \to \Prs^R$ be the corresponding diagram of right adjoints. Then
\begin{enumerate}
\item there is a natural in $Y$ equivalence
$$(\Ran_p^\llax F)(Y) \simeq \rlaxlim\left(\xymatrix{(\fcat{C}^{\op})_Y \ar[r] & \fcat{C}^{\op} \ar[r]^-{F^R} & \Prs^\RR}\right), \;\; Y\in \fcat{D};$$
\item if $f\colon Y\to Z \in \fcat{D}$ is a morphism in $\fcat{D}$, then the functor $(\Ran_p^\llax F)(f)$ is the left adjoint to the canonical functor
$$\rlaxlim\left(\xymatrix{(\fcat{C}^{\op})_{Z} \ar[r] & \fcat{C}^{\op} \ar[r]^-{F^R} & \Prs^\RR}\right) \tto \rlaxlim\left(\xymatrix{(\fcat{C}^{\op})_{Y} \ar[r] & \fcat{C}^{\op} \ar[r]^-{F^R} & \Prs^\RR}\right)$$
induced by the coCartesian lift $\fcat{C}_f\colon \fcat{C}_Y \to \fcat{C}_{Z}$.
\end{enumerate}
\end{cor}

\begin{proof}
By~\cite[Theorem~5.3.6]{haugseng2021lax}, there are equivalences
\begin{align*}
\mathrm{Adj}_{\fcat{C}}\colon \Fun(\fcat{C}, \Prs^\LL)^{\llax} & \xymatrix{\ar[r]^-\sim &} \Fun(\fcat{C}^{\op}, \Prs^\RR)^{\rlax}, \\
\mathrm{Adj}_{\fcat{D}}\colon \Fun(\fcat{D}, \Prs^\LL)^{\llax} & \xymatrix{\ar[r]^-\sim &} \Fun(\fcat{D}^{\op}, \Prs^\RR)^{\rlax}
\end{align*}
sending each diagram of left adjoints $F\colon \fcat{C} \to \Prs^\LL$ (resp.\@ $G\colon \fcat{D} \to \Prs^\LL$) to the corresponding diagram $F^R\colon \fcat{C}^{\op} \to \Prs^\RR$ (resp.\@ $G^R\colon \fcat{D}^{\op} \to \Prs^\RR$) of right adjoints. Under these equivalences, we have
$$\mathrm{Adj}_{\fcat{D}}(\Ran^{\llax}_p F) \simeq \Ran^{\rlax}_{p^{\op}}( \mathrm{Adj}_{\fcat{C}}(F)) \simeq \Ran^{\rlax}_{p^{\op}}(F^R).$$
Since $p^{\op}$ is a Cartesian fibration, \Cref{on_KanExtensions_of_GrothFib} implies the statement.
\end{proof}

\section{Adjoint functors between right lax limits}\label{appendix: right adjoint lax limits}
In this section we will provide formulas for right adjoints between right lax limits of categories in some particular cases. The main results are \cref{lem_adj_functors_to_limit} and \cref{right_adjoint_llax_limits}. First, we consider functors induced by right lax natural transformations.

Let $\eta\colon \fcat{C}_\bullet \to \fcat{D}_\bullet \in \Fun(\fcat{I}, \Cat)^{\rlax}$ be a \emph{right lax} natural transformation between the diagrams of categories. Suppose that each component $\eta_i\colon \fcat{C}_i \to \fcat{D}_i$, $i\in \fcat{I}$ admits a right adjoint $$\eta^R_i\colon \fcat{D}_i \to \fcat{C}_i, \;\; i\in \fcat{I}.$$ Then the functors $\eta^R_i$ form the \emph{right mate} $\eta^R\colon \fcat{D}_\bullet \to \fcat{C}_\bullet$ of $\eta$, i.e.\@ the \emph{left lax} natural transformation
$$\eta^R \colon \fcat{D}_\bullet \to \fcat{C}_\bullet \in \Fun(\fcat{I},\Cat)^{\llax}. $$
We recall that the natural transformation $\eta^R_f\colon \phi_f\circ \eta^R_i \to \eta^R_j \circ \psi_f$, $f\colon i\to j \in \fcat{I}$ is the composite
$$\eta^R_f\colon \phi_f\circ \eta^R_i \to \eta^R_j\circ \eta_j \circ \phi_f \circ \eta^R_i \xrightarrow{\eta^R_j \circ \eta_f \circ \eta^R_i} \eta^R_j\circ \psi_f\circ \eta_i \circ \eta^R_i \to \eta^R_j \circ \psi_f, $$
where $\phi_f\colon \fcat{C}_i \to \fcat{C}_j$ and $\psi_f\colon \fcat{D}_i \to \fcat{D}_j$ are transition functors, see e.g.~\cite[Remark~4.5]{haugseng2021} or~\cite[Theorem~5.3.5]{haugseng2021lax}.

\begin{lem}\label{lemma: right adjoint to oplax limit_lax is strict}
Let $\eta\colon \fcat{C}_\bullet \to \fcat{D}_\bullet \in \Fun(\fcat{I}, \Cat)^{\rlax}$ be a \emph{right lax} natural transformation between the diagrams of categories. Suppose that each component $\eta_i\colon \fcat{C}_i \to \fcat{D}_i$, $i\in \fcat{I}$ is a left adjoint and suppose that the right mate $\eta^R\colon \fcat{D}_\bullet \to \fcat{C}_\bullet$ formed by the right adjoints is \emph{strict}. Then
$$\rlaxlim_{\fcat{I}}(\eta^R) \colon \rlaxlim_{\fcat{I}} \fcat D_\bullet \tto \rlaxlim_{\fcat{I}} \fcat C_\bullet$$
is the right adjoint to the functor $\rlaxlim_{\fcat{I}}(\eta) \colon \rlaxlim_{\fcat{I}}\fcat{C}_\bullet \to \rlaxlim_{\fcat{I}} \fcat D_\bullet$ induced by $\eta$.
\end{lem}

\begin{proof}
Since the right mate $\eta^R\colon \fcat{D}_\bullet \to \fcat{C}_\bullet$ is strict, we can consider $\eta^R$ as a $1$-morphism in the $(\infty,2)$-category  $\Fun(\fcat{I}, \Cat)^\rlax$. Moreover, by \cite[Theorem~4.6(ii)]{haugseng2021}, $\eta^R$ is the right adjoint to $\eta$ in the $(\infty,2)$-category $\Fun(\fcat{I}, \Cat)^\rlax$. Since the right lax limits form the strict $(\infty,2)$-functor
$$\rlaxlim_{\fcat{C}} \colon  \Fun(\fcat{I}, \Cat)^\rlax \to \Cat,$$
the lemma follows.
\end{proof}

\begin{cor}\label{lem_adj_functors_to_limit}
Let $\fcat C_\bullet\colon \fcat{I} \to \Cat$ be a small diagram of categories such that all transition functors $\phi_f\colon \fcat C_i \to \fcat C_j$, $f\colon i \to j$ have left adjoints. Let $\fcat D$ be a category admitting $\fcat{I}^{\op}$-indexed colimits and let $F\colon \fcat D \to \lim_{\fcat{I}} \fcat C_i$ be a functor such that each natural projection 
$$F_i = p_i \circ F \colon \fcat D \to \fcat C_i, \;\; i\in \fcat{I}$$
admits a left adjoint $F_i^L\colon \fcat{C}_i \to \fcat{D}$. Then the functor $F$ admits a left adjoint $F^L$ and there is a natural equivalence
$$F^L(c) \simeq \colim_{\fcat{I}^{\op}} F_i^L(c_i),$$
where $c\in \lim_\fcat{I} \fcat C_i$ and $c_i = p_i(c)$. %Dually, if all the transition functors $\fcat C_i \to \fcat C_j$ admit left adjoints, $\fcat D$ has $\fcat{I}^\op$-indexed colimits, and all functor $F_i$ admit left adjoints $F_i^L$, then $F$ admits left adjoint $F^L$ and
%$$F^L(c) \simeq \colim_{\fcat{I}^\op} F_i^L(c_i), $$
%where $c\in \lim_\fcat{I} \fcat C_i$ and $c_i = p_i(c)$.
\begin{proof}
Let $\const_{\fcat{D}}\in \Fun(\fcat{I},\Cat)$ denote the constant diagram with value in $\fcat{D}$. By assumption, the functors $F_i$ assemble into a strict natural transformation $\eta\colon \const_{\fcat{D}} \to \fcat{C}_\bullet$ and the functor $F$ is the composite 
$$F\colon \fcat{D} \tto \lim_{\fcat{I}} (\const_{\fcat{D}}) \xrightarrow{\lim (\eta)} \lim_{\fcat{I}} \fcat{C}_\bullet.$$
By embedding strict limits into right lax limits, it suffices to construct a left adjoint to the composite
$$\fcat{D} \xrightarrow{\const} \Fun(\fcat{I}^{\op},\fcat{D}) \simeq \rlaxlim_{\fcat{I}} (\const_{\fcat{D}}) \xrightarrow{\rlaxlim(\eta)} \rlaxlim_{\fcat{I}} \fcat{C}_\bullet.$$
By \cref{lemma: right adjoint to oplax limit_lax is strict}, the left adjoint to $\rlaxlim(\eta)$ is induced to its left mate $\eta^L$, which is assembled from the left adjoints $F_i^L, i\in \fcat{I}$. Finally, we note that the left adjoint to the constant diagram functor $\const\colon \fcat{D} \to \Fun(\fcat{I}^{\op},\fcat{D})$ is given by the $\fcat{I}^\op$-indexed colimit.
\end{proof}
\end{cor}

In the rest of section, we consider the functors between right lax limits induced by changing the diagram. We recall first the relative version of Kan extensions.
\begin{defn}\label{defn_relative_Kan}
Let
\[\xymatrix{
\fcat C \ar[r]^-F \ar[d]^p & \fcat E \ar[d]^q \\
\fcat D \ar[r] & \fcat B
}\]
be a commutative diagram of categories. A \emdef{$q$-right} (or \emdef{relative over $\fcat B$ right) Kan extension of $F$ along $p$} is a functor
$$\Ran_p^{\fcat B} F \colon \fcat D \tto \fcat E \in \Fun_{/\fcat{B}}(\fcat{D}, \fcat{E})$$
which is characterized by the universal property
$$\Hom_{\Fun_{/\fcat B}(\fcat D, \fcat E)}(G, \Ran_p^{\fcat B} F) \simeq \Hom_{\Fun_{/\fcat B}(\fcat C, \fcat E)}(G \circ p, F), \;\; G \in {\Fun_{/\fcat B}(\fcat D, \fcat E).}$$

Dually, a \emdef{$q$-left Kan extension of $F$ along $p$} is a functor, if it exists,
$$\Lan_p^{\fcat B} F \colon \fcat D \tto \fcat E \in \Fun_{/\fcat{B}}(\fcat{D}, \fcat{E})$$
characterized by
$$\Hom_{\Fun_{/\fcat B}(\fcat D, \fcat E)}(\Lan_p^{\fcat B} F, G) \simeq \Hom_{\Fun_{/\fcat B}(\fcat C, \fcat E)}(F, G \circ p), \;\; G \in {\Fun_{/\fcat B}(\fcat D, \fcat E).}$$
\end{defn}
\begin{rem}
In the setting above, if the relative right (resp.\@ left) Kan extension along $p$ exists for every $F \in \Fun_{/\fcat B}(\fcat C, \fcat E)$, then the assignment $F \mapsto \Ran_p^{\fcat B} F$ (resp.\@ $\Lan_p^{\fcat B} F$) gives a right (resp.\@ left) adjoint to the precomposition functor
$$p^*\colon \Fun_{/\fcat B}(\fcat D, \fcat E) \tto \Fun_{/\fcat B}(\fcat C, \fcat E).$$
\end{rem}
\begin{ex}\label{ex_adjoint_relative_Kan}
Let
\[\xymatrix{
\fcat C \ar[r]^-F \ar[d]^p & \fcat E \ar[d]^q \\
\fcat D \ar[r]^-\pi & \fcat B
}\]
be a commutative diagram of categories. Assume that $p$ admits a left adjoint $p^L$ over $\fcat B$, i.e.\@$p^L$ is such that the induced natural transformation
$$\xymatrix{\pi \ar[r] & \pi \circ p \circ p^L \ar[r]^-\sim & q \circ F \circ p^L}$$
is an equivalence. Then the adjunction $p^L \dashv p$ induces an adjunction of functor categories
$$p^*\colon \Fun_{/\fcat B}(\fcat D, \fcat E) \xymatrix{\ar@<0.5ex>[r] & \ar@<0.5ex>[l]}  \Fun_{/\fcat B}(\fcat C, \fcat E) \colon p^{L,*},$$
i.e.\@ the relative Kan extension $\Ran^{\fcat B}_p$ exists in this case and it is equivalent to $p^{L,*}$.

Dually, if $p$ admits a right adjoint $p^R$ over $\fcat B$, then the relative left Kan extension $\Lan_p^{\fcat B}$ exists and it is equivalent to $p^{R,*}$.
\end{ex}

Recall that (co)limits may be considered as special cases of Kan extensions. Moreover, the existence of specific (co)limits in the target category provides a sufficient condition for Kan extensions to exist. In the following we review similar results in the relative setting.
\begin{defn}[{Compare to \cite[Definition 4.3.1.1]{Lur_HTT}}]
Let 
\[\xymatrix{
\fcat K \ar[r]^-{s}\ar@{^(->}[d]^i & \fcat E \ar[d]^q \\
\fcat K^{\triangleleft} \ar[r]^-{f} & \fcat B
}\]
be a commutative diagram of categories, where $\fcat K^{\triangleleft}$ denotes the category $\fcat K$ with freely added initial object $\eset$. We call a functor
$$\overline s\colon \fcat K^{\triangleleft} \tto \fcat E$$
a \emdef{$q$-limit diagram} if it has the universal property of the relative right Kan extension $\Ran^{\fcat B}_i s$. In this case the image of the cone point $\overline s(\eset)$ is called a \emdef{$q$-limit of $s$}, and we denote
$$\sideset{}{^{\fcat B}}\lim_{\fcat K} s := \overline s(\eset).$$

As usual, \emdef{$q$-colimit diagram} and \emdef{$q$-colimit $\colim^{\fcat B}_{\fcat K} s$} are defined as a $q$-limit diagram and a $q$-limit of $s$, if they exist, in the opposite categories.
\end{defn}

\begin{ex}\label{ex_rel_lim_initial_obj}
Let 
\[\xymatrix{
\fcat K \ar[r]^-{s}\ar@{^(->}[d]^i & \fcat E \ar[d]^q \\
\fcat K^{\triangleleft} \ar[r]^-{f} & \fcat B
}\]
be a commutative diagram of categories. Assume that $\fcat K$ admits an initial object $\eset_{\fcat K}$ and that the image under $f$ of the natural map
\begin{equation}\label{eq_relative_lim_initial_obj}
\eset \tto i(\eset_{\fcat K})
\end{equation}
is an equivalence in $\fcat B$. Then we claim that the $q$-limit of $s$ exists and is equivalent to $s(\eset_{\fcat K})$. Indeed, the inclusion functor $i$ admits a left adjoint $p$ sending the formal initial object $\eset$ to the actual one $\eset_{\fcat K}$. Moreover, $p$ is a functor over $\fcat B$ if and only if \eqref{eq_relative_lim_initial_obj} is an equivalence. So the assertion follows from \Cref{ex_adjoint_relative_Kan}.
\end{ex}

Relative Kan extensions were studied thoroughly in \cite[Section 4.3.2]{Lur_HTT} provided the functor $p$ is fully faithful.
\begin{prop}[{\cite[Lemma 4.3.2.13 and Proposition 4.3.2.17]{Lur_HTT}}]\label{proposition: relative kan extensions and relative limits}
Let
\[\xymatrix{
\fcat C \ar[r]^-F\ar@{^(->}[d]^p & \fcat E \ar[d]^q \\
\fcat D \ar[r] & \fcat B
}\]
be a commutative diagram of categories such that the functor $p$ is fully faithful. Suppose that the relative limit $\lim_{\fcat C_{d/}}^{\fcat B} F$ (resp.\@ relative colimit $\colim_{\fcat C_{/d}}^{\fcat B} F$) exists for every $d\in \fcat D$. Then the relative right (resp.\@ left) Kan extension $\Ran_p^{\fcat B} F$ exists and there is a natural equivalence
$$\Ran_p^{\fcat B} F(d) \simeq \sideset{}{^{\fcat B}}\lim_{c \in \fcat C_{d/}} F(c) \qquad \text{resp.} \ \Lan_p^{\fcat B} F(d) \simeq \sideset{}{^{\fcat B}}\colim_{c \in \fcat C_{/d}} F(c)$$
for any $d \in \fcat D$. \qed
\end{prop}
\begin{rem}
We present the theory of relative Kan extensions in a slightly different order than~\cite{Lur_HTT}. Namely, in \cite[Section~4.3.2]{Lur_HTT}, relative Kan extensions are defined in terms of the relative limits as in \cref{proposition: relative kan extensions and relative limits} and then \cite[Proposition 4.3.2.17]{Lur_HTT} guarantees that the definition in loc. cit. coincides with our \Cref{defn_relative_Kan}. 
\end{rem}
Using the argument similar to the one in \cite[Section 2.2.1]{QuigleyShah_parametrizedTate}, we can assume the functor $p$ is not necessarily fully faithful in \cref{proposition: relative kan extensions and relative limits}.
\begin{prop}\label{relative_Ran_lim_formula}
Let
\[\xymatrix{
\fcat C \ar[r]^-F \ar[d]^p & \fcat E \ar[d]^q \\
\fcat D \ar[r] & \fcat B
}\]
be a commutative diagram of categories such that the relative limit $\lim_{\fcat C_{d/}}^{\fcat B} F$ (resp.\@ relative colimit $\colim_{\fcat C_{/d}}^{\fcat B} F$) exists for every $d\in \fcat D$. Then the relative right (resp.\@ left) Kan extension $\Ran_p^{\fcat B} F$ exists and there is a natural equivalence
$$\Ran_p^{\fcat B} F(d) \simeq \sideset{}{^{\fcat B}}\lim_{c \in \fcat C_{d/}} F(c) \qquad \text{resp.}\ \Lan_p^{\fcat B} F(d) \simeq \sideset{}{^{\fcat B}}\colim_{c \in \fcat C_{/d}} F(c)$$
for any $d \in \fcat D$.

\begin{proof}
The assertions about left and right Kan extensions are formally equivalent by passing to the opposite categories, so we only treat the latter case. Let $\fcat M \to (\Delta^1)^\op$ be the Cartesian fibration corresponding to the functor $p$. There are fully faithful embeddings $i\colon \fcat C \inj \fcat M$ and $j\colon \fcat D \inj \fcat M$ and $j$ admits the right adjoint $\overline p \colon \fcat M \to \fcat D$ such that $p \simeq \overline p \circ i$. Let us consider $\fcat M$ as a category over $\fcat B$ via the composite $\fcat M \to \fcat D \to \fcat B$. Then the adjunction $j \dashv \overline p$ induces an adjunction
$$\overline p^* \colon \Fun_{/\fcat B}(\fcat D, \fcat E) \xymatrix{\ar@<0.5ex>[r] & \ar@<0.5ex>[l]} \Fun_{/\fcat B}(\fcat M, \fcat E) : j^*.$$
In particular, $j^* \simeq \Ran^{\fcat B}_{\overline p}$.

Note that, for an object $d \in \fcat D$, there is a natural equivalence $\fcat C_{d/} \simeq \fcat C_{j(d)/}$. By assumption, the relative limit $\lim_{\fcat C_{m/}}^{\fcat B} F$ exists for any $m \simeq j(d) \in \fcat M$, $d \in \fcat{D}$. Moreover, the relative limit $\lim_{\fcat C_{i(c)/}}^{\fcat B} F$, $c\in \fcat C$ also exists since the diagram $\fcat C_{i(c)/} \simeq \fcat{C}_{c/}$ admits an initial object, see \Cref{ex_rel_lim_initial_obj}. By \cref{proposition: relative kan extensions and relative limits}, $\Ran_i^{\fcat B} F$ exists and for $m \in \fcat M$ and
$$\Ran_i^{\fcat B} F(m) \simeq \sideset{}{^{\fcat B}}\lim_{c \in \fcat C_{m/}} F(c).$$
It follows that $\Ran_p^{\fcat B} F$ exists and we have
$$\Ran_p^{\fcat B} F(d) \simeq \Ran_i^{\fcat B} (j(d)) \simeq \sideset{}{^{\fcat B}} \lim_{c \in \fcat C_{d/}} F(c)$$
for every $d \in \fcat D$.
\end{proof}
\end{prop}

\Cref{relative_Ran_lim_formula} reduces the computation of the relative Kan extensions to the computation of the relative limits. Under additional assumptions on the functor $q$, the latter can be described in terms of the usual limits.
\begin{prop}\label{rel_limits_along_Cart_fib_preliminary}
Let $\fcat K$ be a small diagram and let 
\[\xymatrix{
\fcat K \ar[r]^-{s}\ar@{^(->}[d] & \fcat E \ar[d]^q \\
\fcat K^{\triangleleft} \ar[r]^-{f} & \fcat B
}\]
be a commutative diagram of categories such that $q\colon \fcat E \to \fcat B$ is a Cartesian fibration and for each $b' \in \fcat B$ the fiber $\fcat E_{b'}$ admits limits for all diagrams indexed by $\fcat K$. Let $S\colon \fcat K \times \Delta^1 \to \fcat{E}$ be a Cartesian lifting of $\fcat K \times \Delta^1 \to \fcat K^{\triangleleft} \xrightarrow{s} \fcat{E}$. Suppose that the canonical morphism
$$\phi_e\left(\lim_{\fcat K} S|_{\fcat{K}\times \{0\}}\right) \tto \lim_{\fcat K} \left(\phi_e\circ S|_{\fcat{K}\times \{0\}}\right) \in \fcat{E}_{b} $$
is an equivalence for all transition functors $\phi_e\colon \fcat{E}_{b} \to \fcat{E}_{b'}$ along all edges $e\colon b'\to b$, where $b$ is the image of the cone point in $\fcat K^{\triangleleft}$ under $f$. Then $\fcat{B}$-relative limit $\lim_{\fcat K}^{\fcat B}(s) \in \fcat{E}_{b}$ exists and it is equivalent to
$$\lim_{\fcat K} S|_{\fcat{K}\times \{0\}}  \in \fcat E_b.$$

\begin{proof}
By~\cite[Proposition 4.3.1.10]{Lur_HTT}, the absolute limit $\lim_{\fcat K} S|_{\fcat{K}\times \{0\}}  \in \fcat E_b$ is also a $q$-relative limit. This is the only input used in the proof of~\cite[Corollary 4.3.1.11]{Lur_HTT}, which implies the assertion.
\end{proof}
\end{prop}

\begin{rem}\label{remark: informal relative limits for cartesian}
The statement of \cref{rel_limits_along_Cart_fib_preliminary} may seem rather technical, so here is an informal idea. We observe that the objects $s(k), k\in \fcat{K}$ may belong to different fibers of the Cartesian fibration $q\colon \fcat{E} \to \fcat{B}$, whereas the relative limit $\lim_{\fcat K}^{\fcat B}(s)$ is an object of the fiber $\fcat{E}_{b}$, where $b=f(\eset)$. To resolve this asymmetry, one can transport the images $s(k), k\in \fcat{K}$ to $\fcat{E}_b$ using the assumption that $q$ is a Cartesian fibration. In this way, one get the diagram $S|_{\fcat{K}\times \{0\}}\colon \fcat{K} \to \fcat{E}_b$ in the fiber $\fcat{E}_b$ and the claim is that the (absolute) limit of $S|_{\fcat{K}\times \{0\}}$ in $\fcat{E}_b$ is exactly the relative limit $\lim_{\fcat K}^{\fcat B}(s)$.
\end{rem}

\begin{cor}[{\cite[Corollary 4.3.1.11]{Lur_HTT} and its proof}]\label{rel_limits_along_Cart_fib}
Let $q\colon \fcat E \to \fcat B$ be a Cartesian fibration and let $\fcat K$ be a small diagram. Suppose that
\begin{enumerate}
\item
for each $b \in \fcat B$, the fiber $\fcat E_b$ admits limits for all diagrams indexed by $\fcat K$;
\item
for each edge $e\colon b \to b^\prime$, the transition functor $\phi_{e} \colon \fcat E_{b^\prime} \to \fcat E_b$ preserves limits of diagrams indexed by $\fcat K$. 
\end{enumerate}
Then for any commutative diagram
\[\xymatrix{
\fcat K \ar[r]^-{s}\ar@{^(->}[d] & \fcat E \ar[d]^q \\
\fcat K^{\triangleleft} \ar[r]^-{f} & \fcat B
}\]
the $\fcat{B}$-relative limit $\lim_{\fcat K}^{\fcat B}(s) \in \fcat{E}_{b}$ exists, where $b$ is the image of the cone point in $\fcat K^{\triangleleft}$ under $f$. Moreover, if $S\colon \fcat K \times \Delta^1 \to \fcat{E}$ is a Cartesian lifting of $\fcat K \times \Delta^1 \to \fcat K^{\triangleleft} \xrightarrow{s} \fcat{E}$, then $\lim_{\fcat K}^{\fcat B}(s) \in \fcat{E}_{b}$ is equivalent to
$$\lim_{\fcat K} S|_{\fcat{K}\times \{0\}}  \in \fcat E_b. \eqno\hbox{\qedsymbol}$$
\end{cor}

\begin{rem}\label{remark: cartesian transport}
In the setting of \cref{rel_limits_along_Cart_fib}, we note that $S|_{\fcat{K}\times \{0\}}(k) \simeq \phi_{f_k}(s(k))$, where $\phi_{f_k}\colon \fcat{E}_k \to \fcat{E}_b$ is the transition functor along an essentially unique map $f_k\colon f(k)\to b$ in the image of $f\colon \fcat{K}^{\triangleleft}\to \fcat{B}$. 
\end{rem}

\begin{rem}
By passing to the opposite categories in \Cref{rel_limits_along_Cart_fib_preliminary} and \Cref{rel_limits_along_Cart_fib}, one obtains similar results about relative colimits provided $q$ is a coCartesian fibration.
\end{rem}

Finally, we compute the right adjoint between right lax limits of categories.

\begin{cor}\label{right_adjoint_llax_limits}
Let $p\colon \fcat J \to \fcat I$ be a morphism of small categories and let $\fcat C_\bullet \colon \fcat I \to \Cat$ be a diagrams of categories indexed by $\fcat I$. Suppose that 
\begin{enumerate}
\item for each pair $i,i^\prime \in \fcat I$, the category $\fcat C_{i'}$ admits small limits indexed by $(\fcat{J}_{/i})^{\op}$; 
\item for each map $f\colon i \to i'$, the transition functor $\phi_{f}\colon \fcat C_i \to \fcat C_{i'}$ preserves small limits indexed by $(\fcat{J}_{/i})^{\op}$. 
\end{enumerate} 
Let
$$F\colon \rlaxlim_{\fcat I} \fcat C_\bullet \tto \rlaxlim_{\fcat J} \fcat C_\bullet$$
denote the induced functor on right lax limits. Then the functor $F$ admits a right adjoint 
$$G \colon \rlaxlim_{\fcat J} \fcat C_\bullet \tto \rlaxlim_{\fcat I} \fcat C_\bullet $$
such that there is a natural equivalence
$$G(c)_i \simeq \lim_{f\colon j\to i \in (\fcat{J}_{/i})^{\op}} \phi_{f}(c_{j})$$
for $c \in \rlaxlim_{\fcat I} \fcat C_\bullet$ and $i\in \fcat I$.

\begin{proof}
Let $q\colon \fcat E \to \fcat I^{\op}$ denote the Cartesian fibration corresponding to the diagram $\fcat C_\bullet\colon \fcat{I} \to \Cat$. Under this identification, we have
$$\rlaxlim_{\fcat I} \fcat C_\bullet \simeq \Fun_{/\fcat I^{\op}}(\fcat I^{\op}, \fcat E), \qquad \rlaxlim_{\fcat J} \fcat C_\bullet \simeq \Fun_{/\fcat I^{\op}}(\fcat J^{\op}, \fcat E),$$
and the functor $F$ is given by the precomposition with $p$. By \cref{defn_relative_Kan}, the right adjoint to
$$F \simeq p^*\colon \Fun_{/\fcat I^{\op}}(\fcat I^{\op}, \fcat E) \tto \Fun_{/\fcat I^{\op}}(\fcat J^{\op}, \fcat E)$$ 
is the relative right Kan extension $\Ran_p^{\fcat I^{\op}}(-)$.

Let $s\colon \fcat J^{\op} \to \fcat E$ be the section over $\fcat I$ corresponding to the object $c=\{c_j\} \in \rlaxlim_{\fcat J} \fcat C_\bullet$. By \Cref{relative_Ran_lim_formula}, we have
$$(\Ran_p^{\fcat I^{\op}} s)(i) \simeq \sideset{}{^{\fcat I^{\op}}}\lim_{\fcat J^{\op}_{i/}}\left(\fcat J^{\op}_{i/} \to \fcat{J}^{\op} \xrightarrow{s} \fcat{E}\right) $$
for $i \in \fcat I$. Finally, by \Cref{rel_limits_along_Cart_fib}, the relative limit $\lim^{\fcat I^{\op}}_{\fcat J_{i/}^{\op}}$ can be computed as the usual limit over the induced diagram in the fiber over $i$:
\[\sideset{}{^{\fcat I^{\op}}}\lim \left(\fcat J^{\op}_{i/} \to \fcat{J}^{\op} \xrightarrow{s} \fcat{E}\right) \simeq \lim_{f\colon j\to i\in \fcat J^{\op}_{i/}} \phi_{f}(s(j)) = \lim_{f\colon j\to i \in (\fcat{J}_{/i})^{\op}} \phi_{f}(c_j).\qedhere\]
\end{proof}
\end{cor}

\begin{cor}\label{corollary: right adjoint to evaluation on the initial object}
Let $\fcat{C}_\bullet \colon \fcat{I} \to \Cat$ be a diagram of categories indexed by $\fcat{I}$. Suppose that $\varnothing \in \fcat{I}$ is an initial object. Then there exists the adjoint pair
$$ev_{\varnothing} \colon \rlaxlim_{\fcat{I}} \fcat{C}_\bullet \xymatrix{\ar@<+0.5ex>[r] & \ar@<+0.5ex>[l]} \fcat{C}_{\varnothing} \colon Q^{\varnothing},$$
where $ev_{\varnothing}$ is the projection on the component indexed by $\varnothing \in \fcat{I}$. Moreover, there are equivalences $$Q^{\varnothing}(c)_i \areq \phi_{\varnothing, i}(c) \in \fcat{C}_i, \;\; c\in \fcat{C}_\varnothing, \;\; i\in\fcat{C},$$ where $\phi_{\varnothing,i}\colon \fcat{C}_{\varnothing}\to \fcat{C}_i$ is the transition functor along a unique map $\varnothing \to i$. 
\end{cor} 

\begin{proof}
Apply \cref{right_adjoint_llax_limits} to the functor $p\colon \fcat{J}:=\{*\} \to \fcat{I}$, $p(*)=\varnothing$. Indeed, in this case $\fcat{J}_{/i} \simeq *$ for any $i\in \fcat{I}$. So, the assumptions from \cref{right_adjoint_llax_limits} are vacuous.
\end{proof}

\begin{cor}\label{corollary: right adjoint to evaluation on any object}
Let $\fcat{C}_\bullet \colon \fcat{I} \to \Prs^\LL$ be a diagram of presentable categories and left adjoint functors. Then there exists the adjoint pair
$$ev_{X} \colon \rlaxlim_{\fcat{I}} \fcat{C}_\bullet \xymatrix{\ar@<+0.5ex>[r] & \ar@<+0.5ex>[l]} \fcat{C}_X \colon Q^{X},$$
where $ev_{X}$ is the projection on the component indexed by an object $X \in \fcat{I}$. Moreover, if for an object $c\in \fcat{C}_X$, the natural morphism
\begin{equation}\label{equation: right adjoint to evaluation assumption}
\phi_{g}\left(\lim_{f\in \Hom(X,i)}\phi_{f}(c)\right) \tto \lim_{f\in \Hom(X,i)}\phi_{g\circ f}(c)
\end{equation}
is an equivalence for any morphism $g\colon i\to j \in \fcat{I}$, then there are equivalences
$$Q^X(c)_i \areq \lim_{f\in \Hom(X,i)}\phi_{f}(c), \;\; i\in \fcat{I}.$$
\end{cor}

\begin{proof}
The right adjoint $Q^X$ exists by the adjoint functor theorem. By \cref{defn_relative_Kan}, the right adjoint $Q^X$ is given by the relative right Kan extension $\Ran_p^{\fcat I^{\op}}(-)$, where $p\colon \{*\} \to \fcat{I}$ is a functor such that $p(*)=X$. Finally, we compute $Q^X(c)\simeq \Ran_p^{\fcat I^{\op}}(c)$ by \cref{relative_Ran_lim_formula} and \cref{rel_limits_along_Cart_fib_preliminary}.
\end{proof}

\section{Object inversion}\label{section: object inversion}
This appendix is devoted to some technical aspects, which arise in the construction of $\otimes$-inversion of a collection of objects in a symmetric monoidal category. The content seems to be folklore, but we were not able to locate certain assertions in the literature. The main result is \cref{lemma: X-linearization of sigma-omega} which presents any object in the category $\fcat{V}[X^{-1}]$ with inverted $X$ as a canonical colimit of objects from $\fcat{V}$.
%\todo{The main text requires only \cref{definition: object inversion}, \cref{proposition: inversion symmetric}, \cref{corollary: right adjoint to inversion many objects}, and \cref{lemma: X-linearization of sigma-omega}. }

\begin{defn}[Proposition~2.9(1) in \cite{Robalo15}]\label{definition: object inversion}
Let $\fcat V$ be a presentably symmetric monoidal category and let $\{V_\alpha\}\in \fcat{V}$ be a set of objects in $\fcat V$. A symmetric monoidal functor 
$$\Sigma^\infty_{\fcat{V}}\colon \fcat V \tto \fcat V[\{V_\alpha^{-1}\}]$$
is the \emdef{formal inversion} of $\{V_\alpha\}$ in $\fcat{V}$ if, for every presentably symmetric monoidal category $\fcat{D}\in \Prs^\LL$, the induced functor
$$\Fun^{\LL,\otimes}(\fcat{V}[\{V_\alpha^{-1}\}], \fcat{D}) \tto \Fun^{\LL,\otimes}(\fcat{V}, \fcat{D}), \;\; F\mapsto F\circ \Sigma^{\infty}_\fcat{V}$$
is \emph{fully faithful} and the essential image is spanned by the symmetric monoidal functors $F\colon \fcat{V} \to \fcat{D}$ such that the objects $F(V_{\alpha})\in\fcat{D}$ are \emph{$\otimes$-invertible}.
\end{defn}

\begin{prop}\label{proposition: object inversion modules}
Let $\Sigma^\infty_{\fcat{V}}\colon \fcat V \tto \fcat V[\{V_\alpha^{-1}\}]$ be a formal inversion. Then the forgetful functor
$$\Mod_{\fcat{V}[\{V_\alpha^{-1}\}]}(\Prs^\LL) \tto \Mod_{\fcat{V}}(\Prs^\LL)$$
is fully faithful and the essential image is spanned by all presentable categories $\fcat M \in\Mod_{\fcat{V}}(\Prs^\LL)$ such that the objects $V_\alpha \in \fcat{V}$ act on $\fcat M$ by equivalences.
\end{prop}

\begin{proof}
See~\cite[Proposition~2.9(2)]{Robalo15}.
\end{proof}

\begin{defn}[Definition~2.16 in \cite{Robalo15}]\label{definition: symmetric object}
An object $X\in \fcat{V}$ of a symmetric monoidal category is \emdef{symmetric} if the cyclic permutation $\sigma\colon X^{\otimes 3} \to X^{\otimes 3}$ is equivalent to the identity map $\Id\colon X^{\otimes 3} \to X^{\otimes 3}$.
\end{defn}

Let $\fcat{M} \in \Mod_{\fcat{V}}(\Prs^\LL)$ be a module over a presentably symmetric monoidal category $\fcat{V}$. Let $\{V_\alpha\}\in \fcat{V}$ be a set of objects. We write $\fcat{M}[\{V_\alpha^{-1}\}]\in \Mod_{\fcat{V}[\{V_\alpha^{-1}\}]}(\Prs^\LL)$ for the tensor product $\fcat{M}[\{V_\alpha^{-1}\}]=\fcat{M}\otimes_{\fcat{V}} \fcat{V}[\{V_\alpha^{-1}\}]$. Note that the assignment
$$\fcat{M} \xymatrix{\ar@{|->}[r] &} \fcat{M}[\{V_\alpha^{-1}\}] $$
is the left adjoint to the fully faithful embedding of~\cref{proposition: object inversion modules}.

\begin{prop}\label{proposition: inversion symmetric}
Let $\fcat V$ be a presentably symmetric monoidal category and let $X\in \fcat{V}$ be a symmetric object. Then there are equivalences
\begin{align*}
\fcat{M}[X^{-1}] &\simeq \colim \left(\fcat{M} \xymatrix{\ar[r]^-{-\otimes X} &} \fcat{M} \xymatrix{\ar[r]^-{-\otimes X} &} \cdots \right) \\
&\simeq \lim \left(\cdots \xymatrix{\ar[r]^-{[X,-]_{\fcat{M}}} &} \fcat{M} \xymatrix{\ar[r]^-{[X,-]_{\fcat{M}}} &} \fcat{M} \right)
\end{align*}
in the category $\Mod_{\fcat{V}}(\Prs^\LL)$ of $\fcat{V}$-modules, where (co)limits above are computed in the category $\Prs^\LL$. In particular, if $\fcat{M}$ is stable, then $\fcat{M}[X^{-1}]$ is stable as well.
\end{prop}

\begin{proof}
See~\cite[Corollary 2.22 and Corollary 2.23]{Robalo15}.
\end{proof}

The formal inversion $\Sigma^\infty_{\fcat{V}}\colon \fcat{V} \to \fcat{V}[\{V_\alpha^{-1}\}]$ is $\fcat{V}$-linear. Therefore it induces a colimit-preserving functor
$$\Sigma^\infty_{\fcat{M}}\colon \fcat{M} \tto \fcat{M}[\{V_\alpha^{-1}\}] $$
of presentable categories. We write
$$\Omega^\infty_{\fcat{M}}\colon \fcat{M}[\{V_\alpha^{-1}\}] \tto \fcat{M}$$
for the right adjoint to $\Sigma^\infty_{\fcat{M}}$ which exists by the adjoint functor theorem \cite[Corollary~5.5.2.9]{Lur_HTT}.

\begin{prop}\label{proposition: right adjoint to inversion}
Let $\fcat{M}\in \Mod_{\fcat{V}}(\Prs^\LL)$ be a module over a presentably symmetric monoidal category $\fcat{V}$ and let $X\in \fcat{V}$ be a symmetric object. Suppose that the endofunctor
$$[X,-]_{\fcat{M}}\colon \fcat{M} \tto \fcat{M}$$
preserves filtered colimits. Then the natural transformation
$$\colim_{n\in \mathbb{N}} [X^{\otimes n}, X^{\otimes n}\otimes -]_{\fcat{M}} \tto \Omega^\infty_{\fcat{M}} \Sigma^\infty_{\fcat{M}}(-)$$
is an equivalence.
\end{prop}

\begin{proof}
By \cref{proposition: inversion symmetric}, the right adjoint functor $\Omega^\infty_{\fcat{M}} \colon \fcat{M}[X^{-1}] \to \fcat{M}$ is the composite
$$\Omega^{\infty}_{\fcat{M}}\colon \fcat{M}[X^{-1}] \simeq \lim_{\mathbb{N}^{\op}_{\geq}} \left(\cdots \xymatrix{\ar[r]^-{[X,-]_{\fcat{M}}} &} \fcat{M} \xymatrix{\ar[r]^-{[X,-]_{\fcat{M}}} &} \fcat{M} \right) \xhookrightarrow{\iota} \rlaxlim_{\mathbb{N}^{\op}_{\geq}} \left(\cdots \xymatrix{\ar[r]^-{[X,-]_{\fcat{M}}} &} \fcat{M} \xymatrix{\ar[r]^-{[X,-]_{\fcat{M}}} &} \fcat{M} \right) \xrightarrow{ev_0} \fcat{M}.
$$
Informally, an object of the middle category is a sequence $c=(c_n)_{n\in\mathbb{N}} \in \fcat{M}$ together with the structure maps $f_n\colon c_n \to [X,c_{n+1}]_{\fcat{M}}$. Then, $ev_0(c)=c_0$ and $\fcat{M}[X^{-1}]$ is the full subcategory spanned by the sequences such that each structure map $f_n$ is an equivalence for $n\geq 0$.

By~\cite[Theorem 5.3.6]{haugseng2021lax}, we have
\begin{align*}
\rlaxlim_{\mathbb{N}^{\op}_{\geq}} \left(\cdots \xymatrix{\ar[r]^-{[X,-]_{\fcat{M}}} &} \fcat{M} \xymatrix{\ar[r]^-{[X,-]_{\fcat{M}}} &} \fcat{M} \right) &\simeq \llaxlim_{\mathbb{N}_{\geq}} \left(\fcat{M} \xymatrix{\ar[r]^-{-\otimes X} &} \fcat{M} \xymatrix{\ar[r]^-{-\otimes X} &} \cdots \right)\\
&\simeq \left(\rlaxlim_{\mathbb{N}_{\geq}} \left(\fcat{M}^{\op} \xymatrix{\ar[r]^-{-\otimes X} &} \fcat{M}^{\op} \xymatrix{\ar[r]^-{-\otimes X} &} \cdots \right)\right)^{\op}.
\end{align*}
Together with \cref{corollary: right adjoint to evaluation on the initial object}, this implies that $\ev_0$ admits a left adjoint 
$$\Sigma^{\infty}_{\mathrm{pre}} \colon \fcat{M} \tto \rlaxlim_{\mathbb{N}^{\op}_{\geq}} \left(\cdots \xymatrix{\ar[r]^-{[X,-]_{\fcat{M}}} &} \fcat{M} \xymatrix{\ar[r]^-{[X,-]_{\fcat{M}}} &} \fcat{M} \right)$$
such that $(\Sigma^{\infty}_{\mathrm{pre}}(c))_n\simeq X^{\otimes n}\otimes c$, $n\geq 0$.

Note that the fully faithful embedding $\iota$ preserves arbitrary limits and filtered colimits, therefore $\iota$ admits a left adjoint
$$L\colon \rlaxlim_{\mathbb{N}^{\op}_{\geq}} \left(\cdots \xymatrix{\ar[r]^-{[X,-]_{\fcat{M}}} &} \fcat{M} \xymatrix{\ar[r]^-{[X,-]_{\fcat{M}}} &} \fcat{M} \right) \tto \lim_{\mathbb{N}^{\op}_{\geq}} \left(\cdots \xymatrix{\ar[r]^-{[X,-]_{\fcat{M}}} &} \fcat{M} \xymatrix{\ar[r]^-{[X,-]_{\fcat{M}}} &} \fcat{M} \right). $$
We claim that $$L(c)_k \simeq \colim_{n\in \mathbb{N}} \left( \cdots \to [X^{\otimes n}, c_{n+k}]_{\fcat{M}} \xrightarrow{\phi_{n}} [X^{\otimes n+1},c_{n+k+1}]_{\fcat{M}} \to \cdots \right)$$
for $c=(c_n)\in \rlaxlim_{\mathbb{N}^{\op}_{\geq}} \left(\cdots \to \fcat{M} \to \fcat{M} \right)$ and $k\geq 0$. Here, the map $\phi_n$ is the composite
$$\phi_n\colon [X^{\otimes n}, c_{n+k}]_{\fcat{M}} \xrightarrow{[X^{\otimes n},f_{n+k}]} [X^{\otimes n}, [X,c_{n+k+1}]_{\fcat{M}}]_{\fcat{M}} \simeq [X^{\otimes n+1},c_{n+k+1}]_{\fcat{M}}.$$ 
Indeed, this formula defines an augmented endomorphism $\eta\colon \Id \to L$ of the right lax limit. By the construction, $\eta\colon c \to L(c)$ is an equivalence if $c\in \fcat{M}[X^{-1}]$. Therefore, it suffices to show that the essential image of $L$ lies in $\fcat{M}[X^{-1}]$. This follows from the assumption that the endofunctor $[X,-]_{\fcat{M}}$ preserves filtered colimits.

Finally, the assertion follows from the natural equivalence \[\Omega^{\infty}_{\fcat{M}}\Sigma^{\infty}_{\fcat{M}}(c) \simeq (L\Sigma^{\infty}_{\mathrm{pre}}(c))_0, \;\; c\in \fcat{M}.\qedhere\]
\end{proof}

\begin{cor}\label{corollary: right adjoint to inversion many objects}
Let $\fcat{M}\in \Mod_{\fcat{V}}(\Prs^\LL)$ be a module over a presentably symmetric monoidal category $\fcat{V}$ and let $\{V_\alpha \}_{\alpha\in I} \in \fcat{V}$ be a set of symmetric objects. Suppose that the endofunctors
$$[V_\alpha ,-]_{\fcat{M}}\colon \fcat{M} \tto \fcat{M}, \;\; \alpha \in I$$
preserve filtered colimits. Then the natural transformation
$$\colim_{\substack{J \subset I, \\ |J|<\infty, \\ n\in\mathbb{N}}} [V_J^{\otimes n}, V_J^{\otimes n}\otimes -]_{\fcat{M}} \tto \Omega^\infty_{\fcat{M}} \Sigma^\infty_{\fcat{M}}(-), \;\; V_J = \bigotimes_{\alpha\in J} V_{\alpha}$$
is an equivalence. \qed
\end{cor}

Next, given a $\fcat{V}$-linear category $\fcat{N}$ and a symmetric object $X \in \fcat{V}$ with some additional assumptions, we present any $c\in \fcat{N}[X^{-1}]$ as a \emph{canonical} colimit of suspension objects. To do so, we will invoke the following general construction. %\todo{repetition of additional}

\begin{notation}\label{definition: right lax linear}
We will write $\Fun_{\fcat{V}}^{\rlax}(\fcat{M}, \fcat{N})$ for the category of \emdef{right lax $\fcat{V}$-linear} functors from $\fcat{M}$ to $\fcat{N}$.
\end{notation}

\begin{construction}\label{construction: X-linearization}
Let $\fcat V$ be a presentably symmetric monoidal category, let $F\colon \fcat{M} \to \fcat{N}$ be a right lax $\fcat{V}$-linear functor, and let $X\in \fcat{V}$ be an object. Suppose that the endofunctor $$[X,-]_{\fcat{N}}\colon \fcat{N} \tto \fcat{N} $$
is continuous. We define a new right lax $\fcat{V}$-linear functor $P_XF\colon \fcat{M} \to \fcat{N}$ by the following formula
$$P_XF(-)=\colim_{n\in \mathbb{N}}\left(\cdots \to [X^{\otimes n}, F(- \otimes X^{\otimes n})]_{\fcat{N}} \xrightarrow{\mu_n} [X^{\otimes n+1}, F(-\otimes X^{\otimes n+1})]_{\fcat{N}} \to \cdots \right), $$
where the natural transformation $\mu_n$ is the composite
$$\mu_n \colon [X^{\otimes n}, F(- \otimes X^{\otimes n})]_{\fcat{N}} \to [X^{\otimes n+1}, F(- \otimes X^{\otimes n}) \otimes X]_{\fcat{N}} \to [X^{\otimes n+1}, F(-\otimes X^{\otimes n+1})]_{\fcat{N}}. $$
By construction, the endomorphism $P_X$ of the category $\Fun^{\rlax}_{\fcat{V}}(\fcat{M},\fcat{N})$ has the augmentation $$p_X\colon \Id \to P_X.$$
\end{construction}

\begin{rem}\label{remark: idea of PX}
Assume additionally that $X\otimes (-)$ is an autoequivalence of both $\fcat{M}$ and $\fcat{N}$. Then one can check that the functor $P_XF$ is a right lax $\fcat{V}$-linear functor which commutes $X\otimes -$, i.e. the natural transformation 
$$X\otimes P_XF(- )\tto P_XF(X\otimes -)$$
is an equivalence. Moreover, $P_XF$ is the universal functor among such with a natural transformation from $F$.

At the moment of writing, we are not aware if $P_XF$ has any interesting universal properties in the general case.
\end{rem}

We show that the construction $F\mapsto P_XF$ has some formal properties.

\begin{prop}\label{proposition: properties of X-linealizarions}
Let $\fcat V$ be a presentably symmetric monoidal category, let $F\colon \fcat{M} \to \fcat{N}$ and $G\colon \fcat{L} \to \fcat{M}$ be right lax $\fcat{V}$-linear functors, and let $X\in \fcat{V}$ be an object. Suppose that the endofunctors $$[X,-]_{\fcat{M}}\colon \fcat{M} \tto \fcat{M},\;\; [X,-]_{\fcat{N}}\colon \fcat{N} \tto \fcat{N}$$
are continuous. Then
\begin{enumerate}
\item the natural transformation $P_XF \to P_X(P_XF)$ is an equivalence;
\item the natural transformation $p^G\colon P_X(F\circ G) \to P_X(F\circ P_X(G))$ is an equivalence;
\item the natural transformation $p^F\colon P_X(F\circ G) \to P_X(P_X(F)\circ G)$ is an equivalence.
\end{enumerate}
\end{prop}

\begin{proof}
The statement and the proof are inspired by \cite[Proposition~1.3.1]{AroneChing11}. For the first part, we note that
\begin{align*}
P_X(P_XF(-))&\simeq \colim_{n,m\in \mathbb{N}}[X^{\otimes n}, [X^{\otimes m}, F(- \otimes X^{\otimes m}\otimes X^{\otimes n})]_{\fcat{N}}]_{\fcat{N}}\\
&\simeq \colim_{n,m\in \mathbb{N}}[X^{\otimes (n + m)}, F(- \otimes X^{\otimes (m+n)})]_{\fcat{N}}\\
&\simeq P_XF(-).
\end{align*}

For the second part, we observe that there is a natural transformation 
$$v^F\colon F\circ P_X(G) \tto P_X(F\circ G)$$
such that the composite
$$P_X(F\circ G) \xrightarrow{p^G} P_X(F \circ P_X(G)) \xrightarrow{P_X(v^F)} P_X(F\circ G)$$
is equivalent to the identity. Indeed, the natural transformation $v^F$ is induced by the compatible transformations 
$$F([X^{\otimes n}, G(- \otimes X^{\otimes n})]_{\fcat{M}}) \tto [X^{\otimes n}, FG(- \otimes X^{\otimes n})]_{\fcat{N}}, \;\; n\geq 0. $$
Therefore, the natural transformation $p^G$ is a retract of the natural equivalence 
$$p^{P_XG} \colon P_X(F\circ P_X(G)) \to P_X(F\circ P_X(P_XG)).$$
This implies that $p^G$ is an equivalence.

For the last part, we observe that there is a natural transformation 
$$u^G\colon P_X(F)\circ G \tto P_X(F\circ G)$$
such that the composite
$$P_X(F\circ G) \xrightarrow{p^F} P_X(P_X(F) \circ G) \xrightarrow{P_X(u^G)} P_X(F\circ G)$$
is equivalent to the identity. Indeed, the natural transformation $u^G$ is induced by the compatible transformations 
$$[X^{\otimes n}, F(G(-) \otimes X^{\otimes n})]_{\fcat{N}} \tto [X^{\otimes n}, FG(- \otimes X^{\otimes n})]_{\fcat{N}}, \;\; n\geq 0. $$
As before, this implies that $p^F$ is an equivalence.
\end{proof}

\begin{lem}\label{lemma: X-linearization of sigma-omega}
Let $\fcat{V}\in\Prs^\LL$ be a presentably symmetric monoidal stable category and let $\fcat{N} \in \Prs^\LL$ be a $\fcat{V}$-module category. Let $X\in \fcat{V}$ be a symmetric object. Suppose that the endofunctor $$[X,-]_{\fcat{N}}\colon \fcat{N} \tto \fcat{N} $$
is continuous and conservative. Then the natural morphism
$$ P_X(\Sigma^{\infty}_{\fcat{N}}\Omega^\infty_{\fcat{N}}) \tto P_X(\Id_{\fcat{N}[X^{-1}]}) \simeq \Id_{\fcat{N}[X^{-1}]} $$
induced by the counit is an equivalence.
\end{lem}

\begin{proof}
By \cref{proposition: inversion symmetric}, the right adjoint functor $\Omega^{\infty}_{\fcat{N}} \colon \fcat{N}[X^{-1}] \to \fcat{N}$ is continuous and conservative. Therefore, $\Omega^{\infty}_{\fcat{N}}$ is a right lax $\fcat{V}$-linear functor and the adjunction counit 
$$\Sigma^{\infty}_{\fcat{N}}\Omega^\infty_{\fcat{N}} \to \Id$$
is a natural transformation of right lax $\fcat{V}$-linear functors.

We will show that $P_X(\Sigma^{\infty}_{\fcat{N}}\Omega^\infty_{\fcat{N}})(c) \simeq c$ for any $c \in \fcat{N}[X^{-1}]$. Since $\Omega^{\infty}_{\fcat{N}}\colon \fcat{N}[X^{-1}] \to \fcat{N}$ is a conservative functor, the category $\fcat{N}[X^{-1}]$ is generated under colimits by the image of the left adjoint $\Sigma^{\infty}_{\fcat{N}}$. Therefore, it suffices to show that 
\begin{equation}\label{equation: sigma-omega, eq1}
P_X(\Sigma^{\infty}_{\fcat{N}}\Omega^\infty_{\fcat{N}})(\Sigma^{\infty}_{\fcat{N}}(d)) \simeq P_X(\Sigma^{\infty}_{\fcat{N}}\Omega^\infty_{\fcat{N}}\Sigma^{\infty}_{\fcat{N}})(d)\tto P_X(\Sigma^{\infty}_{\fcat{N}})(d) \simeq \Sigma^{\infty}_{\fcat{N}}(d)
\end{equation}
is an equivalence for all $d\in \fcat{N}$.

The map~\eqref{equation: sigma-omega, eq1} admits a section 
\begin{equation}\label{equation: sigma-omega, eq2}
P_X(\Sigma^{\infty}_{\fcat{N}}) \tto P_X(\Sigma^{\infty}_{\fcat{N}}\Omega^\infty_{\fcat{N}}\Sigma^{\infty}_{\fcat{N}})
\end{equation}
induced by the adjunction unit $\Id_{\fcat{N}} \to \Omega^\infty_{\fcat{N}}\Sigma^{\infty}_{\fcat{N}}$. By~\cref{proposition: right adjoint to inversion}, we have $P_X(\Id_{\fcat{N}})\simeq \Omega^\infty_{\fcat{N}}\Sigma^{\infty}_{\fcat{N}}$. Therefore, by \Cref{proposition: properties of X-linealizarions}, 
$$P_X(\Sigma^{\infty}_{\fcat{N}}\Omega^\infty_{\fcat{N}}\Sigma^{\infty}_{\fcat{N}}) \simeq P_X(\Sigma^{\infty}_{\fcat{N}}P_X(\Id_{\fcat{N}}))\simeq P_X(\Sigma^{\infty}_{\fcat{N}})$$
and the map~\eqref{equation: sigma-omega, eq2} is an equivalence.
\end{proof}

\printbibliography

\bigskip

\noindent Nikolai~Konovalov, {\sc Department of Mathematics, University of Chicago, Eckhart Hall, 5734 S University Ave, Chicago, IL 60637, USA}
\href{mailto:nkonovalov@uchicago.edu}{nkonovalov@uchicago.edu}

\smallskip

\noindent 
Artem~Prikhodko, {\sc HSE University, Laboratory of Mirror Symmetry, 6 Usacheva str., Moscow, Russia, 119048}
\href{mailto:artem.n.prihodko@gmail.com}{artem.n.prihodko@gmail.com}

\end{document}